\documentclass[12pt]{book}
\usepackage{amsthm,amssymb,amsmath}
\usepackage[all,cmtip]{xy}

\newdir{((}{{\lhook\kern-.1em}}                                              
\newdir{))}{{\rhook\kern-.1em}}                                              
\newdir{ >}{{}*!/-5pt/@{>}}                                                  

\newenvironment{pf}{{\bf Proof.}}{\hfill \qed}
\makeatletter
\newcommand{\subsectionr}
{\@startsection{subsection}{3}{0pt}{\baselineskip}
 {-\fontdimen2\font}{\normalfont\bfseries}}
\newcommand{\subsubsectionr}
{\@startsection{subsubsection}{3}{0pt}{\baselineskip}
 {-\fontdimen2\font}{\normalfont\bfseries}}
\renewcommand{\subsubsection}
{\@startsection{subsubsection}{3}{0pt}{\baselineskip}
 {0.05\baselineskip}{\normalfont\bfseries}}
\makeatother

\newcommand{\la}{\langle}
\newcommand{\ra}{\rangle}

\newtheorem{theorem}{\bf Theorem}[section]
\newtheorem{lemma}[theorem]{\bf Lemma}
\newtheorem{prop}[theorem]{\bf Proposition}
\newtheorem{corollary}[theorem]{\bf Corollary}

{
\newtheorem{definition}[theorem]{\bf Definition}
\newtheorem{remark}[theorem]{\bf Remark}
\newtheorem{example}[theorem]{\bf Example}

\newtheorem{condition}[theorem]{\bf Condition}}

\newcommand{\be}{\begin{equation}}             
\newcommand{\ee}{\end{equation}}        

\newfont{\bfc}{cmbsy10 scaled 1200}  
\newfont{\dr}{msbm10 scaled \magstep1}  
\newfont{\sdr}{msbm8}  
\newfont{\gl}{eufm10 scaled \magstep1}  

\DeclareFontFamily{OT1}{rsfs}{}
\DeclareFontShape{OT1}{rsfs}{n}{it}{<->rsfs10}{}
\DeclareMathAlphabet{\curly}{OT1}{rsfs}{n}{it}

\newcommand{\DGML}{\mathbf{DGMod}_{\Lambda}}
\newcommand{\DGLM}{\mathbf{DG}_{\Lambda}\mathbf{Mod}}
\newcommand{\DGMG}{\mathbf{DGMod}_{\Gamma}}
\newcommand{\DGGM}{\mathbf{DG}_{\Gamma}\mathbf{Mod}}
\newcommand{\DGAlg}{\mathbf{DGAlg}}
\newcommand{\ModZ}{\mathbf{Mod}_{\ZZ}}
\newcommand{\Mor}{\operatorname{Mor}}

\renewcommand{\AA}{{\Bbb A}}
\newcommand{\CC}{{\Bbb C}}
\newcommand{\CP}{{\Bbb CP}}
\newcommand{\DD}{{\Bbb D}}
\newcommand{\EE}{{\Bbb E}}

\newcommand{\KK}{{\mathbf K}}
\newcommand{\LL}{{\Bbb L}}
\newcommand{\NN}{{\Bbb N}}
\newcommand{\PP}{{\Bbb P}}
\newcommand{\QQ}{{\Bbb Q}}
\newcommand{\RR}{{\Bbb R}}
\renewcommand{\SS}{{\mathbf S}}
\newcommand{\TT}{{\mathbf T}}
\newcommand{\UU}{{\mathbf U}}
\newcommand{\VV}{{\Bbb V}}
\newcommand{\WW}{{\Bbb W}}
\newcommand{\ZZ}{{\Bbb Z}}

\newcommand{\bU}{{\mathbf U}}
\newcommand{\bV}{{\mathbf V}}
              
\newcommand{\nau}{\nabla_{\bU}}
\newcommand{\nav}{\nabla_{\bV}}
\newcommand{\nuu}{\nau\bU}

\newcommand{\nuv}{\nau\bV}
\newcommand{\nvu}{\nav\bU}

\newcommand{\ov}{\overline}

\newcommand{\blie}{{\frak b}}
\newcommand{\glie}{{\frak g}}
\newcommand{\hlie}{{\frak h}}
\newcommand{\klie}{{\frak k}}
\newcommand{\llie}{{\frak l}}
\newcommand{\plie}{{\frak p}}

\newcommand{\tlie}{{\frak t}}
\newcommand{\ulie}{{\frak u}}
\newcommand{\zlie}{{\frak z}}

\newcommand{\kE}{\klie_E}
\newcommand{\gE}{\glie_E}

\newcommand{\fW}{{\frak W}}
\newcommand{\fV}{{\frak V}}
\newcommand{\fX}{{\curly X}}
\newcommand{\fY}{{\curly Y}}
\newcommand{\fZ}{\operatorname{Z}}

\newcommand{\Ad}{\operatorname{Ad}}
\newcommand{\Aut}{\operatorname{Aut}}
\newcommand{\bas}{\operatorname{bas}}
\newcommand{\hor}{\operatorname{hor}}
\newcommand{\Chern}{{\mathbf C}}
\newcommand{\codim}{\operatorname{codim}}
\newcommand{\Com}{{\cal C}^{\bullet}}
\newcommand{\Coker}{\operatorname{Coker}}
\newcommand{\cusp}{\operatorname{c}}

\newcommand{\diag}{\operatorname{diag}}
\newcommand{\diam}{\operatorname{diam}}
\newcommand{\dist}{\operatorname{dist}}
\newcommand{\Diff}{\operatorname{Diff}}

\newcommand{\End}{\operatorname{End}}
\newcommand{\ev}{\operatorname{ev}}
\newcommand{\fibr}{\operatorname{fibr}}
\newcommand{\GL}{\operatorname{GL}}
\newcommand{\Gr}{\operatorname{Gr}}
\newcommand{\Hom}{\operatorname{Hom}}
\newcommand{\Id}{\operatorname{Id}}
\newcommand{\id}{\operatorname{id}}
\newcommand{\inter}{\operatorname{int}}
\newcommand{\Ind}{\operatorname{Ind}}
\newcommand{\Jac}{\operatorname{Jac}}
\newcommand{\Ker}{\operatorname{Ker}}
\newcommand{\Lie}{\operatorname{Lie}}
\newcommand{\lin}{\operatorname{lin}}
\newcommand{\loc}{\operatorname{loc}}
\newcommand{\Map}{\operatorname{Map}}
\newcommand{\pt}{\operatorname{pt}}
\newcommand{\reg}{\operatorname{reg}}
\newcommand{\rk}{\operatorname{rk}}
\newcommand{\PSL}{\operatorname{PSL}}
\newcommand{\SL}{\operatorname{SL}}
\newcommand{\Stab}{\operatorname{Stab}}
\newcommand{\SU}{\operatorname{SU}}

\newcommand{\Symp}{\operatorname{Symp}}
\renewcommand{\top}{\operatorname{top}}
\newcommand{\Tor}{\operatorname{Tor}}
\newcommand{\Tr}{\operatorname{Tr}}
\newcommand{\U}{\operatorname{U}}

\newcommand{\Vol}{\operatorname{Vol}}
\renewcommand{\exp}{\operatorname{exp}}

\newcommand{\AAA}{{\curly A}}
\newcommand{\BBB}{{\curly B}}
\newcommand{\CCC}{{\curly C}}
\newcommand{\DDD}{{\cal D}}

\newcommand{\FFF}{{\cal F}}
\newcommand{\GGG}{{\curly G}}
\newcommand{\HHH}{{\cal H}}
\newcommand{\III}{{\curly I}}
\newcommand{\LLL}{{\curly L}}
\newcommand{\MMM}{{\cal M}}
\newcommand{\NNN}{{\cal N}}
\newcommand{\OOO}{{\curly O}}
\newcommand{\PPP}{{\curly P}}
\newcommand{\RRR}{{\cal R}}
\newcommand{\SSS}{{\curly S}}

\newcommand{\VVV}{{\cal V}}
\newcommand{\WWW}{{\cal W}}
\newcommand{\YMH}{{\cal YMH}}
\newcommand{\Met}{{\cal M}et^p_2}
\newcommand{\MetB}{{\cal M}et^p_{2,B}}
\newcommand{\imag}{{\mathbf i}}
\newcommand{\qu}{/\kern-.7ex/}
\newcommand{\exh}{\to\kern-1.8ex\to}
\newcommand{\nsubset}{\subset\kern-1.8ex/}
\newcommand{\oev}{\ov{\ev}}
\newcommand{\oEE}{\ov{\EE}}
\newcommand{\oLL}{\ov{\LL}}
\newcommand{\bM}{\widetilde{\MMM}}
\newcommand{\bR}{\widetilde{\RRR}}
\newcommand{\bAA}{\widetilde{\AA}}

\newcommand{\plf}{{\pi_{\FFF}^{\LLL}}}
\newcommand{\plx}{{\pi_X{\LLL}}} 
\newcommand{\pvf}{{\pi_{\FFF}^{\VVV}}}
\newcommand{\pvx}{{\pi_X^{\VVV}}}
\newcommand{\pfx}{{\pi_X^{\FFF}}}
\newcommand{\nv}{\nabla^{\VVV}}
\newcommand{\nf}{\nabla^{\FFF}}
\newcommand{\nx}{\nabla^{X}}

\renewcommand{\Im}{\operatorname{Im}}

\newcommand{\Ej}{\EE_J}    
\newcommand{\Ez}{\EE_0}
\newcommand{\EJj}{\EE_{J,J}}

\newcommand{\FE}{\FFF^{\EE}}

\newcommand{\FpqE}{\FFF^{p+q,\EE}}

\newcommand{\FpqEz}{\FFF^{p+q,\EE_0}}

\newcommand{\FEj}{\FFF^{\EE_J}}
\newcommand{\FEJj}{\FFF^{\EE_J}_J}

\newcommand{\FpqEj}{\FFF^{p+q,\EE_J}}
\newcommand{\FpqEJj}{\FFF^{p+q,\EE_J}_J}

\newcommand{\oPhi}{\ov{\Phi}}

\newcommand{\ophi}{\ov{\phi}}
\newcommand{\omu}{\ov{\mu}}
\newcommand{\onu}{\ov{\nu}}

\begin{document}
\thispagestyle{empty}
{\sc \large
\noindent Universidad Aut{\'o}noma de Madrid \\
Facultad de Ciencias \\
Departamento de Matem{\'a}ticas}

\vskip 6cm
\begin{center}
{\bf \Huge Teor{\'\i}a de Yang--Mills--Higgs
                                      
para fibraciones simpl{\'e}cticas}
\vskip 1cm 
{\large Ignasi Mundet i Riera}
\end{center}

\vskip 6cm

\hfill
\begin{minipage}[b]{.5\linewidth}
\noindent Memoria para optar al t{\'\i}tulo de
doctor en ciencias matem{\'a}ticas.\\
Director: Oscar Garc{\'\i}a--Prada.

\vskip 1cm

\noindent Madrid, abril de 1999.
\end{minipage}

\newpage{\thispagestyle{empty}
\cleardoublepage
\thispagestyle{empty}

\vspace*{8cm}
\begin{flushright}
{\it Als meus pares}
\end{flushright}}
\newpage

\newpage{\thispagestyle{empty}
\cleardoublepage
\thispagestyle{empty}

\noindent{\bf\Large Agradecimientos}

\hfill
\begin{minipage}[b]{.6\linewidth}
{\it L'aigua era freda i em vaig enrecordar que el dia abans, al mat{\'\i},
a l'hora del casament, havia plogut fort i vaig pensar que a la tarda, 
quan aniria al parc com sempre, potser encara trobaria un toll d'aigua 
pels caminets... i dintre de cada toll, per petit que fos, hi
hauria el cel...}

\begin{flushright}
M. Rodoreda, {\sc La Pla{\c c}a del Diamant}
\end{flushright}
\end{minipage}

\vskip 1cm
Quiero agradecer en primer lugar al director de esta tesis, 
Oscar Garc{\'\i}a--Prada, la ayuda que me ha prestado y el haber
estado siembre dispuesto a escuchar. Sin su constante apoyo (y exigencia) 
esta tesis no habr{\'\i}a visto nunca la luz. 

Agradezco tambi{\'e}n a Ignacio Sols y a Vicente Mu{\~n}oz 
tantas interesantes conversaciones y caf{\'e}s. Mis 
agradecimientos a Luis {\'A}lvarez por todas esas discusiones
sobre f{\'\i}sica, matem{\'a}ticas y teor{\'\i}as gauge. Gracias tambi{\'e}n a todas las 
personas que me animaron a estudiar geometr{\'\i}a y teor{\'\i}as gauge. 
Especialmente, a Vicen{\c c} Navarro y a Sebasti{\`a} Xamb{\'o}.

Quiero agradecer a todos los miembros del Departamento de Matem{\'a}ticas
de la U.A.M. lo bien acogido que me he sentido en todo momento.
A todos los participantes del curso de representaci{\'o}n de grupos
y otros seminarios, con cuyos comentarios he aprendido tantas cosas. 
Especialmente a Rafael Hern{\'a}ndez, Jes{\'u}s Gonzalo y Jos{\'e} Manuel Marco. 
Gracias tambi{\'e}n a mis compa{\~n}eros de despacho por los buenos momentos
y por la paciencia que han tenido (particularmente a J. Guerrero por 
haber aguantado estoicamente tama{\~n}o bombardeo de cuestiones filol{\'o}gicas).

Finalment, voldria agrair als meus pares, als meus germans i a la Lidia 
la seva estimaci{\'o}, el seu recolzament, i l'haver-me fet costat en tot
moment durant tots aquests anys.}

\renewcommand{\thechapter}{\Roman{chapter}}
\tableofcontents

\chapter{Introduction}

Our aim in this thesis is to study a system of equations which generalises
at the same time the vortex equations of Yang--Mills--Higgs theory and
the holomorphicity equation in Gromov theory of pseudoholomorphic curves.
In this work we extend some results and definitions from both theories to  
a common
setting. We introduce a functional generalising Yang--Mills--Higgs functional,
whose minima coincide with the solutions to our equations. We prove a 
Hitchin--Kobayashi correspondence allowing to study the solutions 
of the equations in the Kaehler case. We give a structure of 
smooth manifold to the set of (gauge equivalence classes of) solutions to
(a perturbation of) the equations (the so-called moduli space).
We give a compactification of the moduli space, generalising 
Gromov's compactification
of the moduli of holomorphic curves. Finally, we use the moduli
space to define (under certain conditions) invariants of compact symplectic 
manifolds with a Hamiltonian almost free action of $S^1$.

In this chapter we first introduce the equations studied in this thesis.
Then we briefly recall some of the main features of Yang--Mills--Higgs and 
Gromov theories and finally we explain the contents of each of the
subsequent chapters.

\section{The equations}
\label{etheequations}
Let $K$ be a real compact and connected Lie group, and let $\klie=\Lie(K)$
be its Lie algebra. Let $X$ be a compact connected Kaehler manifold
of complex dimension $n$. Let $\omega_X$ be the symplectic structure
of $X$ and $\Lambda:\Omega^*(X)\to\Omega^{*-2}$ the adjoint of
the exterior product by $\omega_X$. Let also $E\to X$ be a $K$ principal
bundle, with the action of $K$ on the right. Finally, let $F$
be a symplectic manifold with a Hamiltonian action of $K$ on the left.
Let us call $\omega_F$ the symplectic form of $F$ and $\mu:F\to\klie^*$
the moment map. We recall that $\mu$ satisfies two properties: (1)
for any $s\in\klie$ we have $d\mu(s)=\iota_{\fX_s}\omega_F$,
where $\fX_s\in\Gamma(TF)$ is the vector field generated by the
infinitesimal action of $s$ on $F$ and (2) it is equivariant with respect
to the coadjoint action of $K$ on $\klie^*$. Let
$\pi_F:\FFF\to X$ be the fibred product $E\times_K F$.

Let $\AAA=\AAA^E$ be the space of connections on $E$ and
$\SSS=\SSS^E=\Gamma(\FFF)$ the space of sections of $\FFF$.

Let us fix a complex structure $I_F$ of $F$ which is compatible with
$\omega_F$ and invariant under the action of $K$ (such structures
always exist: see lemma 5.49 in \cite{McDS2}). When $F$ is
a complex vector space (and $I_F$ is the standard complex structure)
there is a standard way to define an operator
$\ov{\partial}_A:\Omega^0(\FFF)\to\Omega^{0,1}(\FFF)$
out of a connection $A\in\AAA$, and the sections
$\Phi\in\SSS=\Omega^0(\FFF)$ such that $\ov{\partial}_A\Phi=0$ are usually
called (pseudo)holomorphic\footnote{The prefix {\it pseudo} refers
to the fact that the complex structure defined by $A$ need not be
integrable, similarly to what happens in the theory of (pseudo)holomorphic
curves. In this thesis, however, we will use the word {\it holomorphic}
regardless the integrability of the complex structure.}
with respect to $A$. This notion can be extended to our setting as follows.
Let $T\FFF_v=\Ker d\pi_F$ be the bundle of vertical tangent vectors.
A connection $A\in\AAA$ induces a splitting
$T\FFF\simeq \pi_F^*TX\oplus T\FFF_v$, which allows to define a projection
$\alpha:T\FFF\to T\FFF_v$. On the other hand, the complex structure $I_F$
induces a complex structure on the bundle $T\FFF_v$
(here we need $I_F$ to be $K$-invariant). Then we define the 
covariant derivative (with respect to $A$) of a section $\Phi\in\SSS$
to be
\begin{equation}
d_A\Phi=\alpha\circ d\Phi\in\Omega^1(\Phi^*T\FFF_v),
\label{edefda}
\end{equation}
and the antiholomorphic part of $d_A\Phi$ to be
\begin{equation}
\ov{\partial}_A\Phi=\pi^{0,1}d_A\Phi,
\label{edefdba}
\end{equation}
where $\pi^{0,1}:\Omega^1(\Phi^*T\FFF_v)\to\Omega^{0,1}(\Phi^*T\FFF_v)$
is the projection. These definitions coincide with the classical ones
when $F$ is a vector space (in this case $\FFF\simeq T\FFF_v$
canonically).

Let us take on $\klie$ a biinvariant metric (that is, a metric invariant
under the adjoint action of $K$). Using this metric we get an equivariant
isomorphism $\klie\simeq\klie^*$ which extends to an isomorphisms of
vector bundles
\begin{equation}
E\times_{\Ad}\klie\simeq E\times_{\Ad}\klie^*
\label{eisofibr}
\end{equation}
(observe that we denote with the same symbol $\Ad$ both the adjoint
representation on $\klie$ and the coadjoint representation on $\klie^*$).
Let finally $c\in\klie$ be a central element. The equations which 
we study in this thesis are
\begin{equation}
\left\{\begin{array}{l}
\ov{\partial}_A\Phi=0,\\
\Lambda F_A+\mu(\Phi)=c, \\
F_A^{0,2}=0,\end{array}\right.
\label{eiequs}
\end{equation}
where $A\in\AAA$ is a connection,
$F_A\in\Omega^2(E\times_{Ad}\klie)$ is the curvature of $A$, and
$\Phi\in\SSS$ is a section. Observe that
$\Lambda F_A\in\Omega^0(E\times_{Ad}\klie)$ and that
$\mu(\Phi)\in\Omega^0(E\times_{Ad}\klie^*)$, so to give a sense
to the second equation we need to use the isomorphism (\ref{eisofibr}).
Hence, the equations (\ref{eiequs}) depend on the biinvariant
metric taken on $\klie$. In the third equation
$F_A^{0,2}$ refers to the piece in $\Omega^{0,2}(\gE)$ of the curvature
$$F_A\in\Omega^2(E\times_{Ad}\klie)\subset
\Omega^2(\gE)=\Omega^{2,0}(\gE)\oplus\Omega^{1,1}(\gE)
\oplus\Omega^{0,2}(\gE),$$
where $\gE=E\times_{\Ad}\glie$ and $\glie=\klie\otimes\CC$ is the
complexification of $\klie$. The condition $F_A^{0,2}=0$ is equivalent to
$F_A\in\Omega^{1,1}(\gE)$. We will write $\AAA^{1,1}$ for the set
of connections $A\in\AAA$ such that $F_A^{0,2}=0$. 
By the theorem of Newlander and Niremberg \cite{NewNi} these are the
connections which define integrable complex structures in the complex
vector bundles associated to $E$.

\subsection{Yang--Mills--Higgs theory}
Let us suppose that $F$ is a complex vector space with a Hermitian metric.
The imaginary part of the metric with reversed sign gives a symplectic
form on $F$ compatible with the complex structure. Hence $F$ is a Kaehler
manifold. Let us suppose that the action of $K$ on $F$ is given by
a linear representation $\rho:K\to \U(F)$. The equations
(\ref{eiequs}) become in this situation the equations of Yang--Mills--Higgs
theory. These equations appeared for the first time in physics, in the context
of field theory, and they have been used to model different
phaenomena (such as superconductivity when $F=\CC$, $K=S^1$ and
$\rho$ is the fundamental representation, strong and electroweak forces
when $F=\CC^n$, $K=\SU(n)$ and $\rho$ is a representation depending on
the particles coupled to the gauge fields, etc.). From a mathematical
viewpoint Yang--Mills--Higgs equations have played a prominent role in the
evolution of geometry during the last thirty years. In particular,
the study of the set of solutions (the so--called moduli space) has
been specially fruitful.

\subsubsection{Some examples}
\label{algunsexemples}
When $F=\{\pt\}$ the first equation is unnecessary, and the second one is
Hermite--Einstein equation. If $X$ is a Riemann surface, this 
equation reduces to the condition on a connection $A$ of 
being projectively flat. In general, when $X$ is a Kaehler manifold,
the Hermite--Einstein equation is related to 
the notion of stability of vector bundles coming
from Geometric Invariant Theory (see below section \ref{eihiko}).
The moduli space arising in this situation has been a central object of 
study in geometry: many people have studied its topological properties, 
its properties as a Riemannian variety or the properties of its structure 
of algebraic variety. When $X$ is a compact Riemannian four manifold,
the notion of anti--self--duality generalises the conditions $\Lambda F_A=0$
and $F_A^{0,2}=0$ and the resulting moduli space was used by Donaldson 
to define his celebrated invariants (see \cite{DoKr}).

When $F=\CC^n$, $K=\U(n;\CC)$ and $\rho$ is the fundamental representation,
we get the vortex equations, studied by Jaffe and Taubes
\cite{JT}, Bradlow \cite{Br1, Br2}, Garc{\'\i}a--Prada \cite{GP1,GP2,GP3}
and others. If $X$ is a compact Riemannian four manifod and $n=1$,
these equations generalise to Seiberg-Witten equations (see for example 
\cite{Do4, GP4}). When $F=\klie$ and $\rho$ is the adjoint representation
we get (after twisting the vector bundle $\FFF$ with the cotangent bundle
of $X$) the Higgs bundle equations, studied by 
Hitchin \cite{Hi}, Simpson \cite{Si} and Corlette \cite{Co}, and whose
study has led to important developements in Kaehler geometry
(specially in understanding the fundamental groups of compact Kaehler
manifolds, see \cite{ABCKT}). Another interesting case
is $F=\Hom(W_1,W_2)$ and $K=\U(W_1)\times\U(W_2)$, where $W_i$ are
Hermitian vector spaces. This leads to the equations for holomorphic 
triples, introduced by Garc{\'\i}a--Prada \cite{GP3} and studied also by 
Garc{\'\i}a--Prada and Bradlow in \cite{BrGP3}.

\subsubsection{The Yang--Mills--Higgs functional}
The Yang--Mills--Higgs functional evaluated at a connection $A$ and a
section $\Phi$ is defined as
\begin{equation}
\YMH_c(A,\Phi)=
\|F_A\|_{L^2}^2+\|d_A\Phi\|_{L^2}^2+\|\mu(\Phi)-c\|_{L^2}^2.
\label{edefymh}
\end{equation}
A basic result in the theory is that one can rewrite
the Yang--Mills--Higgs functional as
\begin{align}
\YMH_c(A,\Phi) &= \|\Lambda F_A+\mu(\Phi)-c\|_{L^2}^2+
2\|\ov{\partial}_A\Phi\|_{L^2}^2+4\|F_A^{0,2}\|_{L^2}^2 \notag \\
&+2\int_X\la\Lambda F_A,c\ra\omega^{[n]}
+\int_X \Tr F_A\wedge F_A\wedge \omega^{[n-2]},
\label{eymh2}
\end{align}
where $\omega=\omega_X$ and
$\omega^{[k]}=\omega^k/k!$. (In the integrals of functions 
on $X$ appearing in the sequel we will implicitly use the volume
form $\omega^{[n]}$.)
From this we deduce that, if we fix $E$ and $c$, the pairs $(A,\Phi)$
which minimise $\YMH_c$ are precisely the solutions to equations
(\ref{eiequs}). (Indeed, the terms $\int_X\la\Lambda F_A,c\ra$
and $\int_X \Tr F_A\wedge F_A\wedge \omega^{[n-2]}$ only depend on $E$, 
the representation $\rho$ and $c$.) The equality (\ref{eymh2}) can be proved
using the Kaehler identities (see for example \cite{Br1}).
This equality allows to find $L^2$ bounds on the curvature $F_A$ and the 
covariant derivative $d_A\Phi$ when $(A,\Phi)$ satisfies (\ref{eiequs}).

\subsubsection{The Hitchin--Kobayashi correspondence}
\label{eihiko}
Let $\GGG_K=\Gamma(E\times_{\Ad}K)$ be the gauge group of $E$.
Let $G=K^{\CC}$ be the complexification of $K$, and let 
$\GGG_G=\Gamma(E\times_{Ad}G)$ be the complex gauge group
($\GGG_G$ is the complexification of $\GGG_K$). The group $\GGG_K$
acts on the space of connections $\AAA$ and on the space of sections 
$\SSS=\Omega^0(\FFF)$. On the other hand, both actions
of $\GGG_K$ extend to holomorphic (with respect to a certain natural
complex structure on $\AAA\times\SSS$) actions of $\GGG_G$.

The first equation of (\ref{eiequs}) is invariant under the action of
$\GGG_G$. That is, for any $g\in\GGG_G$ and $(A,\Phi)\in\AAA\times\SSS$
such that $\ov{\partial}_A\Phi=0$ we have $\ov{\partial}_{g(A)}g(\Phi)=0$.
The third equation is also $\GGG_G$ invariant.
The second equation, however, is only invariant under the action of
$\GGG_K$. This suggest the following question: given a pair
$(A,\Phi)\in\AAA\times\SSS$, how can we know whether there exists
a gauge transformation $g\in\GGG_G$ such that $(g(A),g(\Phi))$ satisfies
$$\Lambda F_{g(A)}+\mu(g(\Phi))=c?$$
The so--called Hitchin--Kobayashi correspondence answers this question
by giving a ne\-cessary and sufficient condition for this transformation
to exist. This condition involves certain coherent subsheaves of $\FFF$ 
and coincides (when $X$ is a Riemann surface), in all the cases
studied in the literature, with the condition of stability
arising in the construction of the algebraic moduli space of
pairs $(A,\Phi)$ using Geometric Invariant Theory.
On the other hand one proves that, if such a $g$ exists, then 
it is unique modulo the action of $\GGG_K$ on $\GGG_G$ on the left.

For the case $F=\{\pt\}$ and $K=\U(n)$ (Hermite-Einstein equations)  
Narasimhan and Seshadri \cite{NSe} gave in the 60's a proof
of the Hitchin--Kobayashi correspodence on Riemann surfaces.
The method used in \cite{NSe} is of algebro-geometric nature.
Soon after the appearance of the landmark paper \cite{AB}, 
Donaldson \cite{Do1} gave a proof of the same result
using techniques from gauge theories. In \cite{Do2} Donaldson extended the 
correspondence on algebraic surfaces. Uhlenbeck and Yau \cite{UY} proved the 
correspondence on any compact Kaehler manifold (see \cite{Do3} 
for a proof valid for projective manifolds). Finally, Bartolomeis and Tian 
\cite{BarTi} gave a generalisation of the correspondence to almost complex 
compact manifolds. The case $X$ a Riemann surface, $F=\{\pt\}$ and 
arbitrary $K$ was studied by Ramanathan and Subramanian in \cite{RS}, using
the results of Ramanathan in \cite{R1} on stability of principal bundles.

Several cases of the Hitchin--Kobayashi correspondence for different choices
of vector space $F$, group $K$ and representation $\rho$ have appeared 
in the literature (see the references in section
\ref{algunsexemples}). In 1996 Banfield \cite{Ba}
proved the correspondence for any compact group $K$ and
any representation $\rho:K\to\U(F)$, where $F$ is a Hermitian
vector space. This result generalises all the aforementioned
ones, with the exception of Bartolomeis and Tian result.

\subsection{Gromov theory}
When the group $K$ acting on $F$ is trivial, the second equation in
(\ref{eiequs}) disappears. The section $\Phi$ can be seen as a map
$\Phi:X\to F$ (here we make a little abuse of notation, writing the
map with the same symbol as the section) and the first equation in 
(\ref{eiequs}) can be written $$\ov{\partial}\Phi=0.$$ This is the 
holomorphicity condition (with respect to the complex structures $I_X$ 
and $I_F$ on $X$ and $F$). A relevant situation arises when $X$ is a
Riemann surface. In this case the third equation $F_A^{0,2}=0$
is always satisfied (since $\Omega^{0,2}(X)=0$).
The holomorphic maps $\Phi$ from $X$ 
to $F$ are called (pseudo)holomorphic curves, and the study of their
moduli is the central idea of Gromov theory.

Let $X$ be a compact Riemann surface. In his celebrated paper
\cite{Gr} Gromov uses the moduli space $\MMM=\MMM(A)$ of holomorphic
maps $\Phi:X\to F$ such that $\Phi_*[X]=A$ (where $A\in H_2(F;\ZZ)$)
to the study of the symplectic topology of $F$. One of the most important
results in \cite{Gr} is a natural compactification $\ov{\MMM}$ of the
moduli $\MMM$. Some consequences of the work of Gromov is the 
non-squeezing theorem for symplectic balls or the fact that the group 
of symplectomorhisms $\Symp(M)$ of a compact symplectic manifold $M$
is closed in the group of diffeomorphisms $\Diff(M)$ with respect to
the $C^0$ topology (see \cite{AuLa,McDS1,McDS2}).

Following ideas of Witten \cite{Wi} and Kontsevich and Manin \cite{KoMa}
the moduli of curves $\MMM$ has been used to define invariants of the
symplectic structure of $F$. The idea of the definition of these invariants
consists of using the evaluation map $$\ev:X\times\MMM\to\FFF,$$
which sends $(x,\Phi)$ to $\Phi(x)$, to pullback cohomology classes
from $H^*(\FFF)$ to $H^*(X\times\MMM)$ and obtain classes in $H^*(\MMM)$ 
by means of the slant product. If $\MMM$ is a compact smooth manifold
then there is a fundamental class $[\MMM]\in H_{\top}(\MMM)$. Multiplying 
the obtained cohomology classes and pairing with $[\MMM]$ we get one of the 
Gromov-Witten invariants. In general $\MMM$ is not compact, but 
Gromov compactification allows to extend the evaluation map to
$\ov{\ev}:X\times\ov{\MMM}\to\FFF$. However, $\ov{\MMM}$ is not 
in general a smooth manifold, so the existence of a fundamental class
is not clear. Therefore, to get a rigorous definition 
of the invariants some extra work is needed. During the last
ten years the problem of giving such a rigorous definition has been
intensively studied.

A first step in this line was given in the works of Ruan \cite{Ru}, 
Ruan and Tian \cite{RuTi, RuTi2}, and McDuff and Salamon \cite{McDS1}, 
in which a rigorous definition of Gromov-Witten invariants was given
for semipositive compact symplectic manifolds. In \cite{Ru} these invariants
are used to distinguish two deformation classes of symplectic structures
on a compact differentiable manifold of real dimension 6.

More recently, however, these works have been improved, and there is 
a definition of the invariants valid for any compact symplectic manifold.
This definition has been given independently by Fukaya and Ono \cite{FuOn}, 
Li and Tian \cite{LiTi}, Ruan \cite{Ru2} and Siebert \cite{Sie}. 
There exists also a definition of Gromov-Witten invariants for
projective manifolds in terms of algebraic geometry. This definition
was given independently by Behrend \cite{Beh} and by Li and Tian
\cite{LiTi2}. Finally, Siebert \cite{Sie2} and Li and Tian \cite{LiTi3}
proved that both definitions coincide for projective manifolds.

The consequences of this developement are wide ranging.
For example, the theory of Gromov-Witten invariants is very much
related to enumerative geometry. In this way, the properties of
the invariants have been used to obtain new results such as the
number of rational curves of fixed degree and genus passing
through a certain number of points in the projective space $\PP^n$
(see \cite{KoMa,RuTi}). Others fields to which Gromow-Witten
invariants are related are mirror symmetry and the theory of
integrable systems. For a survey on these and other interesting
applications of Gromov-Witten invariants, see \cite{Ru3}.

\section{Contents of the thesis}

\subsection{Chapter 1: The equations and Yang--Mills--Higgs functional}
In the first part of this chapter we introduce the equations which will be 
studied throughout this thesis. The contents of this part coincides with
what was explained in section \ref{etheequations} of this introduction.
In the second part we define the Yang--Mills--Higgs functional evaluated
at $(A,\Phi)\in\AAA\times\SSS$:
\begin{equation}
\YMH_c(A,\Phi)=
\|F_A\|_{L^2}^2+\|d_A\Phi\|_{L^2}^2+\|\mu(\Phi)-c\|_{L^2}^2,
\label{eiestr}
\end{equation}                                           
Here $d_A\Phi$ is the covariant derivative defined in (\ref{edefda}).
Of course, if $F$ is a vector space and $K$ acts linearly on $F$, 
then the Yang--Mills--Higgs functional defined above coincides with
the classical one given in \ref{edefymh}.
We then prove a formula generalising (\ref{eymh2}) which allows to
write the Yang--Mills--Higgs functional as
$$\YMH_c(A,\Phi)=\|\Lambda F_A+\mu(\Phi)-c\|_{L^2}^2+
2\|\ov{\partial}_A\Phi\|_{L^2}^2+4\|F_A^{0,2}\|_{L^2}^2+T,$$
where $T$ is a constant depending only on the topology of $E$, $F$ and
the section $\Phi$ (and which consequently is invariant under deformations
of $A$ and $\Phi$). This formula will be used to obtain bounds on the
$L^2$ norms of $F_A$ and $d_A\Phi$ for $(A,\Phi)$ solving equations
(\ref{eiequs}) and for $\Phi$ satisfying some fixed homological constraints.

\subsection{Chapter 2: Hitchin--Kobayashi correspondence}
In this chapter we prove a Hitchin--Kobayashi correspondence for the second
equation in (\ref{eiequs}) when $F$ is Kaehler. This correspondence
generalises the results explained in \ref{eihiko}.

When $F$ is Kaehler, the action of $K$ on $F$ extends to a unique holomorphic
action of the complexification $G=K^{\CC}$ of $K$. This allows to extend
the action of the gauge group $\GGG_K=\Gamma(E\times_{\Ad}K)$ on 
$\AAA\times\SSS$ to an action of the complex gauge group
$\GGG_G=\Gamma(E\times_{\Ad}G)$. In this situation we may ask ourselves
the same question as in \ref{eihiko}: 
which orbits in $\AAA\times\SSS$ of the action of $\GGG_G$ contain solutions
to the equation
\begin{equation}
\Lambda F_{A}+\mu(\Phi)=c?
\label{epregunta}
\end{equation}
And, how many $\GGG_K$ orbits of solutions to (\ref{epregunta}) can
contain at most a $\GGG_G$ orbit in $\AAA^{1,1}\times\SSS$? 

The main result of this chapter is a theorem which partially answers 
these questions. We define the notions of simple and $c$-stable pair, 
and we prove
that a simple pair $(A',\Phi')$ is $c$-stable if and only if there exists
$g\in\GGG_G$ such that $(A,\Phi)=g(A',\Phi')$ satisfies
(\ref{epregunta}). Furthermore, this $g$
is unique up to the action of $\GGG_K$ on $\GGG_G$ on the left.
In this chapter we do not ask the pair $(A,\Phi)$ to satisfy the first
equation $\ov{\partial}_A\Phi=0$; but we restrict ourselves to connections
satisfying the third equation $F_A^{0,2}=0$, that is, to connections 
belonging to $\AAA^{1,1}$.
This is a technical condition which probably may be relaxed. On the other
hand, observe that if $X$ is a Riemann surface then 
$\Omega^{1,1}(X)=\Omega^2(X)$, so $\AAA^{1,1}=\AAA$.

To prove the result of this chapter we construct and study a certain 
functional on $\AAA^{1,1}\times\SSS\times\GGG_G$ whose critical points
are exactly the points $(A,\Phi,g)$ which satisfy
$\Lambda F_{g(A)}+\mu(g(\Phi))=c$. The construction of this functional
is rather general, and we call it the integral of the moment map.

We finish the chapter with three examples of the correspondence.
In the first one we take $F$ to be a Hermitian vector space with
a linear unitary action of $K$, and we obtain Banfield's theorem.
In the second example we take $F=\PP(W)$, where $W$ is a Hermitian
vector space. Finally, in the third example we study the case
of $F$ being a Grassmannian or, more generaly, a flag manifold
(this case includes, of course, the one studied in the second
example).

\subsection{Chapter 3: The moduli space}
In this chapter we make the first steps towards a definition 
of invariants of the symplectic manifold $F$ and the Hamiltonian
action of $K$, by constructing certain spaces of solutions to
a perturbation of equations (\ref{eiequs}). From now on
we will suppose that $X$ is a Riemann surface and that $F$ is compact. 
Later we will make more assumptions on our data. 

Let $\GGG$ be the gauge group of $E$. In order to use certain
results such as the implicit function theorem, we extend our
configuration space by completing them with respect to some
$L^p_k$ Sobolev norms. 

Let $\sigma=(\sigma_1,\sigma_2)$ be a $\GGG$ invariant element of 
$\Hom^{0,1}(\pi_F^*TX,T\FFF_v)\oplus\Omega^0(E\times_{\Ad}\klie)$,
where $\pi_F:\FFF\to X$ is the projection.
We consider the following perturbed equations
\begin{equation}
\left\{\begin{array}{l}
\ov{\partial}_A\Phi=\sigma_1,\\
\Lambda F_A+\mu(\Phi)=\sigma_2+c,\end{array}\right.
\label{esiequs}
\end{equation}
where $A$ and $\Phi$ lie in the chosen Sobolev completions of $\AAA$
and $\SSS$, and we call the pairs $(A,\Phi)$ satisfying them
{\bf $\sigma$-twisted holomorphic curves over $X$}
({\bf $\sigma$-THCs} for short). (Recall that, since $X$ is a Riemann
surface, the integrability condition $F_A^{0,2}=0$ is always
satisfied.)

To any homology class $B\in H_2(F_K)$ (where $F_K=EK\times_K F$ is
the Borel construction) we associate a certain $\GGG$-invariant set
$\bM_{\sigma}(B,c)=\bM_{\sigma}^{F,K}(B,c)$ of solutions to 
equations (\ref{esiequs})
and we define the {\bf moduli space of $\sigma$-THCs}
(resp. the {\bf extended moduli space of $\sigma$-THCs}) to be
$\MMM_{\sigma}(B,c)=\MMM_{\sigma}^{F,K}(B,c)=\bM_{\sigma}(B,c)/\GGG$
(resp. $\NNN_{\sigma}(B,c)=\NNN_{\sigma}^{F,K}(B,c)=
\bM_{\sigma}(B,c)/\GGG_0$, where
$\GGG_0=\{g\in\GGG\mid g(x_0)=\id\}$, for a fixed $x_0\in X$).

Let us suppose that $K=S^1$ and that the action on $F$ is semi free
(this means that the action on the complementary of the fixed point
set $F^{S^1}$ is free).
We prove that there exists a discrete subset $C_0\subset\imag\RR$
such that if $c\in\imag\RR\setminus C_0$,
then for a generic perturbation $\sigma\in\Sigma_c(E)$
(where $\Sigma_c(E)$ is a non--empty set of perturbations depending on $c$)
the moduli space $\MMM_{\sigma}(B,c)$ is a smooth
manifold, and we compute its dimension (the point of restricting
to $c\in\imag\RR\setminus C_0$ is that then the isotropy subgroup 
in $\GGG$ of any solution to (\ref{esiequs}) is trivial). 
We also prove that for two different choices of generic perturbation 
$\sigma$ the resulting moduli spaces are cobordant. A similar
result is proved for the extended moduli space. Finally,
we also have for generic perturbation that the projection
$\NNN_{\sigma}(B,c)\to\MMM_{\sigma}(B,c)$
is a $S^1$ principal bundle.

\subsection{Chapter 4: Compactification of the moduli}

In this chapter we prove two basic results. The first one is a regularity
result, which says that the solutions to equations (\ref{eiequs})
are (gauge equivalent to) smooth pairs, and the second one is a theorem
which gives a compactification of the moduli space
$\MMM_{\sigma}(B,c)$ (and of $\NNN_{\sigma}(B,c)$).
In contrast with the preceeding chapter and with the next ones, in
this chapter $K$ can be any compact connected Lie group, and its
action on $F$ is only assumed to be smooth. On the other hand,
$X$ is, as always, a Riemann surface and $F$ is compact.

In the preceeding chapter we gave a definition of the set 
$\bM_{\sigma}(B,c)$ as a subset of a Sobolev completion of 
$\AAA\times\SSS$. A priori it is not clear whether the elements of
$\bM_{\sigma}(B,c)$ are smooth or whether the moduli space
$\MMM_{\sigma}(B,c)$ depends on the chosen Sobolev completions.
In this section we clarify the situation, proving that all the elements
in the moduli are smooth. More concretely we prove that
if $[A,\Phi]\in \MMM_{\sigma}(B,c)$, then there exists
a transformation $g\in\GGG$ such that $g(A)$ and $g(\Phi)$
are smooth. With this we see that the moduli space, as a set, is
intrinsic: it does not depend on the Sobolev norms. In fact, the
structure of the moduli as a smooth differential manifold is also
unique, since it can be given locally in terms of Kuranishi models
and these, by elliptic regularity, are independent of the Sobolev
norms.

Given a sequence $\{(A_k,\Phi_k)\}\subset\MMM_{\sigma}(B,c)$,
it may perfectly happen that there is no convergent subsequence
in the $C^0$ topology. This is similar to what happens in the theory
of pseudoholomorphic curves. The reason for this non compactness is
a phaenomenon called bubbling, which is a consequence of the impossibility
of finding bounds $\|d_A\Phi\|_{C^0}\leq C$ for
$(A,\Phi)\in\MMM_{\sigma}(B,c)$, where $C$ only depends
on $B$ and $c$.
In this section we define the notion of {\bf cusp $\sigma$-THC} and
we prove that for any sequence 
$\{(A_k,\Phi_k)\}\subset\MMM_{\sigma}(B,c)$ there exists a subsequence
converging (in a suitable sense) to a cusp $\sigma$-THC. With this
we obtain a compactification of $\MMM_{\sigma}(B,c)$ which
will allow to define invariants of $F$. This compactification generalises
Gromov's theorem for pseudoholomorphic curves (see \cite{AuLa}).
In fact, the notion of cusp $\sigma$-THC is also a generalisation of 
Gromov's cusp curves.

The next step after compactifying the moduli
$\MMM_{\sigma}(B,c)$ is to study to what extent the compactification
has a smooth structure. In order to define invariants, we would like
to have a fundamental class in the homology of the compactification.
An ideal situation would be that in which the compactification
admited a natural structure of smooth oriented manifold. In this
case we would indeed have a fundamental class. Unfortunately, this
will not happen in general. But we can express the compactification
of the moduli as the union of the moduli plus a countable family
of smooth manifolds. Finally, if the dimensions of these extra manifolds
are lower than that of the moduli minus one, then we will be able
to rigorously define invariants.

\subsection{Chapter 5: The choice of the complex structure}
So far we have not put any restriction on the complex structure
$I_F$ (appart from chapter \ref{hiko}, where we assumed that
$F$ is Kaehler). We only asked $I_F$ to be $K$-invariant.
In this chapter we assume that $K=S^1$ and we prove that
for a generic $S^1$-invariant complex structure $I_F$ the moduli
of simple holomorphic curves is the union of a countable family 
of smooth manifolds, and we compute its dimension. 
(Note that if we did not ask our complex structure to be $S^1$-invariant
then a stronger result could be proved: it is well known that for
a generic complex structure the moduli of simple holomorphic curves
is a smooth manifold; see e.g. \cite{McDS1}.) 
This result will be used when constructing the moduli of cusp $\sigma$-THCs.

\subsection{Chapter 6: The invariants}
In this chapter we assume that $K=S^1$ and that the action on $F$
is semi free. Using the results of the preceeding chapters we define,
under certain conditions, invariants of the symplectic manifold $F$
and the action of $S^1$. 

We define two invariants: the invariant $\Phi$, using the extended moduli
$\NNN$, and the invariant $\oPhi$, using the moduli $\MMM$.
The idea used to define them is very similar to that of Gromov-Witten
invariants. Fix a class $B\in H_2(F_{S^1})$, an element 
$c\in\imag\RR\setminus C_0$ and a generic perturbation 
$\sigma\in\Sigma_c(E)$. Let $\NNN=\NNN_{\sigma}(B,c)$
and $\MMM=\MMM_{\sigma}(B,c)$.
We have maps $\mu_i:H^*_{S^1}(F)\to H^*(\NNN)$,
$\omu_i:H^*_{S^1}(F)\to H^*(\MMM)$, $\nu:H^*(\AAA/\GGG_0)\to H^*(\NNN)$
and $\onu:H^*(\AAA/\GGG_0)\to H^*(\MMM)$. Formaly, the invariant 
$\Phi=\Phi^{X,F}_{B,c}$ (resp. $\oPhi=\oPhi^{X,F}_{B,c}$) is obtained 
by sending classes from $H^*_{S^1}(F)$ and $H^*(\AAA/\GGG_0)$
to $H^*(\NNN)$ (resp. $H^*(\MMM)$) using the maps
$\mu_i$ and $\nu$ (resp. $\omu_i$ and $\onu$), multiplying them
and then pairing the result with the fundamental class
$[\NNN]\in H_{\top}(\NNN)$ (resp. $[\MMM]\in H_{\top}(\MMM)$).
As we said before, in general we cannot prove the existence of
the fundamental classes used above. In this chapter we give a
rigorous definition of the invariants in certain conditions,
bypassing the question of fundamental classes.

Finally, we give an example of a nonzero invariant.

\section{Some questions}
To finish this introduction we list some problems which we would like to
study as a continuation of this thesis.
\begin{itemize}
\item {\bf Computation of the invariants.}
When $F$ is Kaehler, the Hitchin--Kobayashi correspondence proved in this
thesis allows to describe the moduli of THCs. This might be used to
make computations of the invariants in the Kaehler case. Another strategy 
would be to search relations among the invariants (using gluing, as in 
Gromov-Witten theory; see below) which might simplify the computations.

\item{\bf Improve the definition.}
The conditions which were needed to give a rigorous definition of the 
invariants are rather technical and very much restrictive. It would be
nice to get rid of them. To do that, the works
\cite{FuOn, LiTi, Ru2, Sie} (in which a rigorous definition of
Gromov-Witten invariants for any compact symplectic manifold is given)
should be an important source of inspiration.

\item{\bf Equivariant quantum cohomology.}
This would consist of codifying the Hamiltonian Gromov-Witten invariantrs
in a deformation of the ring structure of the equivariant cohomology
of $F$, exactly as is done with Gromov-Witten invariants. A central
question in this problem would be to prove associativity of the resulting
product. We expect that this should give non-trivial relation among
the invariants (as happens with Gromov-Witten invariants).

\item{\bf Gluing.} An interesting question is the following: given
two THCs, one with base $X_1$ and the other with base $X_2$,
how can we obtain a THC with base the connected sum $X_1\sharp X_2$?
A similar question arises in Donaldson theory and also in Gromov-Witten
theory, and the techniques needed in our situation should probably be the
same ones that apear in these theories. This is in our opinion the most
interesting question to be studied after this thesis. Very much likely,
a good understanding of a gluing construction of THCs will be crucial
in developping the latter two questions.

\item{\bf Interpretation of the invariants.}
We would like to find a description of the invariants in terms of enumerative
geometry, just as is done with Gromov-Witten invariants.

\item{\bf Wall crossing.} In this thesis we define invariants using
the moduli spaces $\NNN_{\sigma}(B,c)$ and 
$\MMM_{\sigma}(B,c)$ when $c\in\imag\RR\setminus C_0$. It would
be interesting to study the relation between the invariants
obtained when $c$ belongs to different connected components of
$\imag\RR\setminus C_0$. Presumably we could prove a result similar
to Thaddeus theorem in \cite{Th}.
\end{itemize}

\section{Notations}
All the manifolds, bundles and morphisms in this thesis will be smooth
unless otherwise stated. All vector spaces will be finite dimensional.
We will use the following convention. Finite dimensional manifolds
and vector bundles will be denoted using roman fonts: $X$, $F$, $V$.
General fibre bundles will be denoted using calligraphic fonts:
$\FFF$. Finally, infinite dimensional manifolds will be denoted 
using {\it curly} fonts: $\AAA$, $\GGG$, $\SSS$.
For any real vector space $V$, we will denote $\la,\ra_V$ the natural
pairing $V\times V^*\to\RR$. 

These are the most often used symbols in this thesis.
               
\begin{itemize}
\item $\NN$, $\ZZ$, $\QQ$, $\RR$ and $\CC$ denote as usual the
sets of natural, integer, rational, real and complex numbers.
We denote by $\imag$ the square root $\sqrt{-1}$.
\item $K$ is a compact connected real Lie group
(when defining the invariants we will assume that $K=S^1$);
\item $G$ is the complexification of $K$.
\item $\klie$ and $\glie$ are the Lie algebras of $K$ and $G$.
\item $X$ is a Kaehler manifold with symplectic structure $\omega_X$
and with complex structure $I_X$ (when defining the invariants it will 
be supposed to be a Riemann surface);
\item $F$ is a symplectic manifold with symplectic structure $I_F$
and with a Hamiltonian symplectic action of $K$;
we denote $\mu:F\to\klie^*$ its moment map; we take a
$K$-invariant complex structure $I_F$ compatible with $\omega_F$;
when studying the Hitchin-Kobayashi correspondence
$F$ will be Kaehler and when defining the invariants
it will be compact.
\item $EK\to BK$ is the universal principal $K$-bundle (the action of $K$ 
on any $K$ principal bundle is by definition, as usual, on the right);
\item $E\to X$ is a $K$ principal bundle; $E_G=E\times_{K}G$ be the
$G$ principal bundle associated to $E$ (we take the action of $K$
on $G$ given by left multiplication);
\item $\FFF=\FFF^E=E\times_K F$ is the associated bundle with fibre $F$.
\item $\GGG=\GGG^E=\Gamma(E\times_{\Ad}K)$ is the gauge group of $E$;
in the chapter on Hitchin--Kobayashi correspondence we will denote this
group by $\GGG_K$;                                         
\item $\GGG_G=\Gamma(E_G\times_{\Ad}G)=\Gamma(E\times_{\Ad}G)$
is the gauge group of $E_G$, and is the complexification of $\GGG_K$;
\item $\AAA=\AAA^E$ is the set of connections on $E$;
\item $\SSS=\SSS^E=\Gamma(\FFF)$ is the space of sections of $\FFF$;
\item $\MMM$ denotes moduli spaces in general (there will be
several of them appearing in the thesis);
\item if the manifold $M$ supports an action of $K$ then for any
$s\in\klie$ we denote $\fX^M_s$ the vector field on $M$ generated
by the infinitesimal action of $s$; when the manifold $M$ is clear
from the context, we just write $\fX_s$; finally,
$M^K$ is the set of fixed points.
\end{itemize}

\setcounter{chapter}{0}
\renewcommand{\thechapter}{\arabic{chapter}}


\chapter[Equations and Yang--Mills--Higgs functional]
{The equations and the Yang--Mills--Higgs functional}
\label{eqymh}

\section{The equations}
\label{theequations}

\subsectionr{}
Let $K$ be a compact connected real Lie group, and let $\klie=\Lie(K)$
be its Lie algebra. Let us take on $K$ a metric invariant by the adjoint 
action of $K$ (such metrics are called biinvariant). This metric allows
to identify $\klie\simeq \klie^*$ in a $K$-equivariant way.

\subsectionr{}
Let $F$ be a symplectic manifold with a symplectic left action of $K$.
Let $\omega_F$ be the symplectic structure of $F$. Let $I_F$ be
a $K$-invariant complex structure on $F$ compatible with $\omega_F$,
that is, such that $g_F(\cdot,\cdot):=\omega_F(\cdot,I_F\cdot)$
is a Riemannian metric on $F$. Such a complex structure always exists
(see lemma 5.49 in \cite{McDS2}).
             
\begin{definition}
\label{momentmap}
A {\bf moment map} for the action of $K$ on $F$ is a map 
$\mu:F\to\klie^*$ which satisfies the following two conditions:

(C1) for any $s\in\klie$, $d\mu(s)=\iota_{\fX_s}\omega_F$ (where
$\fX_s\in\Gamma(TF)$ is the vector field generated by the infinitesimal
action of $s$ on $F$) and

(C2) $\mu$ is equivariant with respect to the actions
of $K$ on $F$ and the coadjoint action on $\klie^*$. 
This means that for any $h\in K$,
any $x\in F$, and any $v\in\klie$,
$\la\mu(hx),\Ad(h)v\ra_{\klie}=\la\mu(x),v\ra_{\klie}.$
\end{definition}

In the sequel we will assume that there exists a moment map $\mu:F\to\klie^*$
for the action of $K$ on $F$. (Note that if a moment map exists, then it is 
unique up to addition of a central element in $\klie^*$.)

\subsectionr{}
Let $X$ be a connected Kaehler manifold of dimension $n$. 
Let us write $\omega_X$, $I_X$ and $g_X$
the symplectic form, the complex structure and the Kaehler metric of $X$.
They are related by $g_X(\cdot,\cdot)=\omega_X(\cdot,I_X\cdot)$.
In this chapter we will often write $\omega$ for $\omega_X$.
We will denote $\omega^{[k]}=\omega^k/k!$, and
we will use in the integrals of functions on $X$ the volume element 
$\omega^{[n]}$ (most of the times we will not write it).
Let $\Lambda:\Omega^*(X)\to\Omega^{*-2}(X)$ be the adjoint of
wedging with $\omega_X$. 

\subsectionr{}
Let $E\to X$ be a $K$-principal bundle on $E$. 
As usual, we take the action of $K$ on $E$ to be on the right.
Consider the associated bundle $\pi_F:\FFF=E\times_KF\to X$. 
Let $\AAA$ be the space of connections on $E$ and
let $\SSS$ be the space $\Gamma(\FFF)$ of sections of $\FFF$.

Let $T\FFF_v=\Ker d\pi_F$. A connection $A\in\AAA$ induces a projection
$\alpha:T\FFF\to T\FFF_v$. Let $\pi^{0,1}:\Omega^1(X)\to\Omega^{0,1}(X)$ 
denote the projection and let $\pi^{1,0}=1-\pi^{0,1}$. 
Let $\Phi\in\SSS$ be a section of $\FFF$. We define
the covariant derivative of $\Phi$ (with respect to $A$) to be
$$d_A\Phi=\alpha\circ d\Phi\in\Omega^1(\Phi^*T\FFF_v),$$
and the $\partial$ and $\ov{\partial}$ operators of $A$ acting on $\Phi$
to be 
$$\partial_A\Phi=\pi^{1,0}\alpha(d\Phi)\in\Omega^{1,0}(\Phi^*T\FFF_v)
\qquad\mbox{and}\qquad
\ov{\partial}_A\Phi=\pi^{0,1}\alpha(d\Phi)\in\Omega^{0,1}(\Phi^*T\FFF_v).$$
           
\subsectionr{}
Using the biinvariant metric on $\klie$ we get an equivariant
isomorphism $\klie\simeq\klie^*$ which extends to an isomorphisms of
vector bundles
\begin{equation}
E\times_{\Ad}\klie\simeq E\times_{\Ad}\klie^*
\label{isofibr}                        
\end{equation}
(observe that we denote with the same symbol $\Ad$ both the adjoint
representation on $\klie$ and the coadjoint representation on $\klie^*$).
Let finally $c\in\klie$ be a central element. The equations which 
we study in this thesis are
\begin{equation}
\left\{\begin{array}{l}
\ov{\partial}_A\Phi=0,\\
\Lambda F_A+\mu(\Phi)=c, \\
F_A^{0,2}=0,\end{array}\right.
\label{equs}
\end{equation}
where $A\in\AAA$ is a connection,
$F_A\in\Omega^2(E\times_{Ad}\klie)$ is the curvature of $A$, and
$\Phi\in\SSS$ is a section. Observe that
$\Lambda F_A\in\Omega^0(E\times_{Ad}\klie)$ and that
$\mu(\Phi)\in\Omega^0(E\times_{Ad}\klie^*)$, so to give a sense
to the second equation we need to use the isomorphism (\ref{isofibr}).
Hence, the equations (\ref{eiequs}) depend on the biinvariant
metric taken on $\klie$. In the third equation
$F_A^{0,2}$ refers to the piece in $\Omega^{0,2}(\gE)$ of the curvature
$$F_A\in\Omega^2(E\times_{Ad}\klie)\subset
\Omega^2(\gE)=\Omega^{2,0}(\gE)\oplus\Omega^{1,1}(\gE)
\oplus\Omega^{0,2}(\gE),$$
where $\gE=E\times_{\Ad}\glie$ and $\glie=\klie\otimes\CC$ is the
complexification of $\klie$. The condition $F_A^{0,2}=0$ is equivalent to
$F_A\in\Omega^{1,1}(\gE)$. We will write $\AAA^{1,1}$ for the set
of connections $A\in\AAA$ such that $F_A^{0,2}=0$. 
By the theorem of Newlander and Niremberg \cite{NewNi} these are the
connections which define integrable complex structures in the complex
vector bundles associated to $E$.

\subsectionr{}
\label{THCshol}
Take a connection $A\in\AAA$. This connection induces a splitting
$T\FFF\simeq T\FFF_v\oplus \pfx^*TX.$
Using this splitting, we define a metric $g(A)=g_F\oplus g_X$ and
a complex structure $I(A)=I_F\oplus I_X$ on $T\FFF$.           
Then $I(A)$ and $g(A)$ provide $\FFF$ with an almost Kaehler structure.

\begin{lemma}
Let $A\in\AAA$ be any connection on $E$.
A section $\Phi\in\SSS$ is holomorphic with respect
to $I(A)$ as a map from $X$ to $\FFF$ if and only if
$\overline{\partial}_A\Phi=0$.
\label{compara}
\end{lemma}
\begin{pf} This follows from the formula
$\ov{\partial}_A\Phi=(d\Phi+I(A)\circ d\Phi\circ I_X)/2$.
\end{pf}

\section{The Yang--Mills--Higgs functional}
\label{YMHs}

The biinvariant metric on $\klie$ induces a norm $\|\cdot\|:\klie\to\RR$.
Thanks to $K$-equivariance, we may combine $\|\cdot\|$ with the
volume form $\omega_X^{[n]}$ to obtain an $L^2$ norm
$\|\cdot\|_{L^2}:\Omega^0(E\times_{\Ad}\klie)\to\RR$ and
similarly on $\Omega^0(E\times_{\Ad}\klie^*)$. We can do the same
thing with the $K$-invariant metric $g_F$ on $F$, thus getting
$\|\cdot\|_{L^2}:\Omega^0(T\FFF_v)\to\RR$.

\begin{definition}
Fix a central element $c\in\klie$. The {\bf Yang--Mills--Higgs}
functional $\YMH_c:\AAA\times\SSS\to\RR$ is defined as
$$\YMH_c(A,\Phi)=\|F_A\|_{L^2}^2+\|d_A\Phi\|_{L^2}^2
+\|c-\mu(\Phi)\|_{L^2}^2,$$
where $\Phi\in\SSS$ is a section and $A\in\AAA$ a connection on $E$.
\end{definition}

We will say that two sections $\Phi_0,\Phi_1\in\SSS$ are
{\bf homotopic} iff there exists a map $H_{\Phi}:X\times[0,1]\to\FFF$
such that $\Phi_0=H_{\Phi}|_{X\times\{0\}}$, 
$\Phi_1=H_{\Phi}|_{X\times\{1\}}$
and such that, for any $t\in[0,1]$,
$H_{\Phi}|_{X\times\{t\}}$ is a section, that is,
$\pi_F\circ H_{\Phi}|_{X\times\{t\}}=\Id_X$. Such an homotopy
will be called a {\bf homotopy of sections}.
The relation of homotopy of sections is an equivalence relation. 
For any section $\Phi\in\SSS$, $[\Phi]$ will
denote the homotopy class of $\Phi$ as a section.

\begin{theorem}
\label{minimitza}
Fix a section $\Phi_0\in\SSS$.
The pairs $(A,\Phi)\in\AAA\times\SSS$ which minimize
the functional $\YMH_c$ among the pairs whose section
is homotopic to $\Phi_0$ are those which satisfy equations
\ref{equs}.
\label{idminim}
\end{theorem}

We recall that, just as the equations (\ref{equs}), the Yang--Mills--Higgs 
functional does depend on the biinvariant metric on $\klie$.
A proof of theorem \ref{idminim} will be given in the next section.

\subsection{A weak Kaehler identity}
\label{secnonlinka}

The following result will be used in the proof of theorem
\ref{idminim}.
       
\begin{prop}[Weak Kaehler identity]
For any section $\Phi\in\FFF$ and for any connection $A\in\AAA$, 
the following equality holds:
$$\int_X\la \Lambda F_A,\mu(\Phi)\ra_{\klie}=
\frac{1}{2}(\|\partial_A\Phi\|^2_{L^2}
-\|\overline{\partial}_A\Phi\|^2_{L^2})-C_{[\Phi]},$$
where the constant $C_{[\Phi]}$ depends only on the topological type of 
$E$ and on the homotopy class of sections of $\Phi$. 
\label{nonlinka}
\end{prop}

We give two proofs of this result. The first one works only
when $X$ is a Riemann surface, and follows from a direct computation.
The second proof uses the Chern-Weil map in equivariant cohomology,
and works for any Kaehler manifold $X$. Furthermore, the second
proof gives a geometrical interpretation of the constant
$C_{[\Phi]}$. This interpretation allows to give bounds on 
the $L^2$ norms of $F_A$ and $d_A\Phi$ when $(A,\Phi)$ solve
(\ref{equs}) in terms of homological data (see theorem \ref{boundL2}).

To motivate the name of weak Kaehler identity, consider the case $F=\CC^n$.
Take on $F$ a Hermitian metric $h$. Its imaginary part with reversed
sign gives a symplectic form $\omega_F$ compatible with the complex
structure. Let $K=U(n)$ act on $F$ respecting 
$h$ (and consequently $\omega_F$). In this situation
the moment map is $\mu=-\frac{\imag}{2}x\otimes x^*$, so that
$\mu(\Phi)=-\frac{\imag}{2}\Phi\otimes\Phi^*$. Using the Kaehler 
identities for unitary connections on hermitian bundles 
$\imag[\Lambda,\overline{\partial}_A] = \partial_A^*$ and
$-\imag[\Lambda,\partial_A] = \overline{\partial}_A^*$ we compute
\begin{align}
\int_X\la \Lambda F_A,\mu(\Phi)\ra_{\klie}
&= \int_X \la \Lambda F_A,-\frac{\imag}{2}\Phi\otimes\Phi^*\ra_{\klie} 
= \int_X \frac{1}{2}h(\Phi,\imag\Lambda F_A\Phi) \notag \\
&= \int_X \frac{1}{2}h(\Phi,\imag\Lambda(\overline{\partial}_A\partial_A+
\partial_A\overline{\partial}_A)\Phi)
= \int_X \frac{1}{2}h(\Phi,\partial_A^*\partial_A\Phi-
\overline{\partial}_A^*\overline{\partial}_A\Phi) \notag \\
&= \frac{1}{2}(\|\partial_A\Phi\|^2_{L^2}
-\|\overline{\partial}_A\Phi\|^2_{L^2}). \notag
\end{align}
(Observe that in this case the constant $C_{[\Phi]}$ is equal to zero, no
matter what the topological type of $E$ is.)

\subsection{Proof of proposition \ref{nonlinka} for $X$ a Riemann surface}
\label{Riemsurf}
In this section the manifold $X$ will be a compact Riemann surface 
with a fixed Kaehler structure. Consider first of all two
pairs $(A^1,\Phi^1)$ and $(A^2,\Phi^2)$ which are equal in
the complement of an open set $V\subset X$ and such that 
$[\Phi^1]=[\Phi^2]$. Suppose $V$ is
small enough so that there exists an open set $U\subset X$ containing
$\ov{V}$, a holomorphic chart
$\psi:U\to\CC\simeq\RR^2$ with $[0,1]\times[0,1]\subset\psi(U)$, 
and a trivialisation $E|_U\simeq U\times K$ in such a way that
$\psi(\overline{V})$ is contained in $(0,1)\times (0,1)$.
Write $x$ and $y$ the usual coordinates of $\RR^2$, so that
the complex structure $I\in\End(T\RR^2)$ sends 
$\partial/\partial x$ to $\partial/\partial y$
and $\partial/\partial y$ to $-\partial/\partial x$.

We will write for convenience $S=[0,1]\times[0,1]$.
On $S$ we will consider either the volume form 
$(\psi^{-1})^*\omega$, which we will also call $\omega$,
or $dx\wedge dy$. Define the function $f:S\to\RR$ by 
$\omega=f(dx\wedge dy)$. Take on $S$ the metric coming from the one 
on $F$. From $\omega(\partial/\partial x,\partial/\partial y)=f$,
$|\partial/\partial x|=|\partial/\partial y|$ and the fact that
$\partial/\partial x$ is orthogonal to $\partial/\partial y$
(recall that $I_F$ is an isometry) we deduce that
\begin{equation}
|\partial/\partial x|=|\partial/\partial y|=\sqrt{f}.
\label{normpar}
\end{equation}

Write, for $i=1,2$, $d_{A^i}=d+A^i_x dx+A^i_y dy$ in the chosen
trivialisation of $E|_U$, where $A^i_x$ and $A^i_y$ take values
in $\klie$. The corresponding curvatures in our trivialisation are
$$F_{A^i}=\left(\frac{\partial A^i_y}{\partial x}
-\frac{\partial A^i_x}{\partial y}+[A^i_x,A^i_y]\right)dx\wedge dy.$$
We will consider the restriction of the sections on $U$ as 
maps $\Phi^i:U\to F$. Now we compute 
$$\int_{S}\la\Lambda F_{A^i},\mu(\Phi^i)\ra_{\klie}\omega=
\int_{S}\left\la 
\frac{\partial A^i_y}{\partial x}
-\frac{\partial A^i_x}{\partial y}+[A^i_x,A^i_y],\mu(\Phi^i)
\right\ra_{\klie}dx\wedge dy.$$
Integrating by parts we see that this is equal to 
\begin{align}
&B_i+
\int_{S}(-\la\la d\mu(\Phi),\partial\Phi^i/\partial x\ra_{TF},
A^i_y\ra_{\klie}+
\la\la d\mu(\Phi),\partial\Phi^i/\partial y\ra_{TF},A^i_x\ra_{\klie}
\notag \\
&+\la[A^i_x,A^i_y],\mu(\Phi^i)\ra_{\klie})dx\wedge dy,\notag
\end{align}
where $B_i$ is a boundary term. 
We now use the properties of the moment map to deduce that this equals
\begin{align*}
&B_i+\int_{S}\left(-\omega_F(\fX_{A^i_y},\partial\Phi^i/\partial x)
+\omega_F(\fX_{A^i_x},\partial\Phi^i/\partial y)-
\omega_F(\fX_{A^i_y},\fX_{A^i_x})\right)dx\wedge dy\\
&=B_i+
\int_{S}\left(\omega_F(\partial\Phi^i/\partial x+
\fX_{A^i_x},\partial\Phi^i/\partial y+\fX_{A^i_y})
-\omega_F(\partial\Phi^i/\partial x,\partial\Phi^i/\partial y)
\right)dx\wedge dy
\end{align*}
(recall that $\fX_s$ denotes the vector field given by the action
of $s\in\klie$ on $F$).
By definition $\partial\Phi^i/\partial x+\fX_{A^i_x}$ 
(resp. $\partial\Phi^i/\partial y+\fX_{A^i_y}$)
is $d_{A^i}\Phi^i(\partial/\partial x)$
(resp. $d_{A^i}\Phi^i(\partial/\partial y)$). 
Now, using lemma \ref{trucu} and formula (\ref{normpar}) we deduce
$$\omega_F\left(d_{A^i}\Phi^i\left(\partial/\partial x\right),
d_{A^i}\Phi^i\left(\partial/\partial y\right)\right)=
f\frac{1}{2}
(|\partial_{A^i}\Phi^i|^2-|\overline{\partial}_{A^i}\Phi^i|^2).$$
Using this formula we get
\begin{align}
\int_{S}\la\Lambda F_{A^i},\mu(\Phi^i)\ra_{\klie}
&=B_i+\int_{S}
\frac{1}{2}\left(|\partial_{A^i}\Phi^i|^2-|\overline{\partial}_{A^i}\Phi^i|^2
\right)\omega \notag\\
&-\omega_F(\partial\Phi^i/\partial x,\partial\Phi^i/\partial y)dx\wedge dy.
\notag\end{align}

Since $(A^i,\Phi^i)$ are equal outside $V$ and $\overline{V}\subset
(0,1)\times(0,1)$, the two boundary terms $B_1$ and $B_2$
are the same. Moreover, since $[\Phi^1]=[\Phi^2]$ and $\omega_F$ is closed,
the integral
$\int_{S}\omega_F(\partial\Phi^i/\partial x,\partial\Phi^i/\partial y)
dx\wedge dy$
has the same value for $i=1,2$. Hence we obtain
\begin{align}
&\int_X
\la\Lambda F_{A^1},\mu(\Phi^1)\ra_{\klie}
-\la\Lambda F_{A^2},\mu(\Phi^2)\ra_{\klie}
=\int_{S}
\la\Lambda F_{A^1},\mu(\Phi^1)\ra_{\klie}
-\la\Lambda F_{A^2},\mu(\Phi^2)\ra_{\klie}\notag\\
&=\int_{S}
\frac{1}{2}\left(|\partial_{A^1}\Phi^1|^2-|\overline{\partial}_{A^1}\Phi^1|^2
\right)-\frac{1}{2}\left(|\partial_{A^2}\Phi^2|^2-
|\overline{\partial}_{A^2}\Phi^2|^2\right)\notag\\
&=\int_X\frac{1}{2}\left(
|\partial_{A^1}\Phi^1|^2-|\overline{\partial}_{A^1}\Phi^1|^2\right)
-\frac{1}{2}\left(
|\partial_{A^2}\Phi^2|^2-|\overline{\partial}_{A^2}\Phi^2|^2\right)\notag\\
&=\frac{1}{2}\left(
\|\partial_{A^1}\Phi^1\|^2_{L^2}-\|\overline{\partial}_{A^1}\Phi^1\|^2_{L^2}
\right)
-\frac{1}{2}\left(\|\partial_{A^2}\Phi^2\|^2_{L^2}-
\|\overline{\partial}_{A^2}\Phi^2\|^2_{L^2}\right).\notag
\end{align}
(In all these integrals we omit the volume form, which is $\omega$.)
To finish the proof, observe that given two
pairs $(A,\Phi)$ and $(A',\Phi')$, where $[\Phi]=[\Phi']$,
one can always find a sequence of pairs $(A^i,\Phi^i)$, $i=1,\dots,k$,
such that the homotopy classes $[\Phi^i]$ are all equal to $[\Phi]$
and such that any two consecutive pairs in 
$$\{(A,\Phi),(A^1,\Phi^1),\dots,(A^k,\Phi^k),(A',\Phi')\}$$
coincide outside a small enough set $V\subset X$ so
that we can apply the preceeding reasoning. This implies that
$$\int_X
\la\Lambda F_{A},\mu(\Phi)\ra_{\klie}
-\la\Lambda F_{A'},\mu(\Phi')\ra_{\klie}
=\frac{1}{2}\left(
\|\partial_{A}\Phi\|^2_{L^2}-\|\overline{\partial}_{A}\Phi\|^2_{L^2}\right)
-\frac{1}{2}\left(\|\partial_{A'}\Phi'\|^2_{L^2}-
\|\overline{\partial}_{A'}\Phi'\|^2_{L^2}\right),$$
which is what we wanted to prove.                                                                      
                
\subsection{Proof of proposition \ref{nonlinka} for any Kaehler manifold $X$}
Instead of directly proving proposition \ref{nonlinka} we will prove a 
slightly more general result. In the course of the proof we will find a
geometrical interpretation of the constant $C_{[\Phi]}$.

The symplectic form $\omega_F$ gives an element
of $\Omega^0(\Lambda^2 (T\FFF_v)^*)$, since the action of $K$
leaves $\omega_F$ invariant. On the other hand, the connection $A$ on $E$ 
induces a projection $$\alpha:T\FFF\exh T\FFF_v$$ 
onto the subbundle of vertical tangent vectors. From this we obtain a map 
$\alpha^*:\Lambda^2 (T\FFF_v)^*\to \Lambda^2 T^*\FFF$, and we 
set $\tilde{\omega}_F^A=\alpha^*(\omega_F)\in\Omega^0(\Lambda^2 T^*\FFF)
=\Omega^2(\FFF)$. This 2-form is not in general closed. Consider the 2-form 
$\omega_F^A=\tilde{\omega}_F^A-\la \pi_F^*F_A,\mu\ra_{\klie}$. 
                                 
\begin{prop}
The 2-form $\omega_F^A\in\Omega^2(\FFF)$ is closed, and the cohomology class 
it represents is independent of the connection $A$.
\label{nonlinka2}
\end{prop}

Let us show that proposition \ref{nonlinka2} implies proposition 
\ref{nonlinka}.
(In fact proposition \ref{nonlinka2} is slightly stronger
than \ref{nonlinka}.) We will use the following elementary lemma.

\begin{lemma}
Let $V$ and $W$ be two Euclidean vector spaces with scalar products
$\la,\ra_V$ and $\la,\ra_W$. Suppose that there are complex structures
$I_V\in\End(V)$, $I_W\in\End(W)$ and symplectic forms 
$\omega_V\in\Lambda^2V^*$, $\omega_W\in\Lambda^2W^*$ which satisfy
the following: $\la\cdot,\cdot\ra_V=\omega_V(\cdot,I_V\cdot)$
and $\la\cdot,\cdot\ra_W=\omega_W(\cdot,I_W\cdot)$
(in other words, $V$ and $W$ are Kaehler vector spaces).
Take a linear map $f:V\to W$ and let
$f^{1,0}$ (resp. $f^{0,1}$) be $(f+I_W\circ f\circ I_V)/2$
(resp. $(f-I_W\circ f\circ I_V)/2$). 
Let $2n=\dim_{\RR}V$. Then
$$f^*\omega_W\wedge\omega_V^{[n-1]}=
\frac{1}{2}(|f^{1,0}|^2-|f^{0,1}|^2)\omega_V^{[n]},$$
where, for any $g\in\Hom(V,W)$, $|g|^2=\Tr g^*g$ and $\omega_V^{[k]}=
\omega_V^k/k!$\ .
\label{trucu}
\end{lemma}

\begin{remark}
Note that under the isomorphism
$$V^*\otimes_{\RR}W\simeq
(V^*\otimes_{\RR}\CC)^{1,0}\otimes_{\CC}W\oplus
(V^*\otimes_{\RR}\CC)^{0,1}\otimes_{\CC}W$$
the element $f\in V^*\otimes_{\RR}W$ corresponds precisely to
$f^{1,0}+f^{0,1}$.
\label{trucu2}
\end{remark}

Now assume that proposition \ref{nonlinka2} is true.
Using lemma \ref{trucu} we have 
$$\int_X\Phi^*\tilde{\omega}_F^A\wedge\omega^{[n-1]}=
\frac{1}{2}(\|\partial_A\Phi\|^2_{L^2}
-\|\overline{\partial}_A\Phi\|^2_{L^2})$$
for any section $\Phi:X\to\FFF$. To apply the lemma we set, for
any $x\in X$, $V=T_xX$ and $W=T_{\Phi(x)}\FFF_v$ with the 
induced Kaehler structures, and $f=d_A\Phi(x)$. With these identifications
$f^{1,0}=\partial_A\Phi(x)$ and $f^{0,1}=\overline{\partial}_A\Phi(x)$
(see remark \ref{trucu2}). As a consequence, 
\begin{align}
C_{[\Phi]}&=\frac{1}{2}(\|\partial_A\Phi\|^2_{L^2}
-\|\overline{\partial}_A\Phi\|^2_{L^2})-
\int_X \la \Lambda F_A,\mu(\Phi)\ra \notag \\
&=\int_X(\Phi^*\tilde{\omega}_F^A
-\Phi^*\la \pi_F^*\Lambda F_A,\mu(\Phi)\ra)\wedge\omega^{[n-1]}=
\int_X\Phi^*\omega_F^A\wedge\omega^{[n-1]},\notag 
\end{align}
which by proposition \ref{nonlinka2} depends only on $[\Phi]$.
This proves proposition \ref{nonlinka}.

\subsubsection{The Cartan complex}
We are now going to prove proposition \ref{nonlinka2}.
We will use some results from \cite{BeGeV}, especially
from chapter 7. Define the graded algebra 
$$\Omega_K(F)=(\CC[\klie]\otimes \Omega(F))^K,$$
where $\Omega(F)$ is the algebra of differential forms on $F$, 
and assign to $P\otimes\alpha\in\Omega_K(F)$ the degree
$\deg(P\otimes\alpha)=2\deg(P)+\deg(\alpha)$. 
As usual $(\cdot)^K$ denotes the $K$ invariant elements under
the action of $K$. (The action of $K$ on
$\CC[\klie]\otimes \Omega(F)$ is by pullback
both on $\CC[\klie]$ and on $\Omega(F)$.)
If $\eta\in\Omega_K(F)$, let
$$d_{\klie}(\eta)(s)=d(\eta(s))+\iota(\fX_s)\eta(s),$$
where $s$ is any element in $\klie$ and $\fX_s$ denotes
the vector field generated by $s$ under the action of $\klie$.
Note that in \cite{BeGeV} the field assigned to $s\in\klie$
is $X_s:=-\fX_s$, so that $[X_s,X_{s'}]=X_{[s,s']}$
(recall that the action of $K$ on $F$ is on the left). This explains
the different sign in our definition.

The map $d_{\klie}$ sends $\Omega^{*}_K(F)$ to
$\Omega^{*+1}_K(F)$. One proves that $d_{\klie}^2=0$, 
so that $(\Omega_K(F),d_{\klie})$ is a complex. 
The complex $(\Omega_K(F),d_{\klie})$ is called the Cartan 
complex of $F$ and the action of $K$ on $F$.

According to our hypothesis the symplectic form $\omega_F\in\Omega^2(F)$
of $F$ is invariant under the action of $K$. However, it is not a  
closed form in $\Omega_K(F)$. This can be remedied by substracting
to it the moment map: the form $$\overline{\omega}_F=\omega_F-\mu$$
is equivariantly closed, that is, $d_{\klie}\overline{\omega}_F=0$.

\subsubsection{The Chern-Weil homomorphism}
We will now define a map of differential graded algebras
$$\phi_A:(\Omega_K(F),d_{\klie})\to(\Omega(\FFF),d)$$ 
with the help of a connection $A\in\AAA$. 
We will follow closely section 7.6 in \cite{BeGeV}.
First we give some definitions.

\begin{definition}                 
A {\bf horizontal differential form} on $\FFF$
is a differential form $\alpha\in\Omega(\FFF)$ such that 
$\iota(\fX)\alpha=0$ for all vertical vector fields $\fX$.
For any vector space $W$ with an action
$\rho:K\to \GL(W)$ let $\Omega(E,W)_{\hor}$ be the set  
of horizontal differential forms of the trivial vector bundle 
$E\times W\to E$. A {\bf basic differential form} on $E$ taking 
values in a linear representation $(W,\rho)$ of $K$ is an invariant 
form $\alpha\in\Omega(E,W)_{\hor}$. We will denote
$\Omega(E,W)_{\bas}$ the set of basic forms.
\end{definition}

This is proposition 1.9 in \cite{BeGeV}:

\begin{prop}
There is a natural isomorphism between $\Omega(E\times_{\rho}W)$
and $\Omega(E,W)_{\bas}$. This isomorphism sends 
$\alpha\in \Omega^q(E,W)_{\bas}$ to the form
$\alpha_X\in\Omega^q(E\times_{\rho}W)$ defined as follows
$$\alpha_X(d\pi X_1,\dots,d\pi X_q)(x)=[e,\alpha(X_1,\dots,X_q)(e)],$$
where $e\in \pi^{-1}(x)$, $x\in X$ and $X_j\in T_eE$, and where
$[a,b]$ denotes the element in $E\times_{\rho}W$ represented by
$(a,b)\in E\times W$.
\label{formequi}
\end{prop}

\begin{remark}
In particular, if $W=\RR$ with the trivial action of $K$ the preceeding
proposition says that $\Omega(X)\simeq\Omega(E)_{\bas}$.
If we apply this isomorphism to the principal $K$ bundle
$E\times F\to \FFF$ then we obtain 
$\Omega(\FFF)\simeq\Omega(E\times F)_{\bas}$.
\label{formquot}
\end{remark}

Fix a connection $A\in\AAA$ and call $\alpha:TE\exh TE_v$ the 
associated projection. Since the diagram
$$\xymatrix{E\times F\ar[r]^{\pi_E}\ar[d]\ar@{}[dr]|{S} 
& E\ar[d]^{\pi} \\ 
\FFF\ar[r]_{\pi_F} & X}$$
is cartesian, in particular $E\times F\simeq \pi_F^*E$ as 
principal $K$ bundles on $\FFF$. So we can pull back the connection
$A$ to a connection $\pi_F^*\alpha$ on the principal $K$ bundle
$E\times F\to\FFF$. This new connection gives a map
$\pi_F^*\alpha:T(E\times F)\exh T(E\times F)_v$ 
which induces a projection $h:\Omega(E\times F)\exh \Omega(E\times F)_{\hor}$ 
to the space of horizontal forms on $E\times F\to\FFF$.
Write $\Omega_A\in \Omega^2(E,\klie)$ the curvature form of
the connection $A$.

Now take an element $\alpha=f\otimes\beta\in\CC[\klie]\otimes\Omega(F)$
and define $\alpha(\Omega)\in\Omega(E)\otimes\Omega(F)$ as
$\alpha(\Omega_A)=f(\Omega_A)\otimes\beta$. Extending linearly we
obtain a map from $\CC[\klie]\otimes\Omega(F)$ to 
$\Omega(E)\otimes\Omega(F)$. 

\begin{definition}
The map $\phi_A:\CC[\klie]\otimes\Omega(F)\to
\Omega(E\times F)_{\hor}$ defined by
$$\phi_A(\alpha)=h(\alpha(\Omega_A))$$
is called the {\bf Chern-Weil homomorphism}.
\end{definition}

The following is theorem 7.34 in \cite{BeGeV}.

\begin{theorem}
\label{phimap}
The restriction of $\phi_A$ to the invariant forms
$(\CC[\klie]\otimes\Omega(F))^K$ has image contained in 
$\Omega(E\times F)_{\bas}$ and induces a homomorphism of 
differential graded algebras
$$\phi_A:(\Omega_K(F),d_{\klie})=
((\CC[\klie]\otimes\Omega(F))^K,d_{\klie})\to(\Omega(E\times F)_{\bas},d)
\simeq (\Omega(\FFF),d).$$
\end{theorem}

\begin{remark}
The isomorphism in the right hand side is given by remark \ref{formquot}. 
Using this isomorphism
we will regard $\phi_A$ as taking values in $(\Omega(\FFF),d)$.
\end{remark}

\begin{lemma}
The map induced by $\phi_A$ from the cohomology of
$((\CC[\klie]\otimes\Omega(F))^K,d_{\klie})$ to that of
$(\Omega(\FFF),d)$ does not depend on the connection $A$.
\label{cohonodepen}
\end{lemma}
\begin{pf}
Take two connections $A_0,A_1\in\AAA$ and any closed form
$\eta\in\Omega_K(F)$. Let $E_I\to X\times I$ be the pullback
$\pi_X^*(E)$, where $I=[0,1]$ and $\pi_X:X\times I\to X$ is
the projection. Consider on $E_I$ a connection $A_I$ whose
restriction on $X\times\{t\}$ is $(1-t)A_0+tA_1$ 
($A_I$ is thus in {\it temporal gauge}). 
Denote $\sigma_t:X\to X\times I$ the map which sends 
$x\in X$ to $\sigma_t(x)=(x,t)$. We have
$\phi_{A_j}(\eta)=\sigma_j^*\phi_{A_I}(\eta)$ for $j=0,1$.

Define a map $h:\Omega^{*}(\FFF\times I)\to\Omega^{*-1}
(\FFF)$ as follows. Any form in $\Omega(\FFF)$ may be written 
as $\alpha\wedge dt+\beta$, where $t\in I$ is a coordinate and in
such a way that $\iota_{\partial/\partial t}\alpha=
\iota_{\partial/\partial t}\beta=0$.
Set $h(\alpha\wedge dt+\beta)=\int_I\alpha$. A simple computation shows
that for any $\gamma\in\Omega(M\times I)$,
$(hd+dh)(\gamma)=\sigma_1^*\gamma-\sigma_0^*\gamma.$
Applying this to $\gamma=\phi_{A_I}(\eta)$ we
get $$dh\phi_{A_I}(\eta)=\sigma_1^*\phi_{A_I}(\eta)-\sigma_0^*\phi_{A_I}(\eta)
=\phi_{A_1}(\eta)-\phi_{A_0}(\eta),$$
since by theorem \ref{phimap} $\phi_{A_I}(\eta)$ is a closed 
element in $\Omega(\FFF\times I)$.
This implies that $\phi_{A_0}(\eta)$ and $\phi_{A_1}(\eta)$ are 
cohomologous.
\end{pf}

\begin{lemma}
$\phi_A(\overline{\omega}_F)=\omega_F^A$.
\label{fiomec}
\end{lemma}
\begin{pf}
By definition $\overline{\omega}_F=1\otimes\omega_F-\mu\otimes 1\in
(\CC[\klie]\otimes\Omega(F))^K$. So
$$\phi_A(\overline{\omega}_F)=h(1\otimes\omega_F)-
h(\la\Omega_A,\mu\ra_{\klie}\otimes 1)\in \Omega(E\times F)_{\bas}.$$
We have $h(1\otimes\omega_F)=1\otimes\alpha^*\omega_F=
1\otimes\tilde{\omega}_F^A$. On the other hand, the form 
$\la\Omega_A,\mu\ra_{\klie}\otimes 1$ is horizontal (because the curvature
$\Omega_A\in\Omega^2(E,\klie)$ is a horizontal form, see proposition 1.13
in \cite{BeGeV}),
so $\phi_A(\overline{\omega}_F)=1\otimes\tilde{\omega}_F^A+
\la\Omega_A,\mu\ra_{\klie}\otimes 1$. Now by lemma \ref{formequi}
this form represents $\omega_F^A\in\Omega^2(\FFF)$.
\end{pf}

This lemma, together with theorem \ref{phimap}, finishes the proof
of proposition \ref{nonlinka2}. We can now restate proposition \ref{nonlinka} 
as follows.

\begin{prop}
For any section $\Phi\in\SSS$ and for any connection $A\in\AAA$
$$\int_X\la \Lambda F_A,\mu(\Phi)\ra_{\klie}=
\frac{1}{2}(\|\partial_A\Phi\|^2_{L^2}
-\|\overline{\partial}_A\Phi\|^2_{L^2})-
\int_X \Phi^*\phi_A(\overline{\omega}_F)\wedge \omega^{[n-1]}.$$
\label{defnonlinka}
\end{prop}

\subsubsection{Proof of theorem \ref{idminim}}
The following computation has its origins in an idea of 
Bogomolov in studying vortex equations on $\RR^2$. Here we
mimic \cite{Br1}, except that where he uses the Kaehler
identities we use proposition \ref{defnonlinka}.

\begin{lemma}
For any section $\Phi\in\SSS$ and any connection $A\in\AAA$
\begin{align*}
\YMH_c(A,\Phi)&=\|\Lambda F_A+\mu(\Phi)-c\|_{L^2}^2+
2\|\overline{\partial}_A\Phi\|_{L^2}^2+4\|F_A^{0,2}\|_{L^2}^2 \\
&+2\int_X\la\Lambda F_A,c\ra
-\int_X B(F_A,F_A)\wedge\omega^{[n-2]}
+2\int_X \Phi^*\phi_A(\overline{\omega}_F)\wedge \omega^{[n-1]},
\end{align*}
where $B:\Omega^2(E\times_{\Ad}\klie)\otimes\Omega^2(E\times_{\Ad}\klie)
\to\Omega^4(X)$ combines the wedge product in $\Omega^*(X)$ with
the biinvariant pairing $\klie\otimes\klie\to\RR$.
\label{reescriuYMH}
\end{lemma}
\begin{pf}
Throughout the proof $\|\cdot\|$ will denote $L^2$ norm (which, recall,
is computed using the volume form $\omega^{[n]}$). 
The following formula is well known (and easily checked)
$$\|F_A\|^2=\|\Lambda F_A\|^2-\int_X B(F_A,F_A)\wedge\omega^{[n-2]}
+4\|F_A^{0,2}\|^2.$$
We develop using the above formula and proposition \ref{defnonlinka}
\begin{align*}
\|\Lambda F_A+\mu(\Phi)-c\|^2&+
2\|\overline{\partial}_A\Phi\|^2+4\|F_A^{0,2}\|^2
+2\int_X\la\Lambda F_A,c\ra \\
&-\int_X B(F_A,F_A)\wedge\omega^{[n-2]}
+2\int_X \Phi^*\phi_A(\overline{\omega}_F)\wedge \omega^{[n-1]} \\
&=\|\Lambda F_A\|^2+4\|F_A^{0,2}\|^2-\int_X B(F_A,F_A)\wedge\omega^{[n-2]}
+\|\mu(\Phi)-c\|^2 \\
&+2\int_X\la\Lambda F_A,\mu(\Phi)\ra_{\klie}
+2\|\overline{\partial}_A\Phi\|^2
+2\int_X \Phi^*\phi_A(\overline{\omega}_F)\wedge \omega^{[n-1]} \\
&=\|F_A\|^2+\|\partial_A\Phi\|^2+\|\overline{\partial}_A\Phi\|^2
+\|\mu(\Phi)-c\|^2 \\
&=\|F_A\|^2+\|d_A\Phi\|^2+\|\mu(\Phi)-c\|^2.
\end{align*}
\end{pf}

Theorem \ref{minimitza} follows easily from the preceeding lemma. Indeed, 
$$2\int_X\la\Lambda F_A,c\ra+
2\int_X \Phi^*\phi_A(\overline{\omega}_F)\wedge \omega^{[n-1]}
-\int_X B(F_A,F_A)\wedge\omega^{[n-2]}$$
is a topological quantity, that is, it only depends on the homotopy class
of $\Phi$. That this is true for the second summand is clear; as for
the first summand, by Chern-Weil theory one sees that it is equal
to a linear combination whose coefficients depend on $c$ of first
Chern classes of line bundles obtained from $E$ through representations
$K\to S^1$. Finally, the third summand is equal to a linear combination
of degree 4 pieces of Chern characters of bundles associated to $E$
wedged with $\omega^{[n-2]}$ and integrated over $X$
(see p. 209 in \cite{Br2}).

Finally, we obtain from \ref{minimitza} the following corollary
{\it {\`a} la} Bogomolov

\begin{corollary}
Suppose that a pair $(A,\Phi)$ satisfies equations (\ref{equs}). 
Then the following inequality holds 
$$\int_X\la\Lambda F_A,c\ra+
\int_X \Phi^*\phi_A(\overline{\omega}_F)\wedge \omega^{[n-1]}
-\frac{1}{2}\int_X B(F_A,F_A)\wedge\omega^{[n-2]}\geq 0.$$
\label{bogom0}
\end{corollary}

\subsection{$L^2$ bounds for solutions of (\ref{equs})}
Suppose now that $X$ is a Riemann surface. Let
$\pi_F:F_K=EK\times_K F\to BK$, where $EK\to BK$ is the universal
principal $K$-bundle. This fibration is unique only up to homotopy,
and we chose a model of it for which $EK\to BK$ is a smooth principal
$K$-bundle. Recall that the equivariant (co)homology
of $F$ is by definition the (co)homology of $F_K$.                    
Cartan proved that the cohomology of the complex 
$(\Omega_K(F),d_{\klie})$ is isomorphic to the real $K$-equivariant 
cohomology $H_K^*(F;\RR)$ of $F$.

Let $c_E:X\to BK$ be the classifying map of the bundle $E$.
Let us take an isomorphism of $K$ principal bundles $\phi:E\simeq c_E^*EK$
($\phi$ is equivalently given by any $K$ equivariant map $E\to EK$).
Let $\psi:\FFF\to F_K$ be the map induced by $\phi$.
Let $\AA$ be any connection on $F_K\to BK$. Consider the Chern-Weil map
$$\phi_{\AA}:(\Omega_K(F),d_{\klie})\to\Omega(F_K,d)$$
(here we are assuming that the fibration $EK\to BK$ is smooth).
We clearly have $\phi_{\psi^*\AA}=\psi^*\circ \phi_{\AA}$.
Now, using lemma \ref{cohonodepen} we compute for any pair
$(A,\Phi)\in\AAA\times\SSS$
\begin{align*}
\int_X\Phi^*\phi_A(\ov{\omega}_F)&=
\int_X\Phi^*\phi_{\psi^*\AA}(\ov{\omega}_F)=
\int_X(\psi\Phi)^*\phi_{\AA}(\ov{\omega}_F) \\
&= \la (\psi\Phi)_*[X],[\phi_{\AA}(\ov{\omega}_F)]\ra,
\end{align*}
where $[X]\in H_2(X;\ZZ)$ is the fundamental class and where
$[\phi_{\AA}(\ov{\omega}_F)]\in H_2(F_K;\RR)$ denotes the 
cohomology class represented by the form $\phi_{\AA}(\ov{\omega}_F)$.

On the other hand,
$$\int_X\la\Lambda F_A,c\ra=\int_X\la F_A,c\ra
=\int_X\la F_{\psi^*\AA},c\ra=\la {c_E}_*[X],[F_{\AA}]\ra.$$
But ${c_E}_*=(\pi_F\psi\Phi)_*$ so, using lemma \ref{reescriuYMH},
we conclude the following.

\begin{theorem}
Let $(A,\Phi)$ be a solution of equations (\ref{equs}). There is a
constant $C$, depending only on $(\psi\Phi)_*[X]$, such that
$$\|F_A\|^2_{L^2}\leq C\mbox{ and }\|d_A\Phi\|^2_{L^2}\leq C.$$
\label{boundL2}
\end{theorem}

This result will be useful in the sequel, since we will consider
the space of solutions to equations (\ref{equs}) for different
choices of $E$ and homotopy class of $\Phi$. Note, on the other
hand, that lemma \ref{indepframe} in the appendix implies that
the homology class $(\psi\Phi)_*[X]$ does not depend on the
particular map $\phi:E\to EK$ chosen to construct $\psi$.


\chapter{Hitchin--Kobayashi correspondence}
\label{hiko}
Let us suppose that the complex structure $I_F$ on $F$ is integrable,
that is, $F$ is a Kaehler manifold. Then the action of $K$ on $F$
extends to a unique holomorphic action of the complexification $G=K^{\CC}$
(see \cite{GS}). Let $\GGG_K=\Gamma(E\times_{\Ad}K)$ 
(resp. $\GGG_G=\Gamma(E\times_{\Ad}G)$) be the real (resp. complex)
gauge group. The diagonal action of $\GGG_K$ on $\AAA\times\SSS$ extends in a
natural way to a diagonal action of $\GGG_G$ (see \ref{accioGG} for a 
description of the action on $\AAA$; the action on $\SSS$ comes from the 
action of $G$ on $F$). Then, just as in classical 
Yang--Mills--Theory,
the first and third equations in (\ref{equs}) are $\GGG_G$ invariant,
whereas the second one is only $\GGG_K$ invariant. In this chapter
we will study which $\GGG_G$ orbits of pairs 
$(A,\Phi)\in\AAA^{1,1}\times\SSS$ contain solutions to the second equation.
Observe that we restrict to pairs solving the third equation $F_A^{0,2}=0$,
while on the contrary we do not ask the first equation 
$\ov{\partial}_A\Phi=0$ to be satisfied. This is a technical condition,
and one could probably study the second equation on any $\GGG_G$ orbit
in $\AAA\times\SSS$ using essentially the same methods as here.

More concretely, the question adressed in this chapter is the following:
given a pair $(A,\Phi)\in \AAA^{1,1}\times\SSS$, decide whether there exist 
a gauge transformation $g\in\GGG_G$ such that the pair $g(A,\Phi)$ satisfies
\begin{equation}      
\Lambda F_{g(A)}+\mu(g(\Phi))=c.
\label{hkequs}           
\end{equation}
This question will be partially answered. We will define the notion of simple 
pair (definition \ref{defsimple}) and a condition on pairs $(A,\Phi)$ called 
$c$-stability (definition \ref{parella_estable}), and in theorem
\ref{main} we will prove that, if $(A,\Phi)$ is a simple pair,
then there exist a gauge $g\in\GGG_G$ sending
$(A,\Phi)$ to a pair $g(A,\Phi)$ satisfying the equation (\ref{hkequs})
if and only if $(A,\Phi)$ is $c$-stable.
We will also prove in theorem \ref{main}
that in each $\GGG_G$ orbit inside $\AAA^{1,1}\times\SSS$ there
is at most one $\GGG_K$ orbit of pairs which satisfy (\ref{hkequs}).
We call such a characterization of solutions to (\ref{hkequs})
a Hitchin--Kobayashi correspondence because it generalises the correspondence
in Yang--Mills--Higgs theory with this name. We will prove the correspondence
when the biinvariant metric in $\klie$ satisfies a certain
property (see subsection \ref{repauxiliar}). Just as the equation, our
existence criterion will depend on this metric. 

One can look at theorem \ref{main} from two different points of view.
When $X$ consists of a single point, the curvature term vanishes
in equation (\ref{hkequs}), and so our problem reduces to a well
known one in Kaehler geometry. Namely, that of studying which $G$ orbits
inside $F$ contain zeroes of the moment map $\mu$.
More generally, one studies which $G$ orbits have points whose image is
a fixed central element in $\klie^*$ or belongs to a given coadjoint orbit
in $\klie^*$. If $F$ is a projective manifold,
one can answer this question in a very satisfactory way: 
a $G$ orbit contains a zero of the moment map if and only if
it is stable in the sense of Mumford Geometric Invariant Theory 
(GIT for short) \cite{KeNe, MFK, GS}. 
To extend the notion of GIT stability to actions on any Kaehler manifold
$F$, we use the notion of analytic stability (see definition
\ref{estabilitat}).
This notion coincides with that of GIT stability in the case of projective 
manifolds, and characterizes the $G$-orbits 
in which the moment map vanishes somewhere (see theorem \ref{corr}).
This is the content of the so called Kempf-Ness theory. 
So, in this sense, our result can be viewed as a fibrewise
generalisation of Kempf-Ness theory. 

There is, however, another point of view which allows to look at 
theorem \ref{main} as a result {\it {\`a} la} Kempf-Ness in infinite 
dimensions. One can give a Kaehler structure to the configuration space
$\AAA\times\SSS$. Then the action of the gauge group
$\GGG_K$ on $\AAA\times\SSS$ is symplectic and by isometries, 
and the left hand side in
equation (\ref{hkequs}) is the moment map of this action 
(see sections \ref{exemp1}, \ref{exemp2} and \ref{disgressio}).
Finally, $\AAA^{1,1}\times\SSS$ is a $\GGG_K$ invariant complex 
subvariety (with singularities) of $\AAA\times\SSS$. This point 
of view was adopted for the first time in the context of gauge theories
by Atiyah and Bott \cite{AB} in their study of Yang-Mills equations
over Riemann surfaces, which are a particular case of the equations
that we consider. The idea of Atiyah and Bott was used by Donaldson
\cite{Do1} in his proof of the theorem of Narasimhan and Seshadri
(see below), and it has been subsequently often used in studying other 
particular cases of equation (\ref{hkequs}).

\section{Stability and statement of the correspondence}
\label{statement}
\subsection{The Lie group $K$ and its complexification $G$}
\label{repauxiliar}
Let $G$ be the complexification of $K$ (see for example \cite{BtD} for a 
general construction of the complexification of compact Lie groups). 
Let $\glie=\klie\otimes_{\RR}\CC$ be the Lie algebra of $G$.

Let $\rho_a:G\to\GL(W_a)$ be a faithful representation on a finite
dimensional complex vector space $W_a$, and take on $W_a$ a Hermitian
metric such that $\rho_a(K)\subset U(W_a)$. We will call $\rho_a$ the 
{\bf auxiliar representation}. We denote the restriction
of $\rho_a$ to $K$ with the same symbol, and the induced representation 
of $\glie$ as well. Define the following pairing on $\glie$:
$$\begin{array}{rcl}
\la,\ra:\glie\times\glie & \longrightarrow & \CC \\
(u,v) & \mapsto & \la u,v\ra =\Tr(\rho_a(u)\rho_a(v)^*).
\end{array}$$
This pairing is nondegenerate and its restriction to $\klie$ 
gives a biinvariant metric. In this chapter we will assume that
the metric on $\klie$ used to give a sense to our equations
is precisely this one. Note that not all biinvariant metrics
on $\klie$ come from a representation as above.
For example, the biinvariant metrics on $\Lie(S^1)$ are in bijection with
$\RR_+$ and only a discrete subset of it corresponds to
metrics coming from representations. Most likely, however, this condition
on the metric on $\klie$ could be relaxed.

\subsection{The setting}
\label{setting}
\subsubsectionr{}
In this chapter $X$ will be any compact Kaehler manifold, and we will 
denote $\omega$ the symplectic form on $X$ and $I\in\End(TX)$ the complex 
structure. Recall that we have a principal $K$ bundle $\pi:E\to X$.
Let $\AAA$ be the space of connections on $E$. We denote
$\AAA^{1,1}\subset\AAA$ the set of connections which satisfy $F_A^{0,2}=0$. 
By the theorem of Newlander and Niremberg (see \cite{NewNi}) these are the 
connections $A$ such that $\overline{\partial}_A$ defines an integrable 
holomorphic structure on any associated complex vector bundle $V$.

Let $\GGG_K=\Gamma(E\times_{\Ad}K)$ the gauge 
group of $E$. Let $\pi_G:E_G=E\times_K G\to X$ be the associated bundle 
with fibre $G$ (we consider the action of $K$ on $G$ on the left). This is 
a principal $G$ bundle. Let 
$\GGG_G=\Gamma(E_G\times_{\Ad}G)=\Gamma(E\times_{\Ad}G)$ be the gauge 
group of $E_G$. The group $\GGG_G$ is the complexification of $\GGG_K$. 

\subsubsectionr{}
In this chapter the manifold $F$ will be assumed to be Kaehler. We will 
write its symplectic form $\omega_F$ and its complex structure
$I_F\in\End(TF)$ (recall that both structures are invariant under the action
of $K$ by hypothesis). Let $\FFF=E\times_K F$ be the associated 
bundle and let $\SSS=\Gamma(\FFF)$ the the set of sections of $\FFF$.
The gauge group $\GGG_K$ acts on the left on $\FFF$, and hence acts 
also on $\SSS$.
Since $F$ is Kaehler, the action of $K$ on $F$ extends to a unique 
holomorphic action of $G$ (see \cite{GS}). This allows to extend the 
action of $\GGG_K$ on $\SSS$ to an action of $\GGG_G$.

\subsubsectionr{}
\label{accioGG}
Let $\CCC$ be the set of $G$-invariant complex structures on $E_G$
for which the map $d\pi_G:TE_G\to \pi_G^*TX$ is complex.
We define a map $\Chern:\CCC\to\AAA$, called the {\bf Chern map},
as follows. An invariant complex structure $I\in\CCC$ is mapped to the
connection whose horizontal distribuition is
$I(TE)\cap TE\subset TE$ (this makes sense, since the inclusion
$E=E\times_K K\subset E\times_K G$ given by $K\subset G$
induces an inclusion $TE\subset TE_G$). (This distribuition is $K$ invariant
and hence corresponds to a connection because $I$ is $G$-invariant.) 
The map $\Chern$ is a bijection.
Its inverse sends any connection $A\in\AAA$ to the
complex structure $I_{E_G}(A)$ on $E_G$ defined in \ref{THCshol}
(taking $F=G$).

\begin{lemma}
Let $A\in\AAA$ and let $I_{\FFF}(A)$ be the induced complex structure on 
$\FFF$ as in \ref{THCshol}. 
Let $I_{E_G}(A)=\Chern^{-1}(A)$. By $G$-invariance the 
complex structure $I_{E_G}(A)+I_F$ on $E_G\times F$ descends to give a 
complex structure on $\FFF=E_G\times_G F$. This complex structure coincides 
with $I_{\FFF}(A)$.
\label{integr}
\end{lemma}
          
The group $\GGG_G$ acts on $\CCC$ by pullback.                  
Using the bijection $\Chern:\CCC\to\AAA$ we 
transfer the action of $\GGG_G$ on $\CCC$ to an action on $\AAA$. This 
action extends the action of $\GGG_K$ and leaves invariant the subset 
$\AAA^{1,1}\subset\AAA$.

\subsection{Group actions on Kaehler manifolds}
\label{grkaehler}

We will denote $\la,\ra$ the Kaehler metric on $F$. This metric is given 
by $\la u,v\ra=\omega_F(u,Iv)$. Let $s\in\klie$ be any nonzero element. 
Write $\mu_{s}=\la\mu,s\ra_{\klie}:F\to\RR.$

\begin{lemma}
The gradient of $\mu_{s}$ is $I\fX_{s}$.
\label{gradient}
\end{lemma}
\begin{pf}
Let $x\in F$ and take any vector $v\in T_xF$.
Then $\nabla_v(\mu_{s})=\la d\mu_{s},v\ra_{T_xF}=
\omega_F(\fX_{s},v)=\omega_F(I\fX_{s},Iv)=\la I\fX_{s},v\ra,$ 
by the definition of moment map. 
\end{pf}

Consider the gradient flow $\phi^t_{s}:F\to F$ of the function
$\mu_{s}$, which is defined by these properties:
$\phi^0_{s}=\Id$ and $\frac{\partial}{\partial t}\phi_s^t
=\nabla(\mu_{s})=I\fX_{s}$.
Using the action of $G$ on $F$ we can write $\phi^t_{s}(x)=e^{\imag ts}x$.

\begin{definition}
Let $x\in F$ be any point, and take an element $s\in\klie$.
Define $$\lambda_t(x;s)=\mu_s(e^{\imag ts}x).$$
Define also the {\bf maximal weight $\lambda(x;s)$ of the action 
of $s$ on $x$}
as $$\lambda(x;s)=\lim_{t\to \infty}\lambda_t(x;s)
\in\RR\cup\{\infty\}.$$
\label{defpesmax}
\end{definition}

This limit always exists since by lemma \ref{gradient} the function
$\lambda_t(x;s)$ increases with $t$.
The definition of the maximal weight depends on the
chosen moment map. Since this is not unique, we will sometimes 
write the maximal weight of $s\in\klie$ acting on $x\in F$
with respect to the moment map $\mu$ as $\lambda^{\mu}(x;s)$.

\begin{prop}
The maximal weights satisfy the following properties:
\begin{enumerate}
\item They are $K$-equivariant, that is, for any $k\in K$,
$\lambda(kx;ksk^{-1})=\lambda(x;s)$.
\item For any positive real number $t$ one has
$\lambda(x;ts)=t\lambda(x;s)$.
\end{enumerate}
\label{proppes}
\end{prop}

See sections \ref{banfield} and \ref{filtracions} for explicit
computations of maximal weights in some particular situations.

\subsection{Parabolic subgroups}
\label{parabolic}
A good reference for this material is \cite{R2}.
Let $\glie$ be the Lie algebra of $G$, and split $\glie=\zlie\oplus
\glie^s$ as the sum of the centre plus the semisimple part 
$\glie^s=[\glie,\glie]$ of $\glie$.
Take a Cartan subalgebra $\hlie\subset\glie^s$. 
Let $R\subset\hlie^*$ be the set of roots.
We can decompose
$$\glie=\zlie\oplus\hlie\oplus\bigoplus_{\alpha\in R}\glie_{\alpha},$$
where $\glie_{\alpha}\subset\glie^s$ is the subspace on which
$\hlie$ acts through the character $\alpha\in\hlie^*$.

Fixing a (irrational) linear form on $\hlie^*$, we divide the set
of roots in positive and negative roots: $R=R^+\cup R^-$.
Let us denote the set of simple roots by  $\Delta=(\alpha_1,\dots,\alpha_r)
\subset R^+$. Recall that the set $\Delta$ is characterised by the
following property: any root can be written as 
a linear combination of the elements of $\Delta$ with integer
coefficients all of the same sign. Furthermore, $r$ equals 
$\dim_{\CC}\hlie$, the rank of $G$. The simple coroots are by definition 
$\alpha_j'=2\alpha_j/\la\alpha_j,\alpha_j\ra$, where $1\leq j\leq r$.

We have taken a maximal compact subgroup $K\subset G$.
From now on we will assume that the following relation holds between $K$ and
the Cartan subalgebra $\hlie$: $\zlie\oplus\hlie$ is the
complexification of the Lie algebra $\tlie$ of a maximal torus $T\subset K$.

\begin{lemma}
Chose, for any root $\alpha\in R$, a nonzero element
$g_{\alpha}\in\glie_{\alpha}$ in such a way that $g_{\alpha}$
and $g_{-\alpha}$ satisfy $\la g_{\alpha},g_{-\alpha}\ra=1$.
Let $\RR R^*\subset \hlie$ denote the real span of the 
duals (with respect to the Killing metric) of the roots. 
Assume that $\zlie\oplus\hlie$ is the complexification of the Lie algebra
of a maximal torus $T$ of a maximal compact subgroup $K\subset G$.
Then $\glie^s\cap\klie=\imag\RR R^*\oplus \bigoplus_{\pm\alpha\in R}
\RR(g_{\alpha}+g_{-\alpha})\oplus \RR(\imag g_{\alpha}-\imag g_{-\alpha}).$
\label{hcompact}
\end{lemma}

This lemma (and the following ones in this subsection) can be easily proved 
using basic results on reductive Lie groups (see for example \cite{FH}).

Let $\lambda_1,\dots,\lambda_r$ be the set of fundamental weights,
which belong to $\hlie^*$ and are the duals with respect to the 
Killing metric of the simple coroots.
Let us denote by $\lambda'_1,\dots,\lambda'_r$ the elements in $\hlie$
dual to the fundamental weights through the Killing metric.

To define a parabolic subgroup of $G$, take any subset
$A=\{\alpha_{i_1},\dots,\alpha_{i_s}\}\subset \Delta$. Let
$$D=D_A=\{\alpha\in R | \alpha=\sum_{j=1}^r m_j\alpha_j
\mbox{, where $m_{i_t}\geq 0$ for $1\leq t\leq s$}\}.$$

\begin{definition}
The subalgebra $\plie=\zlie\oplus
\hlie\oplus\bigoplus_{\alpha\in D}\glie_{\alpha}$
will be called the {\bf parabolic subalgebra} of $\glie$ 
with respect to the set $A\subset\Delta$. 
The connected subgroup $P$ of $G$ whose subalgebra is $\plie$ will be called 
the {\bf parabolic subgroup} of $G$ with respect to $A$.
Furthermore, any positive (resp. negative) linear combination of the 
fundamental weights $\lambda_{i_1},\dots,\lambda_{i_s}$ plus an element
of the dual of $\imag(\zlie\cap\klie)$ will be called 
a {\bf dominant} (resp. {\bf antidominant}) {\bf character} on 
$\plie$ (or on $P$).
\end{definition}

\begin{remark}
We will regard $G$ as a parabolic subgroup 
of itself (with respect to the empty set $\emptyset\subset \Delta$).
\end{remark}

Observe that our definition of parabolic subgroup depends upon the 
choice of a Cartan subalgebra $\hlie\subset\glie$ and of a linear
form on $\hlie^*$. In general, any parabolic subgroup
$P\subset G$ obtained from a different choice of Cartan subalgebra
and linear form will be conjugate to a parabolic subgroup obtained
from our data.

Let $\rho:K\to U(W_{\rho})$ be a representation on a Hermitian vector space
$W_{\rho}$.
We will denote its (unique) lift to a holomorphic representation
of the complexification $G$ of $K$
by the same letter $\rho:G\to GL(W_{\rho})$. Take $P\subset G$ to be 
the parabolic subgroup with respect to a set $A=\{\alpha_{i_1},\dots,
\alpha_{i_s}\}\subset\Delta$. Let $\chi$ be the dual of 
an antidominant character of $P$.
Thanks to our conventions (lemma \ref{hcompact}), $\chi$ belongs
to $\imag\klie$. So, since $\rho$ is unitary, $\rho(\chi)$ diagonalises
and has real eigenvalues.
Let $\lambda_1<\dots<\lambda_r$ be the set of different
eigenvalues of $\rho(\chi)$, and let us write $W(\lambda)$ the eigenspace
of eigenvalue $\lambda$. 
Let $W^{\lambda_k}=\bigoplus_{j\leq k}W(\lambda_j)$, 
and let $\fW_{\rho}(\chi)$ be the partial flag
$0\subset W^{\lambda_1}\subset\dots\subset W^{\lambda_r}=W_{\rho}$.

\begin{lemma}
(i) The action of $P$ leaves invariant the partial flag $\fW_{\rho}(\chi)$.
Suppose that the restriction of $\rho$ to the semisimple part
$\plie^s$ of $\plie$ is faithful.
If $\chi=z+\sum_{k=1}^s m_k\lambda'_{i_k}$, where $z\in\zlie$, and, for any
$k$, $m_k<0$, then $P$ is precisely the preimage by $\rho$
of the stabiliser of $\fW_{\rho}(\chi)$.
(ii) Let $\chi\in\imag\klie$ be any element. There is a choice of 
Cartan subalgebra $\hlie\subset\glie$ contained in $\plie$
such that $\chi\in\hlie$ and $\chi$ is antidominant with 
respect to $P$ if and only if the stabiliser of the partial
flag $\fW_{\rho}(\chi)$ contains $P$.
\label{equi1}
\end{lemma}

\begin{lemma}
Let $\chi$ be any element in $\imag\klie$. The preimage by 
$\rho$ of the stabiliser of $\fW_{\rho}(\chi)$ is a parabolic
subgroup $P_{\rho}(\chi)$ of $G$. Moreover, $\chi$ 
is the dual of an antidominant character of $P_{\rho}(\chi)$.
\label{equi3}
\end{lemma}

Let us take now any subspace $W'\subset W_{\rho}$ 
belonging to the filtration $\fW_{\rho}(\chi)$. We define 
$\overline{W}'= G\times_{\rho} W'\to G/P$.
In other words, $\overline{W}=G\times W'/\sim$, where 
$(gp,w)\sim (g,pw)$ for any $(g,w)\in G\times W'$ and $p\in P$. 
This makes sense, since $P$ leaves $W'$ invariant and so
$\rho:P\to GL(W')$. Define also an action of
$G$ on $G\times W'$ by $g'(g,w)=(g'g,g^{-1}g'gw)$. This action is
compatible with the relation $\sim$. Indeed, 
$$g'(gp,w)=(g'gp,p^{-1}g^{-1}g'gpw)\sim(g'g,g^{-1}g'gpw)=g'(g,pw).$$
Repeating this for each subspace in $\fW_{\rho}(\chi)$ we obtain the
following.

\begin{lemma}
The filtration of vector bundles 
$\overline{\fW}_{\rho}(\chi)=G\times_{\rho}\fW_{\rho}(\chi)\to G/P$
is $G$-equivariant and holomorphic.
\label{tautologic}
\end{lemma}

\subsection{Parabolic and maximal compact subgroups}
\label{maxcomp}
Given any parabolic subgroup $P\subset G$ with Lie algebra $\plie$,
we will write $P_K$ (resp. $\plie_K$) the subgroup
$P\cap K$ (resp. the subalgebra $\plie\cap\klie$). 
$P_K$ is a maximal compact subgroup of $P$.

\begin{lemma}
Let $E_G\to X$ be a $G$-principal bundle on any topological space $X$.
If $E_G$ admits reductions of its structure group from 
$G$ to a parabolic subgroup $P$ and to the maximal compact subgroup $K$,
then it also admits a reduction of its structure group from $G$
to $P_K$.
\label{maxcomp1}
\end{lemma}
\begin{pf}                                                       
Consider the surjections $\pi_P:G/P_K\to G/P$ and
$\pi_K:G/P_K\to G/K$. We will prove that, for any pair
$(g_PP,g_KK)\in G/P\times G/K$, the intersection
$\pi_P^{-1}(g_PP)\cap \pi_K^{-1}(g_KK)\subset G/P_K$ consists
of a single point. We can assume, multiplying on the left $g_PP$ and $g_KK$
by $g_P^{-1}$, that $g_PP=P$. So, if $[P]\subset G/P_K$
(resp. $[K]\subset G/P_K$) is $\pi_P^{-1}(P)$ (resp. $\pi_K^{-1}(K)$),
we have to check that for any $g\in G$, $[P]\cap g[K]\subset G/P_K$
is a point. Using intersection theory, it is enough to verify that
$[P]$ and $g[K]$ intersect transversely for any $g$ and that 
$[P]\cap [K]$ consists of a single point 
(indeed, given any point $g\in G$ we can connect $g$ to $1\in G$
with a path, since $G$ is connected; then, if for any point $g$
in the path $[P]$ and $g[K]$ intersect transversely, since $[K]$
is compact and $[P]$ is closed, $\sharp [P]\cap[K]=\sharp [P]\cap g[K]$,
where $\sharp A$ denotes the number of elements in $A$).
Let $\blie=\zlie\oplus\hlie\oplus
\bigoplus_{\alpha\in R^+}\glie_{\alpha}$ (this is a Borel subalgebra
of $\glie$). Thanks to lemma \ref{hcompact}, 
$\blie+\klie=\glie$. Any parabolic subalgebra $\plie$
contains $\blie$ as a subalgebra, so $\plie+\klie=\glie$.
Now, if $g,p\in G$ and $[p]\in [P]\cap g[K]$
($[p]$ denotes the class of $p$ in $G/P_K$), 
then $T_{[p]}[P]=p\plie/\plie_K$ and $T_{[p]}[K]=p\klie/\plie_K$.
So $T_{[p]}[P]+T_{[p]}[P]=T_{[p]}(G/P_K)$, and this means that
the intersection is transverse. On the other hand,
$[P]\cap [K]\subset P_K$ consists of a point by the definition
of $P_K$.

Now suppose that there are reductions $\sigma_P\in\Gamma(E(G/P))$
and $\sigma_K\in\Gamma(E(G/K))$ (here $E(G/P)$ denotes the bundle
$E_G\times_G(G/P)$ associated to $E$ with fibre 
the homogeneous space $G/P$, and similarly $E(G/K)$).
This means that on each point $x\in X$, after identifying 
$(E_G)_x$ with $G$, we have $\sigma_P(x)\in P/G$ and $\sigma_K(x)\in P/K$.
By the preceeding observation, these give a unique point
in $G/P_K$. Doing this in every fibre we get a unique section
$\sigma_{P_K}\in\Gamma(E(G/P_K))$ (smoothness is a consequence
of transversality), which is a reduction of the structure group 
of $E_G$ to $P_K$.
\end{pf}

\begin{lemma}
Let $P$ be a parabolic subgroup with respect to the set
$$A=\{\alpha_{i_1},\dots,\alpha_{i_s}\}\subset\Delta.$$
For any $j\in\{i_1,\dots,i_s\}$, the element $\lambda'_j\in\imag\klie$ 
(dual with respect to the Killing metric of the fundamental weight
$\lambda_j$) is left fixed by the adjoint action of $\plie_K$
on $\glie$.
\label{maxcomp2}
\end{lemma}
\begin{pf}
Fix an $i\in\{i_1,\dots,i_s\}$. We will prove that 
$\lambda'_i$ belongs to the centre of $(\plie\cap\klie)\otimes_{\RR}\CC$.
Let $$D^s_A=\{\alpha\in R|\alpha=\sum_{j=1}^rm_j\alpha_j,
\mbox{ where $m_{i_t}=0$ for $1\leq t\leq s$}\}.$$
Then $(\plie\cap\klie)\otimes_{\RR}\CC=
\zlie\oplus\hlie\oplus\bigoplus_{\alpha\in D^s}\glie_{\alpha}$. 
This stems from lemma \ref{hcompact}.
If $g\in\zlie\oplus\hlie$, then clearly $[\lambda'_i,g]=0$.
And if $g\in\glie_{\alpha}$, for $\alpha\in D^s_A$, then
$[\lambda'_i,g]=\lambda_i(\alpha)=0$.
As a consequence, for any $g\in(\plie\cap\klie)\otimes_{\RR}\CC$,
one has $[\lambda'_i,g]=\lambda_i(g)=0$, which is what we wanted
to prove.
\end{pf}

\subsection{Reductions of the structure group and filtrations}
\label{reduccions}
Let $V=V_{\rho_a}=E\times_{\rho_a}W_a$ be the vector bundle
associated to the auxiliar representation (see section \ref{repauxiliar}).
In this subsection we will see that there is a bijection between the 
reductions of the structure group of $E$ to a parabolic subgroup $P$ 
together with an antidominant character of $P$, and certain filtrations of 
$V$ by subbundles. We denote $E(G/P)$ the bundle $E_G\times_G(G/P)$.
The space of reductions of the structure group of $E_G$ from $G$
to $P$ is $\Gamma(E(G/P))$. 

\subsubsectionr{}
\label{sectfilt}
Fix a parabolic subgroup $P\subset G$ and take a reduction
$\sigma\in\Gamma(E(G/P))$. Let $\chi$ be an antidominant character
for $P$. There is a canonical reduction of
the structure group $G$ of $E_G$ to $K$, since $E_G=E\times_K G$.
By lemma \ref{maxcomp1}, this reduction, together
with $\sigma$, gives a reduction $\sigma_K\in \Gamma(E(G/P_K))$, where
$P_K=P\cap K$. And then, lemma \ref{maxcomp2} implies that
we get a section $g_{\sigma,\chi}\in\Omega^0(E\times_{\Ad}\imag\klie)
=\imag\Lie(\GGG_K)$ which is fibrewise the dual of $\chi$.

With the element $g_{\sigma,\chi}$ we can obtain a filtration of
$V_{\rho}$ as follows. First of all, $\rho(g_{\sigma,\chi})$
has constant real eigenvalues (which are equal to those of
$\rho(\chi)\in\End(W_{\rho})$). Let $\lambda_1<\dots<\lambda_r$ 
be the different eigenvalues, and let $V_{\rho}(\lambda_j)$ 
be the eigenbundle of eigenvalue $\lambda_j$.
Finally, let $V_{\rho}^{\lambda_k}=\bigoplus_{i\leq k}V_{\rho}(\lambda_j)$.
Denote by $\fV_{\rho}(\sigma,\chi)$ the filtration 
$$0\subset V_{\rho}^{\lambda_1}\subset
V_{\rho}^{\lambda_2}\subset\dots\subset V_{\rho}^{\lambda_r}=V_{\rho}.$$

Alternatively, recall that on $G/P$ there is a 
filtration of $G$-equivariant (holomorphic)
vector bundles, $\overline{\fW}_{\rho}(\chi)$ (see lemma \ref{tautologic}). 
$G$-equivariance allows to define the filtration
$\overline{\fV}_{\rho}(\chi)=E\times_G \overline{\fW}_{\rho}(\chi)\to E(G/P)$.
Then $\fV_{\rho}(\sigma,\chi)=\sigma^*\overline{\fV}_{\rho}(\chi)$.

\subsubsectionr{}
\label{filtsect}
Conversely, take $g\in\Omega^0(E\times_{\Ad}\imag\klie)$.
Suppose that $\rho(g)$ has constant eigenvalues, and let
$\lambda_1<\dots<\lambda_r$ be the set of different values
they take. Just as before, we consider the filtration
\begin{equation}
0\subset V_{\rho}^{\lambda_1}\subset
V_{\rho}^{\lambda_2}\subset\dots\subset V_{\rho}^{\lambda_r}=V_{\rho}.
\label{asterix}
\end{equation}
Fix a point $x\in X$. After trivialising the fibre $E_x$ we
can identify $g(x)$ with and element $\chi$ of $\imag\klie$.
Let $P=P_{\rho}(\chi)$ (see lemma \ref{equi3}).
We obtain a reduction $\sigma\in\Gamma(E(G/P))$ as follows.
Let $y\in X$. Trivialise $E_y$ and identify $g(y)$ with
$\chi_y\in\imag\klie$. Let 
$$\sigma(y)=\{g\in G|g(\fW_{\rho}(\chi))=\fW_{\rho}(\chi_y)\}.$$
Then $\sigma(y)$ is invariant under left multiplication by
elements of $P$, and in fact gives a unique point in $G/P$
(here we use lemma \ref{equi3}). Furthermore, the definition of 
$\sigma(y)$ is compatible with change of trivialisation in the
sense that it gives a section $\sigma\in\Gamma(E(G/P))$.

\begin{lemma}
The filtration (\ref{asterix}) is equal to $\fV_{\rho}(\sigma,\chi)$.
\end{lemma}

\subsubsection{Holomorphic reductions of the structure group}
\label{holomred}
Suppose that there is a fixed (integrable) holomorphic structure on $E_G$.
This structure induces a holomorphic structure
on the total space of the associated bundle $E(G/P)$, since
$G/P$ is a complex manifold and the action of $G$ on $G/P$
is holomorphic. 

\begin{definition}
Let $\sigma\in \Gamma(E(G/P))$. A reduction $\sigma$
is {\bf holomorphic} if the map $\sigma:X\to E(G/P)$
is holomorphic.
\end{definition}

One can give an equivalent definition of holomorphicity in terms
of the filtrations induced by the reduction $\sigma$ in the 
associated vector bundles.

\begin{lemma}
Let $\sigma\in \Gamma(E(G/P))$. If the reduction $\sigma$
is holomorphic then, for any antidominant character $\chi$ for $P$ 
and for any representation $\rho:K\to U(W)$, 
the filtration $\fV_{\rho}(\sigma,\chi)$ of $V_{\rho}$ is holomorphic.
Conversely, let $g\in\Omega^0(E\times_{\Ad}\imag\klie)$ have constant 
eigenvalues, and let $P\subset G$, $\sigma\in\Gamma(E(G/P))$, 
$\chi\in\imag\klie$ and $\fV_{\rho}(\sigma,\chi)$ be obtained from it 
as in \ref{filtsect}. If $\fV_{\rho}(\sigma,\chi)$ is holomorphic,
then so is $\sigma$.
\end{lemma}
\begin{pf}         
Since $\fV_{\rho}(\sigma,\chi)=\sigma^*\tilde{\fV}_{\rho}(\chi)$
and $\tilde{\fV}_{\rho}(\chi)\to E(G/P)$ is holomorphic, the
first claim follows.

We now prove the second claim.
Suppose that $\fV_{\rho}(\sigma,\chi)$ is holomorphic. Fix $x\in X$
and take a holomorphic trivialisation $E|_U\simeq U\times G$
on a contractible neighbourhood $U$ of $x$. With this trivialisation,
the restriction of $\sigma$ to $U$ can be viewed as a map from
$U$ to $G/P$. Define a filtration $\fV^U$ of $U\times W_{\rho}$
as $\fV^U(x)=\sigma(x)\fW_{\rho}(\chi)$. Then 
$\fV_{\rho}(\sigma,\chi)|_U$ can be identified with $\fV^U$ and
the holomorphic structure on $V_{\rho}|_U$
corresponds to the trivial $\overline{\partial}$
operator on $U\times W_{\rho}$. Hence if $\fV_{\rho}(\sigma,\chi)|_U$
is holomorphic then $\overline{\partial}$ leaves $\fV^U$
invariant. Since $\rho$ is faithful this is equivalent to
$\overline{\partial}s=0$.
\end{pf}

\subsubsection{Total degree of a reduction of the structure group}
Let $P$ be a parabolic subgroup of $G$ with respect to 
$\{\alpha_{i_1},\dots,\alpha_{i_s}\}\subset \Delta$. Suppose that 
$\sigma\in \Gamma(E(G/P))$ is a reduction. Let $\chi$ be an antidominant
character of $P$.

We begin by defining the degree of the pair $(\sigma,\chi)$.
Let $0\subset V^{\lambda_1}\subset\dots\subset V^{\lambda_r}=V$ be
the filtration $\fV_{\rho_a}(\sigma,\chi)$ of $V$. 
For any vector bundle $V'$ we denote
$$\deg(V')={2\pi}\la c_1(V')\cup [\omega^{[n-1]}],[X]\ra.$$
Here $[\omega^{[n-1]}]$ denotes the cohomology class 
represented by the form $\omega^{[n-1]}$ and $[X]\in H_{2n}(X;\ZZ)$
is the fundamental class of $X$. Then we set
$$\deg(\sigma,\chi)=\lambda_r\deg(V)+
\sum_{k=1}^{r-1}(\lambda_k-\lambda_{k+1})\deg(V^{\lambda_k}).$$

\subsection{Stability, simple pairs and the correspondence}
\label{ssp}
Let $\sigma\in \Gamma(E(G/P))$ be a reduction.
We define the {\bf maximal weight} of $(\sigma,\chi)$ acting on a section
$\Phi\in\SSS=\Gamma(\FFF)$ of the associated bundle $\FFF=E\times_K F$ as
$$\int_{x\in X}\lambda(\Phi(x);-\imag g_{\sigma,\chi}(x)),$$
where $\lambda(\Phi(x);-g_{\sigma,\chi}(x))$ is the maximal
weight of $-g_{\sigma,\chi}(x)$ acting on $\Phi(x)$ as defined
in \ref{defpesmax} (note that here we use the $K$-equivariance of
the maximal weights, as stated in lemma \ref{proppes}).

Finally, given any central element $c\in \zlie\cap\klie$
we define the {\bf $c$-total degree} of the pair $(\sigma,\chi)$
as $$T^c_{\Phi}(\sigma,\chi)=\deg(\sigma,\chi)
+\int_{x\in X}\lambda(\Phi(x);-\imag g_{\sigma,\chi}(x))
+\la\imag\chi,c\ra\Vol(X).$$
Just as the maximal weights, the $c$-total degree is allowed
to be equal to $\infty$.

Now suppose that $X_0\subset X$ has as complement in $X$ a complex 
codimension 2 submanifold. Suppose also that a reduction $\sigma$
is defined only in $X_0$, that is, $\sigma\in\Gamma(X_0;E(G/P))$.
In this case it also makes sense to speak about 
$T^c_{\Phi}(\sigma,\chi)$ for any antidominant character $\chi$.
The only difficulty would be in defining the degree $\deg(\sigma,\chi)$.
However, it is well known that the degree of a vector bundle
can be computed by integrating the Chern-Weil form in the complement 
of a complex codimension 2 variety.

\begin{definition}
A pair $(A,\Phi)\in\AAA^{1,1}\times\SSS$ is {\bf $c$-stable} if
for any $X_0\subset X$ whose complement on $X$ is a complex 
codimension 2 submanifold, for any parabolic subgroup $P$ of $G$, 
for any holomorphic (with respect to the complex structure
$\Chern^{-1}A$ on $E_G$, see lemma \ref{integr}) reduction
$\sigma\in \Gamma(X_0;E(G/P))$ defined on $X_0$, and for any antidominant 
character $\chi$ of $P$ we have $$T^c_{\Phi}(\sigma,\chi)>0.$$
\label{parella_estable}
\end{definition}

We will say that an element $s\in\GGG_G$ is {\bf semisimple} if,
for any $x\in X$, after identifying $(E\times_{\Ad}\glie)_x\simeq \glie$,
$s(x)\in \glie$ is a semisimple element.
(This is independent of the chosen isomorphism  
$(E\times_{\Ad}\glie)_x\simeq \glie$,
because an element of $\glie$ is semisimple if and only if any element 
in its orbit by the adjoint action of $G$ on $\glie$ is semisimple.)

\begin{definition}
A pair $(A,\Phi)$ is {\bf simple} if no semisimple element in 
$\Lie(\GGG_G)$ leaves $(A,\Phi)$ fixed,
that is, for any semisimple $s\in\Lie(\GGG_G)$, 
$\fX^{\AAA\times\SSS}_s(A,\Phi)\neq 0$.
\label{defsimple}
\end{definition}

\begin{remark}
If $(A,\Phi)$ is simple then so is any point in
the $\GGG_G$ orbit through $(A,\Phi)$.
\end{remark}

We are now ready to state the main theorem of this chapter.

\begin{theorem}[Hitchin--Kobayashi correspondence]
Let $(A,\Phi)\in\AAA^{1,1}\times\SSS$ be a simple pair.
There exists a gauge transformation $g\in\GGG_G$ such that
\begin{equation}
\Lambda F_{g(A)}+\mu(g(\Phi))=c 
\label{hk}
\end{equation}
if and only if $(A,\Phi)$ is $c$-stable. Furthermore,
if two different $g,g'\in\GGG_G$ solve equation (\ref{hk}), then
there exists $k\in\GGG_K$ such that $g'=kg$.
\label{main}
\end{theorem}

We briefly explain the idea of the proof of theorem \ref{main}.
We construct on $\AAA^{1,1}\times\SSS\times\GGG_G$ a functional $\Psi$
(that we will call integral of the moment map) whose critical
points give the solutions of equation (\ref{hk}). 
We prove that the pair $(A,\Phi)$ is $c$-stable if and only if the 
functional $\Psi$ is, in a certain sense, proper along the slice 
$\{A\}\times\{\Phi\}\times\GGG_G$. On the other hand we prove that 
the functional being proper along $\{A\}\times\{\Phi\}\times\GGG_G$ is 
equivalent to its having a critical point in 
$\{A\}\times\{\Phi\}\times\GGG_G$, thus proving theorem \ref{main}.

Sections \ref{integral} to \ref{demostracio} are devoted to the proof
of theorem \ref{main}. In section \ref{integral} we explain how
to construct the functional $\Psi$ and prove some basic properties
of it. This is done for any Kaehler action of a Lie group 
(satisfying certain properties which do hold for compact groups and
also for the group $\GGG_K$) on a Kaehler manifold.
In section \ref{esKaehler} we prove that one can apply the results
in section \ref{integral} to the action of $\GGG_K$ on 
$\AAA^{1,1}\times\SSS$. More precisely, we define (using an idea of
Atiyah and Bott \cite{AB}) a Kaehler structure
on $\AAA^{1,1}\times\SSS$ which is respected by the action of 
$\GGG_K$ and such that the action of $\GGG_G$ is holomorphic.
In section \ref{correspondencia} we pause to look at the case
$X=\{\pt\}$ (see the beginning of that section for an explanation).
Finally, in section \ref{demostracio} we give the proof of
theorem \ref{main}.

\subsection{Bogomolov inequality}
In corollary \ref{bogom0} a certain inequality satisfied by all pairs
$(A,\Phi)$ solving equations \ref{equs} is given. Observe, however,
that when $F$ is Kaehler the inequality only depends on the $\GGG_G$
orbit of $(A,\Phi)\in\AAA^{1,1}\times\SSS$. Hence, we may restate
that result as follows, obtaining a necessary topological condition
for existence of solutions to equations (\ref{hk}).

\begin{corollary}
Suppose that a pair $(A,\Phi)\in\AAA^{1,1}\times\SSS$ 
satisfies $\ov{\partial}_A\Phi=0$ and
that there exists a gauge transformation $g\in\GGG_G$ such that 
\ref{hk} is satisfied. Then the following inequality holds 
$$\int_X\la\Lambda F_A,c\ra+
\int_X \Phi^*\phi_A(\overline{\omega}_F)\wedge \omega^{[n-1]}
-\frac{1}{2}\int_X B(F_A,F_A)\wedge\omega^{[n-2]}\geq 0.$$
\label{bogom}
\end{corollary}

\section{The integral of the moment map}
\label{integral}

In this section we consider the following general situation.
Let $H$ be a Lie group which acts on a Kaehler manifold
$M$ respecting the Kaehler structure, and assume that
there exists a moment map $\mu:M\to\hlie^*,$ where
$\hlie=\Lie(H)$. Suppose that there exists the complexification
$L=H^{\CC}$ of $H$, and that the inclusion $H\to L$
induces a surjection $\pi_1(H)\exh \pi_1(L)$. 
Under this assumptions, we construct a functional
$$\Psi:M\times L\to\RR$$
which we call the integral of the moment map $\mu$, and 
which satisfies these two properties:
\begin{itemize}
\item for any $x\in M$, the critical points of the restriction
$\Psi_x$ of $\Psi$ to $\{x\}\times L$ coincide with the points
of the orbit $Lx$ on which the moment map vanishes and
\item the restriction of $\Psi_x$ to lines of the form $\{e^{ts}|t\in\RR\}$,
where $s\in\llie=\Lie(L)$, is convex.
\end{itemize}
If $H$ is compact then $L=H^{\CC}$ always exists and $\pi_1(H)\exh\pi_1(L)$
is always satisfied. But note that we do not need our manifold $M$ or our 
groups $H,\ L$ to be finite dimensional. In fact, we will use this 
construction mainly in the infinite dimensional case
$(M;H,L)=(\AAA^{1,1}\times\SSS;\GGG_K,\GGG_G)$ (in section 
\ref{esKaehler} we will prove that $\AAA^{1,1}\times\SSS$
is a Kaehler manifold, that the action of $\GGG_K$ respects the
Kaehler structure, and we will identify a moment map for this
action). The resulting integral of the moment map will be  
a certain modification of the Donaldson functional, and will
be the key tool to prove theorem \ref{main}.

%
\subsection{Definition of $\Psi$}

Let us fix a point $x\in M$, and let $\phi:L\to M$ be the map which sends
$h\in L$ to $hx\in M$. We define a 1-form on $L$,
$\sigma=\sigma^x\in\Omega^1(L)$, as follows: given
$h\in L$ and $v\in T_hL$, 
$$\sigma_h(v)=\la\mu(hx),-\imag \pi(v)\ra_{\klie},$$
where $\pi:T_hL=\hlie \oplus\imag \hlie \to\imag \hlie $
is the projection to the second summand.

We will use the following formula, which holds for any two 
vector fields $X,Y$ and any 2-form $\omega$ on $M$ 
\begin{equation}
d\omega(X,Y)=L_X(\omega(Y))-L_Y(\omega(X))-\omega([X,Y]).
\label{der}
\end{equation}
Equality (\ref{der}) is a particular case of a formula which
describes the exterior derivative of forms of arbitrary degree
in terms of Lie derivatives (see \cite{BeGeV} p. 18).

\begin{lemma}
The 1-form $\sigma$ is exact.
\end{lemma}
\begin{pf}
Let us first of all prove that $d\sigma=0$. Given $g\in\llie=\Lie(L)$, let
$\fX_g=\fX_g^M$. We will prove that for any pair
$g,g'\in\hlie\cup\imag \hlie$, $d\sigma(g,g')=0$. This implies
by linearity that $d\sigma=0$. We will treat separately three cases,
and will make use of formula (\ref{der}), which in our case reads
$$d\sigma(\fX^L_g,\fX^L_{g'})=\la d(\sigma(\fX^L_{g'})),\fX^L_g\ra_{TL}
-\la d(\sigma(\fX^L_g)),\fX^L_{g'}\ra_{TL}-\sigma([\fX^L_g,\fX^L_{g'}]).$$

Suppose first that $g,g'\in\hlie$. In this case,
$\pi(\fX^L_g)=\pi(\fX^L_{g'})=\pi([\fX^L_g,\fX^L_{g'}])=0$, hence 
by the formula it is clear that $d\sigma(\fX^L_g,\fX^L_{g'})=0$.

Now suppose that $g\in\hlie$ and $g'\in\imag \hlie$. Observe that
$\sigma(\fX^L_g)=0$, so we have to prove that 
$\la d(\sigma(\fX^L_{g'})),\fX^L_g\ra_{TL}-\sigma([\fX^L_g,\fX^L_{g'}])=0$.
Differentiating (C2) in definition \ref{momentmap} 
we have $$\la d\la\mu,v\ra_{\hlie},\fX_g\ra_{TM}+\la\mu,[g,v]\ra_{\hlie}=0.$$
The functoriality of the exterior differentiation $d$
implies that 
$$\la d(\sigma(\fX^L_{g'})),\fX^L_g\ra_{TL}+\sigma(\fX^L_{[g,g']})=0.$$
On the other hand, since the action of $L$ on $M$ is on the left,
$[\fX^L_g,\fX^L_{g'}]=-\fX^L_{[g,g']}$, hence we 
obtain $$\la d(\sigma(\fX^L_{g'})),\fX^L_g\ra_{TL}-
\sigma([\fX^L_g,\fX^L_{g'}])=0,$$
which is what we wanted to prove.
The case $g\in\imag \hlie$ and $g'\in\hlie$ is dealt with in
a very similar way.

Finally, there remains the case $g,g'\in\imag \hlie$.
In this situation $[g,g']\in\hlie$, and so $\sigma([\fX^L_g,\fX^L_{g'}])=0$.
In view of this we have to prove 
$$\la d(\sigma(\fX^L_{g'})),\fX^L_g\ra_{TL}=
\la d(\sigma(\fX^L_g)),\fX^L_{g'}\ra_{TL}.$$
The left hand side is equal to 
$\phi^*(\la d\la\mu,\imag g\ra_{\hlie},\fX_{g'}
\ra_{TM})$ and this, by (C1) in definition \ref{momentmap}, is equal
to $$\phi^*(\omega_M(I\fX_g,\fX_{g'}))=\phi^*(-\la \fX_g,\fX_{g'}\ra),$$
where $\omega_M$ denotes the symplectic form on $M$.
The right hand side is equal to 
$$\phi^*(\omega_M(I\fX_{g'},\fX_g))=\phi^*(-\la \fX_{g'},\fX_g\ra).$$
Both functions are the same by the symmetry of $\la,\ra$.

Once we know that $d\sigma=0$, let us prove that $\sigma$ is exact.
Let $\iota:H\to L$ denote the inclusion. It is clear
that $\iota^*\sigma=0$. On the other hand, by our hypothesis
$\iota_*:\pi_1(H)\to\pi_1(L)$
is exhaustive. 
These two facts imply that $\sigma$ is exact.
Indeed, if it were not exact then we could find a path 
$\gamma:[0,1]\to L$, $\gamma(0)=\gamma(1)=1\in L$ such that 
$$\int_{\gamma}\sigma\neq 0.$$
But then we could deform $\gamma$ to a path $\gamma'\subset H$,
and, since $d\sigma=0$, the value of the integral would not change
and in particular would be nonzero. This is in contradiction with
the fact that $\iota^*\sigma=0$. So $\sigma$ is exact.
\end{pf}

Let $\Psi_x:L\to\RR$ be the unique function such that
$\Psi_x(1)=0$ and such that $d\Psi_x=\sigma^x$.
Define also $\Psi: M\times L \ni (x,g) \mapsto \Psi_x(g)$.
We will call the function $\Psi$ the {\bf integral of the moment map}.

\subsection{Properties of $\Psi$}
\label{propietats}
In this subsection we give the properties of the integral of the
moment map which will be used below.

\begin{prop} Let $x\in M$ be any point, and let $s\in\hlie$. 
\begin{enumerate}
\item $\Psi(x,e^{\imag s})=
\int_0^1 \la\mu(e^{\imag ts}x),s\ra_{\hlie} dt=
\int_0^1 \lambda_t(x;s) dt,$
\item $\frac{\partial\Psi}{\partial t}(x,e^{\imag ts})|_{t=0}
=\la\mu(x),s\ra_{\hlie}=\lambda_0(x;s),$
\item $\forall t_0\in\RR$,
$\frac{\partial^2\Psi}{\partial t^2}(x,e^{\imag ts})|_{t=t_0}\geq 0,$
with equality if and only if $\fX_s(e^{\imag t_0s}x)=0$,
\item $\forall t_0>0$, 
$\Psi(x,e^{\imag ls}x)\geq (l-t_0)\lambda_{t}(x;s)+C_s(x;t_0)$,
where $C_s(x;t_0)$ is a continuous function on $x\in M$,
$s\in\hlie$ and $t_0\in\RR$,
\end{enumerate}
\label{propietats1}
\end{prop}
\begin{pf}
By definition, $\Psi(x,e^{\imag s})=\int_{\gamma}\sigma^x$, where
$\gamma$ is any path in $L$ joining $1\in L$ to $e^{\imag s}$.
If we take $\gamma:[0,1]\ni t\mapsto e^{\imag ts}$, then
the integral reduces to $\int_0^1 \la\mu(e^{\imag ts}x),s\ra_{\hlie} dt$.
This proves (1). 
Property (2) is deduced from (1) differentiating. 
(3) is a consequence of (1) and the fact that $\lambda_t(x;s)$
increases with $t$.
To prove (4), let $C_s(x;t_0)=\int_0^{t_0}\lambda_t(x;s)dt$. Then:
$$\int_0^1 \lambda_t(x;ls)dt=\int_0^l \lambda_t(x;s)dt
\geq (l-t_0)\lambda_{t}(x;s)+C_s(x;t_0);$$
the first equality is obtained making a change of variable and 
using (2) in \ref{proppes}, and the 
inequality comes from the fact that $\lambda_t(x;s)$ increases as 
a function of $t$. 
\end{pf}

\begin{prop} Let $x\in M$ be any point, and let $s\in\hlie$. 
\begin{enumerate}
\item If $g,h\in L$, then
$\Psi(x,g)+\Psi(gx,h)=\Psi(x,hg)$,
\item for any $k\in H$ and $g\in L$, $\Psi(x,kg)=\Psi(x,g)$,
and $\Psi(x,1)=0$,
\item for any $k\in H$ and $g\in L$, $\Psi(kx,h)=\Psi(x,k^{-1}gk)$.
\end{enumerate}
\label{propietats2}
\end{prop}
\begin{pf}
To prove (1), observe that for any $g\in L$, $\sigma^{gx}=R_g^*\sigma^x$,
where $R_g$ denotes right multiplication in $L$ (indeed, 
for any $g'\in L$ one has $\sigma^{gx}(g')=\sigma^{x}(g'g)$
-- as usual, we identify the tangent spaces $T_{g'}(L)$
and $T_{g'g}(L)$ making $L$ act on the right).
This equivalence, together with the requierement that 
$\Psi_{gx}(1)=0$ implies that, for any $h\in L$, 
$\Psi_{gx}(h)=\Psi_x(hg)-\Psi_x(g)$.  
Property (2) is a consequence of (1) together with the fact that, for any 
$x\in M$, $\Psi_x|_H=0$. Finally, to prove (3) we use points (1) and (2):
$\Psi(x;k^{-1}gk)=\Psi(x,gk)+\Psi(gkx,k^{-1})
=\Psi(x,k)+\Psi(kx,g)=\Psi(kx,g).$
\end{pf}

\begin{prop} 
An element $g\in L$ is a critical point of $\Psi_x$ if and only 
if $\mu(gx)=0$.
\label{propietats3}
\end{prop}
\begin{pf} This is a consequence of (2) in \ref{propietats1}
and (1) in \ref{propietats2}.
\end{pf}

Just like maximal weights, the function $\Psi$ depends on the
moment map, which is not unique. When it is not clear from
the context which moment map we consider, we will write
$\Psi^{\mu}$ to mean the integral of the moment map $\mu$.

\subsection{Linear properness}
In this section we restrict to the case
$(M;H,L)=(F;K,G)$. In particular, recall that we have the
auxiliar representation $\rho_a:\glie\to\End(W_a)$
(see section \ref{repauxiliar}).
We define a norm on $\glie$ as follows: for any $s\in\glie$,
$$|s|=\la s,s\ra^{1/2}=\Tr(\rho_a(s)\rho_a(s)^*)^{1/2}.$$
Let $\log_G:G\simeq K\times\exp(\imag\klie)\to\imag\klie$
denote the projection to the second factor of the Cartan decomposition
composed with the logarithm. For any $g\in G$ we will call 
$|g|_{\log}:=|\log_G g|$ the {\bf length} of $g$.

\begin{definition}
We will say that $\Psi_x$ is {\bf linearly proper} if there exist
positive constants $C_1$ and $C_2$ such that for any $g\in G$    
$$|g|_{\log}\leq C_1\Psi_x(g)+C_2.$$
\end{definition}

\begin{prop}
Let $h\in G$ and $x\in F$. If $\Psi_x$ is linearly proper
then $\Psi_{hx}$ is also linearly proper.
\label{tambe}
\end{prop}

Before giving the proof of this proposition we prove the following
technical result.
\begin{lemma} Let $h\in G$. There exists $C\geq 1$ such that
for any $g\in G$
$$N^{-1/2}|gh|_{\log}-\log C
\leq |g|_{\log}\leq N^{1/2}(|gh|_{\log}+\log C).$$
Furthermore, $C$ depends continuously on $h\in G$.
\label{comparison}
\end{lemma}
\begin{pf}
Since the Cartan decomposition commutes with unitary representations,
we may describe the length function as follows. Let $x\in G$ be any element
and write $\rho_a(x)=RS$, where $R\in U(W_a)$ and $S=\exp(u)$,
where $u=u^*$. The matrix $u$ diagonalises and has real
eigenvalues $\lambda_1,\dots,\lambda_N$. So
$|x|_{\log}^2=\sum_{j=1}^N \lambda_j^2.$
Define $\max (x)=\max_{\|v\|=1}|\log \|\rho_a(x)v\||$. Then
we have $\max |\lambda_j|=\max (x)$ and consequently
\begin{equation}
\max(x)\leq |x|_{\log}\leq N^{1/2}\max(x).
\label{ee1}
\end{equation}
Let now $h\in G$. Then there exists $C\geq 1$, depending continuously
on $h$, such that for any $g\in G$ and any $v\in V$,
$C^{-1}\|\rho_a(gh)v\|\leq \|\rho_a(g)v\|\leq 
C\|\rho_a(gh)v\|$, which implies
\begin{equation}
|\max(gh)-\max(g)|\leq \log C.
\label{ee2}
\end{equation}
Putting $x=gh$ in (\ref{ee1}) we obtain
\begin{equation}
N^{-1/2}|gh|_{\log}\leq\max(gh)\leq|gh|_{\log},
\label{ee3}
\end{equation}
and combining (\ref{ee1}) with $x=g$ and (\ref{ee2}) we get
$$\max(gh)-\log C\leq |g|_{\log}\leq N^{1/2}(\max(gh)+\log C).$$
Finally, using (\ref{ee3}) we get
$N^{-1/2}|gh|_{\log}-\log C
\leq |g|_{\log}\leq N^{1/2}(|gh|_{\log}+\log C).$
\end{pf}

\noindent \begin{pf} (Proposition \ref{tambe}.)
Suppose that $\Psi_x$ is linearly proper, that is, for any $g\in G$
$$|g|_{\log}\leq C_1\Psi_x(g)+C_2,$$
where $C_1$ and $C_2$ are positive.
Fix $h\in G$. Let $C\geq 1$ be the constant in lemma \ref{comparison}.
(1) in \ref{propietats2} tells us that
$\Psi_{hx}(g)=\Psi_x(gh)-\Psi_x(h)$, so we get for any $g\in G$
\begin{align}
|g|_{\log} &\leq N^{1/2}(|gh|_{\log}+\log C)
\leq N^{1/2}(C_1\Psi_x(gh)+C_2+\log C)\notag\\
&=N^{1/2}(C_1(\Psi_x(gh)-\Psi_x(h))+C_1\Psi_x(h)+C_2+\log C)\notag\\
&=N^{1/2}(C_1\Psi_{hx}(g)+C_1\Psi_x(h)+C_2+\log C),\notag
\end{align}
so setting $C_1'=N^{1/2}C_1$ and 
$C_2'=\max\{0,N^{1/2}(C_1\Psi_x(h)+C_2+\log C)\}$ then
$C_1'$ and $C_2'$ are positive and 
$|g|_{\log}\leq C_1'\Psi_{hx}(g)+C_2'.$
This proves that $\Psi_{hx}$ is linearly proper.
\end{pf}

\section{The Kaehler structure on $\AAA\times\SSS$}
\label{esKaehler}
In this section we will give, following the classical idea of Atiyah and  
Bott \cite{AB}, a $\GGG_K$-invariant Kaehler structure on 
the manifold $\AAA\times\SSS$. We will identify for this
structure a moment map of the action of $\GGG_K$, the maximal weights 
and the integral of the moment map. 

Recall that the Lie algebras
of the gauge groups are $\Lie(\GGG_K) \simeq \Omega^0(E\times_{\Ad}\klie)$ 
and $\Lie(\GGG_G) \simeq \Omega^0(E\times_{\Ad}\glie)$.
On the other hand, the $K$-equivariance of the Cartan decomposition
implies that $\GGG_G\simeq \GGG_K\times\imag\Lie(\GGG_K)$
(the isomorphism being given by the map from $\GGG_K\times\imag\Lie(\GGG_K)$
to $\GGG_G$ which sends $(g,s)$ to $g\exp(s)$),
and from this fact  
we deduce that $\pi_1(\GGG_K)\to\pi_1(\GGG_G)$ is a surjection
(indeed, $\imag\Lie(\GGG_K)$ is contractible).
As a consequence, the results of section \ref{integral} apply to 
actions of $\GGG_K$ on Kaehler manifolds.
Hence, there is an integral of the moment map 
$\Psi:\AAA\times\SSS\times\GGG_G\to\RR$. This functional
will be the main tool in proving theorem \ref{main}.

\subsection{Unitary connections}
\label{exemp1}

\subsubsection{$\AAA$ is a Kaehler manifold}
Let $\AAA$ be the space of 
$K$-connections on $E$. This is an affine space modelled on 
$\Omega^1(E\times_{\Ad}\klie)$.
We define a complex structure $I_{\AAA}$ on $\AAA$ as follows.
Given any $A\in\AAA$, the tangent space $T_A\AAA$ can be 
canonically identified with 
$\Omega^1(E\times_{\Ad}\klie)=\Omega^0(T^*X\otimes E\times_{\Ad}\klie)$.
Then we set $I_{\AAA}=-I^*\otimes 1$. The complex structure
$I_{\AAA}$ is integrable.
We also define on $\AAA$ a symplectic form $\omega_{\AAA}$.
Let $\Lambda:\Omega^{p,q}(X)\to\Omega^{p-1,q-1}(X)$ be the adjoint of
the map given by wedging with $\omega$. Then, if $A\in\AAA$ and 
$\alpha,\beta\in T_A\AAA\simeq\Omega^1(E\times_{\Ad}\klie)$, we set
$$\omega_{\AAA}(\alpha,\beta)=\int_X \Lambda(B(\alpha,\beta)).$$

Here $B:\Omega^1(E\times_{\Ad}\klie)\otimes\Omega^1(E\times_{\Ad}\klie)
\to \Omega^2(X)$ is the combination of the usual wedge product with
our biinvariant nondegenerate pairing $\la,\ra$ on $\klie$.
It turns out that $\omega_{\AAA}$ is a symplectic form on $\AAA$, and it
is compatible with the complex structure $I_{\AAA}$. Hence $\AAA$ is
a Kaehler manifold. Furthermore, the action of $\GGG_G$ on 
$\AAA$ defined in subsection \ref{setting} is holomorphic and is the 
complexification of the action of $\GGG_K$.

\subsubsection{The moment map}
There exists a moment map for the action of $\GGG_K$ on $\AAA$, which
 takes the following form (see for example \cite{DoKr, Ko}):
$$\begin{array}{rcl}
\mu:\AAA & \longrightarrow & \Lie(\GGG_K)^* \\
A & \mapsto & \Lambda F_A.
\end{array}$$
Here $F_A$ denotes the curvature of $A$. It lies
in $\Omega^2(E\times_{\Ad}\klie)$, so 
$\Lambda F_A\in\Omega^0(E\times_{\Ad}\klie)\subset
\Omega^0(E\times_{\Ad}\klie)^*$, the last inclusion being
given by the integral on $X$ of the pairing $\la,\ra$ on $\klie$.

\subsubsection{Maximal weights}
\label{pesmaxconn}
In the following lemma we compute the $t$-maximal weights $\lambda(A;s)$
for $A\in \AAA$ and $s\in\Omega^0(E\times_{\Ad}\klie)$.
\begin{lemma}
Let $A\in\AAA$ be a connection, and take 
$s\in\Lie(\GGG_K)=\Omega^0(E\times_{\Ad}\klie)$. 
Then 
\begin{equation}
\lambda_t(A;s)=
\int_X\la\Lambda F_A,s\ra + 
\int_0^t \|e^{\imag ls}\overline{\partial}_A(s)e^{-\imag ls}\|^2 dl.
\end{equation}
\label{tpesmaxconn}
\end{lemma}
\begin{pf}                              
Let $\fX^{\AAA}_s\in \Gamma(T\AAA)$ be the
field generated by the action of $s$ on $\AAA$. 
In view of lemma \ref{gradient} we have 
$$\lambda_t(A;s)=\mu_s(A)+\int_0^t \|\fX^{\AAA}_s(e^{lI_{\AAA}s}A)\|^2 dl.$$
We make our computations in $\CCC$, which, as we have seen,
is isomorphic to $\AAA$ as a Kaehler manifold and on which the
action of $\GGG_G$ is easier to deal with.
So let $\overline{\partial}_A=\Chern^{-1}A$. 
By definition,
$e^{lI_{\AAA}s}A=\Chern(e^{lI_{\CCC}s}\overline{\partial}_A)$.
On the other hand, if $\fX^{\CCC}_s$ is the field generated 
by the action of $s$ on $\CCC$, we also have by definition
$\fX^{\AAA}_s(\Chern(\overline{\partial}_{A}))=
D\Chern(\fX^{\CCC}_s(\overline{\partial}_{A})).$
The map $\Chern$ is an isometry, so for any $\overline{\partial}_A\in\CCC$
we have $\|\fX^{\CCC}_s(\overline{\partial}_{A})\|^2=\|\fX^{\AAA}_s(A)\|^2.$
Finally, $\fX^{\CCC}_s(\overline{\partial}_{A})=
-\overline{\partial}_{A}(s)$.
Gathering all these facts together, we conclude that
(we use the $L^2$ norm on $\Omega^0(E\times_{\Ad}\glie)$
induced by the norm $|\cdot|$ on $\glie$):
\begin{align}
\lambda_t(A;s) &= 
\mu_s(A)+\int_0^t 
\|(e^{lI_{\CCC}s}\overline{\partial}_{A})(s)\|^2 dl 
=\mu_s(A)+\int_0^t 
\|(e^{\imag ls}\overline{\partial}_{A})(s)\|^2 dl \notag\\
&= \mu_s(A)+\int_0^t +
\|(e^{\imag ls}\circ\overline{\partial}_A\circ e^{-\imag ls})(s)\|^2 dl 
= \mu_s(A)+\int_0^t +
\|(e^{\imag ls}\circ\overline{\partial}_A)(s)\|^2 dl \notag \\
&= \int_X\la\Lambda F_A,s\ra + \int_0^t
\|e^{\imag ls}\overline{\partial}_A(s)e^{-\imag ls}\|^2 dl.
\label{pesmax}
\end{align}
We have used the fact that 
$e^{\imag ls}$ and $s$ commute -- the action of the gauge group
on $\Lie (\GGG_K)$ is by conjugation! 
\end{pf}

When $s\in L^2_1(E\times_{\Ad}\klie)$ and $A\in\AAA^{1,1}$ 
the maximal weight is given by exactly the same formula. But to prove it
one needs to use a technical theorem of Uhlenbeck and Yau \cite{UY}
which allows to regard $s$ as a genuine smooth section 
of $E\times_{\Ad}\klie$ at the complementary of a complex codimension
two subvariety of $X$, and to check that the integrals appearing
in lemma \ref{tpesmaxconn} converge.

\subsubsection{The integral of the moment map}

The results of section \ref{integral} apply in our case, so
there is an integral of the moment map 
$\Psi^{\AAA}$ which satisfies all the properties given in section
\ref{propietats}. Fix now a connection $A\in\AAA$. By the results of 
subsection \ref{pesmaxconn} and using (1) in proposition \ref{propietats1} 
we see that
\begin{align}
\Psi^{\AAA}_A(e^{\imag s})&=\int_0^1\lambda_t(A,s)
=\int_X\la\Lambda F_A,s\ra+                             
\int_0^1\left(
\int_0^t \|e^{\imag ls}\overline{\partial}_A(s)e^{-\imag ls}\|^2 dl 
\right)dt\notag\\
&=\int_X\la\Lambda F_A,s\ra+\int_0^1
(1-l)\|e^{\imag ls}\overline{\partial}_A(s)e^{-\imag ls}\|^2 dl.
\end{align}

Then, by (2) in \ref{propietats2}, the function $\Psi^{\AAA}_A$ 
factors through $$\Psi^{\AAA}_A:\GGG_G/\GGG_K\to\RR.$$
The resulting functional may be seen as a {\it modified
Donaldson functional}. In fact, when $F=\{\pt\}$,
it coincides (up to a multiplicative constant) with the
Donaldson functional. To see this, one only has to check 
that the Donaldson functional satisfies property (2) in 
\ref{propietats1} (see lemma 3.3.2 in \cite{Br2} for the case
$F=\CC^n$).

\subsubsection{Maximal weights for $A\in\AAA^{1,1}$}
\label{pesmaxconn2}
Note that since $\AAA^{1,1}\subset\AAA$ is a $\GGG_G$ invariant
subvariety (with singularities), the moment map, the maximal 
weights and the integral
of the moment map of the action of $\GGG_K$ on $\AAA^{1,1}$
are the restrictions of their counterparts in $\AAA$.

Recall that $V=E\times_{\rho_a}W_a\to X$ is the vector bundle associated to
the auxiliar representation $\rho_a$. For any $s\in\Lie(\GGG_K)$ we can
view $\rho_a(s)$ as a section of $E\times_{\Ad(\rho_a)}\End(W_a)$.
Take a connection $A\in\AAA^{1,1}$, and consider on $V$ the holomorphic
structure induced by $\overline{\partial}_A$. 
Using lemma \ref{tpesmaxconn} one can prove the following.

\begin{lemma}
Let $s\in\Lie(\GGG_K)$. If $\lambda(A;s)<\infty$, then 
the eigenvalues of $\rho_a(s)$ are constant.
Let $\lambda_1<\dots<\lambda_r$ be the different eigenvalues
of $\imag\rho_a(s)$,
and let $V(\lambda_j)\subset V$ be the eigenbundle of eigenvalue
$\lambda_j$. Put
$V^{\lambda_k}=\bigoplus_{j\leq k} V(\lambda_j)$.
Then, for any $k$, $V^{\lambda_k}$ is a holomorphic subbundle
of $V$. Furthermore
$$\lambda(A;s)=
\lambda_r\deg(V)+\sum_{k=1}^{r-1}(\lambda_k-\lambda_{k+1})
\deg(V^{\lambda_k}).$$
\label{redhol}
\end{lemma}
\begin{pf}
Suppose that $\lambda(A;s)<\infty$. 
Using the previous lemma with $U=\imag\rho_a(s)$ and 
$V=\overline{\partial}_A(\rho_a(s))$, we
get for any $k\geq 1$ and $t\geq 0$ 
\begin{align}
\int_X |\Tr(\rho_a(s)^k\overline{\partial}_A(\rho_a(s)))|
& \leq \int_X \|s\|^k 
\|e^{\imag ts}\overline{\partial}_A(s)e^{-\imag ts}\|
\notag \\
& \leq \left(\int_X \|s\|^{2k}\right)^{1/2}
\left(\int_X
\|e^{\imag ts}\overline{\partial}_A(s)e^{-\imag ts}\|^2\right)^{1/2}
\notag \\
&\leq \|s\|^k_{L^{2k}}
\|e^{\imag ts}\overline{\partial}_A(s)e^{-\imag ts}\|.
\notag
\end{align}
This, together with formula (\ref{pesmax}) implies that 
$\int_X |\Tr(\rho_a(s)^k\overline{\partial}_A(\rho_a(s)))|=0$, so
$$\Tr(\rho_a(s)^k\overline{\partial}_A(\rho_a(s)))=0.$$
On the other hand, for any $p+q=k$, $\Tr(U^kV)=Tr(U^pVU^q)$, so 
$$\Tr(\rho_a(s)^p\overline{\partial}_A(\rho_a(s))\rho_a(s)^q)=0$$ 
as well.
Finally, 
$$\overline{\partial} \Tr(\rho_a(s)^{k+1})=
\Tr(\overline{\partial}_A(\rho_a(s)^{k+1}))=
\Tr(\sum_{p+q=k} \rho_a(s)^p\overline{\partial}_A(\rho_a(s))\rho_a(s)^q)=0.$$
Since $X$ is compact this implies that $\Tr(\rho_a(s)^{k+1})$ is constant
for any $k$. Making $k=1,\dots,n$ we see that the eigenvalues of
$s$ must be constant.
                              
We now prove the second claim.
Using the splitting $V=V(\lambda_1)\oplus\dots\oplus V(\lambda_r)$ 
we can write
$$\overline{\partial}_A=
\left(\begin{array}{cccc}
\overline{\partial}_{A_1} & A_{12} & \dots & A_{1r} \\
A_{21} & \overline{\partial}_{A_2} & \dots & A_{2r} \\
\vdots & \vdots & \ddots & \vdots \\
A_{r1} & A_{r2} & \dots  & \overline{\partial}{A_r}
\end{array}\right),$$
where $A_{ij}\in\Omega^{0,1}(V(\lambda_i)\otimes V(\lambda_j)^*)$.
Let $\pi_k:V\simeq
V^{\lambda_k}\oplus V(\lambda_{k+1})\oplus\dots\oplus
V(\lambda_r)\to V^{\lambda_k}$ denote the projection. Then we can write 
$$u=-\imag\rho_a(s)=
\left(\begin{array}{cccc}
\lambda_1 & 0 & \dots & 0 \\
0 & \lambda_2 & \dots & 0 \\
\vdots & \vdots & \ddots & \vdots \\
0 & 0 & \dots & \lambda_r
\end{array}\right)=
\lambda_r\Id+\sum_{k=1}^{r-1}(\lambda_k-\lambda_{k+1})\pi_k.$$
We compute
$$\overline{\partial}_A(u)=\overline{\partial}_A\circ u
-u\circ\overline{\partial}_A=
((\lambda_j-\lambda_i)A_{ij})_{1\leq i,j\leq r},$$
where $i$ denotes the row and $j$ the column. On the other hand,
\begin{align}
\|e^{\imag ls}\overline{\partial}_A(s)e^{-\imag ls}\|^2
&=\|e^{\imag ls}\overline{\partial}_A(u)e^{-\imag ls}\|^2 \notag \\
&=\left\|((\lambda_j-\lambda_i)e^{l(\lambda_i-\lambda_j)}
A_{ij})_{1\leq i,j\leq r} \right\|^2\notag \\
&=\sum_{1\leq i\neq j\leq r} (\lambda_j-\lambda_i)^2
e^{2l(\lambda_i-\lambda_j)}\|A_{ij}\|^2.\notag
\end{align}
From the fact that
$$\int_0^{\infty} 
\|e^{\imag ls}\overline{\partial}_A(s)e^{-\imag ls}\|^2 dl$$
is finite, we deduce that $A_{ij}=0$ for $i>j$. So $V^{\lambda_k}$ is
holomorphic for any $k$. 
We also deduce that $e^{-\imag ls}\overline{\partial}_A(s)e^{\imag ls}$
converges in the $C^{\infty}$ norm as $s\to\infty$ to
$$\overline{\partial}^{\infty}_A=
\left(\begin{array}{cccc}
\overline{\partial}_{A_1} & 0 & \dots & 0 \\
0 & \overline{\partial}_{A_2} & \dots & 0 \\
\vdots & \vdots & \ddots & \vdots \\
0 & 0 & \dots & \overline{\partial}_{A_r}
\end{array}\right),$$
and this implies that the maximal weight $\lambda(A;s)$
is equal to 
$$\int_X\la\Lambda F_{A^{\infty}},s\ra=
\int_X\Tr(\rho_a(\Lambda F_{A^\infty})\rho_a(s)^*)=
\int_X\Tr(\rho_a(\imag\Lambda F_{A^\infty}) u),$$ 
where $A^{\infty}$ is $\Chern(\overline{\partial}^{\infty}_A)$.
Let also $A_k=\Chern(\overline{\partial}_{A_k})$. Then we
use the formula $\int_X \imag\Tr\Lambda F_{A_k}=\deg(V(\lambda_k))$
to deduce that
$$\lambda(A;s)=\lambda_r\deg(V)+\sum_{k=1}^{r-1}(\lambda_k-\lambda_{k+1})
\deg(V^{\lambda_k}).$$
This finishes the proof.
\end{pf}

\subsubsectionr{}
\label{fipesmaxconn}
If we consider more generally $s\in L^2_1(E\times_{\Ad}\klie)$,
then $\lambda(A;s)<\infty$ leads to a filtration of the locally
free sheaf associated to $V$ by reflexive (coherent) subsheaves,
and not only holomorphic subbundles of $V$ as in the smooth
case. To prove this one uses a theorem of Uhlenbeck and Yau (see \cite{UY}
and section 3.11 in \cite{Br2}).

\subsection{Sections of the associated bundle}
\label{exemp2}

\subsubsection{$\SSS$ is a Kaehler manifold}
Here we will define a symplectic form $\omega_{\SSS}$ and a 
compatible complex structure $I_{\SSS}$ on $T\SSS$,
and we will prove that both are integrable. To do that, consider a section
$\sigma\in\SSS$. Then $T_{\sigma}\SSS=\Gamma(\sigma^* T\FFF_v)$,
where $T\FFF_v\subset T\FFF$ is the subbundle of vertical tangent vectors
of $\FFF$, that is, $T\FFF_v=\Ker(D\pi_F)$. To define the complex
structure, let $\alpha\in\Gamma(\sigma^* T\FFF_v)$. Then 
$I_{\SSS}(\alpha)=I_F\alpha$. This makes sense, since the $K$ invariance
of $I_F$ implies that $T\FFF_v$ inherits the complex structure of $F$.
Now let $\alpha,\beta\in\Gamma(\sigma^* T\FFF_v)$. We
define the symplectic form $\omega_{\SSS}$ as                             
$$\omega_{\SSS}(\alpha,\beta)=\int_X \omega_F(\alpha,\beta).$$
Two things are clear: $\omega_{\SSS}$ is nondegenerate (this 
is a consequence of the nondegeneracy of $\omega_F$) and 
$\omega_{\SSS}$ and $I_{\SSS}$ are compatible, that is,
$\la\alpha,\beta\ra=\omega_{\SSS}(\alpha,I_{\SSS}\beta)$
is a Riemannian pairing. We have to prove that the two
structures are integrable.

Consider the complex structure $I_{\SSS}$. 
First of all, observe that $\FFF$ is a complex manifold.
It is well known that this implies that $\Map(X,\FFF)$ 
is also a complex manifold, with the complex structure induced
by that of $\FFF$. Since
$$\SSS=\{\phi\in\Map(X,\FFF)|\pi_F\circ\phi=\Id\}$$
and the equation $\pi_F\circ\phi=\Id$ is complex,
the set $\SSS$ is complex, considering the restriction
of the complex structure of $\Map(X,\FFF)$. But this restriction
is equal to $I_{\SSS}$.

Let us show now that the 2-form $\omega_{\SSS}$ is closed.
Fix a section $\sigma\in\SSS$. Let $\{U_{\alpha}\}_{\alpha\in A}$ be a
finite covering of $X$ trivialising $E$, with transition functions
$\{\phi_{\alpha\beta}:U_{\alpha}\cap U_{\beta}\to K\}$.
Take a partition of unity $\psi_{\alpha}$ subordinated to 
the covering. The section $\sigma$ translates into a family of sections
$\sigma_{\alpha}:U_{\alpha}\to F$ satisfying the compatibility
condition $\sigma_{\beta}=\phi_{\alpha\beta}^{-1}\sigma_{\alpha}$.
Using Darboux theorem, and possibly refining the covering 
$\{U_{\alpha}\}_{\alpha\in A}$, we can assume that for any
$\alpha$ there exists an open subset $V_{\alpha}\subset F$
symplectomorphic to a neighbourhood of zero of $\RR^{2m}$
with the standard symplectic structure 
$\omega_0=\sum_{i=1}^m dx_i\wedge dx_{i+m}$
and such that $\sigma_{\alpha}(\overline{U_{\alpha}})\subset V_{\alpha}$.
In view of this it is a trivial fact that on 
$\Map(\overline{U_{\alpha}},V_{\alpha})$ the form
$$\omega_{\alpha}=\int_{\overline{U_{\alpha}}}\psi_{\alpha}\omega_0$$
is closed (here we consider the closure $\overline{U_{\alpha}}$
of $U_{\alpha}$ to avoid problems with convergence).
So on $$\prod_{\alpha\in A}\Map(\overline{U_{\alpha}},V_{\alpha})$$ the 
form $\omega_A=\sum \pi^*_{\alpha}\omega_{\alpha}$ is closed ($\pi_{\alpha}$
denotes the projection to the factor 
$\Map(\overline{U_{\alpha}},V_{\alpha})$).
But we can see a neigbourhood $\SSS_{\sigma}$ 
of $\sigma\in\SSS$ as a submanifold of it, namely,
$$\SSS_{\sigma}=\{(\sigma_{\alpha})\in
\prod\Map(\overline{U_{\alpha}},V_{\alpha})|
\sigma_{\beta}=\phi_{\alpha\beta}^{-1}\sigma_{\alpha}
\mbox{ for any $\alpha,\beta\in A$}\}.$$
The form $\omega_{\SSS}$ restricts on $\SSS_{\sigma}$ precisely
to $\omega_A$, which is closed as we have seen. This is true in 
a neighbourhood of $\sigma$ for any $\sigma\in\SSS$, so 
definitively $\omega_{\SSS}$ is closed.

\subsubsection{The actions of $\GGG_K$ and $\GGG_G$ and the moment map} 
Both groups $\GGG_K$ and $\GGG_G$ act on the space
of sections $\SSS=\Gamma(\FFF)$, and the action of
$\GGG_G$ is the complexification of the action of $\GGG_K$.
On the other hand, $\GGG_K$ acts by isometries and respecting the 
symplectic form, and there exists a moment map $\mu_{\SSS}$, which 
is equal fibrewise to $\mu$ (the moment map of the action of $K$ on $F$). 
As such, it is a section of $\Omega^0(E\times_{\Ad}\klie)^*$.

\subsubsection{Maximal weights}
\label{pesmaxsec}
The maximal weight of $s\in\Lie(\GGG_K)=\Omega^0(E\times_{\Ad}\klie)$
acting on a section $\Phi\in\SSS$ is given by the integral of
the maximal weight in each fibre:
$$\int_{x\in X}\lambda(\Phi(x);s(x)).$$
This makes sense due to the $K$ equivariance of $\lambda$
(see (1) in lemma \ref{proppes}).

\subsubsection{The integral of the moment map} Once more, the results
in section \ref{integral} imply that there exists an integral
$\Psi^{\SSS}$ of the moment map of the action of $\GGG_K$ on $\SSS$. If 
$\Psi:F\times G\to\RR$ is the integral of the moment map of the
action of $K$ on $F$, then, for any section $\sigma\in\SSS$ and
gauge transformation $g\in\GGG_G$
$$\Psi^{\SSS}(\sigma,g)=\int_{x\in X} \Psi(\sigma(x),g(x)).$$
This makes sense due to the $K$-equivariance of $\Psi$:
see (3) in \ref{propietats2}. 

\subsection{Symplectic point of view}
\label{disgressio}

We saw that both $\AAA^{1,1}$ and 
and $\SSS$ are Kaehler manifolds, with symplectic forms $\omega_{\AAA}$
and $\omega_{\SSS}$ and with actions of $\GGG_K$ extending
to actions of the complexification $\GGG_G$. Hence
$\AAA^{1,1}\times\SSS$ is also a Kaehler manifold, with symplectic
form $\omega_{\AAA}+\omega_{\SSS}$ (we omit the pullbacks). The moment map
$\mu_{\AAA\times\SSS}$
of the action of $\GGG_K$ on $\AAA\times\SSS$ will simply be the moment
map of the action on $\AAA$ plus that of the action on $\SSS$.
That is, $$\mu_{\AAA\times\SSS}(A,\Phi)=\Lambda F_A+\mu(\Phi).$$
So equation (\ref{hk}) can be written as 
$\mu_{\AAA\times\SSS}=c$, where $c$ denotes the central element
in $(\Lie(\GGG_K))^*=\Omega^0(E\times_{\Ad}\klie)^*$ which is
fibrewise equal to a central element $c\in\klie^*$. 
Furthermore, we have the following result.
\begin{lemma} $T^c_{\Phi}(\sigma,\chi)=
\lambda^{\Lambda F_A+\mu(\Phi)-c}((A,\Phi);-\imag g_{\sigma,\chi}).$
\label{equivalencia}
\end{lemma}
\begin{pf}
Combine subsections \ref{pesmaxconn2} and \ref{pesmaxsec}.
\end{pf}

\section[Analytic stability]
{Analytic stability and vanishing of the moment map in finite dimension}
\label{correspondencia}
In this section we will pause to prove theorem \ref{main} in the case
$X=\{\pt\}$. This is done for two reasons. First of all, this 
particular case has some interest {\it per se} and its proof is 
considerably easier that that of the general case (specially because 
there is no connection and the analysis is elementary). The second reason 
is that theorem \ref{main} can be viewed as an infinite dimensional 
generalisation of the result in this section. 

The results in this section (at least for the case in which 
$F$ is projective) have been known for many years: see \cite{KeNe, Ki}.
That they are related with Hitchin--Kobayashi correspondence
was also known since the first cases of the correspondence were
studied. Our intention here is to make more concrete this relation and 
to stress on the similarities between the
{\it finite dimensional situation} $X=\{\pt\}$ and the general one considered
in theorem \ref{main} (which corresponds to the situation in which
$F=\AAA^{1,1}\times\SSS$ with the actions of $\GGG_K$ and $\GGG_G$). 
For example, the different versions of
Donaldson functional used in the literature are in fact particular
instances of a construction which works for a wide class of Kaehler
actions of Lie groups on Kaehler manifolds (namely, what we have
called the integral of the moment map). Moreover, the $c$-stability
condition is also a particular case of a general notion of stability
for group actions on Kaehler manifolds (the so-called analytic
stability). And the very correspondence coincides almost word by word
with theorem \ref{corr} given in this section. The proof which we
give here works only for Kaehler actions of compact groups,
and so it can not be used in the general situation (in which the group
is $\GGG_K$). Nevertheless, the scheme of the proof will be the same
in the general situation.

Let us write $\Psi:F\times G\to\RR$ for the integral of the moment
map $\mu:F\to\klie^*$.

\begin{definition}
Let $x\in F$. We will say that $x$ is {\bf analytically stable} 
if for any $s\in\klie$ the maximal weight of $s$ acting on $x$ 
is strictly positive: $$\lambda(x;s)>0.$$
\label{estabilitat}
\end{definition}

\begin{lemma} 
A point $x\in F$ is analytically stable if and only if 
$\Psi_x$ is linearly proper.
\label{fita}
\end{lemma}
\begin{pf}
Suppose first that $x$ is analytically stable.
We have to prove that there exist two positive constants 
$C_1, C_2\in\RR$ such that,
for any $s\in\klie$, $\|s\|\leq C_1\Psi_x(e^{\imag s})+C_2$.
Assume that there are not such constants. Then, we can find sequences
$\{s_j\}\subset \klie$ and $\{C_j\}\subset\RR$ such that 
$\|s_j\|\to\infty$, $C_j\to\infty$ and, for any $j$, 
$\|s_j\|\geq C_j\Psi_x(e^{\imag s_j})$. Let $u_j=s_j/\|s_j\|$.
After passing to a subsequence, we can assume that 
$u_j\to s$. Take now any $t>0$. By our hypothesis, and making use
of (4) in proposition \ref{propietats1}, 
$$\frac{1}{C_j}\geq\frac{\Psi_x(e^{\imag s_j})}{\|s_j\|}
\geq \frac{(\|s_j\|-t)}{\|s_j\|}
\lambda_{t}(x;u_j)+\frac{C_{u_j}(x;t)}{\|s_j\|}.$$
Now, making $j\to\infty$, we obtain $0\geq \lambda_{t}(x;s)$,
since, by the compactness of $B_{\klie}(1)=\{s\in\klie|\ \|s\|=1\}$,
$C_{u_j}(x;t)$ is uniformly bounded.
This is true for any $t>0$, so passing to the limit $t\to\infty$
we get $$0\geq \lambda(x;s),$$
which contradicts analytic stability.

Now suppose that there exist positive $C_1$ and $C_2$ such that
for any $s\in\klie$ 
\begin{equation}
\|s\|\leq C_1\Psi_x(e^{\imag s})+C_2.
\label{supo}
\end{equation}
We have to prove that $x$ is analytically stable.
So take $s\in\klie$ and assume that $\lambda(x;s)\leq 0$.
In this case, for any $t\geq 0$, 
$\Psi_x(e^{\imag ts})=\int_0^t\lambda_l(x;u)dl\leq 0$,
which, for $t$ big enough, contradicts (\ref{supo}).
This proves that $x$ is analytically stable.
\end{pf}

\begin{corollary}
Let $x\in F$. Then $x$ is analytically stable if and only if $hx$ is 
analytically stable for any $h\in G$.
\label{tambe2}
\end{corollary}
\begin{pf}
This is a consequence of the preceeding lemma together
with lemma \ref{tambe}.
\end{pf}

\begin{theorem}
Let $x\in F$ be any point. There is at most one $K$ orbit inside
the orbit $Gx\subset F$ on which the moment map vanishes. 
Furthermore, $x$ is analytically stable if and only if:
(1) the stabiliser $G_x$ of $x$ in $G$ is finite and (2) there exists
a $K$ orbit inside $Gx$ on which the moment map vanishes.
\label{corr}
\end{theorem}
\begin{pf}
We first prove uniqueness.
Assume that there are two different $K$ orbits inside a $G$ orbit
on which the moment map vanishes: say, $Kx$ and $Kgx$, where $g\in G$. 
By the polar decomposition we can assume that $g=e^{\imag s}$, 
where $s\in\klie$. Consider the function $\Psi_x:G\to\RR$. 
By proposition \ref{propietats3}, since $\mu(x)=0$, 
both $1,g\in G$ are critical points of $\Psi_x$. Consider now
the path $\gamma(t)=e^{\imag ts}$ connecting $1$ and $g$.
(3) of the proposition tells us that the restriction $\psi$
of $\Psi_x$ to this path has second derivative $\geq 0$. 
Since $0$ and $1$ are critical points of $\psi$, 
the second derivative must vanish at any point between $0$ and $1$.
In particular, 
$\frac{\partial^2\Psi}{\partial t^2}(x,e^{\imag ts})|_{t=0}=0$;
but this implies (again, (3) of the proposition), that the vector field 
$\fX_s(x)=0$, which gives $\fX_{\imag s}(x)=I\fX_s(x)=0$. 
So $e^gx=e^{\imag s}x=x$, and the two orbits $Kgx$ and $Kx$ coincide.

Suppose now that the point $x$ is analytically stable. 
Let us see that there is a $K$ orbit inside $Gx$ on which $\mu$ vanishes. 
By lemma \ref{fita}, the function $\Psi_x$ is linearly proper.
Using (2) in \ref{propietats2}, we
conclude that there must exist a critical point in the $G$
orbit of $x$. Indeed, if $\{s_j\}\subset\klie$ are such that
$e^{\imag s_j}$ is a minimising sequence for $\Psi_x$, then by
the preceeding lemma the set $\{s_j\}$ is bounded; so it
has a subsequence converging to a certain $s\in\klie$, and
$e^{\imag s}$ is a minimum of $\Psi_x$ (of course, here we
use that $\klie$ has finite dimension). At this point (even 
more, at the $K$ orbit through this point) the moment map must vanish.
Let now $y=e^{\imag s}x$. By lemma \ref{tambe2} $y$ is analitically stable.
If the stabiliser $K_y$ of $y$ in $K$ were not finite, 
then, since $K$ is compact, its closure would be a Lie 
subgroup of $K$ of dimension greater than zero. In particular, there 
would exist an $s\in\klie$ such that $\fX_s(y)=0$. But then $e^{ts}y=y$ 
for any $t$, so that the gradient flow $\phi^t_s$ leaves $y$ fixed.
This means that $\lambda(y;s)=-\lambda(y;-s)$, so that either $\lambda(y;s)$
or $\lambda(y;-s)$ (or both) is $\leq 0$. This contradicts 
analytic stability. So $K_y$ is finite. 

Finally, since $\mu(y)$ is invariant under the coadjoint action of 
$K$ in $\klie^*$, it turns out that $G_y$ is the complexification of $K_y$. 
Let us see why (we copy the proof of proposition 1.6 in \cite{Sj}).
One inclusion is easy: $G_y$ contains the complexification of $K_y$.
For the other inclusion,
let $ge^{\imag s}$ be an arbitrary element of $G_x$, where
$g\in K$ and $s\in\klie$. We want to show that $g\in K_x$ and
$s\in\klie_x$ (where $\klie_x$ is the infinitesimal stabliser
of $x$). Using the fact that $\mu$ is $K$-equivariant we have
$$\mu(e^{\imag s}x)=g^{-1}\mu(ge^{\imag s}x)=g^{-1}\mu(x)=\mu(x).$$
Now, lemma \ref{gradient} implies that $s\in\klie_x$, from which we
deduce that $g\in K_x$. This finishes the proof.
So $G_y$ is finite and in consequence $G_x$ is also finite.

To prove the converse, let $x\in F$.
Assume that $G_x$ is finite and that there exists $g\in G$ such 
that $\mu(gx)=0$. Then $G_{gx}$ is also finite and consequently
so is $K_{gx}$. This implies that, for any $s\in\klie$,
$\fX_{\imag s}(gx)\neq 0$, so (lemma \ref{gradient}),
$\lambda(gx;s)>\mu_s(gx)=0$. This means that $gx$ is analytically
stable, hence so is $x$.
\end{pf}

It is an exercise to verify that the property
on analitically stable points of $F$ of being simple 
(see subsection \ref{ssp}) is equivalent to that of having finite 
stabiliser in $G$.

Using the results in this section one can also study the equation
$\mu=c$, where $c\in\klie^*$ is any central element. 
Indeed, $\mu-c$ is a moment map, and so one only has to consider
the maximal weights $\lambda^{\mu-c}$ and the integral 
$\Psi^{\mu-c}$.

\subsection{Kempf-Ness theory}
\label{KN}
Suppose now that $F$ is a projective variety, with polarisation
$\OOO_F(1)\to F$ and such that the action of $G$ on $F$ lifts to an
action on $\OOO_F(1)$. This implies that the action of $G$ on
$F$ extends to an action on the projective space $\PP(W)$ and that
this actions linearises to an action of $G$ on $W$. The following 
definition is due to Mumford:

\begin{definition}
Let $y\in F$ be any point. We will say that
$y$ is {\bf stable} if its stabiliser in $G$ is finite and 
there exists an integer $n\geq 1$ and a $G$ invariant section
$s$ of $\OOO_F(n)$ such that $F_s=\{y'\in F | s(y')\neq 0\}\subset F$ 
is affine, contains $y$ and all the orbits of $G$ in $F_s$ are closed.
\end{definition}

The relevance of this definition comes from this fact.
While in general it is not possible to give an algebraic structure
to the set of orbits $F/G$, if we restrict ourselves to the
set $F^s$ of stable points, then we can give $F^s/G$ a very
natural algebraic structure. (This is the content of Geometric
Invariant Theory; see \cite{MFK}.)

The main point of Kempf-Ness theory is that the condition of
stability defined by Mumford coincides with the condition 
of analytic stability. The link between both definitions
is given by the Hilbert-Mumford numerical criterion, which
allows to decide whether $y\in F$ is stable.

\begin{definition}
A {\bf one parameter subgroups of $G$} (1-PS for short) is a
morphism $\alpha:\CC^*\to G$.
\end{definition}

\begin{lemma}[Hilbert-Mumford]
The point $y\in F$ is stable if and only if
for any 1-PS $\alpha$ and any lift $\hat{y}\in W$,
there exists a weight of $\alpha$ in $\hat{y}$ which is
$>0$, that is: if $\hat{y}=\sum_{n\in\ZZ}y_n$,
where $\alpha(t)y_n=t^ny_n$, there is an integer 
$n>0$ such that $y_n\neq 0$.
\end{lemma}

In (\cite{Ki}, page 107) it is proved that in Hilbert-Mumford's
criterion one only
needs to consider 1-PS which are {\bf compatible with $K$}.
These are the 1-PS which are obtained after complexifying
any group morphism $\alpha_K:S^1\to K$. Such an $\alpha_K$ is completely
determined by its differential at the identity, say $s_{\alpha}\in\klie$.
We will call the elements of the form $s_{\alpha}\in\klie$ 
{\bf integral weights}, and we will write $\klie_{\QQ}$
the set of integral weights. Observe that if $s\in\klie$ is an integral
weight, then the maximal weight of the one parameter subgroup
$\exp(s)$ acting on $y$ is equal to $\lambda(y;-s)$
(see lemma \ref{pesmaxproj}).

\begin{lemma}
Following the notations above, let $x\in\PP(W)$, and suppose that,
for any integral weight $s\in\klie_{\QQ}$, $\lambda(x;s)>0$.
Then, the same inequality is satisfied by any $s\in\klie$.
\end{lemma}
\begin{pf}         
Suppose that the hypothesis of the lemma hold, and take any
$s\in\klie$. The closure 
$$T=\overline{\{\exp(ts)|t\in\RR\}}\subset K$$ 
is a torus. Let $\tlie\subset\klie$ be
its Lie algebra, and let $\tlie_{\CC}=\tlie\otimes_{\RR}\CC$ be
the complexification of $\tlie$.
Then there is a free $\ZZ$-module $\Lambda\subset\tlie$
such that $T=\tlie/\Lambda$. Since $T$ is compact, $\RR\Lambda=\tlie$.
Let $s'\in\QQ\Lambda$. Take $n\in\NN$ such that $ns'\in\Lambda$.
All the elements in $\Lambda$ are integer weights, so
$\lambda(x;s')=\frac{1}{n}\lambda(x;ns')>0$.

Decompose $W=\bigoplus W_{\chi}$, where $\chi\in\tlie_{\CC}^*$
are characters of $T$, in such a way that any $e\in\tlie_{\CC}$
acts on $W_{\chi}$ multiplying by $\chi(e)$. Suppose that
$\lambda(x;s)\leq 0$. This means that if $\chi(s)>0$, then the component
of $x$ in $W_{\chi}$ vanishes: $x_{\chi}=0$. Now, taking into 
account that the weights $\chi$ that appear in the decomposition of
$W$ take rational values when evaluated on $\imag(\QQ\Lambda)$,
we deduce that we can approximate $s$ by an element $s'\in\QQ\Lambda$
such that $\lambda(x;s')\leq 0$. But this is not possible in view
of what we saw in the preceeding paragraph. Therefore, 
for any $s\in\imag\klie$ one has $\lambda(x;s)>0$.
\end{pf}

Together with theorem \ref{corr} this proves the following

\begin{theorem}
A point $x\in F$ is stable in the sense of Mumford if and only
its stabiliser is finite and there exists $g\in G$ such
that $\mu(gx)=0$. 
\end{theorem}

Of course, we have only proved this when the symplectic structure
of $F$ is that induced by the Fubini-Study symplectic form 
through the embedding given by $\OOO_F(1)$. In \cite{Sj}
a much stronger result is proved, which is true even if the
symplectic form in $F$ is different from the one induced
by Fubini-Study.

\subsection{The general case}
In view of lemma \ref{equivalencia}, 
if the results in this section were 
valid for infinite dimensional Lie groups, then theorem \ref{main} would
follow from it. To the best of the author's knowledge, there is no
general result as theorem \ref{corr} valid in infinite dimensions.
However, although the proof of theorem \ref{corr} does not apply
directly to the case of $\GGG_K$ acting on $\AAA^{1,1}\times\SSS$
(since there we make strong use of the
fact that the group acting symplectically is compact), one can
use some of the ideas (with some additional analytic results)
to prove theorem \ref{main}. This will be done in the next section.
The main strategy will be, following the usual approach in proving
the Hitchin--Kobayashi correspondence, to minimize the integral 
of the moment map.

As a final comment, note that so far we have defined the gauge 
group as the space of smooth sections of a certain bundle. Eventually, 
it will be necessary to take a metric on $\GGG_K$ (and $\GGG_G$) and 
complete both spaces with respect to the metric, to assure the convergence
of certain sequences. We will use Sobolev $L^p_2$ and $L^2_1$ norms.

\section{Proof of the correspondence}
\label{demostracio}
\subsection{The length of elements of the gauge group}
There are several ways to extend the notion of length to elements of
the gauge group. We will mainly use these two definitions: if
$g\in\GGG_G$, then $|g|_{{\log},C^0}=\| |g|_{\log}\|_{C^0}$ and
similarly $|g|_{{\log},L^1}=\| |g|_{\log}\|_{L^1}$ (to give 
this a sense we use the $K$ invariance of the length function, which 
is a consequence of the fact that the Cartan decomposition 
$G\simeq K\times\exp(\imag\klie)$ is $K$-equivariant).
Define a norm $\|\cdot\|_{L^p}$ in 
$\Lie(\GGG_G)=\Omega^0(K\times_{\Ad}\glie)$ as 
the $L^p$ norm of $|\cdot|$: if 
$s\in\Omega^0(K\times_{\Ad}\glie)$ then
$$\|s\|_{L^p}=\left(\int_{x\in X} |s(x)|^p\right)^{1/p}.$$
We will usually write $\|\cdot\|$ instead of $\|\cdot\|_{L^2}$.

\subsection{Stability implies existence of solution}
Here we will follow the scheme in section \ref{correspondencia}.
Fix a pair $(A,\Phi)\in\AAA^{1,1}\times\SSS$.
We will make use of the integral of the moment map
$\mu^c(A,\Phi)=\Lambda F_A+\mu(\Phi)-c$,
$\Psi^c=(\Psi^{\AAA\times\SSS})^{\mu^c}_{(A,\Phi)}=
(\Psi^{\AAA})^{\mu^c}_A+(\Psi^{\SSS})^{\mu^c}_{\Phi}$, 
and will see that if the pair $(A,\Phi)$
is simple and $c$-stable, then there exists a $\GGG_K$ orbit inside
the $\GGG_G$ orbit of $(A,\Phi)$ on which $\Psi^c$ attains
its minimum. The main step will be to prove that if the
condition of $c$-stability is satisfied, then the map
$\Psi^c$ satisfies an inequality like that in lemma \ref{fita}.
This method of proof is exactly the same that appears in
\cite{Si, Br2, BrGP1, DaUW} (and in many other places where
similar results are proved), though here we have tried to 
remark the similarities with the finite dimensional case,
so our notation changes a little bit. However, in some steps of the proof
we will only give a sketch, refering to \cite{Br2} for details. 

Recall that on $\glie$ we have a Hermitian
pairing $\la,\ra:\glie\times\glie\to\CC$ and a norm $|\cdot|$,
both obtained by means of the auxiliar representation $\rho_a$.
We will use the following $L^p$ norm on $\Omega^0(E\times_{\Ad}\glie)$:
$$\|s\|_{L^p}=\left(\int_X |s(x)|^p\right)^{1/p},$$
and Sobolev norm
$$\|s\|_{L^p_2}=\|s\|_{L^p}+\|d_A s\|_{L^p}
+\|\nabla d_A s\|_{L^p},$$
where $\nabla:\Omega^0(T^*X\otimes E\times_{\Ad}\glie)\to 
\Omega^1(T^*X\otimes E\times_{\Ad}\glie)$ is
$\nabla_{LC}\otimes d_A$, $\nabla_{LC}$ being the Levi-Civita
connection. 
As usual, $L^p_2(E\times_{\Ad}\glie)$ will denote
the completion of $\Omega^0(E\times_{\Ad}\glie)$ with respect to
the norm $\|\cdot\|_{L^p_2}$.

\subsubsectionr{}
Suppose from now on that $(A,\Phi)$ is simple and $c$-stable. Our aim is to
minimise $\Psi^c$ in $\GGG_G/\GGG_K$. Through the exponential
map we can identify $\GGG_G/\GGG_K$ with 
$\Omega^0(E\times_{\Ad}\imag\klie)$.
Fix from now on $p>2n$ and define
$$\Met=L^p_2(E\times_{\Ad}\imag\klie).$$
The first thing to do is to restrict ourselves to the subset of 
$\Met$ defined as follows:
$$\MetB=\{s\in\Met | \|\mu^c(e^s(A,\Phi))\|^p_{L^p}\leq B\}.$$
Here $B$ is any positive real constant. We prove that if a metric 
minimizes the functional in $\MetB$, then it also minimizes it in $\Met$.
For that it is enough to see that any minimum in $\MetB$
lies away from the boundary of $\MetB$; to verify this claim
one needs the hypothesis that the pair $(A,\Phi)$ is simple. Let us briefly 
explain how this goes (see also \cite{Br2}, Lemma 3.4.2).

Suppose that $s$ minimizes the functional inside $\MetB$. 
Let $B=e^s(A)$, $\Theta=e^s(\Phi)$. Define the differential operator
$L:L^p_2(E\times_{\Ad}\imag\klie)\to L^p(E\times_{\Ad}\imag\klie)$
as $$L(u)=\imag\left.\frac{\partial}{\partial t} 
\mu^c(e^{tu}(B,\Theta))\right|_{t=0}
=\imag\la d\mu^c,u\ra_{T(\AAA\times\SSS)}(B,\Theta).$$
Now, if we can see that there exists an $u$ such that
\begin{equation}
L(u)=-\imag\mu^c(B,\Theta), 
\label{lmh}
\end{equation}
then we can deduce that $\mu^c(B,\Theta)=0$ and, hence, that
$s$ minimizes the functional in the whole space of metrics
$\Met$ (see \cite{Br2}, Lemma 3.4.2 for a proof of this fact).
The operator $L$ is Fredholm and has index zero. Indeed,
modulo a compact operator it is 
$\imag\Lambda\overline{\partial}_{B}\partial_{B}$.
Using the Kaehler identities this is equal to
$\partial_{B}^*\partial_{B}$, which is clearly an
elliptic self adjoint operator.
This implies that if $\Ker(L)=0$ then $L$ is surjective 
and so, in particular, equation (\ref{lmh}) has a solution. 
Assume that $L(u)=0$, where $u\in\Met$. Then, by lemma \ref{gradient}, 
\begin{align}
0&=\la -\imag L(u),-\imag u\ra=
\la \la d\mu^c,u\ra_{T(\AAA\times\SSS)},-\imag u\ra_{\Lie(\GGG_K)}
(B,\Theta)\notag\\
&=\|\fX^{\AAA^{1,1}\times\SSS}_{-\imag u}(B,\Theta)\|^2.
\end{align}
And this implies that $-\imag u$ leaves $(B,\Theta)$ invariant.
Hence if $u\neq 0$ then, since $u$ is semisimple, 
$(B,\Theta)$ is not simple, so neither is
$(A,\Phi)$; and this is a contradiction.

\subsubsectionr{}
The next step is to prove that the functional $\Psi^c$ is
linearly proper with respect to the $C^0$ norm in $\GGG_G$.

\begin{lemma}
There exist positive constants $C_1, C_2$ such that for any $s\in\MetB$ 
one has $\sup |s|\leq C_1 \Psi^c(e^s)+C_2.$
\label{fita2}
\end{lemma}
\begin{remark}
It makes sense to speak about $\sup|s|$ because, since we took
$p>2n$, the Sobolev embedding theorem implies that 
$L^p_2\hookrightarrow C^0$ continuously (in fact, this is
a compact embedding).
\end{remark}
Just as in lemma \ref{fita}, it is here that one uses the stability 
of the pair $(A,\Phi)$. First of all one sees that such a bound is 
equivalent to an $L^1$ bound:
$\|s\|_{L^1}\leq C_1 \Psi^c(e^s)+C_2$
(the constants in both inequalities need not be the same!). One uses that
pointwise
\begin{equation}
|s|\Delta|s|\leq \la\Lambda F_{e^s(A)} - \Lambda F_A,-\imag s\ra.
\label{equ1}
\end{equation}
This is proved in full detail in (\cite{Br2}, Prop. 3.7.1) for
$G=GL(n;\CC)$ and the metric induced by the fundamental representation. 
In our case, we use the auxiliar representation $\rho_a$ to apply
this result to our $G$. 

\begin{lemma} For any point $x\in X$
\begin{equation}
0\leq \la\mu(e^s\Phi(x))-\mu(\Phi(x)),-\imag s(x)\ra_{\klie}.
\label{equ2}
\end{equation}
\end{lemma}
\begin{pf}
The gradient flow of $\mu_{-\imag s}$ is precisely $e^s$
(see lemma \ref{gradient}).
\end{pf}

Summing the inequalities (\ref{equ1}) and (\ref{equ2}), 
using Cauchy-Schwartz, and dividing by $|s|$ we obtain the following
pointwise bound:
$$\Delta|s|\leq |\mu^c(e^s(A,\Phi))-\mu^c(A,\Phi)|.$$
And now, making use of a result of Donaldson (see \cite{Br2}, Lemma 3.7.2), 
this bound allows to relate the $C^0$ and $L^1$ norms
of $s$ provided $s\in\MetB$. More precisely, we conclude that there
exists a constant $C_B$ such that for any $s\in\MetB$ one has
$\|s\|_{C^0}\leq C_B\|s\|_{L^1}$.

\subsubsectionr{}
In order to prove the existence of constants $C_1$ and $C_2$ such that
$\|s\|_{L^1}\leq C_1 \Psi^c(e^s)+C_2$, we suppose the contrary and try
to deduce that in this case the pair $(A,\Phi)$ cannot be $c$-stable.
If there exist not such constants, then we can find a sequence of
real numbers $C_j\to\infty$ and elements $s_j\in\MetB$ with
$\|s_j\|_{L^1}\to\infty$ such that 
$\|s_j\|_{L^1}\geq C_j \Psi^c(e^s)$
(see \cite{Br2}, Lemma 3.8.1). 
Set $l_j=\|s_j\|_{L^1}$, $u_j=l_j^{-1}s_j$ so that
$\|u_j\|_{L^1}=1$ and $\sup|u_j|\leq C$.

\begin{lemma}
After passing to a subsequence, there exists 
$u_{\infty}\in L^2_1(E\times_{\Ad}\imag\klie)$ such that
$u_j\to u_{\infty}$ weakly in $L^2_1(E\times_{\Ad}\imag\klie)$
and such that
$$\lambda((A,\Phi);-\imag u_{\infty})\leq 0.$$
\end{lemma}
\begin{pf}
Just as in lemma \ref{fita}, take $t>0$. Then (4) in proposition
\ref{propietats1} gives
\begin{align}
\frac{1}{C_j}&\geq\frac{\Psi^c(e^{s_j})}{\|s_j\|}\geq
\frac{l_j-t}{l_j}\lambda_t((A,\Phi);-\imag u_j)+\frac{1}{l_j}
\int_0^t\lambda_l((A,\Phi);-\imag u_j)dl\notag\\
&=\frac{l_j-t}{l_j}(\lambda_t(A;-\imag u_j)+\lambda_t(\Phi;-\imag u_j))
\notag\\
&+\frac{1}{l_j}\int_0^t(\lambda_l(A;-\imag u_j)+
\lambda_l(\Phi;-\imag u_j))dl.
\label{desclau}
\end{align}
Now, since $\|u_j\|_{C^0}\leq C_B$, 
and $X$ is compact, $\lambda_t(\Phi;-\imag u_j)$ and 
$\int_0^t\lambda_t(\Phi;-\imag u_j)dl$ are both bounded. Hence, there exists 
$C$ such that for any $j$ $$\frac{l_j-t}{l_j}\lambda_t(A;-\imag u_j)+
\frac{1}{l_j}\int_0^t\lambda_l(A;-\imag u_j)dl<C.$$
Using again the boundedness of $\|u_j\|_{C^0}$ and taking into
account lemma \ref{tpesmaxconn} we obtain 
$$\|\overline{\partial}_A(u_j)\|_{L^2}<C_1.$$
Now, $\overline{u_j}=u_j$ (because the Cartan involution leaves
$\imag\klie$ fixed), and this implies that $\|u_j\|_{L^2_1}$ is also
bounded. So we can take a subsequence (which we again
call $\{u_j\}$) that converges weakly to $u_{\infty}\in L^2_1$. 
We can also assume that there exists the limit
$\lim_{i\to\infty} \lambda_t((A,\Phi);-\imag u_j)$.
On the other hand, since the embedding $L^2_1\hookrightarrow L^2$ 
is compact, we get strong convergence $u_j\to u_{\infty}$ in $L^2$.
$\|u_j\|_{L^1}=1$ and the uniform bound $\|u_j\|_{C^0}\leq C_B$
imply that $\|u_j\|_{L^2}>C_B^{-1}>0$, so $u_{\infty}\neq 0$.
To see that $\lambda_t((A,\Phi);-\imag u_{\infty})\leq
\lim_{i\to\infty} \lambda_t((A,\Phi);-\imag u_j)$ we observe that 
$$u_j\in L^2_{0,C_B}(E\times_{\Ad}\imag\klie)=
\{s\in L^2(E\times_{\Ad}\imag\klie)|\ |s(x)|\leq C_B\mbox{ a.e.}\}.$$
This implies that $u_{\infty}\in L^2_{0,C_B}(E\times_{\Ad}\imag\klie)$,
and this is enough to get the inequality (see \cite{Br2}, proposition 3.2.2). 
Finally, making $j\to\infty$ in formula (\ref{desclau}) we obtain
$$\lim_{i\to\infty} \lambda_t((A,\Phi);-\imag u_j)\leq 0,$$ so in particular
$\lambda_t((A,\Phi);-\imag u_{\infty})\leq 0$. Since this is true for any
$t>0$, we get $\lambda((A,\Phi);-\imag u_{\infty})\leq 0$.
\end{pf}

The next steps are rather standard. One can prove that
$\rho_a(u_{\infty})$ has almost everywhere constant eigenvalues and that it
defines a filtration of $V$ by holomorphic subbundles in the complement
of a complex codimension 2 subvariety of $X$. 
This follows exactly the same lines as sections 3.9 and 3.10 in \cite{Br2},
the main technical point being the use of a theorem of Uhlenbeck and
Yau \cite{UY} on {\it weak subbundles} of vector bundles
(see section 3.11 in \cite{Br2}).
The filtration of $V$ on $X_0$ and the gauge transformation 
$u_{\infty}$ lead to a reduction of the structure group
$\sigma\in\Gamma(X_0;E(G/P))$ defined on $X_0$ by \ref{filtsect}
which will be holomorphic thanks to the results in
subsection \ref{holomred}, and
an antidominant character $\chi$ of $P$. The degree of the pair 
$(\sigma,\chi)$ equals $\lambda((A,\Phi);-\imag u_{\infty})\leq 0$. 
And this contradicts the stability condition, thus finishing the proof
of lemma \ref{fita2}.

\subsubsectionr{} With the inequality of lemma \ref{fita2} in our hands,
we finish the proof of existence of solution to the equations
exactly as is done in \cite{Br2}, section 3.14. This consists
of two steps: the first one is to verify that there exists
an element $s\in\MetB$ minimising $\Psi^c$ and the second
one is to prove the smoothness of this solution $s$.

\subsection{Existence of solutions implies stability}

The method we will follow in this section will be exactly the
same as in the finite dimensional case in section \ref{correspondencia}.
Let us take a simple pair $(A,\Phi)\in\AAA^{1,1}\times\SSS$.
Suppose that there exists a gauge transformation $h\in\GGG_G$ 
such that $h(A,\Phi)$ satisfies equation (\ref{hk}).
Then the pair $h(A,\Phi)$ is analitically stable, 
by exactly the same reasoning as in section \ref{correspondencia}.
The key step is to prove that this implies that $(A,\Phi)$
is also analitically stable.

Take $X_0\subset X$ with complement of complex codimension 2, 
$P\subset G$ parabolic, $\chi$ an antidominant character of $P$ 
and fix a reduction $\sigma\in\Gamma(X_0;E(G/P))$. 
By \ref{sectfilt} we get a section 
$g_{\sigma,\chi}\in\Omega^0(X_0;E\times_{\Ad}\imag\klie)$,
and we have to check that $\lambda((A,\Phi);-\imag g_{\sigma,\chi})>0$.

Since $h(A,\Phi)$ is analitically stable, 
given any $B>0$ there exist constants $C_1$ and $C_2$ such that for any
$s\in\MetB$ there is an inequality
\begin{equation}
\sup |s|\leq C_1\Psi^c_{h(A,\Phi)}(e^s)+C_2.
\label{proper}
\end{equation}
This inequality is valid not only for $s\in\MetB$, but also for any 
$$s\in\Met(C_B)=\{s\in L^p_2(X_0;E\times_{\Ad}\imag\klie)
|\ \|s\|_{C^0}\leq C_B\|s\|_{L^1}\},$$ 
as one can see tracing the proof of lemma \ref{fita2}.
(The only property on $X$ that is used in the proof of the 
inequality for $s\in\Met(C_B)$ is, besides having finite volume, 
that it has no nonconstant holomorphic functions; and this also 
happens in $X_0$, by Hartog theorem.) This proves the following.
\begin{lemma}
Fix a positive constant $C_B$. There exist positive constants
$C_1, C_2$ such that the following holds.
Let $g\in\GGG_G(X_0)=\Omega^0(X_0;E\times_{\Ad}G)$ be
such that $|g|_{{\log},C^0}\leq C_B |g|_{{\log},L^1}<\infty$. 
Then
$$|g|_{{\log},C^0}\leq C_1\Psi^c_{h(A,\Phi)}(g)+C_2.$$
\label{fita3}
\end{lemma}
If we take $C_B=\Vol(X)^{-1}$ then 
$|g_{\sigma,\chi}|_{C^0}\leq C_B|g_{\sigma,\chi}|_{L^1}$
(in fact we have equality, since $|g_{\sigma,\chi}|$ is constant).
So reasoning exactly like in lemma \ref{fita}, the preceeding lemma 
implies that $\lambda(h(A,\Phi);-\imag g_{\sigma,\chi})>0$. 
To deduce that $\lambda((A,\Phi);-\imag g_{\sigma,\chi})>0$ as well, 
we need the following lemma.

\begin{lemma}
There a positive constant $C_B'$ such that 
for any $g_{\sigma,\chi}$ and $h$ and
for big enough (depending on $g_{\sigma,\chi}$ and $h$) $t>0$,
$$|e^{t g_{\sigma,\chi}}h^{-1}|_{{\log},C^0}\leq
C_B' |e^{t g_{\sigma,\chi}}h^{-1}|_{{\log},L^1}.$$
\end{lemma}
\begin{pf}
This is a consequence of lemma \ref{comparison} and the fact
that $X$ is compact (so $|h|$ and $|h^{-1}|$ are bounded functions
on $X$).
\end{pf}

Now we set $C_B=C_B'$ in lemma \ref{fita3}, and proceeding
as in the finite dimensional case (lemma \ref{tambe}) we deduce
that there exist positive constants $C_1''$ and $C_2''$ such
that for any $t>0$
$$\sup |t g_{\sigma,\chi}|=t \sup |g_{\sigma,\chi}|
\leq C_1''\Psi^c_{h(A,\Phi)}(e^{t g_{\sigma,\chi}})+C_2''.$$
This implies that $\lambda((A;\Phi);-\imag g_{\sigma,\chi})>0$.

By lemma \ref{equivalencia} this is equivalent to 
$T^c_{\Phi}(\sigma,\chi)>0$. With this we see that 
$(A,\Phi)$ is $c$-stable.

\subsection{Uniqueness of solutions}
The proof is exactly as in the finite dimensional case (see section 
\ref{correspondencia}): it follows from the convexity of the integral of 
the moment map.

\subsection{Nonsimple pairs}
\label{nosimple}
The Hitchin--Kobayashi correspondence which
we have proved applies only to simple pairs $(A,\Phi)$. This restriction,
however, can often be relaxed. As an example, suppose that there are
elements in the centre $Z=Z(\glie)$ of $G$ which leave $F$ fixed
(trivial example: $F$ equal to a point). Any element 
$z\in Z$ gives an element of the Lie algebra of the gauge group, which
we still denote by $z$. This element is semisimple and for any
$t$ the exponential $\exp(tz)$ fixes all connections in $\AAA$, and
by our assumption fixes also $\Phi$. In this situation, the pair
$(A,\Phi)$ is not simple.

When our group $G$ is $\GL(V)$, there is a standard way to solve
this problem. We assume that the center $Z$ of $GL(V)$ leaves 
$\Phi$ fixed (note that this is {\it not} the case of the 
vortex equations). We {\it split} the equation 
in the $Z$ part and in the $G/Z$ part as follows. 
Define $\GGG_G^0$ to be the set of gauge transformation with
determinant pointwise equal to 1, and suppose that there 
are no semisimple elements in the Lie algebra of $\GGG_G^0$
which leave $(A,\Phi)$ fixed; under this assumption we can
find an element $g\in\GGG_G^0$ so that $g(A,\Phi)$ solves the 
trace-free part of the equation (observe that our proof applies 
to this situation); then Hodge theory gives a central element in 
$\GGG_G$ which, composed with $g$, solves the complete equation.

This idea applies for any reductive Lie group $G$. 
We just need to give a generalisation of the condition
of having determinant pointwise equal to 1 which we imposed
to the elements in $\GGG_G^0$. This is given by the following

\begin{lemma}
Let $G$ be a reductive Lie group. There exists $k\geq 1$
and a morphism $\phi:G\to(\CC^*)^k$ such that $\Ker\phi\cap Z$
is a discrete subgroup of $G$.
\end{lemma}
\begin{pf}
Take a faithful representation $\rho:G\to GL(W)$. Split $W$
in eigenspaces of the roots of $Z$ acting on $W$:
$W=W_1\oplus\dots\oplus W_k$, so that any central element
$z\in Z$ acts on any piece $W_j$ by homotecies.
Then $\rho(G)\subset GL(W_1)\times\dots\times GL(W_k)$, so
that for any $g\in G$ we have $\rho(g)=(g_1,\dots,g_k)$.
Let $\phi:G\to(\CC^*)^k$ be defined as
$\phi(g)=(\det g_1,\dots,\det g_k)$. Now suppose that there exists
$s\in Z(\glie)$ such that, for any $t$, 
$\phi(e^{ts})=(1,\dots,1)$. Since $e^{ts}$ acts
by homotecies on each piece, we must have 
$\rho(e^{ts})\in Z(SL(W_1))\times\dots\times Z(SL(W_k))
\simeq \ZZ/w_1\ZZ\times\dots\times\ZZ/w_k\ZZ$ for any $t$,
where $w_j=\dim W_j$. This implies that $\rho(e^{ts})=(1,\dots,1)$
and, since $\rho$ is faithful, $z=0$. This proves that 
$\Ker\phi\cap Z$ is discrete.
\end{pf}

Suppose for simplicity that the whole center of $G$ leaves
$\Phi$ fixed. 
We then define $\GGG_G^0$ to be the set of gauge transformations
which fibrewise belong to $\Ker\phi$, and proceed as in the case
$G=GL(V)$: we find $g\in\GGG^0_G$ such that the center free part
of the equation is solved and then use Hodge theory to solve the
complete equation.

\section{Example: the theorem of Banfield}
\label{banfield}
Suppose that $F$ is a Hermitian vector space and that $K$ acts on $F$
through a unitary representation $\rho:K\to U(F)$. D. Banfield \cite{Ba} has 
recently proved a general Hitchin--Kobayashi correspondence for this 
situation. The work of Banfield generalises existing results on vortex 
equations, Hitchin equations, and on other equations arising from particular 
choices of $K$ and $\rho$ (see subsection \ref{eihiko}).
In this section we will see how the result of Banfield
can be deduced from theorem \ref{main}. 

\subsection{The stability condition}                 
The first thing we do is to study the maximal weights of elements
in $\klie$ acting on $F$ through $\rho$.
Let $h$ be the Hermitian metric on $F$. The
imaginary part of $h$ with reversed sign defines a symplectic form 
$\omega_F$ compatible with the complex structure and hence a Kaehler
structure. The action of $K$ on $F$ respects
the Kaehler structure and admits a moment map $\mu:F\to\klie^*$
$$\mu(x)=-\frac{\imag}{2}\rho^*(x\otimes x^*).$$ 
In other words, for any $s\in\klie$, 
$\la\mu(x),s\ra_{\klie}=-\frac{\imag}{2} h(x,\rho(s)x)$.
Let $x\in F$ and take an element $s\in\klie$.
Since $\rho(s)\in \ulie(F)$, the endomorphism $\rho(s)$ diagonalises
in a basis $e_1,\dots,e_n$: $\imag\rho(s)e_k=\lambda_k e_k$,
where $\lambda_k$ is a real number for any $k$. 
Write $x=x_1e_1+\dots+x_ne_n$. 
\begin{lemma}
If $\lambda_k\leq 0$ for every $k$ such that $x_k\neq 0$,
then the maximal weight $\lambda(x;s)$ is equal to zero. Otherwise it is 
$\infty$.
\label{pesvect}
\end{lemma}

Let us assume that the representation $\rho$ is contained
in the auxiliar representation $\rho_a$.
Let $E\to X$ be a $G$-principal bundle on a compact
Kaehler manifold $X$. Let $\FFF=E\times_{\rho} F$ be
the vector bundle associated to $E$ through the representation
$\rho$. Take a pair $(A,\Phi)\in\AAA^{1,1}\times\SSS$, and
fix a central element $c\in\klie$. Consider
on $E$ the holomorphic structure given by $\overline{\partial}_A$.
According to definition \ref{parella_estable},
$(A,\Phi)$ is $c$-stable if and only if for any parabolic subgroup
$P\subset G$, for any holomorphic reduction $\sigma\in\Gamma(X_0;E(G/P))$ 
defined on the complement of a complex codimension 2 submanifold $X_0$
of $X$ and for any antidominant character $\chi$ of $P$, the total degree
is positive:
$$T^c_{\Phi}(\sigma,\chi)>0.$$
The total degree is the sum of $\deg(\sigma,\chi)$ plus the
maximal weight of the action of $g_{\sigma,\chi}$ on 
$\Phi$ plus $\la\imag\chi,c\ra\Vol(X)$. The maximal weight is
\begin{equation}
\int_{x\in X}\lambda(\Phi(x);-\imag g_{\sigma,\chi}(x)).
\label{intBan}
\end{equation}
Define now $\FFF^-=\FFF^-(\sigma,\chi)\subset\FFF$ to be 
the subset given by the vectors in $\FFF$ on which $g_{\sigma,\chi}(x)$ 
acts negatively, that is, $v\in\FFF_x$ belongs to $\FFF^-$ if and only if
you can write $v=\sum v_n$ such that
$g_{\sigma,\chi}(x)(v_n)=\lambda_n v_n$ and $\lambda_n\leq 0$.
Since the eigenvalues of $g_{\sigma,\chi}$ are constant, $\FFF^-$
is a subbundle. And since the parabolic reduction is holomorphic,
so is $\FFF^-$.

If $\Phi\subset\FFF^-$, then the maximal weight at each fibre
is equal to zero by lemma \ref{pesvect}, so the stability 
condition reduces to $$\deg(\sigma,\chi)>0.$$
On the other hand, if $\Phi(x)\notin\FFF^-_x$, then there is
an open neighbourhood $U$ of $x$ such that $\Phi(y)\notin\FFF^-_y$
for any $y\in U$. In this situation lemma \ref{pesvect} tells us
that, for any $y\in U$, $\lambda(\Phi(y);-\imag g_{\sigma,\chi}(y))=\infty$.
Since this happens in an open set, the integral (\ref{intBan})
is infinite (since $X$ is compact, $\Phi$ is bounded and so 
$\lambda(\Phi(x);-\imag g_{\sigma,\chi}(x))$ is bounded below).
But the degree $\deg(\sigma,\chi)$ is always a finite number,
so the total degree will be positive (infinite, in fact) in this case. 
To sum up,

\begin{prop}
The pair $(A,\Phi)$ is stable if and only if for any
$P,\sigma,\chi$ as above, if $\Phi$ is contained in 
$\FFF^-(\sigma,\chi)$, then
$$\deg(\sigma,\chi)+\la\imag\chi,c\ra\Vol(X)>0.$$
\end{prop}

And this is precisely Banfield condition.

\subsection{Simple pairs}
To give a characterisation of simple pairs we use the following
definition due to Banfield \cite{Ba}:

\begin{definition}
Suppose that the vector bundle $\FFF$ decomposes into a nontrivial
direct sum $\bigoplus_k \FFF_k$ of holomorphic vector bundles
and that there is a reduction of the structure group of $E$
to $G'\subset G$, compatible with the splitting. 
Suppose further that a central element of the Lie algebra 
$\glie'$ of $G'$ annihilates the section $\Phi$ but acts nontrivially 
on $\FFF$. Then we say $(A,\Phi)$ is a {\bf decomposable pair}. 
If no such splitting exists, the we say that $(A,\Phi)$ is an 
{\bf indecomposable pair}.
\end{definition}

\begin{lemma}
The pair $(A,\Phi)$ is simple if and only if 
it is indecomposable.
\label{bansim}
\end{lemma}
\begin{pf}
Suppose that 
$0\neq s\in\Omega^0(E\times_{\Ad}\glie)$ is semisimple and stabilises 
$(A,\Phi)$. In particular $\fX^{\AAA}_s(A)=0$, and this implies that 
$\overline{\partial}_A(s)=0$. So the eigenvalues of $\rho(s)$ 
are constant, and since $s$ is semisimple $\rho(s)$ diagonalises. 
Let the different eigenvalues of $\rho(s)$ be $\lambda_1<\dots<\lambda_r$, 
and consider the decomposition 
$\FFF=\FFF(\lambda_1)\oplus\dots\oplus\FFF(\lambda_r)$
in eigenbundles, which are holomorphic,
and every $\FFF_k=\FFF(\lambda_k)$ having as structure group
a subgroup $G_k\subset G$. 
Since $s$ leaves $\Phi$ fixed
$\Phi$ must belong to $\FFF(0)$. On the other hand, $0$ in obviously
not the unique eigenvalue of $\rho(s)$, so the decomposition
$$\FFF=\FFF_1\oplus\dots\oplus\FFF_r$$
is not trivial. Finally, the section $s$ provides the central element
killing $\Phi$.

The proof of the converse is similar.
\end{pf}                                                      

\subsection{The equations}
Our equation (\ref{hk}) in the case of linear representations
is the same one given by Banfield (note that Banfield also
considers the holomorphicity condition $\ov{\partial}_A\Phi=0$). 

\section{Example: projective pairs} 

In this section and in the next ones we give some examples in which
the result of Banfield does not apply. 

Suppose that $F=\PP(W)$, where $W$ is a complex vector space with
a Hermitian pairing, and that a compact Lie group $K$ acts on $W$ 
through a representation $\rho:K\to U(W)$. 
Let us remark that not all the actions of compact Lie groups on 
projective spaces arise in this way. More precisely, not always
an action on $\PP(W)$ will lift to an action on $W$.
However, if $K$ acts on $\PP(W)$, one can find a central extension of 
$K$ by $\CC^*$ which does act on $W$, and all the following
discusion adapts easily to this more general situation.

The vector space $W$ is a Kaehler manifold with symplectic form $\omega_W$
equal to the imaginary part of the Hermitian pairing with reversed sign.
Consider the action of $U(1)$ on $W$ given by multiplication. 
This action is symplectic, and
it has a moment map $\mu_{U(1)}(x)=\|x\|^2$. The symplectic quotient,
$\mu_{U(1)}^{-1}(1)/U(1)$, coincides with the projective space
$F=\PP(W)$. So the induced symplectic form $\omega_F$ 
on $F$ comes from the restriction of the symplectic 
form $\omega_W$ of $W$ on $\mu_{U(1)}^{-1}(1)$ (this makes sense, since
the symplectic form is $U(1)$ invariant). The complex structure
on $F=\PP(W)$ is compatible with $\omega_F$, so $F$ is in fact a
Kaehler manifold. 

Let $\mu_K:W\to\klie^*$ be the moment map of the action of $K$ on $W$.
The action of $U(1)$ commutes with that of $K$, so the symplectic 
quotient has an induced action of $K$ (which leaves $\omega_F$ invariant). 
Just as happened with the symplectic form, a moment map $\mu_F$ for this
action on $F$ can be obtained considering the restriction of the moment map 
$\mu_K$ in $\mu_{U(1)}^{-1}(1)$; by the $U(1)$-equivariance of $\mu_K$, this 
descends to the quotient. More explicitly, given any $x\in F$,
one takes any lifting $\hat{x}\in W$ and the moment map at $x$ is
\begin{equation}
\mu_F(x)=-\frac{\imag}{2}\rho^*\left(\frac{\hat{x}\otimes \hat{x}^*}
{\|\hat{x}\|^2}\right).
\label{momproj}
\end{equation}

\subsubsection{Maximal weights} Take a point $x\in F$ and consider
an element $s\in\klie$. We can take a basis
$e_1,\dots,e_n$ of $W$ in which the action of $s$ diagonalizes:
$\imag\rho(s)e_k=\lambda_k e_k$,
where $\lambda_k$ is a real number for any $k$. 
Fix a lifting $\hat{x}\in W$ of $x$ and write 
$\hat{x}=x_1e_1+\dots+x_ne_n$. Then
$$\lambda_t(x;s)=\la \mu(e^{\imag t\rho(s)}x),s\ra_{\klie}=
-\imag \frac{h(e^{\imag t\rho(s)}\hat{x},\rho(s)e^{\imag t\rho(s)}\hat{x})}
{\|\hat{x}\|^2}=\frac{\sum_{k=1}^n \lambda_k e^{2t\lambda_k}|x_k|^2}
{\sum_{k=1}^n e^{2t\lambda_k}|x_k|^2}.$$
\begin{lemma}
The maximal weight of $s$ acting on $x$ is
$$\lambda(x;s)=\max \{\lambda_k|x_k\neq 0\}.$$
\label{pesmaxproj}
\end{lemma}

\subsubsection{The integral of the moment map} 
\label{omegapsi}
The function $\Psi$
takes in this situation the following form: for $x\in F$ and $g\in G$,
$$\Psi(x,g)=\frac{1}{4}\log\frac{\|\rho(g)\hat{x}\|^2}{\|\hat{x}\|^2}.$$
Once again, this is checked by proving that this function satisfies
(2) in proposition \ref{propietats}.

\subsection{The stability condition}
               
The pairs $(A,\Phi)\in\AAA^{1,1}\times\FFF$, where $\FFF=E\times_K\PP(W)$
are called {\bf projective pairs}.
We will give a characterization of stability for projective
pairs very similar to that of Banfield. This characterisation,
however, will only work if we ask $\Phi$ to be a holomorphic
section of $\FFF$ with respect to the holomorphic structure
$\overline{\partial}_A$. (Remark that anywhere else in this chapter
we only wanted it to be smooth.)

As in the preceeding section, we will assume that the representation
$\rho$ is contained inside the auxiliar representation $\rho_a$.
Let $P\subset G$ be a parabolic subgroup and take $\chi$
an antidominant character of $P$. 
Let $\lambda_1<\dots<\lambda_r$ be the set of different
eigenvalues of $\rho(\chi)$, and write $W(\lambda)$ the eigenspace
of eigenvalue $\lambda$. Finally, write $W^{\lambda_k}=\bigoplus_{i\leq k}
W(\lambda_i)$. 

\begin{lemma}
If $\widehat{w}\in W^{\lambda_k}\setminus W^{\lambda_{k-1}}$,
then $\lambda(\widehat{w};-\imag\chi)=\lambda_k$.
\label{equi4}
\end{lemma}

On the other hand, by lemma \ref{equi1}, for any $k$ the subspace
$W^k\subset W$ is invariant by the action of $P$.
So, once we have a holomorphic reduction
of the structure group of $E$ to $P$ and an antidominant character
$\chi$ of $P$, we obtain holomorphic fibrations
$$V^{\lambda_1}\subset\dots\subset V^{\lambda_r}$$ and
$$\PP(V^{\lambda_1})\subset\dots\subset \PP(V^{\lambda_r})=\FFF,$$
defined as $V^{\lambda_k}=E\times_P W^{\lambda_k}$.
(Of course all this may happen to be defined only on the complementary
of a complex codimension 2 submanifold of $X$, but here we will
avoid this technicallity.)
Our section $\Phi$ is assumed to be holomorphic. Hence, if for
some $x\in X$ we have $\Phi(x)\in \PP(V^{\lambda_k})_x\setminus
\PP(V^{\lambda_{k-1}})_x$, then the same happens for almost
any $x\in X$. In consequence, using now lemma \ref{equi4}, for almost
any $x\in X$ the maximal weight $\lambda(\Phi(x);-\imag\chi)$ equals 
$\lambda_k$. This implies that
$$\int_{x\in X} \lambda(\Phi(x);-\imag g_{\sigma,\chi})=\lambda_k.$$

In view of all this it is now easy to prove the following.

\begin{prop}
The pair $(E,\Phi)$ is $c$-stable if and only if for any
pair $(\sigma,\chi)$, if $\Phi\in\PP(V^{\lambda_k})$, then
$$\deg(\sigma,\chi)+\Vol(X)\lambda_k-\int_X\la\chi,c\ra > 0.$$
\end{prop}

On the other hand, the characterisation of simple pair given
in the case of vector pairs works equally well for projective pairs:
the definition of indecomposable pair is valid in the case
of projective pairs, and one can prove that a pair is simple
if and only if it is indecomposable.

Finally, remark that when $F\subset\PP^N$ is a projective variety
the Hitchin--Kobayashi correspondence reduces to that for 
projective pairs. In the following section we will see some
examples of this situation: $F$ will be there either a Grassmannian
or, more generally, a flag manifold.

\section{Example: filtrations of vector bundles}
\label{filtracions}

In this section we study theorem \ref{main} in the particular case
in which $F$ is a Grassmannian or, more generaly, a flag manifold.
We assume, for simplicity, that $X$ is a Riemann
surface and that $\Vol(X)=1$. For the higher dimensional 
case everything that follows
remains valid with the following modification: in the stability condition
one has always to consider reflexive subsheaves, and not only 
subbundles (this is a consequence of the need of
considering reductions of the structure group defined on the
complement of a complex codimension 2 submanifold of $X$ in
the general definition of stability).

%
%
The Lie group $K$ will be $\U(R;\CC)$,
where $R\geq 1$ is an arbitrary integer, and
we will take the standard representation in $\CC^R$ as our
auxiliar representation. 

\subsection{Subbundles}
Let $E\to X$ be a principal $\U(R;\CC)$ bundle on $X$. Consider
the standard representation on $\CC^R$. This provides us with
a vector bundle $V\to X$ of rank $R$. Using theorem \ref{main}, we 
will find a Hitchin--Kobayashi correspondence for subbundles 
$V_0$ of $V$ of fixed rank $0<k<R$. This correspondence has already
been proved in \cite{BrGP1} and in \cite{DaUW}.

Using an idea of \cite{DaUW} we identify the inclusion 
$V_0\hookrightarrow V$ with a section $\Phi$ of the bundle with fibres the 
Grassmannian of $k$-subvectorspaces $\Gr_k(\CC^R)$ associated to
$E$ by the usual action of $\GL(R;\CC)$ on $\Gr_k(\CC^R)$:
$$\FFF=E\times_{\GL(R;\CC)} \Gr_k(\CC^R).$$
The Pl{\"u}cker embedding maps $\Gr_k(\CC^R)$ in a $\GL(R;\CC)$-equivariant 
way into $\PP(\Lambda^k\CC^R)$, and the action 
of $\GL(R;\CC)$ in $\PP(\Lambda^k\CC^R)$ lifts to the obvious
action in $\Lambda^k\CC^R$. So we are in the situation described
at the beginning of this section.
Observe that the centre of $\GL(R;\CC)$ acts trivially on the Grassmannian. 
In consequence, the comments in subsection \ref{nosimple} are relevant in
this situation.

If $\omega$ is the symplectic form in $\Gr_k(\CC^R)$ 
inherited by the Fubini-Study symplectic form on $\PP(\Lambda^k\CC^R)$, 
then $\tau\omega$ also gives $\Gr_k(\CC^R)$ a Kaehler structure when 
$\tau>0$ and everything gets multiplied by $\tau$: the moment map, the
maximal weights and the integral of the moment map.
We fix from now on a constant $\tau>0$ and
we work with the symplectic form $\tau\omega$. The constant
$\tau$ can be identified with the parameter appearing in the
notion of stability and in the equations in \cite{BrGP1,DaUW}.

\subsection{Moment map of $\U(n)$ acting on the Grassmannian}
The action of $\U(n;\CC)$ on $\Gr_k(\CC^R)$ is symplectic. 
Making use of formula (\ref{momproj}) one easily verifies
that if $\pi\in\Gr_k(\CC^R)$, then the moment map of the action 
of $\U(n;\CC)$ at the point $\pi$ is the element in $\ulie(n;\CC)^*$
which sends $\xi\in\ulie(n;\CC)$ to
$$\mu(\pi)(\xi)=-\imag\tau\Tr(\pi\circ\xi),$$
where $\pi$ denotes the orthogonal projection onto $\pi$ (see
\cite{DaUW}, p. 485).

\subsection{Maximal weights of $\U(n)$ acting on the Grassmannian}
Consider the standard action of $\U(n)$ on $\PP(\Lambda^k\CC^R)$.
Take an element $s\in\ulie(n)$.
We now give the maximal weight $\lambda(v;s)$ in the
case when $v=v_1\wedge\dots\wedge v_k\neq 0$, for $v_j\in \CC^R$.
This case is enough for our purposes, since the image of the Grassmanian
$\Gr_k(\CC^R)$ given by the Pl{\"u}cker embedding into $\Lambda^k\CC^R$ 
is precisely the set of points of that form. 

Let $\pi$ be the $k$-subspace of $\CC^R$ spanned by $\{v_j\}$.
Let $\lambda_1<\dots<\lambda_r$ be the eigenvalues of
$\imag s$ acting on $\Lambda^k\CC^R$, 
and for any $1\leq j\leq r$ write 
$E_j=\bigoplus_{i\leq j}\Ker(\imag s-\lambda_k\Id)$.
Set $\alpha_j=\lambda_j-\lambda_{j+1}$. Then
\begin{equation}
\lambda(v;s)=\tau\left(
\dim(\pi)\lambda_r+\sum_{j=1}^{r-1}\dim(\pi\cap E_j)\alpha_j\right).
\label{grass}
\end{equation}
The proof of this formula is an easy exercise which follows
from lemma \ref{pesmaxproj}.

\subsection{Simple extensions}
Reasoning similarly as in lemma \ref{bansim} one can prove this
\begin{lemma}
The pair $(A,\Phi)$ is not simple if and only if one can find
a holomorphic (with respect to $\overline{\partial}_A$)
splitting $V=V'\oplus V''$ such that the subbundle $V_0$
given by the section $\Phi$ is contained in $V'$.
\end{lemma}

\subsection{The stability condition}

Let $c\in\RR$ be a real number. Fix a pair $(A,\Phi)$,
which gives a holomorphic structure on $V$ and an inclusion of bundles
$V_0\subset V.$ 
In this section we will study the $-\imag c\Id$-stability condition
for the pair in terms of $V_0\subset V$.

A (holomorphic) parabolic reduction $\sigma$ of the structure 
group of $E$ is the same as giving a (holomorphic) filtration 
$0\subset V^1\subset\dots\subset V^{r-1}\subset V^r=V,$
and an antidominant character $\chi$ for this reduction is of the form
$$\chi=z\Id+\sum_{j=1}^{r-1}m_j\lambda_{R^j},$$
where $R^j=\rk(V^j)$,
$\lambda_{R^j}=\pi_{\CC^{R^j}}-\frac{R^j}{R}\Id$
($\pi_{\CC^{R^j}}$ is the projection onto $\CC^{R^j}$),
$z$ is any real number and the $m_j$ are real negative numbers.
Taking into account that the auxiliar representation is just
the standard representation of $\GL(n;\CC)$ in $\CC^R$ we deduce that the                   
degree of the pair $(\sigma,\chi)$ is
$$\deg(\sigma,\chi)=z\deg(V)+\sum_{j=1}^{r-1}m_j
\left(\deg(V^j)-\frac{R^j}{R}\deg(V)\right).$$

To calculate the maximal weight of the action of $\chi$ on the section
$\Phi$ we use formula (\ref{grass}). The parameters that appear there
are related to ours as follows: $\alpha_j=m_j$ for any $1\leq j\leq r-1$
and $\lambda_r=z-\sum_{j=1}^{r-1}m_j\frac{R^j}{R}$.
We get, after integration
(recall that the volume of $X$ has been normalized to 1):
\begin{equation}
\int_{x\in X}
\mu(\Phi(x);-g_{\sigma,\chi}(x))=
\rk(V_0)\left(z-\sum_{j=1}^{r-1}m_j\frac{R^j}{R}\right)
+\sum_{j=1}^{r-1}m_j\rk(V_0\cap V^j).
\end{equation}
Hence, the stability notion is as follows: for any filtration
$0\subset V^1\subset\dots\subset V^{r-1}\subset V^r=V$ and any set
of negative weights $\alpha_1,\dots,\alpha_{r-1}$ we must have
\begin{align}
0 &< z\deg(V)+\sum_{j=1}^{r-1}m_j
\left(\deg(V^j)-\frac{R^j}{R}\deg(V)\right) \notag \\
  &+ \tau\left(\rk(V_0)
\left(z-\sum_{j=1}^{r-1}m_j\frac{R^j}{R}\right)
+\sum_{j=1}^{r-1}m_j\rk(V_0\cap V^j)\right)-zc\ R.
\label{st1}
\end{align}
(Observe that thanks to our assumption that $\Vol(X)=1$,
$\la\imag\chi,c\ra\Vol(X)=-zcR$.)
If this is to be satisfied by all possible choices of $z$, then
$$c=\frac{\deg(V)+\tau\rk(V_0)}{R}.$$
So, given the symplectic form $\tau\omega$, there is a unique central 
element $c\in\ulie(n;\CC)$ such that the pair can be
$c$-stable.
Putting the value of the central element inside (\ref{st1}) we get
\begin{align}
0 &< \sum_{j=1}^{r-1} m_j\left(
\deg(V^j)-\frac{R^j}{R}\deg{V}-\tau\rk(V_0)
\frac{R^j}{R}+\tau\rk(V_0\cap V^j)\right) \notag \\
&= \sum_{j=1}^{r-1} m_jR^j\left(
\frac{\deg(V^j)+\tau\rk(V_0\cap V^j)}{R^j}
-\frac{\deg(V)+\tau\rk(V_0)}{R}\right),\notag 
\end{align}
and using the fact that 
the numbers $m_j$ are arbitrary negative numbers, we see that
a necessary and sufficient condition for $(E,\Phi)$ to be stable
is that for any nonzero proper subbundle (in fact, reflexive subsheave)
$V^1\subset V$
$$\frac{\deg(V^1)+\tau\rk(V_0\cap V^1)}{\rk(V^1)}
<\frac{\deg(V)+\tau\rk(V_0)}{R},$$
and this is exactly the same condition that appears in \cite{DaUW, BrGP1}. 

In what concerns the equations, they are exactly those in \cite{DaUW}. 
Instead of writing them in terms of a gauge
transformation, we will put as the variable a metric $h$ in the 
bundle $V$. This is equivalent to our setting, since the relevant
space in our case is the gauge group of complex transformations modulo
unitary gauge transformations, and this coset space can be identified
with the space of metrics. Taking into account the precise form of the
moment map for the action of $\GL(n;\CC)$ in $\Gr_k(\CC^R)$
we can write the equations as follows:
$$\Lambda F_A-\imag\tau\pi^h_{V_0}=-\imag c\Id,$$
where $\pi^h_{V_0}$ is the $h$-orthogonal projection onto
$V_0$. 
The equations considered in \cite{BrGP1}
are written in a different way, but in \cite{DaUW} it is proved that they
are equivalent to the ones considered here.

\subsection{Filtrations}
Here we generalise the preceeding results to
the case of filtrations\footnote{The results of this subsection
were first proved by Luis {\'A}lvarez C{\'o}nsul \cite{Al,AlGP}.}.
Our trick is to identify a filtration
$0\subset V_1\subset\dots\subset V_s\subset V$
with a section $\Phi$ of the associated bundle with fibre the flag 
manifold $F_{i_1,\dots,i_s}$, where $i_k=\rk(V_k)$.
This manifold is embedded in a product of Grassmannians.
The Kaehler structure in the flag manifold is not unique.
We can in fact take as symplectic form any weighted sum 
of the pullbacks of the symplectic forms in the Grassmannians,
provided the weights are positive. So the Kaehler structure
depends on a $s$-uple of positive parameters $\tau=(\tau_1,\dots,\tau_s)$.
We can now work out the stability notion analogously to the case of
extensions, and obtain that (here we write 
$0\subset V_1\subset\dots\subset V_s\subset V$ for the filtration
represented by the section $\Phi$)
\begin{itemize}
\item the equation is 
$\Lambda F_A-\imag\sum\tau_k\pi^h_{V^k}=-\imag c\Id$, 
where $\pi^h_{V^k}$ is the $h$-orthogonal projection onto
$V^k$ and where $c$ is a real constant;
\item the pair $(A,\Phi)$ is simple unless there exists a holomorphic
(with respect to $\ov{\partial}_A$) splitting 
$V=V'\oplus V''$ such that $V_k\subset V'$ for any $k\leq s$;
\item the only value of $c$ for which we can expect our filtration
to be $c$-stable is  
$$c=\frac{\deg(V)+\sum\tau_k\rk(V_k)}{R};$$
\item the stability notion is as follows: for any nonzero proper
reflexive subsheaf $V^1\subset V$,
$$\frac{\deg(V^1)+\sum\tau_k\rk(V_k\cap V^1)}{\rk(V^1)}<
\frac{\deg(V)+\sum\tau_k\rk(V_k)}{R}.$$
\end{itemize}

\subsection{Bogomolov inequality}
In this subsection we state the Bogomolov inequality given in 
corollary \ref{bogom} for the case of filtrations. 
For that we need to compute the cohomology class
$\Phi^*\phi_A(\ov{\omega}_F)$.

We begin with some general observations.
When the cohomology class represented by the symplectic form 
$\omega_F$ of $F$ belongs to $H^2(F;\imag 2\pi\ZZ)$, there
exists a line bundle $L\to F$ with a connection $\nabla$ whose curvature
coincides with $-\imag\omega_F$. Assume that the action of
$K$ on $F$ lifts to a linear action on $L$. Then $\nabla$ can
be assumed to be $K$-equivariant (by just averaging if it is not).
Using the action of $K$ on $L$ we can define a line bundle $\LLL\to\FFF$
as $\LLL=E\times_K L$. Denote $\plx:\LLL\to X$ and
$\plf:\LLL\to \FFF$ the projections. Let $A$ be a connection on $E$. 
The connection $A$ induces a connection on the associated bundle
$\LLL$, which may be seen as a projection 
$\alpha:T\LLL\to\Ker d\plx$. Since $\nabla$ is $K$-equivariant,
we may extend it fibrewise to obtain a projection
$\beta:\Ker d\plx \to\Ker d\plf$. The composition 
$\gamma=\beta\circ\alpha:T\LLL\to \Ker d\plf$ defines a connection
$\nabla^A$ on $\LLL\to\FFF$. It is an exercise to verify that 
$$\phi_A(\ov{\omega}_F)=\imag F_{\nabla^A},$$ 
where $F_{\nabla^A}$ is the curvature of $\nabla^A$.

If $F=\Gr_k(\CC^R)$ is a Grassmannian everything in the preceeding
paragraph works. In particular, the line bundle $L\to F$ can
be identified with the dual of the determinant bundle, that is,
with the line bundle whose fibre on $V\in\Gr_k(\CC^R)$ 
is $\Lambda^k V^*$.

Using this observations, it turns out that in the general case
in which $F=F_{i_1,\dots,i_s}$ and in which $F$ has the Kaehler
structure induced by the parameters $\tau=(\tau_1,\dots,\tau_s)$
(see the preceeding subsection), then for any 
$(A,\Phi)\in\AAA^{1,1}\times\SSS$ we have
$$\int_X\Phi^*\phi_A(\ov{\omega}_F)\wedge\omega^{[n-1]}
=-\sum_{k=1}^s \tau_k\deg(V_k),$$
where $V_1\subset\dots\subset V_s\subset V$ is the filtration
represented by the section $\Phi$.

Finally, one computes
$$\int_X B(F_A,F_A)\wedge\omega^{[n-2]}=
8\pi\la ch_2(V)\cup[\omega^{[n-2]}],[X]\ra,$$
where $ch_2(V)\in H^4(X;\RR)$ is the degree 4 piece of the Chern
character of $V$ (see p. 209 in \cite{Br2}).
So corollary \ref{bogom} takes the following form in this case:

\begin{corollary}
Let $A$ be a connection on $E$, and consider a filtration
$0\subset V_1\subset\dots\subset V_s\subset V$ which is holomorphic
with respect to $\ov{\partial}_A$. Let us write $\Phi$
for the section of $\FFF$ which represents this filtration.
If the pair $(A,\Phi)$ is $\GGG_G$ equivalent to a solution of
$$\Lambda F_A-\imag\sum\tau_k\pi^h_{V^k}=-\imag c\Id,$$
then the following holds
$$\deg(V)\left(\frac{\deg(V)+\sum\tau_k\rk(V_k)}{R}\right)
-\sum_{k=1}^s\tau_k\deg(V_k)
-4\pi\la ch_2(V)\cup[\omega^{[n-2]}],[X]\ra\geq 0.$$
\end{corollary}

\section{A trivial example of stable pair}

Although we have studied some examples to which we can apply our
correspondence, we have still not proved that there exist stable
pairs. This could be achieved by studying a little bit the extensions 
of a stable vector bundle, or, still easier, by taking a
rank two projective bundle on a Riemann surface coming from
a stable rank two bundle and picking a $\tau$ small enough
(in that case the stability of a pair is equivalent to the stability
of the bundle). Here, however, we state a general result concerning
the stability of pairs whose bundle and connection are the trivial ones.

Consider a representation $\rho:G\to GL(W)$ and a $G$-principal
bundle $E\to X$. We will take $F=\PP(W)$. Let $V=E\times_{\rho}W$.
So $\FFF=\PP(E\times_{\rho}W)=\PP(V)$. Denote by $\PP(W)^s\subset \PP(W)$ 
(resp. $\PP(W)^{ss}\subset \PP(W)$) the set of stable (resp. semistable) 
points by the action of $G$. Since a point is stable if and only if so
is any point in its orbit, it makes sense to define
$\PP(V)^s=\amalg_{x\in X} \PP(V)^s_x$
and $\PP(V)^{ss}=\amalg_{x\in X} \PP(V)^{ss}_x$, where 
$\PP(V)_x\simeq\PP(W)$ is the fibre over $x\in X$. 

On the other hand, there is a notion of stability for $G$-principal
bundles due to Ramanathan \cite{R1}. This can be stated using
our notation as follows: $E$ is stable if, for any reduction 
$\sigma\in\Omega^0(E(G/P))$ of the structure group of $E$
to a parabolic subgroup $P\subset G$, and any antidominant character
$\chi$ of $P$, $\deg(\sigma,\chi)>0$. When we only have 
$\deg(\sigma,\chi)\geq 0$, then we say $E$ is semistable.
In fact, our correspondence applies to this case with some
due modifications in the proof, and in particular if
there exists a reduction $h\in\Omega^0(E(G/K))$ of the structure
group to a maximal compact subgroup $K\subset G$ such that
$\Lambda F_h=0$, then $E$ is semistable (see \cite{RS}).
So, for example, the trivial bundle $E=X\times G$ is semistable
(just take a constant section $h\in\Omega^0(E(G/K))$).

Finally, recall that the Kempf-Ness theory (see section \ref{KN})
tells us that $x\in\PP(W)$ is stable (resp. semistable) if and 
only if for any antidominant character $\chi$ of a parabolic subgroup 
of $G$, the maximal weight $\mu(x;\chi)>0$ (resp. $\mu(x;\chi)\geq 0$). 

Putting together all this we obtain the following

\begin{theorem}
Suppose $E$ is semistable and $\Phi\subset\PP(V)^s$. Then the pair
$(E,\Phi)$ is stable for any $\tau>0$. And if 
$E$ is stable and $\Phi\subset\PP(V)^{ss}$, then the pair
$(E,\Phi)$ is also stable for any $\tau>0$.
\end{theorem}

Now, taking $E=X\times G$ and $\rho$ such that
$\PP(W)^s\neq\emptyset$, we can pick a constant section 
$\Phi(x)=w\in \PP(W)^s$ and then the pair $(E,\Phi)$ will
be stable thanks to the preceeding theorem. Furthermore, we
can also chose $w$ such that the pair is simple, by taking
it outside any proper $G$ invariant subspace $W'\subset W$.


\chapter{The moduli space}
\label{moduli}
In this chapter and in all the remaining ones we will assume that $X$
is a Riemann surface (with a fixed Riemannian metric). So from now
on we will forget the third equation $F_A^{0,2}=0$ in (\ref{equs}),
which, as we have already said, is trivially satisfied. We will also
assume henceforth that $F$ is compact.

Our aim in the next chapters is to use the space of solutions to 
(a certain perturbation of) equations (\ref{equs}) to define invariants
of the symplectic manifold $F$ and the action of $S^1$. As a first
step, in this chapter we will construct the moduli space of gauge equivalence
classes of solutions to equations (\ref{equs}). 
The methods used in the construction are rather standard (see for
example \cite{DoKr, FrUh, McDS1}, and consequently
at some steps we will just give a sketch. At some points in our discussion
we will make the assumption that $K=S^1$ and that its action on $F$ is 
almost-free. However, some of the results remain valid in greater
generality.

We begin by fixing Sobolev completions of our ambient space $\AAA\times\SSS$.
This will allow us to use Banach manifold techniques as the implicit function
theorem. Then we define the different moduli spaces appearing in the
thesis. It is important to observe that rather than using equations 
(\ref{equs}), we consider suitable perturbations of them (\ref{sequs}).
This is done because we want to get smooth moduli, and without perturbing 
the equations we can not assure smoothness in general. We compute the 
dimension of the moduli and we prove that the moduli spaces obtained from 
different perturbations are cobordant.

\section{Sobolev completions}

Let $E\to X$ be a principal $K$-bundle and let 
$\kE=E\times_{\Ad}\klie$. Fix a real number $p>2$.
We will consider the completion $\AAA_{L^p_1}$ of the space of connections 
$\AAA=\AAA^E$ on $E$ with respect to the $L^p_1$ norm. This is defined by 
using a fixed smooth connection $A_0\in\AAA$ and then putting 
$\AAA_{L^p_1}=A_0+\Omega^1(\kE)_{L^p_1}$. The space $\AAA_{L^p_1}$ is a 
Banach manifold, which is independent of the particular choice of $A_0$. 
In section \ref{THCshol} we saw how to construct a complex structure $I(A)$  
on $\FFF$ from any connection $A$ on $E$. The same thing can be done
for connections $A$ lying in $\AAA_{L^p_k}$. We obtain the following
result.

\begin{lemma}
Suppose that $X=\DD$ is the unit disk and that
we have a trivialisation $E\simeq K\times \DD\to\DD$. Let us take 
a connection $A=d+\alpha$, where $\alpha\in\Omega^1(\kE)_{L^p_k}$.
Then the complex structure $I(A)$ on $\FFF\simeq F\times\DD$
lies in $\Omega^0(\End T(F\times\DD))_{L^p_k}$.
\label{regAregI}
\end{lemma}

Let $\FFF=\FFF^E=E\times_K F$ be the associated bundle and let 
$\SSS=\SSS^E=\Gamma(\FFF)$. 
Take any embedding $\iota:\FFF\hookrightarrow \RR^N$. 
We define the distance $d_{L^p_1}$ between two sections $\Phi$ and $\Phi'$
to be the sum of the $L^p_1$ norms of the difference of the components
of $\iota\circ\Phi$ and $\iota\circ\Phi'$. This is a metric on $\SSS$.
We consider the completion $\SSS_{L^p_1}$ of $\SSS$ with respect to the 
metric $d_{L^p_1}$. The space $\SSS_{L^p_1}$ is a Banach manifold.
By our choice of $p$ we have a compact embedding
$L^p_1\hookrightarrow C^0$. Consequently, two nearby elements in 
$\SSS$ with respect to $d_{L^p_1}$ are nearby pointwise. This implies
that all the elements in $\SSS_{L^p_1}$ are continuous sections. 
Furthermore, the completion $\SSS_{L^p_1}$ is independent of
the embedding $\iota$. This stems from the fact that
any smooth map $F:V\to W$ of vector bundles over an $n$-dimensional 
manifold which fixes the zero section induces a continuous map
from $\Omega^0(V)_{L^p_k}$ to $\Omega^0(W)_{L^p_k}$ whenever
$pk>n$.

Since there is a Sobolev multiplication $L^p_1\otimes L^p_1\to L^p$, 
for any section $\Phi\in\SSS_{L^p_1}$ and any connection 
$A\in\AAA_{L^p_1}$ the covariant derivative
$d_A\Phi$ lies in $\Omega^0(\Phi^*T\FFF_v)_{L^p}$
and $\overline{\partial}_A\Phi$ lies in $\Omega^{0,1}(\Phi^*T\FFF_v)_{L^p}$.

In a similar way, for any manifold $M$ we define a metric
$d_{L^p_1}$ on $\Map(X,M)$ by using an embedding $M\hookrightarrow\RR^N$
and denote $\Map(X,M)_{L^p_1}$ the completion.

Finally, we consider the completion $\GGG_{L^p_2}$ of the gauge group
$\GGG$ with respect to the                                      
$L^p_2$ norm. The group $\GGG_{L^p_2}$ is a Banach Lie group and
it acts smoothly on $\AAA_{L^p_1}$ and on $\SSS_{L^p_1}$.
Its Lie algebra is
$\Lie(\GGG_{L^p_2})=\Omega^0(\kE)_{L^p_2}$.

\section{The moduli spaces}

We begin introducing some notation. Let $V$ and $W$ be two
complex vector spaces. We will denote $\Hom^{1,0}(V,W)$
(resp. $\Hom^{0,1}(V,W)$) the set of complex linear
(resp. complex antilinear) maps from $V$ to $W$. We obviously
have $$\Hom_{\RR}(V,W)=\Hom^{1,0}(V,W)\oplus \Hom^{0,1}(V,W).$$

\subsection{Moduli of $\sigma$-holomorphic curves}
\label{modcorbes}
Let $M$ be a compact almost Kaehler manifold. Let
$\Sigma_M=\Hom^{0,1}(TX,TM)$ (these are sections
on $X\times M$, and the vector bundles should be taken
to be the pullbacks by the two projections).
Let $B\in H_2(M;\ZZ)$ be any class and let $\sigma\in\Sigma_M$. 
We define $$\MMM_{\sigma}(B)=\MMM_{\sigma}^{M}(B)=\{\Phi\in\Map(X,M)_{L^p_1}
|\ \overline{\partial}\Phi=\sigma,\ \Phi_*[X]=B\}.$$
This is the moduli of $\sigma$-perturbed holomorphic curves on $M$.
Following Ruan we define for any $\sigma\in\Sigma_M$ a
complex structure $I_{\sigma}$ on $X\times M$ as 
$$I_{\sigma}=\left(\begin{array}{cc}I_X & 0 \\ \sigma & I_M\end{array}
\right),$$
where $I_X$ and $I_M$ are the complex structures of $X$ and $M$ and where 
the matrix is 
given with respect to the splitting $T(X\times M)=TX\oplus TM$ (as always, 
we omit the pullbacks). One can prove the following lemma (see lemma 
3.1.1 in \cite{Ru}), which allows to view perturbed holomorphic curves as 
genuine holomorphic curves in $X\times M$.

\begin{lemma}
A map $\Phi:X\to M$ satisfies $\overline{\partial}\Phi=\sigma$
if and only if the map $\Phi^{\id}=(\id,\Phi):X\to X\times M$
is holomorphic with respect to the complex structure $I_{\sigma}$.
\label{genui}
\end{lemma}

When the perturbation $\sigma$ is zero we will usually write
$\MMM^{M}(B)$ instead of $\MMM^{M}_{\sigma}(B)$.

\subsection{Moduli of $\sigma$-twisted holomorphic curves}

Recall that we denote by $EK\to BK$ the universal principal $K$-bundle
and $F_K=EK\times_K F$ the Borel construction of $F$. The equivariant
(co)homology of $F$ is by definition the (co)homology of $F_K$.
Denote $\pi_F:F_K\to BK$ the projection.

\subsubsection{The space of perturbations}
Let $E\to X$ be a principal $K$-bundle and let $\FFF=\FFF^E$.
Let $$\Sigma'(E)=\Hom^{0,1}(\pfx^*TX,T\FFF_v)\oplus\Omega^0(\kE),$$
where $\pfx:\FFF\to X$ is the projection.
The gauge group $\GGG$ of $E$ acts on $\Sigma'(E)$, and we set
$\Sigma(E)$ to be the fixed elements. Eventually, we will consider
the completion of $\Sigma(E)$ with respect to suitable $C^l$ norms.
Observe that if $(\sigma_1,\sigma_2)\in\Sigma(E)$, then
$\sigma_2\in\Omega^0(Z_\klie)$, where $Z_{\klie}$ is the center of $\klie$.

\subsubsection{The moduli space}
\label{themodulispace}
Let $E\to X$ be a principal $K$-bundle, and let $\GGG$ be its gauge group.
Let us fix a pair $\sigma=(\sigma_1,\sigma_2)\in\Sigma(E)$ and a central
element $c\in\klie$. We will consider the space of pairs 
$(A,\Phi)\in\AAA^E\times\SSS^E$ which
satisfy the following two equations
\begin{equation}
\left\{
\begin{array}{l}
\overline{\partial}_A\Phi = \sigma_1  \\
\Lambda F_A+\mu(\Phi) = c+\sigma_2. 
\end{array}\right.
\label{sequs}
\end{equation}
We will call any pair satisfying these equations a
{\bf $\sigma$-twisted holomorphic curve over $X$} 
({\bf $\sigma$-THC} for short). When $\sigma=0$ we will call the solutions
twisted holomorphic curves or THCs.
For any homotopy class of sections $[\Phi_0]\in\SSS^E$ 
we will write the space of $\sigma$-THCs $(A,\Phi)$ such that
$[\Phi]=[\Phi_0]$ as
$$\bM_{\sigma}(E,[\Phi_0],c)=\bM^{F,K}_{\sigma}(E,[\Phi_0],c).$$
Since the sections $(\sigma_1,\sigma_2)$ are gauge 
invariant and the complex structure is also invariant
under the action of $K$, it turns out that the space 
$\bM_{\sigma}(E,[\Phi_0],c)$
is invariant under the action of the identity component of the gauge
group. However, it is not necessarily invariant under the action of
the full gauge group, since in general there may exist a section
$\Phi$ and a gauge transformation $g\in\GGG$ such that
$\Phi$ and $g\Phi$ are not homotopic. (This, of course, does
not happen when $F$ is a vector space.)

Let $P_K(X)$ be the set of (topological isomorphism classes of)
$K$ principal bundles over $X$. Let $\eta:P_K(X)\to H_2(BK;\ZZ)$
be the map which sends a bundle $E\to X$ to ${c_E}_*[X]$, where
$c_E:X\to BK$ is the classifying map of $E$. In lemma \ref{etabij}
of the appendix we prove that $\eta$ is a bijection.
There exists a map ${\rho_E}_*:H_*(\FFF)\to H_*(F_K)$ which is invariant
under the action of the gauge group $\GGG$ on $H_*(\FFF)$ and which
lifts ${c_E}_*$: it is defined by fixing an isomorphism 
$\phi:E\simeq c_E^*EK$, taking the induced isomorphism
$\psi:\FFF\simeq c_E^*F_K$ and putting ${\rho_E}_*=\psi_*$
(see lemma \ref{indepframe} in the appendix for a proof that this
is independent of $\phi$ and that, consequently, this map is
$\GGG$ invariant).

Let $B\in H_2(F_K)$ be a class such that ${\pi_F}_*B=\eta(E)$. We define
\begin{align*}\bM_{\sigma}(B,c)&=
\bM^{F,K}_{\sigma}(B,c)=\{(A,\Phi)\in\AAA\times\SSS\mid   
\mbox{ satisfying (\ref{sequs}) and ${\rho_E}_*\Phi_*[X]=B$}\} \\
&= \coprod_{{\rho_E}_*(\Phi_0)_*[X]=B}\bM^{F,K}_{\sigma}(E,[\Phi_0],c).
\end{align*}
The space $\bM_{\sigma}(B,c)$ is invariant under the action of $\GGG$. 
We define the {\bf moduli space of $\sigma$-THCs} to be the quotient
$$\MMM_{\sigma}(B,c)=\MMM^{F,K}_{\sigma}(B,c)=\bM^{F,K}_{\sigma}(B,c)/\GGG.$$

\subsection{Extended moduli of $\sigma$-twisted holomorphic curves}
Here we keep the notation of the preceeding section. Let us
fix a base point $x_0\in X$, and let $\GGG_0=\{g\in\GGG\mid g(x_0)=1\}$.
Note that $\GGG/\GGG_0=K$ and that $\GGG_0$ acts freely on the
space of connections $\AAA^E$.

For any equivariant homology class $B\in H_2(F_K;\ZZ)$ we define
the {\bf extended moduli space of $\sigma$-THCs} to be
$$\NNN_{\sigma}(B,c)=\NNN^{F,K}_{\sigma}(B,c)=
\bM^{F,K}_{\sigma}(B,c)/\GGG_0.$$

\subsection{The complex structure}
All the moduli spaces that we have defined above depend on the complex structure
of the almost complex manifold $F$. Later it will be convenient to stress this
dependence, and we will specify the complex structure with a subscript.
So for a complex structure $I$ on $F$ we will write
$\MMM_{I,\sigma}(B)$, $\MMM_{I,\sigma}(B,c)$
and $\NNN_{I,\sigma}(B,c)$ (note that in the first moduli
$B$ is a homology class of $F$, whereas in the other ones it is
a homology class of $F_K$).

\section{Local structure}

\subsection{The deformation complex}

Let $[(A,\Phi)]\in\MMM_{\sigma}(B,c)$ be any gauge equivalence class
and consider the sequence of maps
\begin{equation}
\Com_{A,\Phi}:C^0_{A,\Phi}\stackrel{d_1}{\longrightarrow}
C^1_{A,\Phi}\stackrel{d_2}{\longrightarrow}
C^2_{A,\Phi},
\label{defo}
\end{equation}
where (recall that we denote $\kE=E\times_{\Ad}\klie$)
\begin{align*}
C^0_{A,\Phi} &= \Omega^0(\kE)_{L^p_2}, \\
C^1_{A,\Phi} &= \Omega^1(\kE)_{L^p_1}\oplus
\Omega^0(\Phi^*T\FFF_v)_{L^p_1}, \\
C^2_{A,\Phi} &= \Omega^2(\kE)_{L^p}\oplus
\Omega^{0,1}(\Phi^*T\FFF_v)_{L^p}, 
\end{align*}
and where $d_1$ is the infinitesimal action of $\GGG_{L^p_2}$
(recall that $\Lie(\GGG_{L^p_2})=\Omega^0(\kE)_{L^p_2}$), 
and $d_2$ is the linearisation of equations (\ref{sequs}). 
More precisely, for any $\theta\in\Omega^0(\kE)_{L^p_2}$
and for any $(\alpha,\phi)\in
\Omega^1(\kE)_{L^p_1}\oplus\Omega^0(\Phi^*T\FFF_v)_{L^p_1}$
we have
\begin{align}
d_1(\theta) &= -d_A\theta+\theta\cdot\Phi, \notag \\
d_2(\alpha,\phi) &= 
\left(
\begin{array}{c}
\Lambda d_A\alpha+\la d\mu(\Phi),\phi\ra_{T\FFF}\\
\overline{\partial}_{A,\nabla}\phi+\alpha\cdot\phi+
C(A,\Phi,\sigma)(\alpha,\phi)
\end{array}\right).
\label{difeq}
\end{align}
The operator $\overline{\partial}_{A,\nabla}$ is the composition of 
the covariant derivative
$$d_{A,\nabla}:\Omega^0(\Phi^*T\FFF_v)_{L^p_1}
\to\Omega^1(\Phi^*T\FFF_v)_{L^p}$$
(see section \ref{vbfb} in the appendix) with the projection 
$\Omega^1(\Phi^*T\FFF_v)\to \Omega^{0,1}(\Phi^*T\FFF_v)$
(it is a Cauchy-Riemann operator in the vertical direction).

$C$ is a compact operator which depends on the derivative
of $\sigma_1$, on the connection $\nabla$, on $TF$, and
on how we identify a neighbourhood of $\Phi\in\Gamma(\FFF)$
with a neigbourhood of the zero section of $\Omega^0(\Phi^*T\FFF_v)$.
If this identification is made through the exponential map, 
$\sigma_1=0$ and $\nabla$ has torsion equal to $\frac{1}{4}N_J$, where 
$N_J$ is the Nijenhuis operator of $F$, then the dependence of $C$ on 
$\nabla$ can be given in terms of $N_J$ (see \cite{McDS1} p. 28).
We will not give a precise form of it because it is unnecessary for our
purposes. The point is that the operator 
$\overline{\partial}_{A,\nabla}\phi+\alpha\cdot\phi$ has the same
symbol as the (vertical) Cauchy-Riemann equation, and hence is 
elliptic and its index can be computed (see below).

On the other hand, $d\mu(\Phi)$ denotes the section of 
$\Phi^*T\FFF_v\otimes\kE^*$ which arises from extending globally the 
derivative $d\mu\in\Omega^1(F;\klie^*)$ of the moment map (this is possible
thanks to the $K$-equivariance of $\mu$).

Using the fact that $(A,\Phi)$ solves the equations
(\ref{sequs}), one can prove that (\ref{defo}) is a complex, that is, 
$d_2\circ d_1=0$. (Just apply the chain rule to the identity
expressing the gauge invariance of the set of solutions to
(\ref{sequs}).)
The complex (\ref{defo}) is called the deformation complex of 
$\MMM_{\sigma}(B,c)$ at $[(A,\Phi)]$. 
Denote $H^0_{A,\Phi}$, $H^1_{A,\Phi}$ and $H^2_{A,\Phi}$ 
its cohomology groups. 
Using again the fact that the pair $(A,\Phi)$ solves (\ref{sequs})
one can prove that the complex (\ref{defo}) is elliptic,
so $H^0_{A,\Phi}$, $H^1_{A,\Phi}$ and $H^2_{A,\Phi}$ are
finite dimensional vector spaces. We will use them to
give local models of the moduli space $\MMM_{\sigma}(B,c)$.

\subsection{Index of the deformation complex}
\label{locals}
Since the complex (\ref{defo}) is elliptic, the operator 
$$d_1^*+d_2:C_1\to C_0\oplus C_2$$
is Fredholm, and so it has a well defined index
$$\Ind(d_1^*+d_2)=-\Ind(\Com_{A,\Phi})=
\dim H^1_{A,\Phi}-(\dim H^0_{A,\Phi}+\dim H^2_{A,\Phi}).$$
This integer can be computed by means of Atiyah-Singer index theorem.
It is easier, however, to deform the operator $d_1^*+d_2$ by
adding to it a compact operator, and then compute the index of 
the resulting operator (which will coincide with that of
$d_1^*+d_2$). So we take instead of $d_1^*+d_2$ the operator
$$D(\alpha,\phi)=(-{d_A}^*\alpha,\Lambda d_A\alpha,
\overline{\partial}_{A,\nabla}\phi)$$
In other words, we are splitting the complex (\ref{defo}) as the sum
of these two complexes
$$\Com_A:\Omega^0(\kE)_{L^p_2} @>{d_A}>> \Omega^1(\kE)_{L^p_1}
@>{d_A}>> \Omega^2(\kE)_{L^p}$$
and
$$\Com_{\Phi}: 0 @>>> \Omega^0(\Phi^*T\FFF_v)_{L^p_1} 
@>{\overline{\partial}_{A,\nabla}}>> \Omega^{0,1}(\Phi^*T\FFF_v)_{L^p}$$
(note that we have changed the sign of the first map in the complex 
$\Com_A$ and we have omitted the contraction $\Lambda$ in the second
map; this is irrelevant when computing the index of the complex).

\subsubsection{The index of the complex $\Com_A$}
From now on we will often omit the subscripts denoting Sobolev completions,
which will be implicitly assumed. Let $\gE=E\times_{\Ad}\glie$.
Consider the map $f:\Omega^0(\kE)\oplus\Omega^2(\kE)\to\Omega^0(\gE)$
given by $f(\alpha,\beta)=\alpha+\imag * \beta$, where $*$ denotes the Hodge 
star operator, and the map $g:\Omega^{0,1}(\gE)\to\Omega^1(\kE)$
given by $g(\theta)=\theta+\overline{\theta}$.
Both maps are isomorphisms and rend commutative the following
diagram
\begin{equation}
\xymatrix{
\Omega^0(\kE)\oplus\Omega^2(\kE) \ar[r]^-f \ar[d]_{d_A+d_A^*} 
& \Omega^0(\gE) \ar[d]_{\overline{\partial}_A} \\
\Omega^1(\kE) & \Omega^{0,1}(\gE). \ar[l]_g }
\label{isos}
\end{equation}
This implies that the index of the complex $\Com_A$ is equal
to the index of this other complex
$$\Com_{\overline{\partial}_A}:
\Omega^0(\gE)\stackrel{\overline{\partial}_A}
{\longrightarrow}\Omega^{0,1}(\gE).$$
(We have omitted the Sobolev completions; recall that, thanks
to elliptic regularity, the index of the complexes is independent
of the chosen Sobolev completion.)
The complex $\Com_{\overline{\partial}_A}$ is the deformation 
complex for the moduli space of complex
structures on $E_G=E\times_{\Ad}G$ compatible with the complex structure
on $X$. In fact, the commutativity of diagram (\ref{isos}) is 
the infinitesimal version of the isomorphism between the
moduli space of flat connections (whose deformation complex is
$\Com_A$) and that of complex structures on $E_G$ (this isomorphism
is the Chern map; see section \ref{accioGG}).

The index of $\Com_{\overline{\partial}_A}$ is computed using Riemann-Roch, 
and one obtains
$$\Ind(\Com_A)=\Ind(\Com_{\overline{\partial}_A})=
\la c_1(\gE),[X]\ra+\dim_{\RR}\klie(1-g).$$

\subsubsection{The index of the complex $\Com_{\Phi}$}
This is again given by Riemann-Roch, and it is equal to 
$$\Ind(\Com_{\Phi})=-\la c_1(\Phi^*T\FFF_v),[X]\ra-n(1-g),$$
where $2n$ is the real dimension of $F$ (the minus sign
accounts for the fact that the complex $\Com_{\Phi}$ is
the Dolbeaut complex shifted one unit to the right). 

\subsubsection{The index of the complex $\Com_{A,\Phi}$}
We have $\Ind(d_1^*+d_2)=-\Ind(\Com_A)-\Ind(\Com_{\Phi})$. 
Summing up our results we obtain that the index 
$\dim H^1_{A,\Phi}-(\dim H^0_{A,\Phi}+\dim H^2_{A,\Phi})$
(complex dimensions are meant) of the operator $d_1^*+d_2$ is equal to
$$-\Ind(\Com_{A,\Phi})=
\la c_1(\Phi^*T\FFF_v)-c_1(\gE),[X]\ra+(n-\dim_{\RR}\klie)(1-g).$$

Using the functoriality of the Chern classes we can write
$$\la c_1(\Phi^*T\FFF_v),[X]\ra=\la c_1^K(TF),B\ra,$$
where $c_1^K(TF)\in H^2_K(F)$ is the first equivariant Chern
class of the tangent bundle $TF$, and
$$\la c_1(\gE),[X]\ra=\la (\pi^{F_K})^* c_1^K(\glie),B\ra,$$
where $c_1^K(\glie)\in H^*_K(\{\pt\})$ is the first equivariant Chern class
of the bundle $\glie\to\{\pt\}$ viewed as a $K$ bundle using the
adjoint action of $K$ on $\glie$.

So the previous formula for the index of $\Com_{A,\Phi}$ may be 
rewritten as follows
\begin{equation}
-\Ind(\Com_{A,\Phi})=
\la c_1^K(TF)-(\pi^{F_K})^* c_1^K(\glie),B\ra+(n-\dim_{\RR}\klie)(1-g).
\label{vir}
\end{equation}
This number is called the {\bf complex virtual dimension} of the moduli 
space. Under certain transversality conditions to be specified below
it coincides with the actual dimension of the moduli.

\subsection{Local models of the moduli space}
\label{secKuranishi}
In this subsection we give some results on how the deformation complex
allows to model neighbourhoods of given elements in the moduli
of $\sigma$-THCs. The obtained models are called Kuranishi models. These 
results are standard, and they are proved using the implicit function 
theorem for Banach manifolds and the fact that the deformation complex is
elliptic. For details see for example pp. 137--139 in \cite{DoKr}.

\begin{lemma}
Let $[(A,\Phi)]\in\MMM_{\sigma}(B,c)$ and denote
$\Gamma_{A,\Phi}$ the stabiliser of $(A,\Phi)$ in $\GGG$.
There exists a neighbourhood $U$ of $0\in H^1_{A,\Phi}$
and a $\Gamma_{A,\Phi}$ equivariant smooth map
$$f:U\to H^2_{A,\Phi}$$
with vanishing derivative at $0$ such that the quotient
$f^{-1}(0)/\Gamma_{A,\Phi}$ models a neighbourhood of
$[(A,\Phi)]$ in $\MMM_{\sigma}(B,c)$.
\label{Kuranishi}
\end{lemma}

\begin{lemma}
Let $[(A,\Phi)]\in\MMM_{\sigma}(B,c)$. Then
$H^0_{A,\Phi}=\Lie \Gamma_{A,\Phi}$.
\end{lemma}

\begin{corollary}
Let us take any $[(A,\Phi)]\in\MMM_{\sigma}(B,c)$, and let us
suppose that $H^0_{A,\Phi}=H^2_{A,\Phi}=0.$
Then there is a neighbourhood of $[(A,\Phi)]$ in $\MMM_{\sigma}(B,c)$
which is diffeomorphic to $\RR^N$, where $$N=\dim H^1_{A,\Phi}=
2(\la c_1^K(TF)-(\pi^{F_K})^* c_1^K(\glie),B\ra+
(n-\dim_{\RR}\klie)(1-g)).$$
The tangent vector space $T_{A,\Phi}\MMM_{\sigma}(B,c)$
can be canonically identified with $H^1_{A,\Phi}$.
In particular, the dimensional of the moduli space on a neighbourhood
of $[(A,\Phi)]$ coincides with the virtual dimension.
\label{dimensio}
\end{corollary}

\begin{corollary}
Let $[(A,\Phi)]\in\MMM_{\sigma}(B,c)$ be such that
$H^0_{A,\Phi}=H^2_{A,\Phi}=0.$ Then
$T_{A,\Phi}\MMM_{\sigma}(B,c)$ is canonically oriented.
\label{orientacio}
\end{corollary}
\begin{pf}
Since the deformation complex $\Com_{A,\Phi}$ is homotopically
equivalent
to $\Com_A\oplus\Com_{\Phi}$ an orientation of $H^1_{A,\Phi}$
is the same as an orientation of $H^1_A\oplus H^1_{\Phi}$,
where $H^*_A$ and $H^*_{\Phi}$ denote the cohomology of
the complexes $\Com_A$ and $\Com_{\Phi}$ respectively. Now, the complex
$\Com_{\Phi}$ is a complex of modules over $\CC$, so its
cohomology groups are complex vector spaces and hence have
a canonical orientation. So $H^1_{\Phi}$ is canonically
oriented. On the other hand, diagram (\ref{isos}) shows
that $H^1_A=H^1_{\ov{\partial}_A}$ canonically and the complex
$\Com_{\ov{\partial}_A}$ is also one of modules over $\CC$,
so $H^1_{\ov{\partial}_A}$ (and hence $H^1_A$) has a canonical
orientation as well.
\end{pf}

\section{Smoothness of $\MMM_{\sigma}^{F,S^1}(B,c)$ for 
semi-free $S^1$ actions}
\label{smoothness}
In the rest of this chapter we will restrict to the case $K=S^1$
and we will assume that the action of $S^1$ on $F$ is semi-free.
(This means that the action on the complementary
$F\setminus F^{S^1}$ of the fixed point set is free.) Note that some of the
results that follow, however, can be proved in greater generality.

We begin by obtaining some consequences from our assumptions. Since $K=S^1$, 
we have $Z_\klie=\klie=\imag\RR$, $c_1^K(\glie)=0$ and $\dim_{\RR}\klie=1$.
Consequently, formula (\ref{vir}) gives in this case the following 
value for the complex virtual dimension:
\begin{equation}
\la c^{S^1}_1(TF),B\ra+(n-1)(1-g).
\label{vir2}
\end{equation}
On the other hand, for any $S^1$ principal bundle $E$ the associated
bundle $\kE=E\times_{\Ad}\imag\RR$ is the trivial bundle with fibre 
$\imag\RR$.

Fix $B\in H_2(F_K;\ZZ)$ and $c\in \imag\RR$, and write
$\bM_{\sigma}(c)=\bM_{\sigma}(B,c)$ and
$\MMM_{\sigma}(c)=\MMM_{\sigma}(B,c)$. In this section we will
study the smoothness of $\MMM_{\sigma}(c)$ for generic choices 
of $c\in\imag\RR$ and $\sigma\in\Sigma$. Will omit in the sequel the 
subscripts denoting the Sobolev completions of the spaces $\AAA$, $\SSS$
and $\GGG$ (which will be implicitly assumed).

Let $E\to X$ be the unique $S^1$ principal bundle such that 
$\eta(E)={\pi_F}_*(B)$ (see lemma \ref{etabij}). Let $\AAA=\AAA^E$
and $\SSS=\SSS^E$. Then we have $\bM_{\sigma}(c)\subset\AAA\times\SSS$.
Let us write $$\FFF_0=X\times F^{S^1}=E\times_{S^1}F^{S^1}\subset\FFF
\qquad\mbox{ and }\qquad
\SSS^*=\{\Phi\in\SSS\mid \Phi(X)\nsubseteq\FFF_0\}.$$
Let 
$$C_0=\mu(F^{S^1})-2\pi\imag\frac{\deg(E)}{\Vol(X)}\subset \imag\RR.$$
This is a finite subset, since $F^{S^1}$ is a finite union 
of compact connected submanifolds and $\mu$ is locally constant
on $F^{S^1}$.

\begin{lemma}
Let $c\in\imag\RR\setminus C_0$. For any small enough (in $C^0$)
perturbation $\sigma$ we have $\bM_{\sigma}(c)\subset\AAA\times\SSS^*$.
In other words, if $(A,\Phi)$ is supported on $E$ and
satisfies equations (\ref{sequs}) then $\Phi(X)\nsubseteq \FFF_0.$
\label{nodinsfixos}
\end{lemma}
\begin{pf}
Let $\sigma=(\sigma_1,\sigma_2)\in\Sigma(E)$ be a perturbation
and suppose that the pair $(A,\Phi)$ is supported on $E$,
$\Phi(X)\subset\FFF_0$, and that 
$$\Lambda F_A+\mu(\Phi)=c+\sigma_2.$$
Since $X$ is connected and $\Phi(X)\subset\FFF_0$, $\mu(\Phi)$ takes
a constant value $c_0\in\imag\RR$.
Integrating the equation above over $X$ and using Chern-Weil theory we 
deduce 
$$c+\frac{1}{\Vol(X)}\int\sigma_2=c_0-2\pi\imag\frac{\deg(E)}{\Vol(X)},$$
and if $|\sigma_2|_{C^0}<d(c,C_0)$ this is a contradiction.
\end{pf}

From now on we will take $c\in\imag\RR\setminus C_0$.
Let $\Sigma_c(E)=\{(\sigma_1,\sigma_2)|\ |\sigma_2|_{C^0}<d(c,C_0)\}$.

\begin{corollary}
If $\sigma\in\Sigma_c(E)$ then the action of $\GGG$ on $\bM_{\sigma}(c)$
is free.
\label{eslliure}
\end{corollary}
\begin{pf}
The only points in $\AAA\times\SSS$ whose stabiliser is nontrivial
are the pairs $(A,\Phi)$ such that $\Phi(X)\subset\FFF_0$ (these
are fixed by the constant gauge transformations).
\end{pf}

\begin{lemma} 
Let $\sigma\in\Sigma_c(E)$ and $(A,\Phi)\in\bM_{\sigma}(c)$.
Then $\Phi(X)\cap \FFF_0$ is a finite set of points.
\label{esfinit}
\end{lemma}
\begin{pf}
Suppose that $\Phi(X)$ and $\FFF_0$ meet at infinite points.
Let $\theta\in S^1$. Both $\Phi$ and
$\theta\cdot\Phi$ are perturbed holomorphic curves on $\FFF$
and by lemma \ref{genui} $\Phi^{\id}$ and $(\theta\cdot\Phi)^{\id}$
are holomorphic curves on $X\times\FFF$ with respect to the complex 
structure obtained from $\sigma_1$. By our assumption 
$\Phi^{\id}$ and $(\theta\cdot\Phi)^{\id}$ meet in an infinite
number of points. But since both curves are everywhere injective,
a result of McDuff (see lemma \ref{interseccio}) implies that they
have the same image. So $\Phi=\theta\cdot\Phi$ and since 
this is true for any $\theta\in S^1$, we conclude that
$\Phi(X)$ must be included in $\FFF_0$.
\end{pf}

\begin{theorem}
Let $c\in\imag\RR\setminus C_0$. There is a subset
$\Sigma_c^{\reg}(E)\subset\Sigma_c(E)$ of
Baire of the second category\footnote{We recall that a set of Baire
of the second category is by definition any countable intersection of 
dense open subsets of a topological space.} 
such that for any $\sigma\in\Sigma_c^{\reg}(E)$ the moduli
of $\sigma$-THCs $\MMM_{\sigma}(c)$ is a smooth oriented manifold of real
dimension equal to $$2\la c^K_1(TF),B\ra+2(n-1)(1-g).$$
\label{gensigmasmooth}
\end{theorem}
\begin{pf}
Our proof will follow the ideas of similar results in Donaldson
theory and Gromov theory (see for example \cite{FrUh,McDS1}). Let us take
a big 
positive integer $l>0$ (later on we will specify how big $l$ has to be)
and consider the completion $\Sigma_c(E)^l$
of $\Sigma_c(E)$ with respect to the $C^l$ topology.
Let $\BBB^l=\AAA\times\SSS^*\times\Sigma_c(E)^l$. We consider the universal
set of THCs
$$\bM_{\Sigma}(c)^l=\{(A,\Phi,\sigma)\mid \sigma\in\Sigma_c(E)^l
\mbox{, $(A,\Phi)$ is a $\sigma$-THC and }{\rho_E}_*\Phi_*[X]=B\}.$$
(this is a subset of $\BBB^l$).
We prove that $\bM_{\Sigma}(c)^l$ is a smooth Banach manifold as follows.
There is a Banach vector bundle $\WWW\to\BBB^l$ whose fibre over
$(A,\Phi,\sigma)$ is $\Omega^{0,1}(\Phi^*T\FFF_v)\oplus\Omega^1(\imag\RR)$
and a section $\SS=\SS_{c,\sigma}:\BBB^l\to\WWW$ which sends any
$(A,\Phi,\sigma)$ 
to $(\ov{\partial}_A\Phi-\sigma_1,\Lambda F_A+\mu(\Phi)-c-\sigma_2)$.
We have by definition $\bM_{\Sigma}(c)^l=\SS^{-1}(0)$.
To prove that $\bM_{\Sigma}(c)^l$ is smooth it is enough to verify
that $\SS$ is transverse to the zero section.

So let $(A,\Phi,\sigma)\in\bM_{\Sigma}(c)^l$. The tangent space of
$(A,\Phi,\sigma)$ at $\BBB^l$ is
$$T_{(A,\Phi,\sigma)}\BBB^l=\Omega^1(\imag\RR)
\oplus \Omega^0(s^*T\FFF_v)\oplus (\Hom^{0,1}(TX,s^*T\FFF_v)^{S^1}_{C^l}
\oplus \Omega^0(\imag\RR)_{C^l}),$$
where $\Hom^{0,1}(TX,s^*T\FFF_v)^{S^1}_{C^l}$ is the set of $C^l$
sections $s_1\in \Hom^{0,1}(TX,s^*T\FFF_v)$ which satisfy
$s_1(x)\in\Hom^{0,1}(T_xX,T_{\Phi(x)}\FFF_v)^{S^1}$
for any $x\in X$ such that $\Phi_0(x)\in\FFF_0$ (this is so due
to the gauge invariance of the elements of $\Sigma_c(E)^l$).

The differential of the map $\SS$ at the point $(A,\Phi,\sigma)$ is 

$$\begin{array}{rcl}D\SS:T_{(A,\Phi,\sigma)}\BBB^l & @>>> &
\Omega^{0,1}(s^*T\FFF_v)_{L^p}\oplus\Omega^0(\imag\RR)_{L^p}\\
(\alpha,\phi,(s_1,s_2)) & \mapsto &
\left(
\begin{array}{c}
\overline{\partial}_{A,\nabla}\phi+\alpha\cdot\phi+
C(A,\Phi,\sigma)(\alpha,\phi)+s_1\\
\Lambda d_A\alpha+\la d\mu(\Phi),\phi\ra_{T\FFF}+s_2
\end{array}\right).\end{array}$$
(See (\ref{difeq}).)
The image of $D\SS$ is closed because modulo the infinitesimal action
of the Lie algebra of the gauge group $\GGG$ it is an elliptic
operator. So if $D\SS$ were not exhaustive then there would exist
a nonzero element $(\eta_1,\eta_2)\in
\Omega^{0,1}(s^*T\FFF_v)_{L^q}\oplus\Omega^0(\imag\RR)_{L^q}$
($1/p+1/q=1$) such that for any $(\alpha,\phi,(s_1,s_2))$ 
\begin{equation}
\int_X\la\eta_1,\overline{\partial}_{A,\nabla}\phi+\alpha\cdot\phi+
C(A,\Phi,\sigma)(\alpha,\phi)+s_1\ra=0
\label{primeraint}
\end{equation}
and 
\begin{equation}
\int_X\la\eta_2,\Lambda d_A\alpha+\la d\mu(\Phi),\phi\ra_{T\FFF}+s_2\ra=0.
\label{segonaint}
\end{equation}
Let $X_0=\{x\in X\mid \Phi(x)\in\FFF_0\}$. Any section 
$s_1\in \Hom^{0,1}(TX,s^*T\FFF_v)$ whose support is in $X\setminus X_0$
lies inside $\Hom^{0,1}(TX,s^*T\FFF_v)^{S^1}$. This means that 
$\eta_1$ has to be zero in $X\setminus X_0$ because otherwise one
could take $\alpha=\phi=0$ and $s_1$ a suitable bump function which
would make the integral (\ref{primeraint}) nonzero.
Now, by lemma \ref{esfinit} the set $X_0$ is finite. So $\eta_1$
has to vanish identically.
Similarly $\eta_2$ has to be zero
because otherwise one could make $\alpha=\phi=0$ and $s_2$
a bump function making the integral in (\ref{segonaint}) nonzero.

This proves that $\bM_{\Sigma}(c)^l$ is a smooth Banach manifold.
Now, by corollary \ref{eslliure} the action of $\GGG$ on 
$\bM_{\Sigma}(c)^l$ is free. Uhlenbeck's gauge fixing theorem
(see theorem \ref{Uhlenbeck}) implies that there are local slices
for this action. Hence the quotient 
$\MMM_{\Sigma}(c)^l=\bM_{\Sigma}(c)^l/\GGG$
is a smooth Banach manifold. 

We now consider the projection 
$\pi_{\Sigma}:\MMM_{\Sigma}(c)^l\to\Sigma_c(E)^l$. The map
$\pi_{\Sigma}$ is Fredholm and its index depends on the homology class
$B$ (see the computations in section \ref{locals}). 
If $l\geq\Ind(\pi_{\Sigma})+2$
then the Sard-Smale theorem \cite{Sm} tells us that the set of regular
values $\Sigma_c^{\reg}(E)^l\subset\Sigma_c(E)^l$ is of the
second category of Baire. For any $\sigma\in\Sigma_c^{\reg}(E)^l$
and any $(A,\Phi)\in\MMM_{\sigma}(c)$
the second cohomology group $H^2_{A,\Phi}$ vanishes, and hence
we can apply lemma (\ref{Kuranishi}) and obtain that the dimension of
$\MMM_{\sigma}(c)$ is $$\la c^K_1(TF),B\ra+(n-1)(1-g).$$

To finish the argument we deduce from the preceeding reasoning
that $\Sigma_c^{\reg}(E)\subset\Sigma_c(E)$ is of the second
category with respect to the $C^{\infty}$ topology.
The idea is due to Taubes (see p. 36 in \cite{McDS1}) and goes
as follows. One considers for any $K>0$ 
the set $\Sigma_c^{\reg,K}(E)\subset\Sigma_c(E)$ 
of perturbations $\sigma$ such that for any $(A,\Phi)\in\MMM_{\sigma}(c)$
which satisfies $|d_A\Phi|_{C^0}\leq K$ the cohomology group
$H^2_{A,\Phi}$ vanishes (that is, $\MMM_{\sigma}(c)$ is smooth
at $(A,\Phi)$). We obviously have
$$\Sigma_c^{\reg}(E)=\bigcap_{K>0}\Sigma_c^{\reg,K}(E).$$
The set $\Sigma_c^{\reg,K}(E)$ is open for any $K$. Indeed, its
complementary is closed, since for any sequence 
$\sigma_n\in\Sigma_c(E)$ which converges to $\sigma$ and
$(A_n,\Phi_n)\in\MMM_{\sigma_n}(c)$ which satisfy 
$|d_{A_n}\Phi_n|_{C^0}\leq K$ for any $n$, one can take a subsequence
converging to $(A,\Phi)$ (using nonlinear elliptic bootstrapping,
see p. 192 in \cite{McDS1}) and the property of being exhaustive
is open. Then one uses the identity
$$\Sigma_c^{\reg,K}(E)=\Sigma_c(E)\cap \Sigma_c^{\reg,K}(E)^l$$
and the preceeding arguments to deduce that $\Sigma_c^{\reg,K}(E)$
is dense. This concludes the reasoning (see pp. 36-37 in \cite{McDS1}
for more details on this last step).

Now, if $\sigma\in\Sigma^{\reg}_c(E)$, for any pair $(A,\Phi)\in
\bM_{\sigma}(c)$ the infinitesimal stabiliser $H^0_{A,\Phi}=0$
by corollary \ref{eslliure} and the obstruction $H^2_{A,\Phi}=0$
as well. Hence we have a smooth Kuranish model of a neighbourhood
of $[(A,\Phi)]$ in the moduli space $\MMM_{\sigma}(c)$. So we may
apply corollary \ref{dimensio} to compute the dimension of the
moduli (recall that in our case the formula given in corollary
\ref{dimensio} reduces to formula (\ref{vir2})).
Finally, corollary \ref{orientacio} implies that the moduli space
$\MMM_{\sigma}(c)$ has a canonical orientation.
\end{pf}

A similar argument proves that the cobordism class of
$\MMM_{\sigma}(c)$ is independent of the particular choice of perturbation
$\sigma\in\Sigma_c^{\reg}(E)$ and the invariant complex structure
on $F$, and that it only depends on the connected
component of $\imag\RR\setminus C_0$ in which $c$ lies. 
More precisely,

\begin{theorem}
Let $c_0,c_1$ belong to the same connected component of
$\imag\RR\setminus C_0$. For any pair of perturbations
$\sigma_i\in\Sigma^{\reg}_{c_i}(E)$, $i=0,1$, one can find
paths $[0,1]\ni t\mapsto c_t\in\imag\RR\setminus C_0$ and
$[0,1]\ni t\mapsto \sigma_t\in\Sigma_{c_t}(E)$
such that
$$\MMM_{\sigma_{[0,1]}}(c_{[0,1]})
=\bigcup_{t\in [0,1]}\MMM_{\sigma_t}(c_t)$$
is a smooth oriented cobordism between $\MMM_{\sigma_0}(c_0)$
and $\MMM_{\sigma_1}(c_1)$. Likewise, any two moduli
$\MMM_{I_{F,0},\sigma}(c)$ and $\MMM_{I_{F,1},\sigma}(c)$ arising
from different complex structures $I_{F,0},I_{F,1}$ on $F$ are oriented
cobordant.
\end{theorem}                                    

\section{Smoothness of $\NNN_{\sigma}^{F,S^1}(B,c)$
for semi-free $S^1$ actions}
Fix a homology class $B\in H_2(F_K;\ZZ)$ and write
$\NNN_{\sigma}(c)$ for $\NNN_{\sigma}(B,c)$.
Using exactly the same methods as in the preceeding section
one can prove the following.
                            
\begin{theorem}
Let $c\in\imag\RR\setminus C_0$. For any $\sigma\in\Sigma_c^{\reg}(E)$
the extended moduli space $\NNN_{\sigma}(c)$ is a smooth oriented
manifold of real dimension equal to 
$$2\la c_1^{S^1}(TF),B\ra+2(n-1)(1-g)+1$$
and the natural map $\NNN_{\sigma}(c)\to\MMM_{\sigma}(c)$ is a principal
$S^1$ bundle.
Furthermore, if $c_0,c_1$ belong to the same connected component
of $\imag\RR\setminus C_0$ and $\sigma_i\in\Sigma_{c_i}^{\reg}(E)$,
$i=0,1$, then there exist paths
$[0,1]\ni t\mapsto (c_t,\sigma_t)$ such that 
$\sigma_t\in\Sigma_{c_t}(E)$ for any $t$ and
$$\NNN_{\sigma_{[0,1]}}(c_{[0,1]})
=\bigcup_{t\in [0,1]}\NNN_{\sigma_t}(c_t)$$
is a smooth oriented cobordism between $\NNN_{\sigma_0}(c_0)$
and $\NNN_{\sigma_1}(c_1)$. Likewise, any two extended moduli
$\NNN_{I_{F,0},\sigma}(c)$ and $\NNN_{I_{F,1},\sigma}(c)$ arising from
two different complex structures $I_{F,0},I_{F,1}$ on $F$ are oriented
cobordant.
\label{gensigmasmoothN}
\end{theorem}


\chapter{Compactification of the moduli}
\label{compact}

In the previous chapter we defined the moduli space 
$\MMM_{I,\sigma}(B,c)$ of THCs as 
a subset of $\AAA_{L^p_1}\times\SSS_{L^p_1}/\GGG_{L^p_2}$. So {\it a priori}
the elements $[(A,\Phi)]\in\MMM_{I,\sigma}(B,c)$ are not necessarily
smooth. In particular, it is not clear to what extent our moduli 
space depends on the Sobolev
norms we have chosed to complete $\AAA\times\SSS$. In this chapter
we will clarify the situation by proving that any pair
$(A,\Phi)\in \AAA_{L^p_1}\times\SSS_{L^p_1}/\GGG_{L^p_2}$
which satisfies equations (\ref{equs})
is smooth. Moreover, the same thing can be proved if we take 
Sobolev norms different from the ones we have chosen, provided 
they are in a certain {\it reasonable} rank. The conclusion is that 
at least as a set, the moduli space does not depend on the Sobolev 
norm. This automatically implies that the structure of the
moduli space as a differentiable manifold is intrinsic as well.
Indeed, the deformation complex at any point of the moduli is
elliptic. Hence, its cohomology groups (which model neighbourhoods
of the moduli) are the same for any Sobolev norm we may 
use to complete the modules appearing in the complex.

The other thing we will do in this chapter is to give a
compactification of the moduli space of twisted holomorphic curves.
This compactification is inspired in Gromov's compactification theorem
for holomorphic curves and makes use of Uhlenbeck's gauge fixing
theorem. Note that the construction works for any compact Lie group $K$
and without requiering any condition on the action on $F$.
                                                              
\section{Regularity of THCs}

\subsection{Preliminary results}
We will use the following result on regularity of holomorphic
curves (see \cite{McDS1} pp. 179--194).

\begin{lemma}
Let $p>2$ and $k\geq 1$. 
Let $\Sigma\subset\CC$ be an open set with smooth boundary.
Let $M$ be a compact manifold with a complex structure
$I\in\Omega^0(\End TM)_{L^p_k}$ and let 
$\phi\in\Map(\Sigma,M)_{L^p_{1,\loc}}$ satisfy $\ov{\partial}\phi=\eta$,
where $\eta\in\Omega^{0,1}(\phi^*TM)_{L^p_k}$.
Then $$\phi\in\Map(\Sigma,M)_{L^p_{{k+1},\loc}}.$$
Moreover, for any compact subset $Q\subset\Sigma$ there is
a bound
$$\|\phi\|_{L^p_{k+1}}\leq c\left(\|\phi\|_{L^p_k}
+\|\eta\|_{L^p_k}\right)$$
where $c$ depends only on $p$, $Q$, $\Sigma$ and $\|I\|_{L^p_k}$.
\label{fitacorb}
\end{lemma}

We will also make use of Uhlenbeck's theorem, which says that
any connection on the trivial bundle on the unit disk $\DD$ is
gauge equivalent to a connection in Coulomb gauge \cite{Uh}
(this is property {\it i)} in theorem \ref{Uhlenbeck}).
(Note that Uhlenbeck proves in \cite{Uh} a theorem valid in
any dimension.) As before, we write $\kE=E\times_{\Ad}\klie$.

\begin{theorem}
Let $p\geq 1$ and
consider the trivial principal $K$-bundle $E=K\times\DD\to\DD$
on the unit disk $\DD\subset\CC$. There exist constants
$\delta>0$ and $\epsilon>0$ such that any connection $A'=d+\alpha'$ on $P$, 
where $\alpha'\in\Omega^1(\kE)_{L^p_1}$, whose curvature $F_{A'}$
satisfies $\|F_{A'}\|_{L^1}\leq\delta$, is gauge equivalent by
an element $s\in \Omega^0(\DD,K)_{L^p_2}$ to a connection 
$A=d+\alpha$ where $\alpha$ satisfies
\begin{align*}
i)&\ d^*\alpha=0,\\
ii)&\ x\cdot \alpha=0\mbox{ for any }x\in\partial\DD,\\
iii)&\ \|\alpha\|_{L^p_1}\leq\epsilon\|F\|_{L^p_0}.
\end{align*}
\label{Uhlenbeck}
\end{theorem}

\begin{lemma}
Let $A=d+\alpha$ be a connection on the trivial principal bundle
$E=K\times\DD\to\DD$ which satisfies conditions
{\it i)-iii)} in theorem \ref{Uhlenbeck}.
Suppose that the curvature $F_A=d\alpha+\frac{1}{2}[\alpha,\alpha]$
lies in $\Omega^2(\kE)_{L^p_k}$, where $k\geq 1$.
Then $\alpha\in\Omega^1(\kE)_{L^p_{k+1}}$ and 
there is a bound
$$\|\alpha\|_{L^p_{k+1}}\leq C(\|F_A\|_{L^p_k}+\|\alpha\|_{L^p_1}),$$
where $C$ depends on $k$.
\label{fitaconn}
\end{lemma}
\begin{pf}
For $k\geq 1$ we have in real dimension $2$ Sobolev multiplication
$L^p_k\otimes L^p_k\to L^p_k$. On the other hand, the operator $d+d^*$ 
with the boundary condition {\it ii)} in theorem \ref{Uhlenbeck} is elliptic 
of degree $1$. So, using the generalisation of G{\aa}rding's inequality to 
$L^p$ spaces (see \cite{GiTr}) we get for any $k$
$$\|\alpha\|_{L^p_{k+1}}\leq C_k(\|(d+d^*)\alpha\|_{L^p_k}
+\|\alpha\|_{L^p_k})=C_k(\|d\alpha\|_{L^p_k}+\|\alpha\|_{L^p_k}.$$
Combining this inequality with Sobolev multiplication and using induction
the desired bound follows.
\end{pf}

In our case we will have natural upper bounds for the $L^2$ norm of the 
curvature of connections $A$ appearing in pairs 
$[(A,\Phi)]\in\MMM_{I,\sigma}(B,c)$ (this bounds are obvious from the 
equations when $F$ is compact, since then $\mu$ is bounded; for noncompact 
$F$ the Yang--Mills--Higgs functional provides the upper bound, see
theorem \ref{boundL2}). Since over 
a compact domain such as $X$ we have a continuous embedding 
$L^2\hookrightarrow L^1$, we will be able to use the preceeding theorem. 

\subsection{Rescaling}
Concerning the $L^2$ norm of the curvature, there is a crucial point
which will be used several times in our discussion. 
In real dimension 2 the $L^2$ norm of the curvature of connections is not 
conformally invariant (this is in contrast with the situation of real 
dimension 4). The situation is even better for our purposes, 
since the curvature {\it transforms in the good way} under 
conformal maps, as the following lemma says.

\begin{lemma} 
Let $0<r<1$ be a real number. Let $A=d+\alpha$ be a connection on 
the trivial principal bundle $K\times\DD_r\to\DD_r$ over the disk 
$\DD_r\subset\CC$ of radius $r$. Consider the homotecy
$\lambda_r:\DD\to\DD_r$ which sends any $z\in\DD$ to $rz\in\DD_r$.
Take on $\DD$ and $\DD_r$ the metric induced by the canonical one on $\CC$.
Then
$$\Lambda F_{\lambda_r^*A}=r^2\Lambda F_A\mbox{ and }
\|F_{\lambda_r^*A}\|^2_{L^2(\DD)}=r^2\|F_A\|^2_{L^2(\DD_r)}.$$
\label{noncomf}
\end{lemma}
\begin{pf}
We have $F_{\lambda_r^*A}=\lambda_r^*F_A$.
So for any $x\in\DD$ and any $v_1,v_2\in T_x\DD$ we have
$$F_{\lambda_r^*A}(x)(v_1,v_2)=F_A(rx)(rv_1,rv_2)=r^2F_A(rx)(v_1,v_2),$$
and this proves the first formula. To prove the second one observe that 
$|F_{\lambda_r^*A}(x)|^2=r^4|F_A(rx)|^2$ and use the change of variables 
formula.
\end{pf}

We will also use the following lemma, whose claim is analogous
to the fact that the energy is conformally invariant in dimension 2.
Its proof is an easy exercise.
        
\begin{lemma}
Using the same notation as in the above lemma, let
$\Phi:\DD_r\to F$ be a $L^p_1$ section. Then
$$\|d_{\lambda_r^*A}(\lambda_r^*\Phi)\|^2_{L^2(\DD)}=
\|d_A \Phi\|^2_{L^2(\DD_r)}.$$
\end{lemma}

These two lemmae imply the following. Let $(A,\Phi)$ be a pair
consisting of a connection one form $A$ on the trivial bundle
$K\times\DD_r\to\DD_r$ and a section $\Phi:\DD_r\to F\times\DD_r$.
Suppose that $|\mu|<C$. Then
$$\YMH_c(\lambda_r^*(A,\Phi))\leq \YMH_c(A,\Phi)+\pi (C+|c|)^2,$$
the $\pi(C+|c|)^2$ summand accounting for the term $\|\mu-c\|^2_{L^2}$
in Yang--Mills--Higgs functional.

Another consequence is that if $(A,\Phi)$ is a THC on $\DD_r$,
then $\lambda_r^*(A,\Phi)$ is no longer a THC. This justifies
the following definition. For any $\lambda\geq 0$ we will say that
$(A,\Phi)$ is a $\lambda$-THC if the following two equations 
are satisfied
$$\left\{\begin{array}{l} \overline{\partial}_A\Phi = 0  \\
\Lambda F_A+\lambda \mu(\Phi) = \lambda c.
\end{array}\right.$$
Now, if $(A,\Phi)$ is a THC on $\DD_r$, then, by lemma \ref{noncomf},
$\lambda_r^*(A,\Phi)$ is a $r^2$-THC on $\DD$.

For any $x\in X$ and any $r>0$ lower than
the injectivity radius of $X$, $\psi_{x,r}:\DD\to D(x;r)$ will denote
in the sequel the diffeomorphism 
$$\psi_{x,r}(z)=\exp_x(rz)$$
between the unit disk $\DD$
(viewed as a subset of $T_xX$) and the geodesic disk of radius $r$
centered at $x$. The map $\psi_{x,r}$ will be called the
geodesic chart of $D(x;r)$.

\subsection{Regularity}

We are now ready to prove the main result of this section.

\begin{theorem}
Suppose that the pair $(A,\Phi)\in\AAA_{L^p_1}\times\SSS_{L^p_1}$
satisfies the following equations
\begin{equation}
\left\{
\begin{array}{l}
\overline{\partial}_A\Phi = \sigma_1  \\
\Lambda F_A+\mu(\Phi) = c+\sigma_2, 
\end{array}\right.
\label{equsb}
\end{equation}
where $(\sigma_1,\sigma_2)\in\Sigma$ (so that $\sigma_1,\sigma_2$
are smooth). Then there is a gauge transformation $g\in\GGG$
such that both $g^*A$ and $g\Phi$ smooth.
\label{regular}
\end{theorem}
\begin{pf}
Let $x\in X$ be any point. Consider a geodesic
disk $D(x;r)$ centered at $x$ of radius $r>0$, and
write $\lambda$ for the geodesic chart $\psi_{x,r}$.
Denote $(A',\Phi')$ the pullback $(\lambda^*A,\lambda^*\Phi)$.
Fix a smooth trivialisation $\lambda^*E\simeq K\times\DD\to\DD$.
If $r$ has been chosen small enough the metric $\lambda^*g(X)$
is very near to the flat metric on $\DD$.
Taking this into account, lemma \ref{noncomf} implies that for small enough
$r$ the $L^2$ norm of the curvature of $A'=\lambda^*A$, computed
with the flat metric on $\DD$, is smaller than the constant
$\delta$ in theorem \ref{Uhlenbeck}:
$$\|F_{A'}\|_{L^2(\DD)}<\delta.$$
So we can take a gauge transformation $s\in\Omega^0(\DD;K)_{L^p_2}$
such that $s^*A'$ satisfies {\it i)-iii)} in \ref{Uhlenbeck}.

Let us write $(A'',\Phi'')=s(A',\Phi')$ for simplicity.
We prove the smoothness of $(A'',\Phi'')$ by a two-step bootstrapping,
using first lemma \ref{fitacorb} and then lemma \ref{fitaconn}.
Since the connection $A''$ lies in $\AAA_{L^p_1(\DD)}$, lemma
\ref{regAregI} implies that $I(A'')$ is of class $L^p_1$.
Now lemma \ref{fitacorb} and the first equation in
(\ref{equsb}) gives $\Phi''\in\SSS_{L^p_2(\DD)}$.
We now go to the second equation in (\ref{equsb}).
Since $\Phi''\in\SSS_{L^p_2(\DD)}$ we get that $F_{A''}$ has finite
$L^p_2$ norm, and so by lemma \ref{fitaconn} the connection
$A''\in\AAA_{L^p_3(\DD)}$.
Now, by lemma \ref{regAregI}, $I(A'')$ is of class $L^p_3$.
We use lemma \ref{fitacorb} again, and so on. At the end we deduce
that $(A'',\Phi'')$ is smooth. Note that, since lemma \ref{fitacorb}
gives only interior regularity, we may have to stretch our disk a little
bit each step. In any case, we end up with smoothness at any disk smaller 
that $\DD$, say $\DD_{1/2}$.

Observe also that if $s'$ is any smooth gauge transformation
then $((s's)^*A',(s's)^*\Phi')$ is still smooth. In particular,
we can take $s'$ near $s^{-1}$ in the $C^0$ norm (since
$s$ is in $L^p_2$ and $L^p_2\hookrightarrow C^0$). This
means that the gauge transformation sending $(A',\Phi')$
to a smooth pair may be taken arbitrarily small in the
$C^0$ norm.

So far we have proved that for any $x\in X$ there is 
a ball $D(x)$ centered at $x$ and a gauge transformation 
$g_x:D(x)\to K$ defined on $D(x)$ which sends $(A,\Phi)$ to a 
smooth pair (over $D(x)$). At the intersection of two overlapping balls
$D(x)$ and $D(x')$ the transformations $g_x$ and $g_{x'}$ will
differ by a smooth map $g:D(x)\cap D(x')\to K$
(recall that we used smooth trivialisations of the restriction
of $E$ to the balls $D(x)$). A standard gluing argument proves
then that there exists the desired global gauge transformation.

We show for example how to glue $g_x$ and $g_{x'}$ to a gauge
$g_{x,x'}$. Take a smooth
function $\psi$ defined in a neighbourhood of $D(x)\cup D(x')$
such that $\psi|_{D(x)\setminus D(x')}=1$ and
$\psi|_{D(x')\setminus D(x)}=0$ (we may have to restrict to
the complementary of a small tubular neighbourhood of the boundary of
$D(x)\cup D(x')$ for this $\psi$ to exist). By a previous comment
we may assume, without lose of generality, that $|g_x|$ and $|g_{x'}|$
are small enough so that at $D(x)\cap D(x')$ we can write
$g_xg_{x'}^{-1}=\exp(h)$, for $h$ smooth. We then define
$$g_{x,x'}=\left\{\begin{array}{ll}
g_x & \mbox{ at } D(x)\setminus D(x') \\
\exp(\psi h)g_{x'} & \mbox{ at }D(x)\cap D(x') \\
g_{x'} & \mbox{ at } D(x')\setminus D(x)\end{array}\right.$$
\end{pf}

\section{Equivariant Gromov-Schwartz lemma}

Take a connection $\nabla_0^F$ on $TF$ which is compatible with
the complex structure $I_F$ and whose $(1,1)$-torsion vanishes
(that such a connection exists is proved, for example, in \cite{Wo}).
Let $\nabla^F$ be the $K$-invariant connection obtained from
averaging $\nabla_0^F$. Since $I_F$ is $K$-invariant, it turns out
that the connection $\nabla^F$ still has vanishing $(1,1)$-torsion.

\begin{lemma}
Let $c\in\klie$ be a central element. There exist a constant
$\epsilon>0$ with the following property.
Let $X=\DD_r\subset\CC$ be the disk of radius $r>0$ with the standard
metric. Consider a pair $(A,\Phi)$, where $A\in\Omega^1(X;\klie)$ is a
connection one form on the trivial bundle $K\times X\to X$ and
$\Phi:X\to F\times X$ is a section of the associated bundle.
Suppose that the pair $(A,\Phi)$ is a $\lambda$-THC, that is, it
satisfies the equations           
\begin{equation}
\left\{\begin{array}{l}
\overline{\partial}_A\Phi = 0  \\
\Lambda F_A+\lambda \mu(\Phi) = \lambda c,
\end{array}\right.
\label{eqlthc}
\end{equation}
where $\lambda\in[0,1]$ is a real number. Suppose that 
$\YMH_c(A,\Phi)<\epsilon$. Then 
$$|d_A\Phi(0)|^2<C(1+r^{-2})(1+\|d_A\Phi\|^2_{L^2(X)}),$$
where $C$ is a constant independent of $A$, $\Phi$ and $\lambda$.
\label{schwartz}
\end{lemma}
\begin{pf}
Throughout the proof $C$ will denote a constant independent of $A$, $\Phi$
and $\lambda$, and will not always be necessarily the same.
Let $r'=\min(r,\sqrt{\epsilon/2\pi})$. We will prove that
\begin{equation}
|d_A\Phi(0)|^2<C{r'}^{-2}((1+\|d_A\Phi\|^2_{L^2(D_{r'})}).
\label{ufff}
\end{equation}
This clearly implies the lemma.

Take on $TX$ the trivial connection $\nx$. This connection is obviously
compatible with the canonical complex structure $I_X$ on $X\subset\CC$
and has zero torsion. Let $g_X$ be the canonical metric on $X$.
Let $\pfx:\FFF=X\times F\to X$. The connection $A$ 
gives a splitting $T\FFF\simeq \pfx^*TX\oplus T\FFF_v$, which allows to
define a complex structure $I(A)=I_X\oplus I_F$ 
and a metric $g(A)=g_X\oplus g_F$ 
on $T\FFF$. On the other hand, since $\nabla^F$ is $K$-invariant,
we can extend it fibrewise to obtain a connection on $T\FFF_v$
(see section \ref{vbfb}).
Summing the resulting connection with $\pfx^*\nx$ (here we use the
splitting $T\FFF\simeq \pfx^*TX\oplus T\FFF_v$ given by $A$) we get a
connection $\nf=\nf(A,\nabla^F,\nx)$ on $T\FFF$.
The connection $\nf$ is compatible with $I(A)$ and $g(A)$, that is,
$\nf I(A)=\nf g(A)=0$. 

\begin{lemma}
Take a metric on $\klie$. With respect to the metrics $g_X$, $g_F$ and $g(A)$
we have $|F_{\nf}| \leq |F_{\nabla^F}|+C|F_A|$, where $C>0$ is a constant.
\label{fitacurv}
\end{lemma}
\begin{pf}
Since $F$ is compact there is a real number $C>0$ such that, for any 
$x\in F$ and for any $s\in\klie$, $|\fX_s^F(x)|\leq C|s|$ and 
$|\Omega^{TF}(s)(x)|\leq C|s|$. Now lemma \ref{curvatura} in the appendix
gives the result.
\end{pf}
                                                                         
\begin{lemma}
Suppose that the $(1,1)$ parts of the torsions of $\nabla^F$ and $\nx$ 
vanish. Then the $(1,1)$ part of the torsion of the connection $\nv$ is
$-d\pfx^*F_{A^{\FFF}}^{1,1}.$
\label{torsiouu}
\end{lemma}
\begin{pf} This follows from lemma \ref{torsio} in the appendix.\end{pf}

Let $(u,v)$ be the standard coordinates in $\CC\simeq\RR^2$.
Let $\bU=\partial\Phi/\partial u$ and $\bV=\partial\Phi/\partial v$.
In the rest of the proof we will write $\nabla=\nf$ and $J=I(A)$.

As defined, the fields $\bU$ and $\bV$ are sections 
of $\Phi^*T\FFF$. To give a rigorous sense to the following computations,
we extend smoothly $\bU$ and $\bV$ to get fields on $\FFF$; however,
we are only interested on the values they take on
$\Phi(X)$, and all the equalities between fields on $\FFF$ should
be understood in this proof as refering only to their restrictions to
$\Phi(X)$. For example, we may write $J\bU=\bV$ and $J\bV=-\bU$
(no matter how $\bU$ and $\bV$ have been defined outside 
$\Phi(X)$). We also have $|\bU|=|\bV|$.

We have $[\bU,\bV]=0$ (this is always true on $\Phi(X)$, regardless
the way we have extended $\bU$ and $\bV$ out of $\Phi(X)$).
By lemma \ref{torsiouu} we have
$\nuv-\nvu=\Tor_{\nf}(\bU,\bV)=\fX^F_{\Lambda F_A}$.
Using the second equation in (\ref{eqlthc})
we deduce from this the following inequalities
\begin{align}
|\nau(\nuv-\nvu)|&<C|\bU|,\label{ineq1} \\
|\nav(\nuv-\nvu)|&<C|\bU|,\label{ineq2}
\end{align}
where the constant
$C$ depends only on the (norm of the) differential $d\mu$ of the moment map
$\mu$ and on the constant $C$ in the proof of lemma \ref{fitacurv}.

Using $\nabla J=0$ we compute
\begin{align*}
\frac{1}{2}\Delta|\bU|^2 &= 
|\nuu|^2+|\nvu|^2+\la\bU,\nau\nuu\ra+\la\bU,\nav\nvu\ra \\
&\geq \la\bU,\nau\nuu\ra+\la\bU,\nav\nvu\ra \\
&= -\la\bU,J\nau\nuv\ra+\la\bU,\nav\nuv\ra-\la\bU,\nav(\nuv-\nvu)\ra \\
&= -\la\bU,J\nau\nvu\ra
+\la\bU,J\nav\nuu\ra \\
&+\la\bU,J\nau(\nvu-\nuv)\ra
-\la\bU,\nav(\nuv-\nvu)\ra.
\end{align*}
The first and the second summands in the last expression 
are equal to $-\la\bU,J F_{\nabla}(\bU,\bV)\bU\ra$, which can be
bounded by $C|\bU|^4$ thanks to lemma \ref{fitacurv}
(recall that we write $\nabla=\nf$). Finally, the
third and the fourth summands can be bounded by
$C|\bU|^2$ using the inequalities
(\ref{ineq1}) and (\ref{ineq2}).
So we obtain
\begin{equation}
\frac{1}{2}\Delta|\bU|^2\geq-C|\bU|^2-C|\bU|^4.
\label{Laplace}
\end{equation}
On the other hand, we have $|d\Phi|^2=2|\bU|^2$ and 
$|d\Phi|^2=2+|d_A\Phi|^2$, so from the bounds
$\|d_A\Phi\|_{L^2(D_r)}^2\leq\YMH_c(A,\Phi)<\epsilon$ and
$r'\leq\sqrt{\epsilon/2\pi}$ we deduce $\|d\Phi\|^2_{L^2(D_{r'})}<2\epsilon$.
Now, as in theorem 2.3 in \cite{PaWo}, we deduce from the inequality
(\ref{Laplace}) that if $\epsilon$ is small
enough (with respect to the constants appearing in (\ref{Laplace}))
then $$|d\Phi|^2\leq C{r'}^{-2}\|d\Phi\|_{L^2(D_{r'})}^2,$$
from which (\ref{ufff}) easily follows.
\end{pf}

\section{Removability of singularities}

Let $r\in\RR^+\setminus\NN$, and write $r=l+\alpha$, where
$k\in\NN$ and $\alpha\in (0,1)$. Let $C^r$ denote the H{\"o}lder 
$C^{l,\alpha}$ norm. Recall that for any compact $n$-dimensional
manifold $M$ we have a continuous embedding
$L^q_k(M)\hookrightarrow C^r(M)$
whenever $r\leq k-\frac{n}{q}$ (here $C^r(M)$ denotes the completion
of $C^{\infty}(M)$ with the $C^r$ norm). Since we are assuming that
$p>2$, it turns out that there exists $r>0$ and  
\begin{equation}
L^p_1(M)\hookrightarrow C^r(M)
\label{Sobolev}
\end{equation}
for any Riemann surface $M$.

We will use the following theorem on removal of singularities
for holomorphic curves (see p. 180 in \cite{Sik}).

\begin{theorem}
Let $M$ be a compact manifold with a complex structure $I$
of H{\"o}lder class $C^r$, where $r>0$. Let $f_0:\DD\setminus\{0\}\to M$
be a $I$-holomorphic map of finite area.
Then $f_0$ can be extended to a $I$-holomorphic map
$f:D\to M$ of class $C^{r+1}$.
\label{remcorb}
\end{theorem}

The following theorem generalises the preceeding result to THCs.

\begin{theorem}
Let $c\in\klie$ be a central element. Let $X=\DD\subset\CC$
be the unit disk centered at $0$, and let $X^*=X\setminus
\{0\}$. Let $\alpha_0\in\Omega^1(X;\klie)_{L^p_1}$ be a connection 
one form and take a section $\Phi_0:X^*\to
F\times X^*$ of the bundle $F\times X^*\to X^*$. Let $A_0=d+\alpha_0$
be the connection induced by $\alpha_0$. Suppose that on $X^*$
\begin{equation}
\left\{\begin{array}{l}
\overline{\partial}_{A_0}\Phi_0 = 0 \\
\Lambda F_{A_0}+\lambda \mu(\Phi_0) = \lambda c
\end{array}\right.
\label{equss}
\end{equation}
where $\lambda\in[0,1]$, 
and that $\YMH_c(A_0,\Phi_0)<\infty$. Then there exists
a gauge transformation $g:X\to K$ such that $g(A_0,\Phi_0)$ 
can be extended to a $\lambda$-THC $(A,\Phi)$ defined on $X$. 
\label{removal}
\end{theorem}

\begin{pf}
The proof is similar to that of theorem \ref{regular}: we use
already existing results alternatively with the first equation and then
with the second one.
We can assume (maybe after reducing to a smaller disk $\DD_r$, $r<0$,
and rescalling), that $\|F_{A_0}\|_{L^2(X^*)}<\delta$, where
$\delta>0$ is the constant in theorem \ref{Uhlenbeck}. Then
using theorem \ref{Uhlenbeck}, we may gauge the pair $(A_0,\Phi_0)$
to a pair $(A,\Phi)$ such that $A=d+\alpha$ is in Coulomb gauge: 
$d^*\alpha=0$. 
Consider the complex structure $I=I(A)$ on $F\times X$. By the embedding 
(\ref{Sobolev}) the connection $A$ is of class $C^r$, for $r>0$. The same 
applies to the complex structure $I$, so we may use theorem \ref{remcorb}.
Indeed, by the first equation in (\ref{equss}), $\Phi$ is
holomorphic with respect to $I$ and, since $\YMH_c(A,\Phi)<\infty$, 
the energy $\|d\Phi\|_{L^2}<\infty$.
We thus obtain an extension $\Phi:X\to F\times X$ of $\Phi_0$
of class $C^{r+1}$. We now turn to the second equation in 
(\ref{equss}) and deduce that $F_A$ is of class $L^p_1$.
Since we are in Coulomb gauge this gives a $L^p_2$ bound on $A$
and hence on $I$. And so on.
\end{pf}

\section{Compactness}

Let $\sigma\in\Sigma(E)$ be any perturbation.
                                   
\begin{definition} A {\bf cusp $\sigma$-THC} is the following set of data.
\begin{enumerate}
\item A connected singular curve $X^{\cusp}$ with only nodal 
singularities, of the form $X^{\cusp}=X_0\cup X_1\cup\dots \cup X_K$,
where $X_0=X$, and where the other components are rational curves $\CP^1$
and are called {\bf bubbles}; furthermore, two different components
$X_i$ and $X_j$ meet at most at one point. We call $X_0$ the
{\bf principal component} of the cusp curve $X^{\cusp}$.
\item A $S^1$-principal bundle $E\to X_0$, a connection $A$ on $E$,
a section $\Phi_0:X_0\to \FFF=\FFF^E=E\times_{S^1}F$ and an element
$c\in\imag\RR$ satisfying the equations
$$\left\{\begin{array}{l}
\ov{\partial}_A\Phi_0=\sigma_1 \\
\Lambda F_A+\mu(\Phi_0)=c+\sigma_2. \end{array}\right.$$
\item For any $k\neq 0$, a holomorphic map $\Phi_k:X_k\to\FFF$ whose 
image is inside a unique fibre $\FFF_{x_k}$ of $\FFF\to X$ (note $\Phi_k$
is holomorphic with respect to the complex structure on $F$). The maps
$\Phi_0,\Phi_1,\dots,\Phi_K$ are required to glue together to give
a map $\Phi:X^{\cusp}\to\FFF$. 
\end{enumerate}
\label{defcuspTHC}
\end{definition}
We denote cusps with tuples of the form $(E,X^{\cusp},A,\Phi,c)$.

\begin{theorem}
Let $B\in H_2(F_K;\ZZ)$ be any homology class. Let $E\to X$ be the unique 
$K$ principal bundle such that $\eta(E)=(\pi_F)_*B$, where $\pi_F:F_K\to BK$ 
is the projection (see lemma \ref{etabij} in the appendix).
Let $\AAA$ be the space of connections on $E$, $\SSS=\Gamma(E\times_K F)$, 
and $\GGG$ the gauge group of $E$. Consider a sequence of $\sigma$-THCs 
$(A_j,\Phi_j)\in\MMM_{\sigma}(B,c_j)$, where the elements 
$c_j$ lie in a bounded subset of $\klie$.

Then, after passing to a subsequence, there exists a cusp $\sigma$-THC 
$(E,X^{\cusp},A,\Phi,c)$ and gauge transformations $g_j\in\GGG$ 
such that if $(A_j',\Phi_j')=g_j(A_j,\Phi_j)$ we have
\begin{enumerate}
\item $c_j\to c$,
\item $A_j'\to A$ in $C^{\infty}$,
\item ${\rho_E}_*({\Phi_0}_*[X_0]+\sum_{k=1}^K{\Phi_k}_*[X_k])=B$,
where the map ${\rho_E}_*$ is given in subsection \ref{themodulispace}
(see also \ref{pullingback} in the appendix),
\item the images $\Phi_j'(X)\subset\FFF$ converge pointwise to
$\Phi(X^{\cusp})$, that is, for any sequence $x_j\in X_j$
there exists $x\in X^{\cusp}$ such that $\Phi_j'(x_j)\to \Phi(x)$.
\end{enumerate}
\label{compactificacio}
\end{theorem}
\begin{pf} We prove the theorem in several steps.
\subsubsectionr{} We may assume that $\sigma_1=0$. If it is not,
then we substitute $F$ by $X\times F$ with the complex structure
$I_{\sigma_1}$ as defined in subsection \ref{modcorbes} and with the
trivial action of $S^1$ on the first factor. In the rest of the
argument we will assume that $\sigma_2=0$ just to simplify the
notation (the general case is not more difficult, but longer to write).

\subsubsectionr{}
Since $(A_j,\Phi_j)\in\MMM_{\sigma}(B,c)$, theorem \ref{boundL2}
gives bounds
\begin{equation}
\|F_{A_j}\|_{L^2}^2\leq C
\mbox{ and }
\|d_{A_j}\Phi_j\|_{L^2}^2\leq C,
\label{bounds}
\end{equation}
where $C>0$ is independent of $j$. (Note that theorem \ref{boundL2}
refers to THCs and not to $\sigma$-THCs. The needed modifications
to deal with $\sigma$-THCs are easy and left to the reader.)

\subsubsectionr{} We next show that there is a finite subset
$\{x_1,\dots,x_l\}\subset X$, a subsequence of $\{(A_i,\Phi_i)\}$ 
(which we call again $\{(A_i,\Phi_i)\}$) and gauge transformations 
$\{g_i\}$ such that $g_i(A_i,\Phi_i)$ converge uniformly in the 
$C^{\infty}$ topology on $X\setminus\{x_1,\dots,x_l\}$ to a THC 
$(A_{0},\Phi_{0})$.

Fix $\epsilon'>0$ lower than the $\epsilon$ in lemma \ref{schwartz}.
For any $k,n\in\NN$ let
\begin{equation}
C_{k,n}=\{x\in X|\ \YMH_{c_k}(\psi_{x,2^{-n}}(A_k,\Phi_k))>\epsilon'\}.
\label{defC}
\end{equation}

\begin{lemma}
$C_{k,n}$ can be covered by $N$ balls
$\{\psi_{x_l^{k,n},3.2^{-n}}(\DD)\}_{1\leq l\leq N}$,
where $N$ depends only on $\epsilon$.
\end{lemma}
\begin{pf}
Fix $k\in\NN$. Take any point $x_1^{k,n}\in C_{k,n}$ and repeat the
following process.
Once $X_m=\{x_1^{k,n},x_2^{k,n},\dots,x_m^{k,n}\}$ have been selected,
take $x_{m+1}^{k,n}$ to be any point in $C_{k,n}$ at distance
$\geq 2.2^{-n}$ from any point in $X_m$. This process must finish at
some step. Indeed, if $\{x_1^{k,n},\dots,x_N^{k,n}\}$ have been selected
in this way, then the disks $\{D(x_j^{k,n};2^{-n})\}_{1\leq j\leq N}$
are disjoint, and hence by the definition (\ref{defC}) of $C_{k,n}$
we must have $N\epsilon'\leq\YMH_{c_k}(A_k,\Phi_k)=C$.
Finally, if $\{x_1^{k,n},\dots,x_N^{k,n}\}$ is a maximal set,
then $C_{k,n}$ is contained in 
$$\{D(x_j^{k,n};3.2^{-n})\}_{1\leq j\leq N}.$$
\end{pf}

Now, by taking a subsequence of $(A_k,\Phi_k)$ we can assume that
the elements in $$\{x_1^{k,n},\dots,x_N^{k,n}\}$$ converge
as $k\to\infty$ to $\{x_1^n,\dots,x_N^n\}$. Repeating this 
for any $n\in\NN$ and using the diagonal argument, we may assume
that $\{x_1^n,\dots,x_N^n\}=\{x_1,\dots,x_N\}$ for all $n$.
Let $$X_0=X\setminus\{x_1,\dots,x_N\}.$$
Gromov-Schwartz lemma \ref{schwartz} implies the following. If
$x\in X\setminus C_{k,n}$ then
$$|d_{A_k}\Phi_k(x)|<1+2^n C,$$
where $C$ is independent of $k$ and $n$. As a consequence, for any
$x\in X_0$ we have a bound for any $k$
\begin{equation}
|d_{A_k}\Phi_k(x)|<C(d(x,\{x_1,\dots,x_N\})),
\label{bdist}
\end{equation}
where $C:\RR_+\to\RR_+$ is a continuous function.

From the bound on the curvature in (\ref{bounds})
and lemma \ref{noncomf} we deduce 
that for any $\epsilon>0$, there exists $r=r(\epsilon)$ such that for
any $x\in X$ and $0\leq r'\leq r$
$$\|F_{\psi_{x,r}^*A}\|_{L^2(\DD)}\leq\epsilon.$$
Taking this into account, we define $r_0=r(\delta)$, where $\delta$
is the constant in Uhlenbeck's theorem \ref{Uhlenbeck}.

Cover $X$ with open balls of radius $r_0$ and take a finite
cover $\{D_1,\dots,D_s\}$ such that the collection $\{D_1',\dots,D_s'\}$
of disks concentric with $\{D_1,\dots,D_s\}$ and with radius $r_0/2$
still covers $X$.
Using Uhlenbeck's theorem \ref{Uhlenbeck} we get for any $k,l$
trivialisations
$$\sigma_{k,l}:E|_{D_l}\to D_l\times K$$
for which the connections $A_k|_{D_l}=d+\alpha_{k,l}$ are in Coloumb gauge.
That is, $\alpha_{k,l}\in\Omega^1(D_l;\klie)$ has bounded $L^p_1$ norm:
\begin{equation}
\|\alpha_{k,l}\|_{L^p_1(D_l)}\leq C.
\label{boundconn}
\end{equation}
On any nonempty intersection $D_i\cap D_j$ we have transition functions
$g_{k,i,j}:D_i\cap D_j\to K$ such that
$\sigma_{k,j}=g_{k,i,j}\sigma_{k,i}$ and which satisfy the cocycle condition
$g_{k,i,l}=g_{k,i,j} g_{k,j,l}.$
On $D_i\cap D_j$ we have $\alpha_{k,j}=g_{k,i,j}^*\alpha_{k,i}$.
This is equivalent to
$dg_{k,i,j}=g_{k,i,j}^{-1}\alpha_{k,j}-\alpha_{k,i}g_{k,i,j}^{-1}.$
The transition functions $g_{k,i,j}$ take values in the
compact group $K$. So using the bound (\ref{boundconn}) and 
the Sobolev multiplication theorems we get uniform bounds
\begin{equation}
\|g_{k,i,j}\|_{L^p_2(D_i\cap D_j)}\leq C'.
\label{boundtrans}
\end{equation}

Note that the gauge equivalence class of the pair $(E,A_k)$
is completely determined by the collection of $1$-forms
$\{\alpha_{k,l}\}\in\Omega^1(D_l;\klie)$ and by the transition
functions $g_{k,i,j}$.

Using the bounds (\ref{boundconn}) and (\ref{boundtrans}) we deduce
that there exists a collection of $1$-forms
$\{\alpha_l\in\Omega^1(D_l;\klie)$ and a set of maps
$\{g_{i,j}:D_i\cap D_j\to K\}$ and that,
after possibly restricting to a subsequence of $\{A_k\}$,
we have weak convergences for any $i,j,k,l$
$$\alpha_{k,l} \rightharpoonup \alpha_l\mbox{ in }L^p_1\mbox{ and }
g_{k,i,j} \rightharpoonup g_{i,j}\mbox{ in }L^p_2.$$

This implies strong convergence $g_{k,i,j}\to g_{i,j}$ in $L^p_1$.
In particular, making $k\to\infty$ in the cocycle equation 
$g_{k,i,l}=g_{k,i,j} g_{k,j,l}$
we deduce that the set $\{g_{i,j}\}$ satisfies the cocycle condition
$$g_{i,l}=g_{i,j} g_{j,l}$$
and hence defines a $K$-principle bundle isomorphic to $E$
(since $L^p_1\hookrightarrow C^0$ and consequently the convergence
$g_{k,i,j}\to g_{i,j}$ is uniform; see proposition 3.2 and
corollary 3.3 in \cite{Uh}). On the other hand, by the 
Sobolev multiplications theorems we can make $k$ go to $\infty$ in the
equation $\alpha_{k,j}=g_{k,i,j}^*\alpha_{k,i}$ and deduce that
the $1$-forms $\alpha_i$ satisfy
$$\alpha_j=g_{i,j}^*\alpha_i,$$
and so they define a connection $A\in \AAA_{L^p_1}$ on $E$.
Let $A_{0}=A|_{X_0}$. Our construction in terms of local trivialisations 
implies the existence of a sequence of global gauge transformations 
$g_k$ such that $g_k^*A_k\rightharpoonup A$ weakly in $L^p_1$.

Let now $x\in X_0$ and suppose that $x\in D_l'$. 
Take a compact disk $K=D(x;r)$, where 
$$r=\min\{r_0/2,d(x,\{x_1,\dots,x_r\}/2)\}.$$
We get from (\ref{bdist}) a uniform bound
$$|d_{A_k}\Phi_k|_{C^0(K)}\leq C''$$
which only depends on $d(x,\{x_1,\dots,x_r\}/2)$.

Next, we prove that the restriction of $\Phi_k$ to $K$ and $\alpha_{k,l}$
are bounded in the $C^{\infty}$ topology. For that we use the bounds
in lemmae \ref{fitacorb} and \ref{fitaconn} alternatively, exactly
as in the proof of theorem \ref{regular}. So by Ascoli-Arzela
and after passing to a subsequence, there exists 
$\Phi_{0}^K:K\to F$ such that $\Phi_k|_K\to\Phi_{0}^K$, and the 
convergence $\alpha_{k,l}\to\alpha_l$ is in $C^{\infty}$. 
Hence on $K$ the pair $(d+\alpha_l,\Phi_{0}^K)$ is a THC.
Doing this for any $x\in X_0$ we get a smooth section $\Phi_{0}'$
defined on $X_0$ such that $(A_{0},\Phi_{0}')$
is a THC on $X_0$.
             
\subsubsectionr{} Using theorem \ref{removal} on removability of 
singularities we see that the pair $(A_{0},\Phi_{0}')$
extends to a pair $(A,\Phi_{0})$ defined on the whole
$X$. Indeed, from construction we know that $A_{0}$ extends
to a connection $A$ with $\|A\|_{L^p_1}\leq\infty$.
                       
\subsubsectionr{} The total energy $\YMH_{c_{0}}(A_{0},\Phi_{0})$ may be 
lower than $\YMH_{c_k}(A_k,\Phi_k)$ if bubbling occurs. Equivalently,
the cohomology class ${\rho_E}_*{\Phi_0}_*[X_0]$
need not be $B$. The last step is to
consider this possibility and to study the possible bubbling.

For that one does the following. Fix $x\in\{x_1,\dots,x_N\}$.
We consider sequences $\{t_k\}\in\RR$ and $\{z_k\}\in\CC$
(on which later certain conditions will be imposed) and the maps
$\psi_k:\CC\simeq T_xX\to X$ which send any $z\in\CC$ to
$$\psi_k(z)=\exp(t_k^{-1}(z-z_k)).$$
(Note that in general $\psi_k$ may fail to be defined in the complementary
of a disk in $\CC$.)
Suppose that $t_k\to\infty$ and that $t_k z_k^{-1}$ remains bounded.
Then for any compact $K\subset\CC$ and big enough $k$ (with respect
to $K$) the THC
$$(A_k'',\Phi_k'')=\psi_k(A_k,\Phi_k)$$
is defined on $K$. Furthermore, the energy $\YMH_{c_k}(A_k'',\Phi_k'')$
remains bounded. So we can apply the same reasoning as in the beginning
and deduce that out of a set $y_1,\dots,y_M\subset\CC$ the
THCs $(A_k'',\Phi_k'')$ converge (after conveniently regauging) to
a THC $(A_\infty'',\Phi_\infty'')$ defined on $\CC\setminus
\{y_1,\dots,y_M\}$. 

The connection $A_\infty''$ is (gauge equivalent to) the trivial one. 
This follows from the fact that, since $t_k\to\infty$, lemma \ref{noncomf} 
implies that $\|F_{A_k''}\|_{L^2}\to 0$. So we end up with a genuine
holomorphic curve 
$\Phi_{0}'':\CC\setminus\{y_1,\dots,y_M\}\to \FFF_x=F$
of finite energy. Hence we may use removal of singularities to get
a map from $\CP^1$ to $F$. And now we may repeat the process again.

Now we have to specify how the sequences $\{t_k\}$ and $\{z_k\}$
are to be chosen in all the process. The following two things are necessary:
\begin{enumerate}
\item the recursive process ends after a finite number of steps
(bounded above by a funtion of $\YMH_{c_k}(A_k,\Phi_k)$ or 
equivalently by $B$) and
\item there is no lose of energy, that is, after we finish
the process we get a set of bubbles $X_1,\dots,X_R$ such that
$$X^{\cusp}=X_0\cup X_1\cup\dots\cup X_R$$ is an admissible curve 
and a cusp $\sigma$-THC $(E,X^{\cusp},A,\Phi,c)$ whose energy 
$\YMH_c(A,\Phi_0)+\sum\|d\Phi_k\|_{L^2}^2$ is equal to
$\YMH_c(A_j,\Phi_j)$ for any $j$ and such that
${\rho_E}_*({\Phi_0}_*[X_0]+\sum{\Phi_k}_*[X_k])=B$
\end{enumerate}

In \cite{PaWo} a precise algorithm is given to chose the sequences
$\{t_k\}$ and $\{z_k\}$ in the case of holomorphic curves. 
This is done by chosing $z_k$ so that the energy $\psi_k^*\Phi$
is centered in the south pole (here we view $\CP^1=\CC\cup\{\infty\}$ and we
identify the south pole with $0\in\CC$) and such that the amount
of energy concentrated in the north hemisphere is equal to a certain
suitable constant $C_0$. We now give a sketch of their construction
adapted to our situation. See \cite{PaWo} pp. 83-85 for more details. 
\begin{enumerate}
\item Let $S^2_x\subset\RR^3$ be the unit sphere centered at $(0,0,0)$,
and let $p_+=(0,0,1)$ and $p_-=(0,0,-1)$ be the north and south pole.
Fix $\epsilon>0$ smaller than the injectivity radius of $F$.
Let $\sigma:S^2_x\to T_xX$ be the stereographic projection which
maps $\sigma(p_-)$ to $0$ and $\sigma(p_+)$ to $\infty$.
\item Let us fix a retraction $\rho:\FFF|_{D(x;\epsilon)}\to\FFF_x\simeq F$ 
(for example by radial parallel transport). Identifying a neighbourhood
of $0$ in $T_xX$ with a neighbourhood of $x$ in $X$ by means of the
exponential map we may define
$$f_j=\sigma^*(\rho\circ\Phi_j):D(\epsilon)\to F,$$
where $D(\epsilon)\subset S^2_x$. Let $T_j:S^2_x\to S^2_x$
be a conformal map such that if $f_j'=f_j\circ T_j$ then the
center of mass of the measure $||df_j'|^2-|df_0'|^2|$ lies in
the $z$-axis.
\item For a fixed suitable constant $C_0$ (see p. 83 in \cite{PaWo})
we consider a radial dilation $d_j:S_x^2\to S_x^2$ which keeps
the north and south poles fixed and such that if
$D_j=d_j^{-1}(D(\epsilon))$, $f_j''=f_j'\circ d_j$ and $H_-$ is the
south hemisphere in $S_x^2$, then
$$\int_{D_j\setminus H_-}||df_j''|^2-|df_0''|^2|=C_0.$$
\item Then $R_j=\sigma\circ T_j\circ d_j:S_x^2\to S_x^2$
keeps the two poles fixed and its restriction to $\CC=S_x^2\setminus\{p_+\}$
is $R_j(x)=t_jx+z_j$, and $t_j\to\infty$. So by our previous reasoning
this gives a bubble $B(\epsilon)$ concentrated in $\FFF_x$.
Now we make $\epsilon\to 0$ and we take as our final bubble the
limit bubble (again this exists by Gromov compactness theorem).
\end{enumerate}
Translating word by word the argument in \cite{PaWo} we deduce
that if we follow this renormalisation method then there is no
loss of energy, and consequently the third claim in the theorem
follows.

For the last claim we invoque corollary 6.4 in \cite{PaWo}, where
pointwise convergence is proved for holomorphic curves. Again,
the argument applies to our situation.

Another very good reference for the bubble tree compactification
is \cite{Pa}, where harmonic maps are studied instead of holomorphic
curves (which are particular cases of harmonic maps).
\end{pf}


\chapter{The choice of the complex structure}
\label{cxstr}

\setcounter{theorem}{0}

In this chapter we will assume $K=S^1$ and we will call 
$\omega=\omega_F$ the symplectic form on $F$.
In chapter \ref{moduli} we constructed the moduli space of
$\sigma$-THCs and proved that for generic $\sigma$ the moduli
was a smooth manifold. In chapter \ref{compact} we gave a
compactification of the moduli by adding to it what we called
cusp $\sigma$-THCs. We may view this compactification as adding
some pieces to the moduli in order to get a compact stratified
topological space. Our ultimate goal is to use the compactification
of the moduli to define invariants, and for that we shall need
the strata to be smooth. 
In this chapter we prove that for a generic invariant complex structure 
on $F$ certain moduli spaces of holomorphic curves are smooth manifolds.
This will be used in chapter \ref{invar} to prove that
the moduli of cusp $\sigma$-THCs is a stratified variety with 
strata admitting ramified coverings by smooth manifolds.

Let $\III_{\omega}\subset\End(TF)$ be the set of complex 
structures on $F$ compatible with $\omega$ (that is, for any 
$I\in\III_{\omega}$, $\omega(\cdot,I\cdot)$ is a Riemannian metric).
Our aim is to study compatible complex structures $I\in\III_{\omega}$ 
which are invariant under the action of $S^1$ (recall that invariance
was requiered in order to define the equations, see section 
\ref{theequations} in chapter \ref{eqymh}). A first existence result,
which was already mentioned in chapter \ref{eqymh} is the following (see 
lemma 5.49 in \cite{McDS2}).

\begin{lemma}
Let $\III_{\omega,S^1}\subset\III_{\omega}\subset\End(TF)$ 
be the set of $S^1$-invariant complex structures on $F$ which are compatible 
with $\omega$. The space $\III_{\omega,S^1}$ is nonempty and contractible.
\end{lemma}

We are only interested on some elements of $\III_{\omega,S^1}$, namely,
those for which certain spaces of holomorphic curves are smooth.
If we forget the action of $S^1$ and we take all compatible complex 
structures, then for a generic structure we get a smooth moduli of simple
holomorphic curves. Multicovered curves (that is, those which are not
simple) cannot be included in the moduli without possibly giving rise to
singular points, no matter what complex structure we consider.
Due to $S^1$-invariance, in our setting there will appear more 
{\it problematic} types of curves, appart from multicovered ones.
In fact, in order to get smooth moduli spaces of simple curves
we will have to restrict ourselves to curves with a fixed isotropy pair
(see subsection \ref{comentaris}).

\section{Isotropy pairs}

\begin{definition}
Let $s:\Sigma\to F$ be any smooth map. We define the {\bf isotropy pair} 
of $s$, $(L(s),H(s))$, to be the pair of closed subgroups 
$L(s)\subset H(s)\subset S^1$ defined as follows
\begin{align*}
H(s) &:= \{\theta\in S^1|\ \theta\cdot s(\Sigma)=s(\Sigma)\} \\
L(s) &:= \{\theta\in H(s)|\ \theta|_{s(\Sigma)}=\Id\}.
\end{align*}
\end{definition}

Let $s:\Sigma\to F$ be a simple map, and let $g\in H(s)$ be any element.
Let $\Sigma_i$ be the set of injective points of $s$, that is,
$$\Sigma_i=\{x\in\Sigma|\ ds(x)\neq 0,\ \sharp s^{-1}s(x)=1\}.$$
The action of $g$ on $s(\Sigma)$ induces a holomorphic bijection
$\gamma_i(g):\Sigma_i\to\Sigma_i$
which can be extended to a homeomorphism
$\gamma(g):\Sigma\to\Sigma.$
Now, since the map $s$ is simple, 
the noninjective points $\Sigma\setminus\Sigma_i$ can
only accumulate at a finite set of points (namely, the critical
points $\Ker ds$), and hence the map $\gamma(g)$ is holomorphic
by standard removability of singularities. This way we have defined
a map 
$\gamma:H(s)\to\Aut(\Sigma).$
Obviously, $\Ker \gamma=L(s)$. Let $\Gamma=\gamma(H(s))$.
Since $L$ and $H$ are closed subgroups of $S^1$, the extension
\begin{equation}
1\to L\to H\to \Gamma\to 1
\label{salvat}
\end{equation}
is uniquely determined from $L$ and $\Gamma$.

On the other hand, note that $s(\Sigma)$ is contained by definition
in $F^{L(s)}$, the fixed point set of the action of $L(s)$ on $F$.
We define
$$\Map^L(\Sigma,F)=\{s\in\Map(\Sigma,F^L)|L(s)=L\}\subset
\Map(\Sigma,F^L).$$

This inclusion is not an equality in general, since there might be some
$s\in\Map(\Sigma,F^L)$ such that $L(s)$ strictly contains $L$.                                                       
The following technical results will be crucial in proving the
smoothness of certain moduli spaces for generic choices of
invariant complex structure $I\in\III_{\omega,S^1}$.

\begin{theorem}
Let us take a 
complex structure $I\in\III_{\omega,S^1}$, and let $s:\Sigma\to F$ 
be a simple holomorphic map. Let $H=H(s)$. Then there exists
a disk $D\subset\Sigma$ such that
$$S^1\cdot s(D)\cap s(\Sigma)=H\cdot s(D).$$
\label{stabilH}
\end{theorem}

\begin{theorem}
Let us take a complex structure $I\in\III_{\omega,S^1}$,
and let $s\in\Map^L(\Sigma,F^L)$ be a simple holomorphic map.
Then the set $\{x\in\Sigma|\ L\neq \Stab_{S^1}s(x)\}$ is finite.
\label{stabilL}
\end{theorem}

\subsection{Proof of theorem \ref{stabilH}}
In order to prove the theorem we will use the following result
on holomorphic curves (see lemma 2.2.3 in \cite{McDS1}).

\begin{lemma}
Let $I\in\III_{\omega}$ be any complex structure. Let $s_1,s_2:\Sigma\to F$ 
be two simple holomorphic maps with respect to $I$. Let 
$K\subset\Sigma$ be a closed subset of noninjective points for both $s_1$ and 
$s_2$. If the intersection $s_1(K)\cap s_2(K)$ contains infinite points, then
$s_1=s_2$.
\label{interseccio}
\end{lemma}

From now on we fix a complex structure $I\in\III_{\omega,S^1}$, and
we take on $F$ the metric $\omega(\cdot,I\cdot)$. This metric is
$S^1$-invariant because of the invariance of $\omega$ and $I$.

\begin{definition}
Let $\fX\in \Gamma(TF)$ be the vector field generated
by $\imag\in\imag\RR=\Lie(S^1)$. For any $x\in\Sigma$ and
any smooth map $s:\Sigma\to F$ we define
$$\theta_s(x):=\dist(\fX(x),ds(T\Sigma)).$$
\end{definition}

We assume for simplicity that $H=\{1\}$,
and at the end we will say some words about the general case.
We proceed as follows: we assume that the claim of the theorem is not true
and we show that this leads to a contradiction. So we suppose
that for any open set $U\subset\Sigma$ there exists a point 
$x\in U$ and $\alpha\in S^1$ such that $\alpha\cdot s(x)\in s(\Sigma)$.

\subsubsectionr{}
Let $Z=s^{-1}(\{s(z)|z\in\Sigma,\ ds(z)=0\})$ be the set of critical
points. This is a finite set (see lemma 2.2.1 in \cite{McDS1}).
So $s(\Sigma)$ is not contained in $S^1\cdot s(Z)$, since the
latter is a disjoint union of points and circles. Let $T$
be a $S^1$-invariant tubular neighbourhood of $S^1\cdot s(Z)$.
Put $\Sigma'=s^{-1}(F\setminus T)$.

The set of noninjective points $Z'=\{z\in\Sigma|\sharp s^{-1}s(z)>1\}$
can only accumulate at critical points (see the comment before
lemma 2.2.3 in \cite{McDS1}), so $Z''=Z'\cap\Sigma'$ is finite.
Hence, $s(\Sigma')$ is not contained in $S^1\cdot s(Z'')$ and
so we may take a $S^1$-invariant tubular neighbourhood $T''$
of $S^1\cdot s(Z'')$ so that $\Sigma''=s^{-1}(F\setminus T'')\cap\Sigma'$ 
has nonempty interior. 

Let $Y=\{z\in\Sigma|\theta_s(z)=0\}$. This is a closed set.
If the interior $\inter Y\neq\emptyset$ then for $\alpha\in S^1$ near the
identity $\alpha\cdot s(\Sigma)$ and $s(\Sigma)$ meet at
an open set and hence, by lemma \ref{interseccio}, they coincide.
But this is implies that $H(s)\neq\{1\}$, and this finishes the argument.
So we may suppose that there is a small open disk $D_a\subset\Sigma''$
such that $\inf \theta_s|_{D_a}=a>0$. Suppose also that
$S^1\cdot s(D_a)\subset W\subset F$, where $W$ is open and 
$S^1$-invariant, and all points in $W$ have the same stabiliser,
so that $W/S^1$ is a smooth manifold.

\subsubsectionr{}
The composition $$D_a\stackrel{s}{\longrightarrow} W
\stackrel{\pi}{\longrightarrow} W/S^1$$
is an embedding if $D_a$ is chosen small enough. Let 
$N\subset W/S^1$ be an open neighbourhood of
$\pi s(D_a)$ with a submersion $p:N\to\pi s(D_a)$ which is
a left inverse for the inclusion $\pi s(D_a)\hookrightarrow N$.
Let $Y_N=Y\cap (\pi s)^{-1}(N)$. The critical points of
$$\Sigma\cap (\pi s)^{-1}(N)\stackrel{s}{\longrightarrow}
N\stackrel{p}{\longrightarrow}\pi s(D_a)$$
contain $Y_N$. Hence, by Sard's theorem $\pi s(Y_N)\subset \pi s(D_a)$
has measure zero. Since $\pi s(Y_N)$ is closed, its complementary
contains a closed disk $\Sigma_0$. Furthermore, there exists
$b>0$ such that for any $x\in \Sigma_0$ and $\alpha\in S^1$
if $\alpha\cdot s(x)=s(y)\in s(\Sigma)$, then
$\theta_s(y)\geq b$.

From the construction of $\Sigma_0$ we deduce the following result. 

\begin{prop}
\label{chart}
There exist real positive numbers $r,\ \eta,\ \epsilon$ such that
for any $x\in\Sigma_0$ and $\alpha\in S^1$ if
$z=\alpha\cdot s(x)\in s(\Sigma)$, $s^{-1}(z)$ has a unique element
$y\in\Sigma$ and if $D_y=D(y;r)$ is the disk centered at $y$
of radius $r$,
\begin{enumerate}
\item If $w\in s(\Sigma)$ and $d(w,y)<\eta$, then $w\in s(D_y)$.
\item There exists an open neighbourhood $V\subset F$ of
$s(y)$ containing $s(D_y)$ and a chart
$\phi=(\phi_1,\dots,\phi_{2n}):V\to\RR^{2n}$ with 
$\phi(s(y))=0$ such that
\begin{itemize}
\item[2.a.] For any $v\in D_y$, $\phi_3(s(v))=\dots=\phi_{2n}(s(v))$.
\item[2.b.] If $\beta\in[-\epsilon,\epsilon]\subset S^1$, then for
any $v\in D_y$, $\beta\cdot s(v)\in V$ and
$$\phi(\beta\cdot s(v))=\phi(s(v))+(0,0,\beta,0,\dots,0).$$
\end{itemize}
\end{enumerate}
\end{prop}

\subsubsectionr{}
We assume for the rest of the argument that $\diam(s(\Sigma_0))<\eta/2$.
Let us identify $S^1\simeq [0,2\pi)$ so that $0$ is the identity
and consider
\begin{equation}
I=\{(\alpha,x)\in (0,2\pi)\times\Sigma_0|\alpha\cdot s(x)\in s(\Sigma)\}.
\label{defI}
\end{equation}
Thanks to the inequality $\theta_s|_{\Sigma_0}\geq b$ we know that
there exists $\delta>0$ such that
$I\subset [\delta,2\pi-\delta]\times \Sigma_0$. Clearly $I$ is closed.
By our assumption the image of the projection 
$\pi_{\Sigma}:I\to\Sigma_0$ is dense and so (since it is also closed)
coincides with $\Sigma_0$. Let now $[0,\mu]\subset[-\epsilon,\epsilon]$ be 
a subset such that for any $\nu\in[0,\mu]$ and for any 
$x\in F$, $d(x,\nu\cdot x)<\eta/2$.

\subsubsectionr{}
Cover $[\delta,2\pi-\delta]$ with closed intervals $A_1,\dots,A_r$ of
length $<\mu$ and let $I_k=I\cap A_k\times\Sigma_0$.
Since $\pi_{\Sigma}(I_1)\cup\dots\cup\pi_{\Sigma}(I_r)=\Sigma_0$
and $\pi_{\Sigma}(I_l)$ is closed for any $l$, there exists
a $\pi_{\Sigma}(I_k)$ with nonempty interior. Let $D\subset
\inter \pi_{\Sigma}(I_k)$ be a disk, and take $x\in D$.
By assumption there exists $\alpha\in A_k$ such that
$\alpha\cdot s(x)=s(y)$, $y\in\Sigma$. By property 2 in \ref{chart} 
there exists
an open set $V\subset F$ containing $s(D_y)=s(D(y,r))$ and a chart
$$\phi:V\to\RR^{2n}.$$
For any $z\in D$ there exists $\beta\in A_k$ such that 
$\beta\cdot s(Z)\in s(\Sigma)$. On the other hand, since
$d(\alpha\cdot s(z),\alpha\cdot s(x))=d(s(z),s(x))<\eta/2$
and $|\alpha-\beta|<\mu$ we have
$$d(\beta\cdot s(z),s(y))<\eta.$$
Hence, by property 1 in proposition \ref{chart}, 
$\beta\cdot s(z)\in s(D_y)$. 
So by 2.b in \ref{chart}, if $w=s(z)$, then $\phi_3(w)=\alpha-\beta$,
$\phi_4(w)=\dots=\phi_{2n}(w)=0$. This implies that for any
$z\in D$, $\sharp I_k\cap\{z\}\times A_k=1$.
Let $(z,h(z))$ be the unique element of this set. The function
$h:D\to A_k$ is $h(z)=\alpha-\phi_3(s(z))$ and so is continuous.
Hence there exists $c\in A_k$ such that $\sharp h^{-1}(c)=\infty$
(see lemma \ref{tururut} below). From this we see that 
$c\cdot s(\Sigma)\cap s(\Sigma)$ has infinite points which do not
accumulate on critical points (since $s(\Sigma_0)$ is at positive distance
from the $S^1$-orbit of the image of any critical point of $s$). 
Finally, using lemma \ref{interseccio} we deduce that
$c\cdot s(\Sigma)=s(\Sigma)$, in contradiction with the
assumption $H=\{1\}$. This finishes the proof of the case $H=\{1\}$.

\subsubsectionr{}
In the general case we proceed as follows. Since the case $H=S^1$
is trivial, we suppose that $\sharp H<\infty$. We assume that
for any open set $U\subset\Sigma$ there exists $x\in U$ and
$\alpha\in S^1\setminus H$ such that $\alpha\cdot s(x)\in s(\Sigma)$.
We do exactly the
same thing as in the case $H=\{1\}$ to get a subset $\Sigma_0\subset\Sigma$
(note that the function $\theta_s(x)$ is equivariant under the action of
$H$). Now, the set $I$ defined in (\ref{defI}) is at
positive distance from $H\times\Sigma_0$. So the element
$c\in S^1$ found at the end of the reasoning does not belong
to $H$, and hence the fact that $c\cdot s(\Sigma)=s(\Sigma)$ leads to  
a contradiction.

\subsubsectionr{}
The following lemma finishes the proof of the theorem.

\begin{lemma}
Let $I=[0,1]$, and let $h:I^2\to I$ be a continuous map. 
There exists $c\in I$ such that $\sharp h^{-1}(c)=\infty$.
\label{tururut}
\end{lemma}
\begin{pf}
The image $h(I^2)$ is an interval $[a,b]\subset I$. If
$a=b$ we put $c=a=b$, and the proof is finished. Otherwise, let 
$a<c<b$. Assume that $C=h^{-1}(c)$ is a finite set of points. 
Let $x,y\in I^2$ such that $h(x)=a$ and $h(y)=b$.
We can find a path $\gamma:I\to I^2\setminus C$ from $x$ to $y$.
Now the sandwich principle implies that there exists
$m\in I$ such that $h(\gamma(m))=c$, in contradiction with
the definition of $C$. (In fact, $C=h^{-1}(c)$ will neither
be countable. If it were, then $I^2\setminus C$ would 
be arc connected as in the finite case, and the same reasoning
would lead to a contradiction.)
\end{pf}

\subsection{Proof of theorem \ref{stabilL}}
Let $\Sigma'=\{x\in\Sigma|\ L\subsetneqq\Stab_{S^1}s(x)\}$ and
suppose that $\sharp\Sigma'=\infty$.
Since the set $$\{\Stab_{S^1}x|\ x\in F\}$$ 
of stabilizers is finite (see lemma \ref{stabfinit} in the appendix), we 
may assume that there exists a group $L''$ strictly containing $L$ such that
$$\Sigma''=\{x\in\Sigma|\ \Stab_{S^1}s(x)=L''\}$$
has infinite elements. Let now $\theta\in L''\setminus L$.
Then $s(\Sigma)$ and $\theta\cdot s(\Sigma)$ intersect 
at an infinite set $\Sigma''$ of points. Hence by theorem
\ref{interseccio} they coincide, and so $\theta\in H(s)$. 
But now $\gamma(\theta)\in\Aut(\Sigma)$ has infinitely many 
fixed points (all the points in $\Sigma''$), and so it must 
be the identity. So we deduce from the exact sequence (\ref{salvat}) 
that $\theta\in L$, which is a contradiction.

\section{Smoothness of moduli of holomorphic curves}

Let us take a closed subgroup
$L\subset S^1$. Recall that the fixed point set $F^L\subset F$
is a compact symplectic submanifold (with possibly several
connected components of different dimension). 
The action of $S^1$ on $F$ gives an action of the Lie group
$S^1/L$ on $F^L$.

Fix a compact Riemann surface $\Sigma$ and
a group of automorphisms $\Gamma\subset\Aut(\Sigma)$.
Assume that there is an injection $\rho:\Gamma\to S^1/L$ 
with closed image. The morphism $\rho$ allows us to speak about
$\Gamma$-equivariant maps $s:\Sigma\to F^L$. These are simply
the maps $s$ which satisfy $s(g\cdot x)=\rho(g)\cdot s(x)$ 
for any $x\in\Sigma$. We will denote
\begin{equation}
\Map^L(\Sigma,F)_{L^p_1}^{\Gamma,\rho}
\label{mapequi}
\end{equation}
the set of $\Gamma$-equivariant maps $s\in\Map^L(\Sigma,F)_{L^p_1}$.

Let us take a $S^1$-invariant complex structure $I\in\III_{\omega,S^1}$.
Let $p>2$ be any real number. We define the moduli of $(L,\Gamma,\rho)$ 
equivariant curves with respect to $I$ to be
$$\MMM_I^{L,\Gamma,\rho}(\Sigma,F)
=\{s\in\Map^L(\Sigma,F)_{L^p_1}^{\Gamma,\rho}|\ \ov{\partial}_Is=0,\mbox
{ $s$ simple and $\Gamma$-equivariant }\}.$$
For any $A\in H_2(F;\ZZ)$, let also $\MMM_I^{L,\Gamma,\rho}(\Sigma,F;A)=
\{s\in\MMM_I^{L,\Gamma,\rho}(\Sigma,F)\mid s_*[\Sigma]=A\}.$

\begin{theorem}
There is a subset $\III^{L,\Gamma,\rho}\subset\III_{\omega,S^1}$
of Baire second category such that for any $I\in\III^{L,\Gamma,\rho}$
the moduli space $\MMM_I^{L,\Gamma,\rho}(\Sigma,F)$ is smooth and oriented. 
\label{regularitat}
\end{theorem}
\begin{pf}
The proof, with due modifications, is exactly like that of 
theorem \ref{gensigmasmooth} or of theorem 3.1.2 in \cite{McDS1}.
We start considering the completion $\III_{\omega,S^1}^l$ (resp.
$\III_{\omega}^l$) of $\III_{\omega,S^1}$ (resp. $\III_{\omega}^l$) 
in the $C^l$ norm, where $l>0$ is a big
enough integer, and we define, for any $A\in H_2(F;\ZZ)$,
the {\bf $A$-universal moduli space} to be
$$\MMM_{\III^l}^{L,\Gamma,\rho}(\Sigma,F;A)
=\left\{(s,I)\in\Map^L(\Sigma,F)_{L^p_1}^{\Gamma,\rho}
\times \III_{\omega,S^1}^l\Big|
\begin{array}{l}\ov{\partial}_Is=0,\ s_*[\Sigma]=A,\\
\mbox{and $s$ simple }\end{array}\right\}.$$

\begin{prop}
The $A$-universal moduli space $\MMM_{\III^l}^{L,\Gamma,\rho}(\Sigma,F;A)$
is a smooth Banach manifold.
\label{udoscxstr}
\end{prop}
\begin{pf} Let $(s,I)\in 
\MMM_{\III^l}^{L,\Gamma,\rho}(\Sigma,F;A)$
be any point. We will prove that the $A$-universal moduli space
is smooth at $(s,I)$. For that we have to check that the linearisation
$$D\FFF(u,I):\Omega^0(s^*TF^L)^{\Gamma}_{L^p_1}\times T_I\III_{\omega,S^1}^l
\to \Omega^1(s^*TF^L)^{\Gamma}_{L^p}$$
of the equation at $(s,I)$ is surjective. Here we write 
$\Omega^i(s^*TF^L)^{\Gamma}$ the $\Gamma$-invariant sections 
of $\Lambda^{0,i} T\Sigma\otimes_I s^*TF^L$ (this bundle has an
action of $\Gamma$ through the representation $\rho$).
The tangent space $T_I\III_{\omega,S^1}^l\subset T_I\III_{\omega}^l$ 
is equal to the subspace of $\Gamma$-invariant elements in 
$T_I\III_{\omega}^l$. This latter space is the set of $C^l$ sections
of the bundle $\End(TF,I,\omega)$ whose fibre at $x\in F$ is the
space of linear maps $Y:T_xF\to T_xF$ which satisfy
$$YI+IY=0\mbox{ and }\omega(Y\cdot,\cdot)+\omega(\cdot,Y\cdot)=0$$
(see p. 34 in \cite{McDS1}). 

We now follow the notation (and the ideas) of the proof of proposition 
3.4.1 in \cite{McDS1}. We may write the differential
$$D\FFF(s,I)(\xi,Y)=D_s\xi+\frac{1}{2}Y(s)\circ ds\circ j,$$
where $j$ is the complex structure in $\Sigma$ and $D_s$ is a first 
order differential operator whose symbol coincides with that of
Cauchy-Riemann operator. Hence $D_s$ is elliptic and consequently
Fredholm. So if it were not exhaustive there would exist a nonzero element
$\eta\in \Omega^1(s^*TF^L)^{\Gamma}_{L^q}$ (where $1/p+1/q=1$) such
that for any $\xi\in \Omega^0(s^*TF^L)^{\Gamma}$ and for any
$Y\in T_I\III_{\omega,S^1}^l$
\begin{equation}
\int_{\Sigma}\la\eta,D_s\xi\ra=0
\mbox{ and }\int_{\Sigma}\la\eta,Y(s)\circ ds\circ j\ra=0.
\label{conseq}
\end{equation}

We now invoque theorem \ref{stabilH} and obtain a disk $D\subset\Sigma$
such that $$S^1\cdot s(D)\cap s(\Sigma)=H\cdot s(D).$$
Using theorem \ref{stabilL} we deduce that (after possibly
shrinking $D$) all the elements in $s(D)$ have stabiliser equal to $L$.
Then $\eta$ vanishes on an open subset of $D$. For suppose that
$\eta(x)\neq 0$, where $x\in D$. One can always find an endomorphism
$Y_0\in\End(T_{s(x)}F,I_{s(x)},\omega_{s(x)})^L$ such that 
$$\la \eta(x),Y_0\circ ds(x)\circ j(x)\ra\neq 0,$$
since $\eta(x)\in T_{s(x)}F^L$.
We extend $Y_0$ to $S^1\cdot s(x)$ in a $S^1$ equivariant 
way (we can do this because $\Stab_{S^1}s(x)=L$ and we took $Y_0$ to 
be $L$-invariant) and then we use a $S^1$-invariant smooth cutoff 
function to extend $Y_0$ to a small neighbourhood of $S^1\circ s(x)$.
It turns out that this can be done in such a way that,
being $\eta$ $\Gamma$-equivariant, the right
hand side integral in (\ref{conseq}) does not vanish. And this is
a contradiction. 

Consequently $\eta$ vanishes in $D$. Since it also satisfies the left 
hand side equation in (\ref{conseq}), Aronszajn's theorem \cite{Ar} 
(see theorem 2.1.2 in \cite{McDS1}) implies that $\eta$ vanishes identically. 
So the linearisation $D\FFF(s,I)$ must be exhaustive. This finishes the proof.
\end{pf}

Now the proof of theorem \ref{regularitat} is resumed as that of
theorem \ref{gensigmasmooth} or in p. 36 in \cite{McDS1}. One uses the
Sard-Smale theorem (for that $l$ has to be big enough, depending on the 
index of the linearisation $D\FFF$, which on its turn is a function
of $A\in H_2(F;\ZZ)$) to prove the existence of
a subset $(\III^{L,\Gamma,\rho})^{A,l}\subset\III_{\omega,S^1}^l$ 
of the second category
such that for any $I\in (\III^{L,\Gamma,\rho})^{A,l}$ the moduli space
$\MMM_I^{L,\Gamma,\rho}(\Sigma,F;A)$ is smooth. Then a trick of
Taubes allows to deduce from this, making $l\to\infty$, that
there exists a subset $(\III^{L,\Gamma,\rho})^A\subset\III_{\omega,S^1}$
of the second category with the same property, but consisting
of smooth complex structures and not of $C^l$ ones as before.
Since the set of homology classes $A\in H_2(F;\ZZ)$ is
countable, the intersection 
$$\III^{L,\Gamma,\rho}=\bigcap_{A\in H_2(F;\ZZ)}(\III^{L,\Gamma,\rho})^A$$
is again of the second category.

To finish the proof, note that the linearisation of the equations
is, modulo a compact operator, the Cauchy-Riemann operator. Hence
the cohomology groups of the deformation complex carry natural
orientations (because they are complex vector spaces) and 
consequently so does the moduli space (see also section \ref{secKuranishi}
on Kuranishi models).
\end{pf}

Let us write for simplicity 
$$\MMM_I^{L,\Gamma,\rho}(A)=\MMM_I^{L,\Gamma,\rho}(\Sigma,F;A)
\qquad\mbox{ and }\qquad\MMM_{\III}^{L,\Gamma,\rho}(A)=
\MMM_{\III}^{L,\Gamma,\rho}(\Sigma,F;A).$$ We now define 
$$\III^{\reg}_{\omega,S^1}=\bigcap_{L,\Gamma,\rho} \III^{L,\Gamma,\rho},$$
where the intersection is taken for the triples $(L,\Gamma,\rho)$
such that the moduli $\MMM_{\III}^{L,\Gamma,\rho}(A)$ is nonempty
for some $A\in H_2(F;\ZZ)$. Again, this is a Baire set of the
second category.

Using the same methods as above one can prove the following theorem,
which says that the moduli spaces $\MMM_I^{L,\Gamma,\rho}(A)$
obtained from different complex structures $I$ in $\III^{\reg}_{\omega,S^1}$
are smooth oriented cobordant.                                             

\begin{theorem}
Let us fix a triple $(L,\Gamma,\rho)$ and a homology class $A\in H_2(F;\ZZ)$.
Let $I_0,I_1\in \III^{\reg}_{\omega,S^1}$. There exists a path
$[0,1]\ni\lambda\mapsto I_{\lambda}\in\III_{\omega,S^1}$
such that the space
$$\bigcup_{\lambda\in[0,1]}\MMM_{I_\lambda}^{L,\Gamma,\rho}(A)$$
is a smooth oriented cobordism between 
$\MMM_{I_0}^{L,\Gamma,\rho}(A)$ and $\MMM_{I_1}^{L,\Gamma,\rho}(A)$.
\end{theorem}

\section{Dimension of the moduli (case $\Sigma=\CP^1$)}

In this section we will compute the dimensions of the moduli 
$\MMM_I^{L,\Gamma,\rho}(\Sigma,F;A)$ for $I\in \III^{\reg}_{\omega,S^1}$.
We will restrict to the case $\Sigma=\CP^1$ (this will be enough for our
purposes,      
since we will use the moduli spaces $\MMM_I^{L,\Gamma,\rho}(\Sigma,F;A)$
to parametrise bubbles). We will also make the assumption that
the action of $S^1$ on $F$ is almost free. Recall that this means
that the action of $S^1$ on $F\setminus F^{S^1}$ is free and that it implies
that the action of $S^1$ on the normal bundle of $F^{S^1}\subset F$
has weights belonging only to $\{-1,1\}$ (see lemma \ref{nomesmesmenysu} 
in the appendix).

\subsection{Preliminaries}
Let us fix a triple $(L,\Gamma,\rho)$. Recall that $L$ is a subset of $S^1$,
$\Gamma$ is a compact subgroup of $\Aut(\CP^1)$ and $\rho:\Gamma\to S^1/L$
is an injection. This latter condition implies that $\Gamma$ is
abelian. Since $\Aut(\CP^1)=\PSL(2,\CC)$, any element in 
$\Aut(\CP^1)$ which spans a compact subgroup must fix two points
of $\CP^1$. And since $\Gamma$ is abelian, there
must exist two points $x_+$ and $x_-$ which are fixed by all the
elements of $\Gamma$. 

Using one of the fixed points, say $x_+$, we
get an injection $\Gamma\to S^1\subset\CC^*$ by assigning to
any $\gamma\in\Gamma$ the induced endomorphism 
$\iota(\gamma)\in\GL(T_{x_+}\CP^1)$.
In the sequel we will identify $\Gamma$ with its image in $S^1$.
There are two possibilities. Either $\Gamma$
is a finite group or $\Gamma\simeq S^1$. When $\Gamma$ is a finite group,
the map $\iota$ fixes an isomorphism $\Gamma\simeq\ZZ/m\ZZ$, and
when $\Gamma$ is infinite $\iota$ gives an identification with $S^1$.

On the other hand, if $\Gamma\neq\{1\}$ then, for any $\Gamma$-equivariant
map $s:\CP^1\to F$, the fixed points $x_\pm$ are mapped by $s$
to the fixed point set $F^{S^1}$ (because the action on $F\setminus F^{S^1}$
is free). Let $z$ be a holomorphic coordinate in $\CP^1$ centered at
$x_+$. Taking $S^1$-equivariant coordinates in a neighbourhood of $s(x_+)$ 
the map $s$ can be written (see p. 16 in \cite{McDS1})
\begin{equation}
s(z)=az^l+O(|z|^{l+1}),
\label{locally}
\end{equation}
and the constant $a$ can be identified with an element of $T_{s(x_+)}F$.
We call $l$ the {\bf index} of $s$ at $x_+$.
Let $T_{x_\pm}^PF$ (resp. $T_{x_\pm}^ZF$, $T_{x_\pm}^NZ$) be the subspace
of $T_{x_\pm}F$ spanned by vectors of weight $1$ (resp. $0$, $-1$)
under the action of $S^1$. Since there are no more weights, $a$ must
lie in $T_{x_\pm}^PF\cup T_{x_\pm}^ZF\cup T_{x_\pm}^NF$ (otherwise
it would not be invariant under the action of $\Gamma$). Using
(\ref{locally}) we may write for any $\theta\in\Gamma$ and $z$ near $x_+$
$$s(\theta z)=\rho(\theta)\cdot s(z)=\theta^l s(z)$$
modulo $O(|z|^{l+1})$. The $\cdot$ in the second term refers to
the action of $S^1$ on $T_{x_+}F$. From this we deduce the following
\begin{itemize}
\item $a$ cannot belong to $T_{x_\pm}^ZF$,
\item if $a\in T_{x_\pm}^PF$ then $\rho(\theta)=\theta^l$,
\item if $a\in T_{x_\pm}^NF$ then $\rho(\theta)=\theta^{-l}$.
\end{itemize}
In fact, after possibly composing $s$ with the holomorphic
map $r:\CP^1\to\CP^1$ defined $r([x:y])=[y:x]$ in coordinates
for which $x_+=[0:1]$ and $x_-=[1:0]$, we may assume
that $a\in T_{x_\pm}^PF$. Hence $\rho(\theta)=\theta^l$ for
any $\theta\in\Gamma$, where $l$ is a positive integer.
If $\Gamma=S^1$, then $l$ must be $1$, and if $\Gamma=\ZZ/m\ZZ$
then $l$ and $m$ must be coprime and the representation $\rho$
only depends on the class of $l$ modulo $m$. 

In the sequel we will write $\MMM_I^{L,\Gamma,l}(A)$ 
instead of $\MMM_I^{L,\Gamma,\rho}(\Sigma,F;A)$. When $\Gamma=1$
we will write $\MMM_I^{L}(A)$ and when $L=\Gamma=1$
we will write $\MMM_I(A)$.

\subsection{The deformation complex}
Let $A\in H_2(F;\ZZ)$ be any cohomology class, 
take a complex structure $I\in\III^{\reg}_{\omega,S^1}$, and
let $s\in\MMM_I^{L,\Gamma,l}(A)$. The deformation complex
of the moduli at $s$ is
$$D_s^{\Gamma}:\Omega^0(s^*TF^L)^{\Gamma}\to\Omega^1(s^*TF^L)^{\Gamma},$$
where $D_s^{\Gamma}$ is equal to the Cauchy-Riemann operator modulo a compact
operator (see p. 28 in \cite{McDS1}). Since $I\in\III^{\reg}_{\omega,S^1}$,
this operator is exhaustive and consequently the dimension of 
$\MMM_I^{L,\Gamma,l}(A)$ at $s$ is equal to 
$\dim (\Ker D_s^{\Gamma})$. To compute this dimension we consider the 
natural extension of $D_s^{\Gamma}$
\begin{equation}
D_s:\Omega^0(s^*TF^L)\to\Omega^1(s^*TF^L)
\label{defcomg}
\end{equation}
(this is the deformation complex of the moduli of holomorphic
curves in $F^L$). The operator $D_s$ is $\Gamma$-equivariant,
and hence acts on the cohomology groups $H^i_s$ of the complex.
We have $\Ker D_s^{\Gamma}=(H^0_s)^{\Gamma}$
and $\Coker D_s^{\Gamma}=(H^1_s)^{\Gamma}=0$. So the complex
dimension at $s$ is equal to 
\begin{equation}
\dim_\CC T_s\MMM_I^{L,\Gamma,l}(A)=\dim (H^0_s)^{\Gamma}
-\dim (H^1_s)^{\Gamma}.
\label{ladimensio}
\end{equation}
This dimension can be computed putting instead of $D_u$ any Dolbeaut operator 
on $s^*TF^L$, since they have the same symbol. Because the action of $S^1$ 
on $F$ is almost-free, we need only distinguish these possibilities.

\noindent{\bf Case 1.} $L=S^1$, $\Gamma=\{1\}$. 
Let $F^{S^1}=F_1\cup\dots\cup F_r$ be the connected 
components of the fixed point set. Suppose that 
$A\in H_2(F_k;\ZZ)\subset H_2(F;\ZZ)$. Then by Riemann-Roch
the moduli space has dimension
$$\dim_\CC \MMM_I^{S^1}(A)=\la c_1(TF_k),A\ra+\dim_\CC F_k.$$

\noindent{\bf Case 2.} $L=\{1\}$, $\Gamma\neq\{1\}$.
As we said in the preceeding section,
there exist two fixed points $x_{\pm}\in\CP^1$ of the action of $\Gamma$.
Since the map $s$ is $\Gamma$-equivariant, we have a natural lift of 
the action $\rho$ of $\Gamma$ on $\CP^1$ to $E=s^*TF\to\CP^1$. 
Let us write it $$\gamma:\Gamma\to\Aut(E),$$
where $\Aut(E)$ denotes the automorphisms of $E$ as vector bundle.
The map $\gamma$ induces representations $\gamma_{\pm}$ of 
$\Gamma$ on the fibres $E_{x_\pm}$ over $x_{\pm}$. The weights of this 
representation are $l$ times the weights of the representation of $S^1$ 
on $TF_{s(x_\pm)}$. Since the action of $S^1$ is almost-free, the
weights of $S^1$ acting on $TF_{s(x_\pm)}$ belong to $\{-1,0,1\}$.
Let $P_\pm$ (resp. $Z_\pm$, $N_\pm$) be the number of weights
of the representation $\gamma_\pm$ which are equal
to $1$ (resp. $0$, $-1$). Using the Riemann-Roch theorem and the Atiyah-Bott
fixed point formula one can compute the dimension (\ref{ladimensio})
in terms of $P_\pm$, $Z_\pm$, $N_\pm$ and $l$. Let us write the dimension
$$\Ind_{\gamma}(E)=\dim (H^0(E))^{\Gamma}-\dim (H^1(E))^{\Gamma}.$$
We will compute the dimensions $\Ind_{\gamma}(E)$ in the next two sections.
The results obtained may be summarized in the following two theorems,
which correspond respectively to theorems \ref{inds1} and \ref{indtor}.
\begin{theorem}
Suppose that $\Gamma=S^1$.
The degree of $E$ is $\deg(E)=P_++N_--P_--N_+$ and the index is
$$\Ind_{\gamma}(E)=(P_++Z_+)+(N_-+Z_-)-\rk(E).$$
\label{dim1}
\end{theorem}
\begin{theorem}
Suppose that $\Gamma=\ZZ/m\ZZ$ and
take $l'=l+km$, $k\in\ZZ$, such that $1\leq l'\leq m-1$. 
Then $$\Ind_\gamma(E)=\frac{1}{m}(\deg(E)+m\rk(E)-m(P_-+N_+)+
l'(P_-+N_+-P_+-N_-)).$$
\label{dim2}
\end{theorem}
Note that corollary \ref{omesuomenysuozero} in the appendix implies that 
\begin{equation}                          
P_++N_++P_-+N_-\geq 2.
\label{prineq}
\end{equation}

From these theorems we deduce the following result.
\begin{theorem} Let $I\in\III^{\reg}_{\omega,S^1}$.
If $F$ is positive with respect to the complex structure $I$ then
for any $\Gamma\neq\{1\}$ and $l$
$$\dim_\CC T_s\MMM_I^{1,\Gamma,l}(A)\leq\dim_\CC T_s\MMM_I(A)-2.$$
\label{claudim}
\end{theorem}
\begin{pf}
Let $\Gamma\neq\{1\}$ and let $s\in\MMM_I^{1,\Gamma,l}(A)$. Let us write 
$E=s^*TF$. Since $I\in\III^{\reg}_{\omega,S^1}$ the moduli spaces
$\MMM_I^{1,\Gamma,l}(A)$ and $\MMM_I(A)$ are 
smooth, and their dimensions are equal to their virtual dimensions.
Hence, what we have to prove is                                  
$\Ind_{\gamma}(E)\leq\Ind(E)=\deg(E)+\rk(E)$.
Note that, since $F$ is positive, we have by definition 
$\deg(E)=\deg(s^*TF)\geq 1$.

Suppose that $\Gamma=S^1$. In this case we have to prove that
$$(P_++Z_+)+(N_-+Z_-)-\rk(E)\leq \deg(E)+\rk(E)-2.$$
Writing $2\rk(E)=P_++N_++Z_++P_-+N_-+Z_-$
the inequality gets transformed into
$2\leq\deg(E)+(N_++P_-)$.
Now, summing the inequalities (\ref{prineq}) and $\deg(E)\geq 1$
(and using $\deg(E)=P_++N_--P_--N_+$) we obtain 
$3/2\leq\deg(E)+(N_++P_-)$. This implies our inequality because
$\deg(E)+(N_++P_-)$ is an integer.

Now suppose that $\Gamma=\ZZ/m\ZZ$. Let $l'$ be as in theorem \ref{dim2},
so that $1\leq l'\leq m-1$. Since the indices $\Ind_{\gamma}(E)$
and $\Ind(E)$ are integers, it is enough for our purposes to prove
$$\Ind_{\gamma}(E)\leq \deg(E)+\rk(E)-(1+1/m).$$
Writing the value of $\Ind_{\gamma}(E)$ given by theorem \ref{dim2},
multiplying by $m$ and simplifying we arrive at the (equivalent) inequality
$$m+1\leq (m-1)\deg(E)+(m-l')(P_-+N_+)+l'(P_++N_-),$$
which is a consequence of (\ref{prineq}) and $\deg(E)\geq 1$, taking
into account that $m-l'\geq 1$ and $l'\geq 1$.
\end{pf}

\subsection{Comments on the case $\Gamma\neq\{1\}$}
\label{comentaris}
When $\Gamma=\ZZ/m\ZZ\subset S^1$ the topological space $F/\Gamma$ has the 
structure of an orbifold, and $\omega$ together with any 
$I\in\III_{\omega,S^1}$ provide $F/\Gamma$ with the structure of an almost 
Kaehler orbifold. The curves $s:\CP^1\to F$ for which 
$H(s)=\Gamma$ give multicovered curves when composed with the projection
$F\to F/\Gamma$. This {\it explains} why they give problems when trying 
to obtain smooth moduli for invariant complex structures.

One can interpret the curves in
$\MMM_I^{1,S^1,1}(A)$ in the following terms. 
Let $s\in\MMM_I^{1,S^1,1}(A)$, and take coordinates
$[x:y]$ on $\CP^1$ for which $x_+=[1:0]$ and $x_-=[0:1]$. 
Then the action of $S^1$ on $\CP^1$ is given by multiplication of
the second coordinate: $\theta[a:b]=[a:\theta b]$.

Let $C$ be the cilinder $\RR\times S^1$. Consider the embedding
$\iota:C\to\CP^1$ which maps $(t,\alpha)$
to $[1:e^t\alpha]$. This is a conformal map, so that
$s\circ\iota:C\to F$ is holomorphic. Furthermore,
$s\circ\iota$ is $S^1$ equivariant with respect to the $S^1$
action on $\RR\times S^1$ given by multiplication on the second
factor.

Let $\fX^C$ (resp. $\fX^F$) be the vector field
generated by the infinitesimal action of $\imag\in\Lie(S^1)$
on $C$ (resp. on $F$). By equivariance we have
\begin{equation}
d(s\circ\iota)\fX^C=\fX^F.
\label{quasgrad}
\end{equation}
Let $I_C$ (resp. $I_F$) denote the complex structure on $C$
(resp. on $F$). The lines $$l_\alpha=\RR\times\{\alpha\}\subset C$$ are
integral curves of the field $-I_C\fX^C$.
Hence, by (\ref{quasgrad}) and holomorphicity, $(s\circ\iota)l_\alpha$ 
are integral lines of $-I_F\fX^F$. But the field $-I_F\fX^F$ 
is the gradient of the function $f=\mu(\imag)$ with respect
to the metric $g=\omega(\cdot,I_F\cdot)$ (see lemma \ref{gradient}).
So for any $\alpha\in S^1$ the curve $(s\circ\iota)l_{\alpha}$
is a line of steepest descent joining the critical points
$s(x_+)$ and $s(x_-)$. Furthermore, the energy of $s$ is
$$E(s)=2\pi\int_{\RR}|\gamma_{\alpha}'(t)|^2dt=
2\pi|\mu(\imag)(s(x_+))-\mu(\imag)(s(x_-))|,$$
where $\gamma_\alpha(t)=(s\circ\iota)(t,\alpha)$.
Assigning to any $s\in\MMM_I^{1,S^1,1}(A)$ the path $\gamma_1$
we get the following result.

\begin{prop}
Let $F_+$, $F_-$ be two connected components of the fixed point set
of the $S^1$ action on $F$. There is a one to one correspondence between 
holomorphic $S^1$ equivariant maps $s:\CP^1\to F$ with $s(x_\pm)\in F_\pm$
and lines of steepest descent of $\mu$ from $F_+$ to $F_-$ with respect
to $g=\omega(\cdot,I\cdot)$. Furthermore, all these maps have
energy $$2\pi|\mu(\imag)(s(x_+))-\mu(\imag)(s(x_-))|.$$
\end{prop}

Let us take now two connected components $F_+$ and $F_-$ of the fixed point
set $F^{S^1}$. Let $Z_\pm=\dim F_\pm$.
Consider the gradient flow $\phi^t:F\to F$ of the
function $f$, which is defined by $\phi^0=\Id$ and 
$\frac{\partial}{\partial t}\phi^t=\nabla f$. Let $W_+$ be the
stable set of $F_+$
$$W_+=\{x\in F\mid \lim_{t\to\infty}\phi^t(x)\in F_+\}$$
and $W_-$ the unstable set of $F_-$
$$W_-=\{x\in F\mid \lim_{t\to-\infty}\phi^t(x)\in F_-\}.$$
The sets $W_\pm$ are locally closed and their dimensions
are $\dim W_+=Z_++P_+$ and $W_-=Z_-+N_-$, where $P_+$ (resp. $N_-$)
is number of positive (resp. negative) weights of the action 
of $S^1$ on the normal bundle of $F_+$ (resp. $F_-$).
The union of the lines of steepest descent connecting $F_+$ and $F_-$ is
$W_+\cap W_-$. It is well known that for a generic metric the resulting
(un)stable sets intersect transversely, and hence the dimension of
their intersection is
$$\dim W_+\cap W_-=\dim(W_+)+\dim(W_-)-\dim(F)=Z_++P_++Z_-+N_--\dim(F).$$

This number coincides with the dimension of $\MMM^{1,S^1,1}$
given in theorem \ref{dim1}. In fact, if $I\in\III^{\reg}_{\omega,S^1}$ and
we take the metric $g=\omega(\cdot,\cdot I)$, then the (un)stable     
sets $W_\pm$ intersect transversely. Moreover, this metric 
is $S^1$-invariant.

As a final comment, we recall that the function $f=\mu(\imag)$
is an equivariantly perfect Morse function (see for example 
\cite{Ki}). It seems interesting to pursue these ideas relating
lines of steepest descent and $S^1$-invariant holomorphic
curves.

\section{Index computations: $S^1$ actions}
Let us identify the projective line $\CP^1$ with the sphere $S^2$ and
consider the standard embedding $S^2\to\RR^3$. Define the morphism
$$\rho:S^1\to\Aut(\CP^1)$$ by sending any $\theta\in S^1$
to the rotation of angle $\theta$ and axis $x=y=0$.
Let $x_+=(0,0,1)$ and $x_-=(0,0,-1)$ be the north and south pole.
These are the fixed points of the action $\rho$.

Let $E\to\CP^1$ be a holomorphic vector bundle. Suppose that there
is a lift $$\alpha:S^1\to\Aut(E)$$ of the action $\rho$
($\Aut(E)$ is defined to be the automorphisms of $E$ which are
linear on fibres). Then there is an induced action of $S^1$
on $H^0(E)$ and $H^1(E)$. Let us write
$$\Ind_{\alpha}(E)=\dim H^0(E)^{S^1}-\dim H^1(E)^{S^1}.$$
On the other hand, since $x_\pm$ are kept fixed by the action of
$\rho$, the morphism $\alpha$ induces an action of $S^1$ on the
fibres over $x_\pm$. Hence we have two representations
$$\alpha_{\pm}:S^1\to\GL(E_{x_\pm}),$$
where $E_x$ denotes the fibre over $x\in\CP^1$.

In this section we will compute $\Ind_{\alpha}(E)$ in terms of the 
degree and rank of $E$ and the weights of the representations $\alpha_{\pm}$.
Our tool will be Grothendieck theorem on vector bundles
over the projective line. We will prove the following

\begin{theorem}
Let $P_\pm$ (resp. $Z_\pm$, $N_\pm$) be the number of strictly positive
(resp. zero, strictly negative) weights in the representation $\alpha_\pm$.
Then $\deg(E)=P_++N_--P_--N_+$ and
$$\Ind_{\alpha}(E)=(P_++Z_+)+(N_-+Z_-)-\rk(E).$$
\label{inds1}
\end{theorem}

\subsection{Preliminaries}
In this subsection we define some natural actions of $S^1$ on
$\OOO(k)\to\CP^1$ which lift the composition of $\rho$ with
the map $m_2:S^1\ni\theta\mapsto\theta^2\in S^1$.
For this action we will be able to completely determine the $S^1$-module
structure of $H^i(\OOO(k))$.

\subsubsectionr{}
What follows is probably the most simple case of Borel-Weil theorem.
Define $$B=\left\{ \left(
\begin{array}{cc}* & * \\ 0 & *\end{array}\right)\in \SL(2,\CC)\right\}.$$
This is a Borel subgroup of $\SL(2,\CC)$, and the corresponding homogeneous 
space is isomorphic to the projective line:
\begin{equation}
\CP^1\simeq\SL(2,\CC)/B. 
\label{isom1}
\end{equation}
Define also the nilpotent subgroups
$$N_+=\left\{ \left(
\begin{array}{cc}1 & * \\ 0 & 1\end{array}\right)\in \SL(2,\CC)\right\}
\mbox{ and }
N_-=\left\{ \left(
\begin{array}{cc}1 & 0 \\ \mbox{$*$} & 1\end{array}\right)\in \SL(2,\CC)\right\}.$$
The principal $\CC^*$-bundle $P$ associated to the
tautological line bundle $\OOO(-1)\to\CP^1$ is also isomorphic to a
homogeneous space: 
\begin{equation}
P\simeq\SL(2,\CC)/N_+.
\label{isom2}
\end{equation}
Using the isomorphism (\ref{isom1}) we get a left action of
$\SL(2,\CC)$ on $\CP^1$. This action admits a canonical lift to 
$P$ (and hence to $\OOO(-1)$) thanks to isomorphism (\ref{isom2}). 
Tensoring we get canonical lifts to any bundle of the form 
$\OOO(k)$, $k\in\ZZ$. 

\subsubsectionr{}
Let us fix from now on an integer $k\in\ZZ$, and let $L_k=\OOO(k)$.
The cohomology groups $H^0(L_k)$ and $H^1(L_k)$ are acted on
linearly by $\SL(2,\CC)$. They are in fact irreducible representations.
One can prove this as follows. Take the subgroup of diagonal matrices
$$T=\left\{ \left(
\begin{array}{cc}* & 0 \\ 0 & *\end{array}\right)\in \SL(2,\CC)\right\}$$
as a maximal torus of $\SL(2,\CC)$. 
Suppose for example that the representation $H^0(L_k)$ is not irreducible.
Then there are two linearly 
independent elements $v_1,v_2\in H^0(L_k)$ of maximal weight. Since they
are elements of maximal weight, $v_1$ and $v_2$
are invariant under the action of the nilpotent subgroup $N_-$.
The action of $N_-$ on $\CP^1$ has orbits which are Zariski
dense. But the ratio $v_1/v_2$ must
be constant along any such orbit, and since $v_1$ and $v_2$ are 
holomorphic, by density the ratio must be the same everywhere. Hence
$v_1$ and $v_2$ are not linearly independent.

\subsubsectionr{}
Consider the inclusion $\iota:S^1\to T$ defined as
$$\iota(\theta)=\left(\begin{array}{cc}
\theta & 0 \\ 0 & \theta^{-1}\end{array}\right).$$
Let $\rho_2$ (resp. $\tau_k$) be the composition of $\iota$ with the 
restriction to $T$ of the canonical action of $\SL(2,\CC)$ on $\CP^1$ 
(resp. $L_k$).
The dimensions of the representations $H^i(L_k)$ 
are $$\dim H^0(L_k)=\max\{0,k+1\}\mbox{ and } \dim H^1(L_k)=
\max\{0,-k-1\}.$$
An irreducible representation of $\SL(2,\CC)$ of dimension
$n$ splits as a representation of $T$ into a sum of $n$ one-dimensional
representations of weights $\{n,n-2,\dots,-n+2,-n\}$
(see for example \cite{FH}).
Thus we obtain a description of the action of $S^1$ on $H^i(L_k)$ 
for any $i$ and $k$. For example, the weights of $\tau_k$ acting 
on $H^0(L_k)$ for $k\geq 0$ are $\{k, k-2,\dots,-k+2,-k\}$.
From this we deduce the following formula
\begin{equation}
\Ind_{\tau_k}(L_k)=\left\{\begin{array}{l}
1\mbox{ if $k$ is even and $k\geq 0$}\\
-1\mbox{ if $k$ is even and $k\leq -2$}\\
0\mbox{ otherwise.}\end{array}\right.
\label{caslinia}
\end{equation}
Note that if $m_2:S^1\to S^1$ sends $\theta$ to $\theta^2$,
then $\rho_2=\rho\circ m_2$. In fact, the action $\tau_k$ descends
to a lift of $\rho$ only when $k$ is even. Indeed, if it descends 
then $\tau_k(-\Id)=\Id$, so that the representation of 
$\SL(2,\CC)$ on $H^i(L_k)$ descends to a representation of
$\PSL(2,\CC)$. But all the weights of the representations of $\PSL(2,\CC)$
are even, so by the classification of representations of $\SL(2,\CC)$
$k$ has to be even.

Finally observe that restricting $\tau_k$ to the fixed points $x_{\pm}$
of $\rho$ we get representations $(\tau_k)_{\pm}:S^1\to\GL((L_k)_{x_\pm})$. 
In fact, if we write $\lambda:S^1\to\CC^*$ the standard inclusion, then
\begin{equation}
(\tau_k)_{\pm}=\lambda^{\pm k}.
\label{pestau}
\end{equation}

\subsection{Proof of theorem \ref{inds1}}
Using Grothendieck theorem we may decompose
\begin{equation}
E=\bigoplus_{k=1}^n E_k
\mbox{, where }E_k=\bigoplus_{r_k}\OOO(\lambda_k).
\label{grot}
\end{equation}
We sort the summands so that $\lambda_1<\lambda_2<\dots<\lambda_n$.
By means of the splitting (\ref{grot}) 
we may use the maps $\tau_k:S^1\to\Aut(\OOO(k))$ defined in the previous
subsection to define a lift $\tau:S^1\to\Aut(E)$ of $\rho_2$.
Let $\alpha_2=\alpha\circ m_2$. This is another lift of $\rho_2$.
Note that we have $$\Ind_{\alpha}(E)=\Ind_{\alpha_2}(E)$$
and the weights of the representations $(\alpha_2)_{\pm}$ are 
twice the weights of $\alpha_{\pm}$.

\subsubsectionr{}
We define a map $\Delta:S^1\to H^0(\End E)$ as 
$\Delta=\tau^{-1}\circ\alpha_2$. Note that this need not be 
a morphism of groups, since $\tau$ and $\alpha$ do not necessarily
commute. Since $H^0(\OOO(k))=0$ for $k<0$, the map $\Delta$ takes
the following form in terms of the splitting $E=\bigoplus E_k$
(we will have this splitting in mind in all the matrices which we
will write in the sequel)
$$\Delta=\left(
\begin{array}{cccc}
\Delta_1 & * & \dots & * \\
0 & \Delta_2 & \dots & * \\
\vdots & \vdots & \ddots & \vdots \\
0 & 0 & \dots & \Delta_n
\end{array}\right).$$
Furthermore (since $H^0(\OOO)=\CC$), $\Delta_k(\theta)$ is a constant 
matrix for any $k$ and $\theta\in S^1$.
Now let $c(t)$ be the diagonal matrix $\diag(1,t,t^2,\dots,t^n)$.
Then the representations $\alpha_2^t=c(t^{-1})\alpha_2 c(t)$
converge as $t\to 0$ to a diagonal representation 
$\alpha_2'=\diag(A_1,\dots,A_n)$ with constant terms. 
Then $\Delta'=\tau^{-1}\circ\alpha_2'$ is equal to
$\diag(\Delta_1,\Delta_2,\dots,\Delta_n)$. By the deformation invariance
of the (equivariant) index, we have
$$\Ind_{\alpha_2}(E)=\Ind_{\alpha_2'}(E).$$
So we will compute the right hand side. Observe also that 
the weights of the representations $(\alpha_2)_{\pm}$ are the
same as those of $(\alpha_2')_{\pm}$.

\subsubsectionr{}
The representations $(\alpha_2')_+$ and $\tau_+$ commute, and
hence $\Delta'$ is really a morphism of groups. So we may diagonalise
each block $\Delta_k$ and 
deduce the result from the case $\rk(E)=1$.
Let $$(w^{k,1},w^{k,2},\dots,w^{k,r_k})$$
be the weights of the representations $\Delta_k$.
Using formula (\ref{pestau}) we deduce that the weights of the representation
$(A_k)_{\pm}$ are 
$$(2a^{k,1}_\pm,\dots,2a^{k,r_k}_\pm)=
(w^{k,1}\pm \lambda_k,\dots,w^{k,r_k}\pm \lambda_k).$$
Observe that, since $\alpha_2=\alpha\circ m_2$,
the weights of $\alpha_\pm$ are $a^{k,j}_\pm$. 
This proves the equality $\deg(E)=P_++N_--P_--N_+$, since we
obviously have $\deg(E)=\sum r_k\lambda_k$.

\subsubsectionr{}
We prove now the formula giving $\Ind_{\alpha}(E)$.
Let us work out the case $\rk(E)=1$. Suppose that $E=\OOO(k)$. 
The weight of $(\alpha_2)_{\pm}$ is $2a_\pm=w \pm k$, where
$w$ is the weight of $\Delta$. 
So $w$ and $k$ have the same parity. Hence
\begin{equation}
\Ind_{\alpha_2}(L_k)=\left\{\begin{array}{l}
1\mbox{ if $k\geq 0$ and }|w|\leq |k|\\
-1\mbox{ if $k\leq -2$ and }|w|\leq |-k-2|\\
0\mbox{ otherwise.}\end{array}\right.
\end{equation}
Using the formulae $k=a_+-a_-$ and $w=a_++a_-$ we may rewrite
the above conditions as
\begin{equation}
\Ind_{\alpha_2}(L_k)=\left\{\begin{array}{l}
1\mbox{ if $a_+\geq 0$ and $a_-\leq 0$} \\
-1\mbox{ if $a_+\leq -1$ and $a_-\geq 1$} \\
0\mbox{ otherwise.}\end{array}\right.
\end{equation}

\subsubsectionr{}
In the general case, let $\{a^1_\pm,\dots,a^N_\pm\}$ be the
weights of the representations $\alpha_{\pm}$. Suppose that after
writing $E$ as sum of line bundles and diagonalising the action
of $S^1$ the two weights $a^k_+$ and $a^k_-$ correspond to the
same line bundle. Then 
\begin{align*}
\Ind_\alpha(E)=\Ind_{\alpha'_2}(E) &= \sharp\{k|\ a^k_+\geq 0,\ a^k_-\leq 0\}
-\sharp\{k|\ a^k_+\leq -1,\ a^k_-\geq 1\} \\
&=\min(P_++Z_+,N_-+Z_-)-\min(N_+,P_-) \\
&= \rk(E)+\min(-N_+,-P_-)-\min(N_+,P_-)\\
&= \rk(E)-\max(N_+,P_-)-\min(N_+,P_-) \\
&= \rk(E)-(N_++P_-)\\
&=(P_++Z_+)+(N_-+Z_-)-\rk(E),
\end{align*}
which is what we wanted to prove (here we have used several times
that $\rk(E)=P_++Z_++N_+=P_-+Z_-+N_-$).

\section{Index computations: $\ZZ/m\ZZ$ actions}

%
Let us fix a natural number $m\geq 2$. Let $\iota_m:\ZZ/m\ZZ\to S^1$ be 
the morphism which sends $1$ to $e^{2\pi\imag/m}$. 
Let $E\to\CP^1$ be a holomorphic vector bundle and suppose that
there is an action $\beta:\ZZ/m\ZZ\to\Aut(E)$ which lifts
$\rho_m=\rho\circ\iota_m$ (we keep the notation of the preceeding
section). We then get an induced action of $\ZZ/m\ZZ$
on $H^i(E)$ and, just as before, our aim is to relate
$$\Ind_{\beta}(E)=\dim H^0(E)^{\ZZ/m\ZZ}-\dim H^1(E)^{\ZZ/m\ZZ}$$
to the weights of the representations $\beta_\pm:\ZZ/m\ZZ\to\GL(E_{x_\pm})$.

Note that whereas the character ring of $S^1$ is $\ZZ$, that of
the group $\ZZ/m\ZZ$ is isomorphic to $\ZZ/m\ZZ$. However,
in this section we will take the extra assumption that the
weights of the representations $\beta_\pm$ belong to 
$\{-l,0,l\}\in\ZZ/m\ZZ$, for $m\nmid l$. We will prove the following

\begin{theorem}
Let $P_\pm$ (resp. $Z_\pm$, $N_\pm$) be the number of weights
of $\beta_\pm$ which are $l$ (resp. $0$, $-l$). Take a representative
$l'\in\ZZ$ of $l$ such that $1\leq l'\leq m-1$. 
Then $$\Ind_\beta(E)=\frac{1}{m}(\deg(E)+m\rk(E)-m(P_-+N_+)+
l'(P_-+N_+-P_+-N_-)).$$
\label{indtor}
\end{theorem}

\subsection{Proof of theorem \ref{indtor}}
We may write
\begin{align*}
\Ind_{\beta}(E) &=
\frac{1}{m}\left(
\sum_{k\in\ZZ/m\ZZ}\Tr(\beta(k),H^0(E))-\Tr(\beta(k),H^1(E))\right) \\
&=\frac{1}{m}\left(\Ind(E)+
\sum_{k=1}^{m-1}\Tr(\beta(k),H^0(E))-\Tr(\beta(k),H^1(E))\right), 
\end{align*}
where $\Tr(\beta(k),H^i(E))$ denotes the trace of $\beta(k)$ 
acting on $H^i(E)$. The index $\Ind(E)$ is equal to $\deg(E)+\rk(E)$ by 
Riemann-Roch. We will compute the value of
$$\Tr(\beta(k),H^0(E))-\Tr(\beta(k),H^1(E))$$ 
for $1\leq k\leq m-1$ using Atiyah-Bott fixed point theorem 
(see \cite{BeGeV}).

\begin{theorem}[Atiyah-Bott]
Let $M$ be a compact complex manifold and $W\to M$ a holomorphic
vector bundle. Let $\gamma:M\to M$ be a complex diffeomorphism
which lifts to $\gamma:W\to W$. Suppose that the fixed points of
$\gamma$ are isolated. Then
$$\sum_i (-1)^i\Tr(\gamma,H^i(W))=
\sum_{x_0\in M^{\gamma}} \frac{\Tr(\gamma_{x_0}^W)}
{\det_{T_{x_0}^{1,0}M}(1-\gamma_{x_0}^{-1})},$$
where $H^i(W)$ is the $i$-th Dolbeaut cohomology group and 
$\gamma_{x_0}^W:W_{x_0}\to W_{x_0}$ is the linear endomorphism of
the fibres over the fixed points induced by $\gamma$.
\label{AtiBot}
\end{theorem}

In our case we have for any $1\leq k\leq m-1$ a complex diffeomorphism
$\rho_m(k)\in\Aut(\CP^1)$ whose fixed points are $x_\pm$. Let
$\theta=\exp(2\pi\imag/m)$. We then have
$$\det(1-\rho_m(k)_{x_\pm}^{-1})=(1-\theta^{\mp 1}).$$
Let $N=\rk(E)$ and let $b^1_\pm,\dots,b^N_\pm\in\ZZ/m\ZZ$ be
the weights of $\beta_\pm$. Then
$$\Tr(\beta(k)_\pm)=\sum_{j=1}^N\theta^{b^j_{\pm}k}.$$
So using theorem \ref{AtiBot} we conclude that
\begin{equation}
\Ind_{\beta}(E)=\frac{1}{m}
\left(\deg(E)+\rk(E)+\sum_{k=1}^{m-1}\sum_{j=1}^N
\left(\frac{\theta^{b^j_+k}}{1-\theta^{-k}}
+\frac{\theta^{b^j_-k}}{1-\theta^k}\right)\right).
\label{quasifi}
\end{equation}

\begin{lemma} Let $\theta=\exp(2\pi\imag/m)$. Then for $1\leq w\leq m-1$
\begin{align*}
\sum_{k=1}^{m-1}\frac{1}{1-\theta^k}&=
\sum_{k=1}^{m-1}\frac{1}{1-\theta^{-k}}=
\frac{m-1}{2} \\
\sum_{k=1}^{m-1}\frac{\theta^{wk}}{1-\theta^k}&=
\sum_{k=1}^{m-1}\frac{\theta^{-wk}}{1-\theta^{-k}}=
-\frac{m-1}{2}+w-1.
\end{align*}
\end{lemma}
\begin{pf}
Let $f(x)=\prod_{k=1}^{m-1}(x-\theta^k)$. We have
$f(x)=1+x+\dots+x^{m-1}$ and
$$\sum_{k=1}^{m-1}\frac{1}{1-\theta^k}=
\frac{f'(1)}{f(1)}=\frac{m(m-1)/2}{m}=\frac{m-1}{2}.$$
In general, for any $1\leq w\leq m-1$ 
\begin{align*}
\sum_{k=1}^{m-1}\frac{\theta^{wk}}{1-\theta^k} &=
\sum_{k=1}^{m-1}\left(-\frac{1-\theta^{wk}}{1-\theta^k}
+\frac{1}{1-\theta^k}\right)
=\sum_{k=1}^{m-1}-(1+\theta^k+\dots+\theta^{(w-1)k})+\frac{m-1}{2}\\
&=-\frac{m-1}{2}+w-1,
\end{align*}
since, for any $w\in \ZZ$, $\sum_{k=1}^{m-1}\theta^{wk}$
is $m-1$ if $m\mid w$ and $-1$ otherwise.
\end{pf}

Now, combining the above lemma with (\ref{quasifi}) we get 
$$\Ind_\beta(E)=\frac{1}{m}(\deg(E)+m\rk(E)-m(P_-+N_+)+
l'(P_-+N_+-P_+-N_-)),$$
which is what we wanted to prove.

\begin{remark}
If the weights of the representations $\alpha_\pm$ induced by an
action $\alpha:S^1\to\Aut(E)$ which lifts $\rho$
belong to $\{-1,0,1\}$, then theorem \ref{inds1}
may be deduced from theorem \ref{indtor} by taking subgroups
of the form $\ZZ/2^r\ZZ\subset S^1$ and making
$r\to\infty$.
\end{remark}

It would be interesting to give a proof of theorem \ref{indtor}
in the lines of the proof of theorem \ref{inds1}, using only
Grothendieck's theorem.


\chapter{The invariants $\Phi$ and $\oPhi$}
\label{invar}

In this chapter the Lie group $K$ will be $S^1$ and its action
on $F$ will be assumed to be almost-free. We will explain how to 
define invariants of the symplectic manifold $F$ and the Hamiltonian 
action of $S^1$ by means of the moduli of THCs. All the (co)homology 
groups appearing in this chapter will be with coefficients in $\ZZ$.

\subsection{The invariant $\Phi$}
\label{ideaphi}
To define this invariant we will use the extended moduli space
$\NNN$. Let $E\to X$ be a principal $S^1$ bundle, $\AAA$ the space
of connections on $E$, and $\GGG=\Map(X,S^1)$ the gauge group of $E$.
As before, $\GGG_0=\{g:X\to S^1\mid g(x_0)=1\}\subset\GGG$, where $x_0\in X$.
Let $\FFF=E\times_{S^1}F$. We will write $\SSS=\Gamma(\FFF)$,
$\SSS^*=\Gamma(\FFF)\setminus\Gamma(\FFF)^{\GGG}$ and 
$\BBB_0=(\AAA\times\SSS^*)/\GGG_0$.

The group $\GGG_0$ acts freely on $\AAA\times E\to\AAA\times X$
(by that we mean that it acts freely both on the base $\AAA\times X$
and on the total space $\AAA\times E$ in a way compatible with the
fibration) with the diagonal actions (the action on $X$ is the trivial one).
Let $\EE_J=(\AAA\times E)/\GGG_0\to\AAA/\GGG_0\times X$ be the quotient
and $\FEj=\EE_J\times_{S^1}F=(\AAA\times \FFF)/\GGG_0$. Let us consider
the map $\oev_J:\AAA\times\SSS^*\times X\to \AAA\times\FFF$
which sends $(A,\Phi,x)$ to $(A,\Phi(x))$. This map is $\GGG_0$
equivariant, so it descends to give a map $$\ev_J:\BBB_0\times X\to\FEj,$$ 
which we will call the {\bf evaluation map}.

We have a map $\rho_{\EE_J}^*:H^*_{S^1}(F)\to H^*(\FEj)$ defined as follows.
Let $\phi:\EE_J\to ES^1$ be any $S^1$ equivariant map
(so that $\phi$ is a lift of the classifying map
$\AAA/\GGG_0\times X\to BS^1$ of the bundle $\EE_J$). This induces
a map $\psi:\FEj\to F_{S^1}$, and $\rho_{\EE_J}^*$ is the morphism
in cohomology induced by $\psi$. In lemma \ref{indepframe}
in the appendix we prove that $\rho_{\EE_J}^*$ is independent
of the choice of $\phi$. Let then $\mu_i:H^*_{S^1}(F)\to H^*(\BBB_0)$, $i=0,2$,
be the maps defined as
$$\mu_2(\delta):=\ev_J^*\rho_{\EE_J}^*\delta/[X]\qquad\mbox{ and }
\qquad \mu_0(\delta):=\ev_J^*\rho_{\EE_J}^*\delta/[\pt],$$
where $\delta\in H^*_{S^1}(F)$, and $[X]\in H_2(X)$ (resp. $[\pt]\in H_0(X)$)
is the fundamental class of $X$ (resp. the class of a point).
Let finally $$\nu:H^*(\AAA/\GGG_0)\to H^*(\BBB_0)$$
be the map induced by the projection $\BBB_0\to\AAA/\GGG_0$.

Let us take a homology class $B\in H_2(F_{S^1})$ such that
${\pi_F}_*B=\eta(E)$, where $\pi_F:F_{S^1}\to BS^1$ is the projection and
where $\eta$ is the map of lemma \ref{etabij} in the appendix.
Let $c\in\imag\RR\setminus C_0$, $I\in\III_{\omega,S^1}$ and
$\sigma\in\Sigma_c^{\reg}(E)$. Then the extended moduli of
$\sigma$-THCs $$\NNN=\NNN^{F,S^1}_{I,\sigma}(B,c)\subset\BBB_0$$
is a smooth and oriented manifold. Let us assume for the moment that
there exists a fundamental class $[\NNN]\in H_{\top}(\NNN)$.
In this case we define, for cohomology classes
$$\alpha_1,\dots,\alpha_p,\beta_1,\dots,\beta_q
\in H^*_{S^1}(F)\quad\mbox{ and }\quad\gamma\in H^*(\AAA/\GGG_0),$$
the {\bf Hamiltonian Gromov-Witten} invariant $\Phi_{B,c}^{X,F}$ to be 
\begin{align*}
\Phi_{B,c}^{X,F}&(\alpha_1,\dots,\alpha_p\mid\beta_1,\dots,\beta_q\mid\gamma)
\notag\\
&=\la \mu_2(\alpha_1)\cup\dots\cup\mu_2(\alpha_p)\cup
\mu_0(\beta_1)\cup\dots\cup\mu_0(\beta_q)\cup\nu(\gamma),[\NNN]\ra.
\end{align*}
                                                         
Assuming that the cobordisms given by theorem \ref{gensigmasmoothN}
also support fundamental classes, one could prove that the invariant
$\Phi$ only depends on the positive deformation class
of the symplectic form of $F$, of the action of $S^1$, and of
the connected component of $\imag\RR\setminus C_0$ on which $c$
lies. Hence, it is independent of the perturbation $\sigma$
and of the invariant and compatible complex structure $I$.

If the extended moduli space $\NNN$ is compact and so are the cobordisms
between the moduli arising from different perturbations and complex
structure, then the fundamental classes needed above do exist automatically.
However, most of the time the moduli space $\NNN$
will be noncompact, and we shall need the results in chapter \ref{compact}.
Furthermore, the compactification of $\NNN$ obtained by adding 
cusp THCs given by theorem \ref{compactificacio} does not have a priori
so good a structure to necessarily carry a fundamental class. It will be a 
stratified space with strata admitting (for generic complex structure)
a ramified covering by smooth manifolds, and the only thing we will be
able to prove is that the extra strata in the compactification
will have codimension at least 2
when $F$ is a positive manifold and a certain conditions
are satisfied by the fixed point data. In this situation we will be
able to define the Hamiltonian Gromov-Witten invariants.

There are some cases where one can be sure that there
is a fundamental class $[\NNN]$. For example, if there is a complex
structure $I\in\III_{\omega,S^1}$ on $F$ for which there are no rational 
curves, then for any complex structure $I\in\III^{\reg}_{\omega,S^1}$ 
near $I$ the moduli $\NNN_I$ is automatically compact. Indeed, by
theorem \ref{compactificacio} the only cause of noncompactness
is the appearance of bubbles, which are rational curves.
Furthermore, the energy the bubbles is bounded above by the
Yang--Mills--Higgs functional. Finally, the property of having no rational
curve of bounded energy is open in $\III$ thanks to Gromov compactness.
(See remark \ref{onconditions}.)          
Hence in this situation the above definition of the mixed
Hamiltonian Gromow-Witten invariants makes perfect sense.

\subsection{The invariant $\oPhi$}

To define this other invariant we will used the moduli space $\MMM$.
Following the notation above, let $\BBB=(\AAA\times\SSS^*)/\GGG$.
The gauge group $\GGG$ acts freely on 
$\AAA\times\SSS^*\times E\to\AAA\times\SSS^*\times X$.
Let $\EE\to\BBB\times X$ be the quotient and
$\FE=\EE\times_{S^1} F=(\AAA\times\SSS^*\times \FFF)/\GGG$.

Let $\oev:\AAA\times\SSS^*\times X\to\AAA\times\SSS^*\times\FFF$ 
be the universal section, defined as $\oev(A,\Phi,x)=(A,\Phi,\Phi(x))$.
The section $\oev$ is $\GGG$ equivariant, so it descends to give
a section $$\ev:\BBB\times X\to\FE.$$
Just as before we have a map $\rho_{\EE}^*:H_{S^1}^*(F)\to
H^*(\FE)$ and we define $\omu_i:H^*_{S^1}(F)\to H^*(\BBB)$, $i=0,2$, as
$$\omu_2(\delta):=\ev^*\rho_{\EE}^*\delta/[X]\qquad\mbox{ and }
\qquad \omu_0(\delta):=\ev^*\rho_{\EE}^*\delta/[\pt],$$
where $\delta\in H^*_{S^1}(F)$ and $[X]\in H_2(X)$ (resp. $[\pt]\in H_0(X)$)
is the fundamental class of $X$ (resp. the class of a point).
Finally, let $$\onu:H^*(\AAA/\GGG)=H^*(\AAA/\GGG_0)\to H^*(\BBB)$$
be the map induced by the projection $\BBB\to\AAA/\GGG$.

Let us take $B\in H_2(F_{S^1})$ such that ${\pi_F}_*B=\eta(E)$,
$c\in\imag\RR\setminus C_0$, $I\in\III_{\omega,S^1}$ and
$\sigma\in\Sigma_c^{\reg}(E)$. Then
$$\MMM=\MMM^{F,S^1}_{I,\sigma}(B,c)\subset\BBB$$ 
is a smooth oriented manifold. Let us suppose as before that there is
a fundamental class $[\MMM]\in H_{\top}(\MMM)$.
In this situation we define, for cohomology classes
$$\alpha_1,\dots,\alpha_p,\beta_1,\dots,\beta_q
\in H^*_{S^1}(F)\quad\mbox{ and }\quad\gamma\in H^*(\AAA/\GGG),$$
the Hamiltonian Gromov-Witten invariant $\oPhi_{B,c}^{X,F}$ to be
\begin{align*}
\oPhi_{B,c}^{X,F}&(\alpha_1,\dots,\alpha_p\mid\beta_1,\dots,\beta_q\mid\gamma)
\notag\\
&=\la \omu_2(\alpha_1)\cup\dots\cup\omu_2(\alpha_p)\cup
\omu_0(\beta_1)\cup\dots\cup\omu_0(\beta_q)\cup\onu(\gamma),[\MMM]\ra.
\end{align*}
As before, if the cobordisms between the moduli spaces $\MMM$ obtained
from different choices of $I$ and $\sigma$ have fundamental classes,
then the numbers $\oPhi$ are invariant. Here the same comments
as in the preceeding subsection are valid: our assumptions on 
existence of fundamental classes will be satisfied when $\MMM$
is compact, but when we need to compactify $\MMM$, then
we will have to study the added strata.
   
\section{Conditions}
\label{positivity}
In some part of this work we have made several assumptions on $F$
and the action of $S^1$. In order to give a perfect sense to our
invariants we will need to add still some extra conditions on the fixed
point set and on the complex structure of $F$. 
We state here all these conditions. The first one has been assumed
several times before, and the next three are new.

\begin{condition}
$F$ is a compact symplectic manifold. The Hamiltonian action of $S^1$ on $F$
is almost free. As a consequence, the action on the complementary
$F\setminus F^{S^1}$ of the fixed point set is free.
\label{cond0}
\end{condition}
This condition is equivalent to all Marsden-Weinstein quotients
at regular values of the moment map $\mu$ being smooth. That is,
if we do not assume that the action is almost free then there
are quotients which are orbifolds.

There are two cases in which we can define the invariants.

\subsubsection{First case}
In the first situation there is simply no bubbling, and this happens
if the following condition is satisfied.
\begin{condition}
The symplectic form $\omega$ vanishes on the set of spherical homology
classes in $H_2(F;\ZZ)$ (these are the classes which belong to the image
of the map $\pi_2(F)\to H_2(F;\ZZ)$).
\end{condition}
When there is no bubbling the extended moduli spaces is automatically
compact and hence there is fundamental class in the homology of the moduli
spaces $\NNN$, $\MMM$. Therefore, we can define the invariants as in the
introduction.

\subsubsection{Second case}
If we allow bubbling to occur, we need to control the strata added
to the moduli spaces $\NNN$ and $\MMM$ in theorem 
\ref{compactificacio}. The following three conditions will be assumed in the 
sequel, and they will be used when proving that these strata have codimension
greater than two.

\begin{condition}
$(F,I)$ has to be {\bf positive}. This means that for any
$I$-holomorphic map $s:\CP^1\to F$ the pullback $s^*TF$
has positive degree                                                          
$$\deg(s^*TF)>0.$$
\label{cond1}
\end{condition}

This condition is fulfilled for example when $F$ is a Fano manifold,
that is, when the anticanonical bundle $\Lambda^{\top}TF$ is
ample. This means that there exists and integer $k>0$ such that
$(\Lambda^{\top}TF)^k$ is very ample. In this case, for any $I$-holomorphic
map $s:\CP^1\to F$ and different points $x,y\in\CP^1$ there is a section
$\sigma\in H^0((\Lambda^{\top}TF)^k)$ which vanishes on $s(x)$
but not on $s(y)$. So the restriction $s^*\sigma$
is a nonzero holomorphic section of $s^*(\Lambda^{\top}TF)^k$.
Hence, this bundle has positive degree, and consequently the same
happens to $\Lambda^{\top}TF$. As examples of Fano manifolds we have
the projective spaces $\CP^n$ and any hypersurface of $\CP^n$ of degree
$<n$.

\begin{condition}
For any connected component $F_0$ of the fixed point set $F^{S^1}$
and any $I$-holomorphic $s:\CP^1\to F_0$,
$$\deg(s^*TF_0)\leq \deg(s^*TF).$$
\label{cond2}
\end{condition}

\begin{condition}
All the connected components $F_0$ of the fixed point set
$F^{S^1}$ have complex codimension at most 3.
\label{cond3}
\end{condition}

\begin{remark}
\label{onconditions}
Due to our need of taking generic complex structures $I$
to get smooth strata in the compactification of $\NNN_I$, it is 
desirable to know whether conditions \ref{cond1} and \ref{cond2}
are deformation invariant with respect to the complex structure on $F$.
The answer is no. However, something weaker but still enough
for our purposes is true. Suppose we restrict our atention in 
conditions \ref{cond1} and \ref{cond2} to maps $s:\CP^1\to F$
of bounded energy $K>0$. Then the resulting conditions are open
in $\III_{\omega,S^1}$ (this follows from Gromov compactness theorem,
see lemma 5.1.2 in \cite{McDS1}). On the other hand, in the compactifications
of $\NNN_{I,\sigma}(B,c)$ and $\MMM_{I,\sigma}(B,c)$ the energy of bubbles is 
bounded above by the value of $\YMH_c$ on any pair 
$(A,\Phi)\in\bM{I,\sigma}(B,c)$ (which only depends on
$B$ and on $|\sigma|_{C^0}$). So if the conditions are fulfilled by
a complex structure $I_0\in\III_{\omega,S^1}$ then the invariants
obtained from $\NNN_{I,\sigma}(B,c)$ and $\MMM_{I,\sigma}(B,c)$
will be well defined for any complex structure $I$ near enough $I_0$ 
(how near it must be depends on $B\in H_2(F_{S^1})$).

In particular, if $F$ is a Fano manifold, then condition \ref{cond1}
can be granted. 
\end{remark}

\begin{example}
Consider the action of $S^1$ on $\CP^4$ such that an element
$\lambda\in S^1$ maps $[x_0:x_1:x_2:x_3:x_4]$ to
$[\lambda x_0:\lambda x_1:x_2:x_3:x_4]$.
The connected components of the fixed point locus are the following
subspaces
\begin{align*}
F_1 &= \{x_2=x_3=x_4=0\}, \\
F_2 &= \{x_0=x_1=0\}. 
\end{align*}
This action is almost free, the projective space $\CP^4$ is a Fano
manifold and condition \ref{cond2} is easily seen to be satisfied
(considering the standard complex structure on $\CP^4$). Finally, the
codimension condition \ref{cond3} also holds. Hence the invariants
are well defined in this case.
\end{example}

Of course, if we consider the product of $\CP^4$ by any compact positive
symplectic manifold $M$ and take the diagonal action of $S^1$ (with the trivial
action on $M$) then we get a symplectic manifold with a Hamiltonian action
of $S^1$ which also satisfies the conditions. 
 
\section{Moduli of cusp $\sigma$-THCs}
    
\subsection{Evaluation maps are submersions}
\label{submersions}
In this subsection we will generalise the result in section 6.1
of \cite{McDS1} for curves in $\MMM^{L,\Gamma,\rho}(A)$, where
$A\in H_2(F;\ZZ)$. For any $x\in\CP^1$ we have an evaluation map
$$\ev_x:\MMM_{\III}=\MMM^{L,\Gamma,\rho}_{\III}(A)\to F$$
which sends any $s\in\MMM_{\III}$ to $\ev_x(s)=s(x)$. When
$\Gamma=1$ theorem 6.1.1 in \cite{McDS1} says that the map $\ev_x$
is a submersion. When $\Gamma\neq 1$ this need not hold any longer.
In fact, we must distinguish two possibilities. If $x\neq x_\pm$,
then the map $\ev_x:\MMM_{\III}\to F$ is a submersion, and if
$x=x_\pm$ then the evaluation map $\ev_x$ takes values in $F^{S^1}$
and the map $\ev_x:\MMM_{\III}\to F^{S^1}$ is a submersion.
We state this in the following lemma.

\begin{lemma}
Suppose that $\Gamma\neq S^1$. Given $I\in\III_{\omega,S^1}$, a curve
$s\in\MMM^{L,\Gamma,\rho}_I(A)$ and a point $x\in\CP^1$
different from $x_\pm$ (resp. equal to $x_\pm$) ---note that the
second condition only makes sense when $\Gamma\neq 1$---
there exists $\delta>0$ such that for any $v\in T_{s(x)}F^L$
(resp. for any $v\in T_{s(x)}F^{S^1}$) and every 
$0<\rho<r<\delta$ there exists a smooth $\Gamma$-equivariant
vector field $\xi\in\Omega^0(s^*TF^L)^{\Gamma}$ and an infinitesimal
variation of almost complex structure $Y\in T_I\III_{\omega,S^1}$
(see proposition \ref{udoscxstr}) such that the following
holds

i) $D_s\xi+\frac{1}{2}Y(s)\circ ds\circ j=0$ (that is, the pair
$(\xi,Y)$ belongs to $T_{(s,I)}\MMM_{\III}$),

ii) $\xi(x)=v$ and

iii) $\xi$ is supported in $\Gamma\circ B_\delta(x)$ and $Y$ is supported
in and arbitrarily small neighbourhood of 
$s(\Gamma\cdot(B_r(x)\setminus B_{\rho}(x))$.
\end{lemma}
\begin{pf}
Since the proof is almost the same as that of lemma 6.1.2 in
\cite{McDS1}, we will just give a sketch and mention the differences.
The first thing to do is to find a local solution $\xi_0$ of
$D_s\xi=0$ in $B_{\delta}(x)$ satisfying $\xi_0(x)=v$. 
This is done by solving a boundary value problem (see proposition
4.1 in \cite{McD} and the references therein). Then one multiplies
$\xi$ by a cutoff function with support in a neighbourhood of
$B_r(x)\setminus B_{\rho}(x)$ to extend $\xi_0$ to a section
of $s^*TF^L$. One then averages $\xi_0$ by the action of $\Gamma$
and obtains a section $\xi\in\Omega^0(s^*TF^L)^\Gamma$.
Finally, one must modify $I$ by a suitable infinitesimal
$Y\in T_I\III_{\omega,S^1}$ so that {\it i)} is satisfied
(in order to take $Y$ $\Gamma$-equivariant one needs to be careful with
the fixed point locus of the action of $S^1$; this may be done 
using theorem \ref{stabilL}, as was done in the proof
of proposition \ref{udoscxstr}). This $Y$ can be taken
fulfilling property {\it iii)}, repeating the argument in 
\cite{McDS1} but taking into account $\Gamma$-equivariance.
\end{pf}

When $\Gamma=S^1$ the result is even easier. As before, we distinguish
two possibilities. If $x=x_\pm$ then the same result as above holds.
When $x=x_\pm$, then it can be improved in the following sense.
The perturbation $Y$ of the complex structure may be chosen with
support in a neighbourhood of $S^1\cdot x_a\cap S^1\cdot x_b$,
where $d(x_+,S^1\cdot x_a)<d(x_+,S^1\cdot x)<d(x_+,S^1\cdot x_b)$.
This is a straightforward consequence of the interpretation of 
$\MMM_I^{1,S^1,1}$ in terms of lines of steepest descent of the moment
map with respect to the metric $\omega(\cdot,I\cdot)$
(see subsection \ref{comentaris}).

\begin{definition}
We will say that a point $x\in\CP^1$ is {\bf critical} with respect to the
tuple $(L,\Gamma,\rho)$ if either $L=S^1$ or $\Gamma\neq 1$
and $x=x_\pm$.
\label{defcritic}
\end{definition}

\subsection{Framings}

\subsubsectionr{}
Let $q\in\NN$ and fix points $y_1,\dots,y_q\in X$. We will call these
points the {\bf marked points}. Recall that we denote cusps THCs with 
tuples of the form $(E,X^{\cusp},A,\Phi,c)$ (see definition \ref{defcuspTHC}).
For any bubble $X_k\subset X^{\cusp}$ the image $\Phi_k(X_k)\subset\FFF$
is contained in a unique fibre. Let $x_k\in X$ be the base point 
corresponding to this fibre.

The gauge groups $\GGG_0$ and $\GGG$ act on the set of cusp curves.
On the other hand, any bubble $X_k$ has a reparameterisation group $A_k$.
This is the subset of $\Aut(X_k)=\PSL(2;\CC)$ given by automorphisms with 
keep fixed the intersections of $X_k$ with the other irreducible components
of $X^{\cusp}$. We call {\bf moduli of cusp curves} (resp. {\bf extended 
moduli of cusp curves}) the set of orbits of the action of $\GGG$ (resp. 
$\GGG_0$) and the reparameterisation groups $A_k$ of the bubbles.

\subsubsectionr{}
The maps $\Phi_k$ in a cusp curve may be multicovered. This means that
$\Phi_k$ factors as $\Phi_k'\circ r_k$, where $r_k:\CP^1\to\CP^1$
is a ramified covering. If $r_k$ has maximal degree, we will call
$\Phi_k'$ the {\bf simplification} of $\Phi$. This will be
a simple map.

Let us write $F_1,\dots,F_r$ for the connected components of the fixed
point set $F^{S^1}$. We will follow this notation. For any natural number
$K$, $\KK$ will denote the set $\{1,\dots,K\}\subset\NN$ and 
$\KK_0=\{0\}\cup\KK$. 
We define the framing $D=D(E,X^{\cusp},A,\Phi,c)$ of the cusp 
$\sigma$-THC $(E,X^{\cusp},A,\Phi,c)$ to be the following set of 
data: 

\begin{enumerate}
\item The class $B_0={\rho_E}_*{\Phi_0}_*[X_0]\in H_2(F_{S^1};\ZZ)$.
\item The number $K$ of bubbles in $X^{\cusp}$. 
\item Homology classes $B_1,\dots,B_K\in H_2(\FFF;\ZZ)$ describing the
image of the bubbles $X_k$ in $\FFF$.
\item For every $k\in\KK$ the tuple $(L_k,\Gamma_k,\rho_k)$ such that,       
after identifying $\FFF_{x_k}\simeq F$,
the simplification $\Phi'_k$ of $\Phi_k$ belongs to 
$\MMM^{L_k,\Gamma_k,\rho_k}(B_k)$. 
\item For any $k\in\KK$ such that $L_k=S^1$, the number $c(k)$ such that, 
for all $s\in\MMM^{L_k,\Gamma_k,\rho_k}(B_k)$, $s(\CP^1)\subset F_{c(k)}$.
\item For any $k\in\KK$ such that $\Gamma_k\neq 1$, numbers
$c(k)_+$ and $c(k)_-$ such that, for all
$s\in\MMM^{L_k,\Gamma_k,\rho_k}(B_k)$, $s(x_\pm)\in F_{c(k)_\pm}$.
\item A set $C\subset\KK_0^2$ containing the pairs $(i,j)$
such that $i<j$ and $X_i\cap X_j\neq\emptyset$ and the following.
\begin{enumerate}
\item A partition $C=C_{00}\cup C_{01}\cup C_{10}\cup C_{11}$ defined
as follows. For any pair $(i,j)\in C$, let $x=X_i\cap X_j$.
Put $\epsilon(i)$ to be $1$ if $x\in X_i$ is a critical point with respect to
$(L_i,\Gamma_i,\rho_i)$ and $0$ otherwise (see definition \ref{defcritic}), 
and define $\epsilon(j)$ similarly. Then $(i,j)$ belongs to 
$C_{\epsilon(i)\epsilon(j)}$.
\item A set $C''\subset\KK^3$ containing the sets $(i,j,k)$
such that $\Gamma_i(X_i\cap X_j)\cap(X_i\cap X_k)\neq\emptyset$
(this is a subset of $X_i$).
\end{enumerate}
\item An element $c\in\imag\RR$.
\item A set $S\subset\KK^2$ containing the pairs $(i,j)$
such that $\Phi_i(X_i)=\Phi_j(X_j)$.
\item The set of $k\in\KK$ such that $x_k\in\{y_1,\dots,y_q\}$
(that is, such that $\Phi_k(X_k)$ is contained in the fibre over
a marked point) and the marked point $y_{m(k)}=x_k$.
\end{enumerate}

\begin{definition}
We will denote $\MMM_{I,\sigma}'(D)$ the moduli of cusp $\sigma$-THCs
with framing $D$. The extended moduli of cusp $\sigma$-THCs will be
denoted $\NNN_{I,\sigma}'(D)$.
\end{definition}

\subsubsectionr{}
\label{coho}
We define the {\bf total homology class of the frame $D$} 
to be the equivariant homology class $B(D)=B_0+B_1+\dots+B_K\in 
H_2(F_{S^1};\ZZ)$. 
Lemma \ref{triviallocsys} in the appendix says that there is 
a natural map $H_*(F)\to H_*(\FFF)$. In the sequel the cohomology
classes $B_k$ will mean either the element ${\Phi_k}_*[X_k]\in H_*(\FFF)$ or
any of their preimages in $H_*(F)$ by that map. 
For example, we may write for any $k$
\begin{equation}
\la c_1(TF),B_k\ra=\la c_1^{S^1}(TF),B_k\ra,
\label{equinonequi}
\end{equation}
where on the left hand side we view $B_k\in H_2(F)$ and on
right hand side we view $B_k$ as an equivariant homology class.

\subsubsectionr{}
For any cusp $\sigma$-THC with frame $D$ we make the following reduction
process. First we substitute the bubble maps $\Phi_k$ by their 
simplifications $\Phi_k'$, and then we identify bubbles with the same image 
in $\FFF$. Finally, if necessary we forget some intersection points
in order that no two irreducible components of the cusp we have obtained
intersect at more than one point.

After this process we end up with another cusp $\sigma$-THC with
frame $\ov{D}$. We call the resulting cusp a {\bf reduced cusp} and
$\ov{D}$ a {\bf reduced frame}. 
We will denote $\MMM_{I,\sigma}(\ov{D})\subset\MMM_{I,\sigma}'(\ov{D})$
(resp. $\NNN_{I,\sigma}(\ov{D})\subset\NNN_{I,\sigma}'(\ov{D})$)
the moduli (resp. extended moduli) of reduced cusp curves with
framing $\ov{D}$.

Note that the total homology
class of $\ov{D}$ will not necessarily be equal to that of $D$.
If $\ov{B}_0,\ov{B}_1,\dots,\ov{B}_K$ are the homology classes
of $\ov{D}$ we will have
$$B(D)=\ov{B}_0+r_1\ov{B}_1+\dots+r_K\ov{B}_K,$$
where $r_k\geq 1$ are integers. This motivates the following definition.

\begin{definition}
If $B\in H_2(F_{S^1};\ZZ)$
and the homology classes $B_0,B_1,\dots,B_K$ of a frame $D$ satisfy
$B=B_0+r_1B_1+\dots+r_KB_K$ for some integers $r_k\geq 1$, then we will
say that the frame $D$ is {\bf $B$-admisible}.
\end{definition}

We will denote $\DDD(B,c)$ the set of $B$-admisible framings $D$ such
that $c(D)=c$. This is obviously a numerable set. It contains a 
distinguished element $D^T$ which represents the cusp curve with no bubbling.
We will call $D^T$ the {\bf top framing} of $B$.

\subsubsectionr{}
The main result of this subsection is the following theorem.

\begin{theorem}
Let $D$ be a reduced frame.
Suppose that the element $c=c(D)\in\imag\RR$ lies in the complementary
of $C_0$. For generic $\sigma\in\Sigma_c^{\reg}(E)$ and complex structure
$I\in\III^{\reg}_{\omega,S^1}$ the moduli $\MMM_{I,\sigma}(D)$
of reduced cusp $\sigma$-THC with frame $D$ is branchedly covered by 
a smooth manifold. Furthermore, if the positivity conditions 
in \ref{positivity} are satisfied and if $D$ is $B$-admisible for some 
$B\in H_2(F_{S^1};\ZZ)$ then 
$$\dim_{\RR} \MMM_{I,\sigma}(D)\leq \dim_{\RR} \MMM_{I,\sigma}(B,c)-
2(1+b(D)),$$
where $b(D)$ is the number of marked points whose fibre contains
a bubble.
\label{gensigmasmoothcusp}
\end{theorem}
\begin{pf}
Through all the proof $\dim$ and $\codim$ will denote real dimension
and codimension. Let $2n$ be the dimension of $F$.

\subsubsectionr{}
For any tuple $(L,\Gamma,\rho)$ and bundle $\FFF=\FFF^E\to X$ denote 
$\Map^L(\CP^1,\FFF)^{\Gamma,\rho}_{\fibr}$ the set of maps 
$\Phi:\CP^1\to \FFF$ whose image is included in a single fibre
$\FFF_x$ and such that $\Phi\in\Map^L(\CP^1,\FFF_x)^{\Gamma,\rho}$
(see (\ref{mapequi})). We have
$$\Map^L(\CP^1,\FFF)^{\Gamma,\rho}_{\fibr}=
E\times_{S^1}\Map^L(\CP^1,F)^{\Gamma,\rho}.$$
Define for any homology class $B\in H_2(F;\ZZ)$
$$\MMM_{\III,\fibr}^{L,\Gamma,\rho}(B)
=\left\{(\Phi,I)\in\Map^L(\CP^1,F)_{\fibr}^{\Gamma,\rho}
\times \III_{\omega,S^1}^l\Big|
\begin{array}{l}\ov{\partial}_I\Phi=0,\ \Phi_*[\Sigma]=B\\
\Phi\mbox{ simple }\end{array}\right\}.$$
Incidentally, this is the moduli space used to defined
fibrewise and equivariant quantum cohomology by Givental, Kim
and Lu (see \cite{Lu}).
Note that we have 
$\MMM_{\III,\fibr}^{L,\Gamma,\rho}(B)
=E\times_{S^1}\MMM_{\III}^{L,\Gamma,\rho}(B).$
In particular, 
\begin{equation}
\dim\MMM_{\III,\fibr}^{L,\Gamma,\rho}(B)=
\dim\MMM_{\III}^{L,\Gamma,\rho}(B)+2.
\label{dimfibr}
\end{equation}

\subsubsectionr{}
Fix a reduced frame $D$ and
suppose that the element $c=c(D)\in\imag\RR$ lies in the complementary
of $C_0$. 
Let $K$ be the number of bubbles,
$C\subset \{0,1,\dots,K\}^2$ the set of pairs describing which 
irreducible components intersect, $B_0,B_1,\dots,B_K$ the homology
classes of $D$, and $(L_k,\Gamma_k,\rho_k)$ the tuples telling the 
moduli in which $\Phi_k$ sits. We denote $X=X_0,X_1,\dots,X_K$
the irreducible components of the cusps with frame $D$.

\subsubsectionr{}
Define for any $k$ the group $G_k$ to be $\CC^*$ if $\Gamma_k\neq 1$ and 
$\PSL(2;\CC)$ if $\Gamma_k=1$. 
These groups act effectively on $\MMM_I^{L_k,\Gamma_k,\rho_k}(B_k)$
as follows: any $s\in\MMM_I^{L_k,\Gamma_k,\rho_k}(B_k)$ is mapped
by $g\in G_k$ to $g(s):=s\circ g:\CP^1\to F$. Note that when
$\Gamma_k\neq 1$ the group $\CC^*$ acts on $\CP^1$ keeping fixed
$x_\pm$. 

\subsubsectionr{}
Let us write
$$\MMM_{\III}^*(D)=\prod_{k=1}^{K}
\MMM_{\III,\fibr}^{L_k,\Gamma_k,\rho_k}(B_k)\setminus\Delta,$$
where $\Delta$ is the multidiagonal, that is, the set of elements
$(s_1,\dots,s_K)$ such that $s_i=s_j$ for some $i\neq j$. 
$\MMM_{\III}^*(D)$ parameterizes
tuples of $K$ different holomorphic maps $\Phi_k:\CP^1\to \FFF$
whose image is contained in any fibre.
Reasoning exactly like in proposition \ref{udoscxstr} one proves
that $\MMM_{\III}^*(D)$ is a smooth manifold of finite dimension
on which the gauge group $\GGG$ acts smoothly. Note that 
$\dim_\RR\MMM_{\III}^*(D)/\GGG=\dim_\RR\MMM_{\III}^*(D)-1.$
We will prove that the universal moduli of reduced cusp $\sigma$-THCs
with frame $D$
$$\MMM_{\III,\Sigma}(D)=\bigcup_{I\in\III_{\omega,S^1},\sigma\in\Sigma_c(E)}
\MMM_{I,\sigma}(D)$$
is a smooth Banach submanifold of 
$\bM_{\III,\Sigma}(B_0,c)\times\MMM_\III^*(D)/(\GGG\times\prod G_k).$

\subsubsectionr{}
Recall that $\FFF_0=E\times_{S^1}F^{S^1}$.
For any pair $e=(i,j)\in C$ we define
$$\FFF_e=\left\{\begin{array}{ll}
\FFF\times\FFF & \mbox{ if $e\in C_{00}$,}\\
\FFF\times\FFF^0 & \mbox{ if $e\in C_{01}$,}\\
\FFF^0\times\FFF & \mbox{ if $e\in C_{10}$,}\\
\FFF^0\times\FFF^0 & \mbox{ if $e\in C_{11}$,}
\end{array}\right.$$
and we write $\Delta_e\subset\FFF_e$ for the diagonal in $\FFF_e$.
We also define $X_e=X_i^e\times X_j^e$, where:
if $i=0$ and $x_j=y_{m(j)}$ is a marked point then $X_i^e=y_{m(j)}$
and $X_j^e=x_{\pm}^j$ if $X_i\cap X_j=x_{\pm}^j$ and $X_j$ otherwise;
and if $i\neq 0$ then $X_i^e$ is $x_\pm^i\in X_i$ if $X_i\cap X_j=x_\pm^i$ 
and $X_i$ otherwise, and $X_j^e$ is defined similarly. 

We then have an evaluation map
$$\ev_C:\bM_{\III,\Sigma}(B_0,c)\times\MMM_\III^*(D)\times
\prod_{e\in C}X_e\to\prod_{e\in C}\FFF_e.$$
Consider the projection 
$$\Theta_\III:\bM_{\III,\Sigma}(B_0,c)\times\MMM_\III^*(D)\times
\prod_{e\in C}X_e\to\III_{\omega,S^1}^{K+1}$$
and write $\Delta_\III$ for the diagonal in $\III_{\omega,S^1}^{K+1}$.
Since $c\in\imag\RR\setminus C_0$, theorem \ref{gensigmasmooth}
and \ref{regularitat} imply that $\Theta_{\III}^{-1}(\Delta_{\III})$
is a smooth Banach manifold (this is exactly like the proof of lemma 4.9 
in \cite{RuTi}). Let us define
\begin{equation}
\RRR_{\III,\Sigma}(D)=
(\Theta_\III^{-1}(\Delta_\III)\cap\ev_C^{-1}(\prod_{e\in C}\Delta_e)/\GGG)
/\prod G_k.
\label{defR}
\end{equation}
Now, the results in subsection \ref{submersions} imply that 
the evaluation map $\ev_C$ restricted to $\Theta_\III^{-1}(\Delta_\III)$
is a submersion. Furthermore, the action of
$\GGG$ on $\bM_{\III,\Sigma}(B_0,c)\times\MMM_\III^*(D)\times\prod_{e\in C}X_e$ 
is free, so $\RRR_{\III,\Sigma}(D)$ is a smooth Banach manifold.
Finally, the projection 
$$p:\RRR_{\III,\Sigma}(D)\to\PPP=\III_{\omega,S^1}\times\Sigma_c(E)$$
is a Fredholm map. Hence, the theorem of Sard-Smale implies that
there is a Baire set of the second category $\PPP^{\reg}\subset\PPP$
such that for any $(I,\sigma)\in\PPP^{\reg}$ the set
$$p^{-1}(I,\sigma)=\RRR_{I,\sigma}(D)$$
is a smooth manifold. 
The set of cusp curves is equal to
$\MMM_{I,\sigma}(D)=\pi_\MMM(\RRR_{I,\sigma}(D)),$
where $$\pi_\MMM:\bM_{I,\sigma}(B_0,c)\times\MMM_I^*(D)\times
\prod_{e\in C}X_e\to\bM_{I,\sigma}(B_0,c)\times\MMM_I^*(D)$$
is the projection. This is a ramified covering, and the fibre
over any element in $\MMM_{\III,\Sigma}(D)$ is just the set of 
points of $X^{\cusp}$ whose images in $\FFF$ coincide.

\subsubsectionr{}
To prove the claim on the dimension of $\MMM_{I,\sigma}(D)$, we
will prove that $$\dim \RRR_{I,\sigma}(D)\leq\dim\MMM_{I,\sigma}(B,c)-2.$$
Consider the following commutative diagram

$$\xymatrix{\Theta_{\III}^{-1}(\Delta_I)\cap\ev_C^{-1}(\prod\Delta_e)
\ar[r]^{\tilde{\iota}} \ar[d]
& \bM_{I,\sigma}(B_0,c)\times\MMM_I^*(D)\times\prod X_e \ar[d] \\
\bR_{I,\sigma}(D)
:=\Theta_{\III}^{-1}(\Delta_I)\cap\ev_C^{-1}(\prod\Delta_e)/\GGG
\ar[r]^{\iota} 
& \bM_{I,\sigma}(B_0,c)\times\MMM_I^*(D)/\GGG\times\prod X_e}$$           
in which the horizontal maps $\tilde{\iota}$ and $\iota$ are
embeddings and the vertical maps are the projections to the sets
of $\GGG$ orbits.

Since $\Theta_{\III}^{-1}(\Delta_I)\cap\ev_C^{-1}(\prod\Delta_e)$
is $\GGG$ invariant, the codimension of the embeddings $\tilde{\iota}$
and $\iota$ is the same. Now, 
$\Theta_{\III}^{-1}(\Delta_I)=
\bM_{I,\sigma}(B_0,c)\times\MMM_I^*(D)\times\prod X_e$ and the map
$\ev_C$ restricted to $\Theta_{\III}^{-1}(\Delta_I)$ is a submersion. 
Consequently,
$$\codim {\tilde{\iota}}=
\codim_{\Theta_{\III}^{-1}(\Delta_I)}\ev_C^{-1}(\prod\Delta_e)
=\sum\codim_{\FFF_e}\Delta_e.$$
Combining (\ref{defR}) and the above reasoning we conclude that
\begin{align*}
\dim\RRR_{I,\sigma}(D) &=\dim\bR_{I,\sigma}(D)-\sum\dim G_k \\
&=\dim\MMM_{I,\sigma}(B_0,c)+\dim\MMM_I^*(D)/\GGG \\
&+\sum\dim X_e-\sum\codim_{\FFF_e}\Delta_e-\sum\dim G_k\\
&=\dim\MMM_{I,\sigma}(B_0,c)+\dim\MMM_I^*(D)-1 \\
&+\sum\dim X_e-\sum\codim_{\FFF_e}\Delta_e-\sum\dim G_k.
\end{align*}

To bound this dimension we divide the set of bubbles $\KK$ in 
three subsets. Let $\SS$ (resp. $\TT$ and $\UU$) denote the set
of $k\in\KK$ such that $L_k=1$, $\Gamma_k=1$ (resp.
$L_k=1$, $\Gamma_k\neq 1$ and $L_k=S^1$, $\Gamma_k=1$).
Let $S=|\SS|$, $T=|\TT|$ and $U=|\UU|$. Theorems \ref{dim1} and \ref{dim2}, 
and formula (\ref{equinonequi}), imply the following.
\begin{itemize}
\item If $k\in\SS$ then $\dim G_k=6$ and
$$\dim\MMM_{I,\fibr}^{L_k,\Gamma_k,\rho_k}(B_k)=
2+2\la c_1(TF),B_k\ra+2n=2+2\la c_1^{S^1}(TF),B_k\ra+2n,$$
where in the last term we view $B_k$ as an element of $H_*(F_{S^1})$
(see \ref{coho}).
\item If $k\in\TT$ then $\dim G_k=2$ and
$$\dim\MMM_{I,\fibr}^{L_k,\Gamma_k,\rho_k}(B_k)\leq
2+2\la c_1(TF),B_k\ra+2n-4=2\la c_1^{S^1}(TF),B_k\ra+2n-2,$$
by theorem \ref{claudim}.
\item If $k\in\UU$ then $\dim G_k=6$ and
\begin{align*}
\dim\MMM_{I,\fibr}^{L_k,\Gamma_k,\rho_k}(B_k) 
&=2+2\la c_1(TF^{S^1}),B_k\ra+\dim F_{c(k)} \\
&\leq2+2\la c_1(TF),B_k\ra+\dim F_{c(k)} \\
&=2+2\la c_1^{S^1}(TF),B_k\ra+\dim F_{c(k)},
\end{align*}
by condition \ref{cond2}.
\end{itemize}
On the other hand, since $D$ is $B$-admissible and condition
\ref{cond1} is satisfied, we have 
$$\sum_{k=1}^K\la c_1^{S^1}(TF),B_k\ra\leq\la c_1^{S^1}(TF),B\ra.$$
Hence, 
\begin{align*}
\dim\RRR_{I,\sigma}(D) &\leq 2\la c_1^{S^1}(TF),B\ra+2(n-1)(1-g) \\
&+(S+T)2n+\sum_{k\in\UU}\dim F_{c(k)}-4K+\sum_{e\in C}
(\dim X_e-\codim_{\FFF_e}\Delta_e).
\end{align*}
To find an upper bound for the last two terms we proceed as follows.
Suppose to begin with that $b(D)=0$.
Since for any $e$ we have $\dim X_e-\codim_{\FFF_e}\Delta_e\leq 0$, an
upper bound for $\sum_{e\in C'}$ where $C'\subset C$ will also give a
bound on $\dim\RRR_{I,\sigma}(D)$. So we take any subset $C'\subset C$ of $K$
elements with the following property. The graph whose vertices are the 
elements of $\KK_0$ and which has an edge joining $i$ to $j$ if either 
$(i,j)$ or $(j,i)$ belong to $C'$ is connected. This implies that 
$C'\nsubseteq C_{11}$ (because otherwise the vertex $0\in\KK_0$ would be 
disconnected from the rest). Take an injective map
$$v:C'\to\KK$$
which assigns to $(i,j)$ either $i$ or $j$. Let $k\in\KK$ and $e=v^{-1}(k)$.
\begin{itemize}
\item If $k\in\SS$ then $\dim X_e-\codim_{\FFF_e}\Delta_e=-2n+2$.
\item If $k\in\TT$ and $e\notin C_{11}$ then
$\dim X_e-\codim_{\FFF_e}\Delta_e\leq-2n+2$ and if $e\in C_{11}$ then 
$$\dim X_e-\codim_{\FFF_e}\Delta_e\leq
\max\{-\dim F_{c(k)_+},-\dim F_{c(k)_+}\}-2\leq -2n+4,$$
by condition \ref{cond3}.
\item If $k\in\UU$ then 
$\dim X_e-\codim_{\FFF_e}\Delta_e=-\dim F_{c(k)}+2$.
\end{itemize}
Since $C'\nsubseteq C_{11}$, we get
$$\dim\RRR_{I,\sigma}(D)\leq 2\la c_1^{S^1}(TF),B\ra+2(n-1)(1-g)-2=
\dim\MMM_{I,\sigma}(B,c)-2,$$
which is what we wanted to prove. If $b(D)>0$ then the same reasoning
as above works. Just observe that we have to substract at the 
end $2b(D)$ (that is, two units for each $e=(0,j)$ such that
$x_j=y_{m(j)}$).
\end{pf}

\section{Extended moduli of cusp $\sigma$-THCs}
Using the same techniques as in the proof of theorem \ref{gensigmasmoothcusp}
one proves the following theorem on the extended moduli of reduced cusp 
$\sigma$-THCs.

\begin{theorem}
Let $D$ be a reduced frame.
Suppose that the element $c=c(D)\in\imag\RR$ lies in the complementary
of $C_0$. For $\sigma\in\Sigma^{\reg}_c(E)$ and complex structure
$I\in\III^{\reg}_{\omega,S^1}$ the extended moduli $\NNN_{I,\sigma}(D)$
of reduced cusp $\sigma$-THC with frame $D$ is branchedly covered by
a smooth manifold. If the positivity conditions
in \ref{positivity} are satisfied and if $D$ is $B$-admissible for some
$B\in H_2(F_{S^1};\ZZ)$ then 
$$\dim_{\RR} \NNN_{I,\sigma}(D)\leq \dim_{\RR} \NNN_{I,\sigma}(B,c)-
2(1+b(D)),$$
where $b(D)$ is the number of marked points whose fibre contains
the image of some bubble. Furthermore, for any pair 
$I_0,I_1\in \III^{\reg}_{\omega,S^1}$
and $\sigma_0,\sigma_1\in\Sigma^{\reg}_c(E)$ we can find paths
$[0,1]\ni t\mapsto (I_t,\sigma_t)\in\III_{\omega,S^1}\times\Sigma_c(E)$
such that 
$$\bigcup_{t\in [0,1]}\NNN_{I_t,\sigma_t}(D)$$
is a smooth cobordism between $\NNN_{I_0,\sigma_0}(D)$
and $\NNN_{I_1,\sigma_1}(D)$.
\label{gensigmasmoothcuspN}
\end{theorem} 

\section{Definition of the invariants}
 
\subsection{A retraction of $\BBB$ and $\FE$}
Let $E\to X$ be a $S^1$ principal bundle of degree $d$, and let $\AAA$ be the 
set of connections on $E$. Let $\GGG_0\subset\GGG=\Map(X,S^1)$
be the set of gauge transformations of $E$ fixing the fibre
over $x_0\in X$.

Let $$\Jac_d(X)=\{A\in\AAA\mid F_A=-\imag 2\pi d\omega_X/\Vol(X)\}/\GGG_0$$
be the Jacobian of degree $d$
(here as usual we write $\omega_X$ for the symplectic form in $X$). 
This is a torus of real dimension twice the genus of $X$.
We will construct a retraction $\AAA/\GGG_0\to\Jac_d(X)$.
Recall that we have a metric on $X$, which induces metrics on the exterior 
algebra of forms $\Omega^*(X)$. Let $\HHH^j$ be the space of harmonic 
$j$-forms with respect to this metric.

Let $$F:\AAA/\GGG_0\to\Omega^2(\imag\RR)$$
be the map which sends any $[A]$ to
$ F_A+\imag 2\pi d\omega_X/\Vol(X)$.
It is easy to see, using Hodge theory, that the image of $F$ is the 
orthogonal of $\imag\HHH^2(X)$ in $\Omega^2(\imag\RR)$.
The preimage of $0\in\Omega^2(\imag\RR)$ is precisely $\Jac_d(X)$. 
In fact, $F:\AAA/\GGG_0\to\imag\HHH^2(X)^{\bot}$ is a smooth fibration
with fibres diffeomorphic to $\Jac_d(X)$. We will construct a connection
on this fibration by specifying its horizontal distribuition.

Given any $[A]\in\AAA/\GGG_0$, the tangent space $T_{[A]}\AAA/\GGG_0$
can be canonically identified with $\Ker d_1^*$, where
$d_1:\Omega^0(\imag\RR)\to\Omega^1(\imag\RR)$ is the exterior derivation.
Then we set the horizontal space at $[A]$ to be
$$(T_{[A]}\AAA/\GGG_0)_h:=\Ker d_1^*\cap (\Ker d_2)^{\bot},$$
where $d_2:\Omega^1(\imag\RR)\to\Omega^2(\imag\RR)$
is the exterior derivation.
Now, using parallel transport along lines going through 
$-\imag 2\pi d\omega_X/\Vol(X)\in\Omega^2(\imag\RR)$ we get
the desired retraction
$$R:\AAA/\GGG_0\to\Jac_d(X).$$

Let $\GGG^{\CC}=\Map(X,\CC^*)$. This is the complexification of $\GGG$
and it acts on $\AAA$ holomorphically (see \ref{setting}). 
Each $\GGG^{\CC}$ orbit in $\AAA$ contains exactly one $\GGG$ orbit
giving an element in $\Jac_d(X)$ (this is the most simple case of
Hitchin-Kobayashi correspondence).
\begin{lemma} The map $R$ factors through the projection
$\AAA/\GGG_0\to\AAA/\GGG^{\CC}$.
\end{lemma}
\begin{pf}
Suppose for simplicity that $d=0$. The vector field we have used to construct
the retraction is the gradient of the Yang--Mills functional 
$YM(A)=\|F_A\|_{L^2}^2$. Now, $\mu(A)=F_A$ is a moment map for the
action of $\GGG$ on $\AAA$, and so $YM(A)=\|\mu(A)\|_{L^2}^2$. But
it is a general fact that $\nabla\|\mu\|^2=-2\imag\fX_{\mu}$, where
we identify $\Lie(\GGG)\simeq\Lie(\GGG)^*$ to view $\mu\in\Lie(\GGG)$
and where $\fX$ is the field on $\AAA$ generated by $\mu$. Hence
the integral lines of the gradient of $\|F_A\|_{L^2}^2$ are contained
in the orbits of the action of $\GGG^{\CC}$.
\end{pf}

We now show how to lift the above retraction to a retraction of $\FEj$ onto
its restriction on $\Jac_d(X)\times X$. Let $\FEJj=\FEj|_{\Jac_d(X)\times X}$.
The universal connection $\bAA$ on $\pi_X^*E\to\AAA\times X$ (which is
by definition equal to $\pi_X^*A$ on the slice $\{A\}\times X$)  
descends to a connection $\AA$ on $\EE_J$. Using the induced connection on 
$\FEj$ and the connection on $\AAA/\GGG_0\to\imag\HHH^2(X)^{\bot}$
we get a connection on the fibration $\FEj\to\imag\HHH^2(X)^{\bot}$
which allows to define the retraction
$$R^{\FFF}:\FEj\to\FEJj.$$

\subsection{Pseudocycles in smooth manifolds}
The following definitions (with some modifications) are taken from 
chapter 7 in \cite{McDS1}. Let $M$ be any oriented smooth manifold. 

\begin{definition}
Given a manifold $V$ and a continuous map $f:V\to M$, the 
{\bf omega-limit-set} of $f$ is
$$\Omega_f=\bigcap_{K\subset V\\K\text{ compact}}\ov{f(V\setminus K)}.$$
\end{definition}

The omega-limit-set of $f$ is, in a certain sense, the boundary of
$f(V)$. More precisely, it consists of the limit points of sequences
$f(x_n)$, where $\{x_n\}\subset V$ has no convergent subsequence.

\begin{definition}
A {\bf pseudocycle} of real dimension $k$ in $M$ is a pair of maps 
$(f,g):(V,W)\to M$, where $V$ and $W$ are $\sigma$-compact
\footnote{Recall that this means that they can be covered by countably many 
compact sets.} oriented smooth manifolds, with $V$ of real dimension $k$ and
with all the components in $W$ having real dimension
at most $k-2$, such that $\Omega_f\subset g(W).$
\end{definition}

\begin{definition}
Two pseudocycles $(f_i,g_i):(V_i,W_i)\to M$, $i=1,2$, 
are said to be {\bf bordant} if there exists
a pseudocycle $(f_{\VV},g_{\WW}):(\VV,\WW)\to M$ such that
$\partial\VV=V_1-V_2$ and $f_{\VV}|_{V_i}=f_i$ for $i=1,2$.
\end{definition}

In particular, two pseudocycles $(f_i,g_i):(V_i,W_i)\to M$, $i=1,2$, 
such that $V_1=V_2$ and $f_1=f_2$ are trivially bordant. Note that,
since what we are really interested in is bordism classes of pseudocycles,
we could have defined pseudocycles, following \cite{McDS1}, as maps
$f:V\to M$ from an oriented manifold $V$ whose omega-limit-set may
be covered by the image of a map from a manifold of real dimension 
at most that of $V$ minus two to $M$. 

\begin{definition}
Two pseudocycles $(f_i,g_i):(V_i,W_i)\to M$, $i=1,2$, 
are said to be {\bf transverse}
if all the intersections $f_{V_1}(V_1)\cap f_{V_2}(V_2)$,
$f_{V_1}(V_1)\cap f_{W_2}(W_2)$, $f_{W_1}(W_1)\cap f_{V_2}(V_2)$
and $f_{W_1}(W_1)\cap f_{W_2}(W_2)$ are transverse in $M$.
\end{definition}

\begin{lemma}
If two pseudocycles $(f_i,g_i):(V_i,W_i)\to M$, $i=1,2$, are transverse
and of complementary dimension in $M$, then 
$\Omega_{f_i}\cap \ov{g_j(W_j)}=\emptyset$ for $\{i,j\}=\{1,2\}$ and the set
$\{(v_1,v_2)\in V_1\times V_2\mid f_1(v_1)=f_2(v_2)\}$
is finite. Let us define
$$f_1\cdot f_2=\sum_{(v_1,v_2)\in V_1\times V_2\\ f_1(v_1)=f_2(v_2)}
\nu(v_1,v_2),$$
where $\nu(v_1,v_2)$ is the intersection number of $f_1(V_1)$
and $f_2(V_2)$ at the point $f_1(v_1)=f_2(v_2)$ (since the intersection
is transverse, this number is $\pm 1$, depending on the orientations).
Then the number $f_1\cdot f_2$ only depends on the cobordism class
of $f_1$ and $f_2$.
\end{lemma}
\begin{pf}
This is lemma 7.1.3 in \cite{McDS1}.
\end{pf}

\begin{lemma}
Given two pseudocycles $(f_i,g_i):(V_i,W_i)\to M$, $i=1,2$, there is a subset
$\Diff(M,V_1,V_2,W_1,W_2)^{\reg}\subset\Diff(M)$ of Baire
of the second category such that the pseudocycles
$(\phi\circ f_{V_1},\phi\circ f_{W_1})$ and $(f_{V_2},f_{W_2})$
are transverse for any $\phi\in \Diff(M,V_1,V_2,W_1,W_2)^{\reg}$.
\end{lemma}
\begin{pf}
Given two $\sigma$-compact manifolds $V,W$ and two smooth maps 
$f,g:V,W\to M$, it is well known that there is a subset
$\Diff(M,V,W)^{\reg}\subset\Diff(M)$ of Baire
of the second category such that the maps $\phi\circ f$ and $g$
are transverse for any $\phi\in \Diff(M,V,W)^{\reg}$. Then we set
$$\Diff(M,V_1,V_2,W_1,W_2)^{\reg}=\bigcap_{i,j}\Diff(M,V_i,W_j)^{\reg},$$
which is of course of the second category.
\end{pf}

\begin{lemma}
Let $\dim_\RR M=m$, and let $(f,g):(V,W)\to M$ be an $(m-d)$ dimensional
pseudocycle. Any homology class $\beta\in H_d(M;\ZZ)$ can be represented
by a pseudocycle $(f_{\beta},g_{\beta}):(V_{\beta},W_{\beta})\to M$
in the sense that $V_{\beta}$ carries a fundamental class $[V_{\beta}]$
of dimension $d$ and that $(f_{\beta})_*[V_{\beta}]=\beta$. 
Furthermore, the map
$$\Psi_f:H_d(M;\ZZ)\to\ZZ$$
which sends $\beta$ to $f\cdot f_{\beta}$ is well defined (that is,
it does not depend on the chosen pseudocycle representative of $\beta$)
and only depends on the bordism class of $f$.
\end{lemma}
\begin{pf}
This is a consequence of the preceeding lemma together with
remark 7.1.1 and lemma 7.1.5 in \cite{McDS1}.
\end{pf}

For convenience, we extend the map $\Psi_f$ by zero to the rest of 
the homology of $M$.

\subsection{Definitions of the Hamiltonian Gromov-Witten invariants}

\subsubsection{The invariant $\Phi$}

Let $p\in\NN$ and let $c_{p+q}:(F^{p+q})_{S^1}=ES^1\times_{S^1}(F^{p+q})
\to (F_{S^1})^{p+q}=(ES^1\times_{S^1}F)^{p+q}$ be the natural map (we 
consider on $F^{p+q}$ the diagonal action of $S^1$). Given cohomology classes
$\alpha_1,\dots,\alpha_p,\beta_1,\dots,\beta_q\in H^*_{S^1}(F)$
we will write $$c(\alpha_1,\dots,\alpha_p,\beta_1,\dots,\beta_q)=
c_{p+q}^*(\alpha_1\otimes\dots\otimes\alpha_p\otimes\beta_1\otimes
\dots\otimes\beta_q).$$
Let us consider the map $$\oev^{p,q}_J:\AAA\times\SSS^*\times X^p
\to \AAA\times E\times_{S^1}(F^{p+q})$$ which sends 
$(A,\Phi,(x_1,\dots,x_p))$ to 
$(A,\Phi(x_1),\dots,\Phi(x_p),\Phi(y_1),\dots,\Phi(y_q))$.
This map is $\GGG_0$ equivariant, so it descends to give a map
$$\ev^{p,q}_J:\BBB_0\times X^p\to\FpqEj,$$ 
where $\FpqEj=\EE_J\times_{S^1}(F^{p+q})$.
On the other hand, the retraction $R^{\FFF}$ can easily be generalised
to a retraction $R^{p+q,\FFF}:\FpqEj\to\FpqEJj$, where
$\FpqEJj=\FpqEj|_{\Jac_d(X)\times X}$.
Let $f_{\NNN}^{p,q}=R^{p+q,\FFF}\circ\ev^{p,q}_J:\NNN\times X^p\to\FpqEJj$.
\begin{lemma} The map $f_{\NNN}^{p,q}$ is a pseudocycle.
\label{espseudocicle}
\end{lemma}
\begin{pf}
Observe first of all that $\FpqEJj$ is a smooth and oriented manifold.
Let $D\neq D^T$ be a $B$-admissible framing, different from the top one.
For any $1\leq j\leq q$ let us write 
$$T_j=\{y_j\}\amalg\coprod_{x_k=y_j}X_k,$$
that is, $T_j$ is the disjoint union of the marked point $y_j$
and all the bubbles which are mapped to the fibre over $y_j$ by the
cusp $\sigma$-THCs with framing $D$. Note that the Cartesian product
$T_1\times\dots\times T_q$ is a disjoint union of manifolds whose
(complex) dimensions are at most $b(D)$.
Let $B_0\in H_2(F_{S^1})$ be the homology class of the principal component
of a cusp with framing $D$. Consider the evaluation map
$$\oev^{p,q,D}:\bM_{I,\sigma}(B_0,c)\times\MMM_I^*(D)
\times(X_0\amalg\dots\amalg X_k)^p\times(T_1\times\dots\times T_q)
\to\AAA\times E\times_{S^1}(F^{p+q})$$
defined as follows: the point 
$$((A,\Phi_0),(\Phi_1,\dots,\Phi_K),(x_1,\dots,x_p),(z_1,\dots,z_q))$$
is mapped to $(A,\Phi(x_1),\dots,\Phi(x_p),\Phi(z_1),\dots,\Phi(z_q))$,
where $\Phi:X^{\cusp}\to\FFF$ is the map obtained from gluing
the maps $\Phi_0,\Phi_1,\dots,\Phi_K$ (of course, we view
$T_j\subset X_0\amalg\dots\amalg X_K$).
Quotienting the map $\oev^{p,q,D}$ by the reparameterisation groups $G_k$ and
the gauge group $\GGG_0$, restricting to $\RRR_{I,\sigma}(D)$ and
composing with $R^{p+q,\FFF}$ we get a map $f_{\NNN}^{p,q,D}$ with
target $\FpqEJj$. Since $\dim T_1\times\dots\times T_q\leq 2b(D)$, we deduce
(using theorem \ref{gensigmasmoothcusp}) that the domain of
$f_{\NNN}^{p,q,D}$ has dimension $\leq \dim\NNN\times X^p-2$. Finally, the 
compactness theorem \ref{compactificacio} implies that
$$\Omega_{f_{\NNN}^{p,q}}\subset\bigcup_{D\neq D^T}
\Im f_{\NNN}^{p,q,D}$$ 
(here we implicitly use the fact that the image of a cusp $\sigma$-THC in 
$\FFF$ coincides with the image of its reduction).
This proves the lemma.
\end{pf}

Let $\pi:\NNN\to\Jac_d(X)$ be the projection to $\AAA/\GGG_0$ composed
with the retraction $R$. Then $f_{\NNN}^{p,q}\times\pi$ is
also a pseudocycle. Finally, the Hamiltonian Gromov-Witten invariant
$\Phi^{X,F}_{B,c}(\alpha_1,\dots,\alpha_p\mid
\beta_1,\dots,\beta_q\mid\gamma)$ is defined as 
$$\Psi_{f^{p,q}_{\NNN}\times\pi}
(PD(\rho_{\EE_J}^*c(\alpha_1,\dots,\alpha_p,\beta_1,\dots,\beta_q))
\times PD(\gamma)),$$
where $PD$ denotes Poincar{\'e} dual.

This construction gives the invariants sketched in section
\ref{ideaphi}. Indeed, the diagram
$$\xymatrix{\FpqEj\ar[r]^{R^{p+q,\FFF}}\ar[d] & \FpqEJj \ar[d]\\
\AAA/\GGG_0\times X\ar[r] & \Jac_d(X)\times X}$$
is induced by a diagram of $S^1$ principal bundles
$$\xymatrix{\Ej\ar[r]^{R^E}\ar[d] & \EJj \ar[d]\\
\AAA/\GGG_0\times X\ar[r] & \Jac_d(X)\times X,}$$
where $\EJj=\EE_J|_{\Jac_d(X)\times X}$
(this follows from the construction of $R^{p+q,\FFF}$).
Consequently, the map $\rho_{\Ej}^*:H_{S^1}^*(F)\to H^*(\FpqEj)$
is equal to $(R^{p+q,\FFF})^*\rho_{\EJj}^*$, where
$\rho_{\EJj}^*:H_{S^1}^*(F)\to H^*(\FpqEJj)$
(thanks to lemma \ref{apobvi} in the appendix).

\begin{theorem}
The invariant $\Phi^{X,F}_{B,c}$ is well defined, and it only depends on
the manifold $F$, its symplectic structure, the action of $S^1$
and the connected component of $\imag\RR\setminus C_0$ in which
$c$ lies. Furthermore, the invariant 
$$\Phi^{X,F}_{B,c}(\alpha_1,\dots,\alpha_p\mid\beta_1,\dots,\beta_q\mid\nu)$$
is zero unless the following 
relation holds, where $|\alpha|$ denotes the degree of any cohomology class
$\alpha\in H^*_{S^1}(F)$
\begin{equation}
\sum_{j=1}^p|\alpha_j|-2p
\sum_{j=1}^q|\beta_j|+|\nu|=2\la c_1^K(TF),B\ra+2(n-1)(g-1)+1.
\label{conddim}
\end{equation}
\end{theorem}
\begin{pf}
The first claim follows from our discussion on pseudocycles, lemma
\ref{espseudocicle} and theorem \ref{gensigmasmoothN}. The second
claim follows from dimension counting. Just observe that the right hand side 
of formula (\ref{conddim}) is the dimension of the extended moduli space 
$\NNN_{I,\sigma}(B,c)$.
\end{pf}

\subsubsection{The invariant $\oPhi$}
To define the invariant $\Phi$ we have composed the map $\ev^{p,q}_J$
with the retraction $R^{p+q,\FFF}$, thus getting a pseudocycle
$f_{\NNN}:\NNN\times X^p\to\FpqEJj$, where $\FpqEJj$
is a smooth oriented and compact fibration with fibre $F^{p+q}$.
We have then used Poincar{\'e} duality to express the product of 
cohomology classes in terms of intersection of subvarieties.

Let $\oev^{p,q}:\AAA\times\SSS^*\times X^p\to\AAA\times\SSS^*\times
E\times_{S^1}(F^{p+q})$ be the map which sends
$(A,\Phi,x_1,\dots,x_p)$ to $(A,\Phi,\Phi(x_1),\dots,\Phi(x_p),
\Phi(y_1),\dots,\Phi(y_q))$. This map is $\GGG$ equivariant.
The map $\ev^{p,q}=\oev^{p,q}/\GGG$ goes from $\BBB\times X$ and takes
values in $\FpqE=\EE\times_{S^1}(F^{p+q})=
\AAA\times\SSS^*\times (F^{p+q})/\GGG$.
In contrast with what happens with $\FpqEj$, it is not clear whether
$\FpqE$ admits as a retract a smooth oriented compact submanifold.
To use the technique of pseudocycles we will consider the following
construction.

Let us denote to simplify $M=\MMM\times X^p$, $N=\NNN\times X^p$
and $B=\Jac_d(X)$. Let also $F_M=\FpqE$, 
$F_N=\FpqEz=\Ez\times_{S^1}(F^{p+q})$, where 
$\Ez=(\AAA\times\SSS^*\times E)/\GGG_0$ and $F_B=\FpqEJj$. We then have
the following diagram of fibrations
$$\xymatrix{
F_M \ar[d]_{\pi_M} && F_N \ar[ll]_{\pi_1}\ar[d]_{\pi_N}
\ar[rr]^{R^{p+q,\FFF}\pi_3}&& F_B \ar[d]_{\pi_B}\\
M \ar@/_/@{.>}@<-1ex>[u]_{\ev^{p,q}} &&
N \ar@/_/@{.>}@<-1ex>[u]_{\ev^{p,q}_0} \ar[ll]^{\pi_0}\ar[rr]_{R\pi_2}
\ar@{.>}[urr]_{f_{\NNN}^{p,q}} && B,}$$
where $\ev^{p,q}_0=\oev^{p,q}/\GGG_0$,
$\pi_0$ and $\pi_1$ denote the quotients of the action of 
$S^1=\GGG/\GGG_0$, and $\pi_2$ and $\pi_3$ are induced 
(taking quotient of the $\GGG_0$ action) respectively by the
projections $\AAA\times\SSS^*\times X\to\AAA\times X$ and
$\AAA\times\SSS^*\times F\to\AAA\times F$ 
(so that $\ev_J^{p,q}=\pi_3\ev^{p,q}_0$).
Suppose now that we have a smooth oriented and compact manifold $W$
with a free action of $S^1$ and a $S^1$ equivariant map $g:N\to W$.
Then $S^1$ acts on the diagram
$$\xymatrix{F_N\ar[rrr]^{(R^{p+q,\FFF}\pi_3,g\pi_N)}\ar[d] &&& 
F_B\times W\ar[d]\\
N \ar[rrr]^{(R\pi_2,g)} &&& B\times W,}$$
and quotienting we get a diagram of fibrations
$$\xymatrix{F_M\ar[r]^-{R_M}\ar[d] & (F_B\times W)/S^1 \ar[d]\\
M\ar[r] & (B\times W)/S^1.}$$
Now, $(F_B\times W)/S^1$ is a smooth compact oriented fibration with
fibre $F^{p+q}$. Furthermore, the map
$f_{\MMM}^{p,q}=R_M\circ\ev^{p,q}$ is a pseudocycle.
Repeating the construction of $\Phi$ but using $(F_B\times W)/S^1$
instead of $F_B$ we get a rigorous definition of the invariant $\oPhi$.

It remains now to construct $W$ and the $S^1$ equivariant map $g:N\to W$.
Let $P\subset X$ be a finite set of points,
and let $\oev_P:\AAA\times\SSS^*\times X\to\AAA\times F^{|P|}$
be the map which sends $(A,\Phi,x)$ to $(A,\Phi,\prod_{p\in P}\Phi(p))$.
This map is $\GGG_0$ equivariant, so it descends to give a map
$\ev_P:N\to W_0=(\AAA\times F^{|P|})/\GGG_0$.
The group $S^1=\GGG/\GGG_0$ acts on $W_0$.
\begin{lemma}
One can take $P\subset X$ so that $\ov{\Im \ev_P}\cap W_0^{S^1}=\emptyset$.
\end{lemma}
\begin{pf}
For any $\epsilon>0$ we will denote $P_\epsilon\subset X$ any finite subset
such that the union of the disks of radius $\epsilon$
centered at the points $p\in P_\epsilon$ covers $X$. 
Suppose that the claim of the lemma is not true. 
Then there exists a sequence $\epsilon_j\to 0$, sets $P_{\epsilon_j}$
and $\sigma$-THCs $(A_j,\Phi_j)\in\bM_{\sigma,I}(B,c)$ so that for any $j$
the image of the points in $P_{\epsilon_j}$ by the section $\Phi_j$
is contained in $\FFF_0=E\times_{S^1}F^{S^1}$. By the compactification
theorem \ref{compactificacio} one may take a subsequence of $(A_j,\Phi_j)$
which, after suitably regauging, converge pointwise to a cusp
$\sigma$-THC. Now, by construction, the image of the principal component
$X_0$ of this limit cusp must be inside $\FFF_0$. But this
is in contradiction with our assumption that $\sigma\in\Sigma_c(E)$
(see lemma \ref{nodinsfixos}).
\end{pf}

On the other hand, generalizing the construction of $R^{p+q,\FFF}$
we can construct a $S^1$ equivariant retraction 
$R^{W_0}:W_0\to W'={W_0}|_{\Jac_d(X)\times X^p}$.
Then $g'=R^{W_0}\ev_P$ does not meet ${W'}^{S^1}$. Finally, let
$T$ be a $S^1$ invariant tubular neighbourhood of
${W'}^{S^1}\subset W$ whose closure does not meet the closure of
$\Im g$. Then we set $W=(W\setminus T)\cup_{\partial T}-(W\setminus T)$.
This is a smooth compact and oriented manifold with a free action of 
$S^1$ and the map $g'$ induces $g:N\to W$ with the desired properties.

\begin{theorem}
The invariant $\oPhi^{X,F}_{B,c}$ is well defined, and it only depends on
the manifold $F$, its symplectic structure, the action of $S^1$
and the connected component of $\imag\RR\setminus C_0$ in which
$c$ lies. Furthermore, the invariant $$\oPhi^{X,F}_{B,c}
(\alpha_1,\dots,\alpha_p\mid\beta_1,\dots,\beta_q\mid\nu)$$
is zero unless the following 
relation holds, where $|\alpha|$ denotes the degree of any cohomology class
$\alpha\in H^*_{S^1}(F)$
\begin{equation}
\sum_{j=1}^p|\alpha_j|-2p
\sum_{j=1}^q|\beta_j|+|\nu|=2\la c_1^K(TF),B\ra+2(n-1)(g-1).
\label{conddim2}
\end{equation}
\end{theorem}

\section{An example}

In this section we will study a particular case of our construction
in which the invariant $\oPhi$ is nonzero. We will take $F$ to be
the sphere $S^2$ with the action of $S^1$ given by rotation through
a fixed axis. 

\subsection{The data}
Consider on $S^2$ the metric obtained identifying
$S^2$ with the sphere of radius $1$ in $\RR^3$. This metric (as any
other one in $S^2$) is Kaehler. With this choice the symplectic structure
on $F$ is just the induced volume form and the complex structure 
acting on any vector is given by rotation of $90^{\circ}$
counterclockwise (looking towards the center). Suppose that the action 
of $S^1$ is induced by the action of $S^1$ on $\RR^3$ given by rotation 
through the axis $x=y=0$. Then $S^2$ has as fixed points 
$(0,0,1)$ and $(0,0,-1)$. A moment map for this action is
$\mu(x,y,z)=\imag z$.

We can identify $S^2$ with $\CP^1$ in such a way that the action
of $S^1$ is given by $\lambda\cdot [x:y]=[\lambda x:y]$ for any
$\lambda\in S^1\subset\CC^*$. Then $[0:1]$ corresponds to $(0,0,1)$
and $[1:0]$ to $(0,0,-1)$. The action of the complexification $\CC^*$
of $S^1$ takes the same form: any $\lambda\in\CC^*$ sends 
$[x:y]$ to $[\lambda x:y]$. 
The maximal weights are easily seen to be the following.
\begin{equation}
\lambda([x:y];\imag)=\left\{
\begin{array}{l} 1\mbox{ if }y\neq 0 \\
                 -1\mbox{ if }y=0\end{array}\right.
\qquad\mbox{ and }\qquad
\lambda([x:y];-\imag)=\left\{
\begin{array}{l} -1\mbox{ if }x=0 \\
                 1\mbox{ if }x\neq 0\end{array}\right.
\label{mommap}
\end{equation}                 
Let $E\to X$ be a $S^1$ principal bundle.
With our identification $S^2\simeq \CP^1$ we see that
$\FFF=E\times_{S^1}S^2=\PP(L^E\oplus\OOO)$, where $L^E=E\times_{S^1}\CC$
is the line bundle associated to $E$. The action of the gauge
group $\GGG=\Map(X,S^1)$ is given as follows: any section
$\phi=[\phi_0:\phi_1]$, where
$(\phi_0,\phi_1)\in H^0(L^E\otimes K)\oplus H^0(K)$ and
$K\to X$ is a line bundle, is mapped by $g\in\GGG$ to
$\phi=[g\phi_0:\phi_1]$.

\subsection{Some topology}
We begin recalling the Leray-Hirsch theorem: let $V\to X$ be a 
vector bundle of rank $n+1$, where $X$ is any topological space. 
Then, as a ring,
$$H^*(\PP(V))=H^*(X)[t]/(t^{n+1}+t^nc_1(V)+\dots+tc_n(V)+c_{n+1}(V)),$$
where $c_j(V)\in H^{2j}(X)$ is the $j$-th Chern class of $V$.
Furthermore, $t$ has degree $2$ and is the first Chern class
of the line bundle $\OOO_{\PP(V)}(-1)\to\PP(V)$. 

This allows to compute the equivariant cohomology of $S^2$.
Indeed, we have $S^2_{S^1}=ES^1\times_{S^1}S^2=\PP(\OOO(-1)\oplus\OOO)$,
where $\OOO(-1)\to \CP^{\infty}=BS^1$ is the tautological
bundle. We have $H^*(\CP^{\infty})=\ZZ[a]$, where $a=c_1(\OOO(-1))$ so, 
by the Leray-Hirsch theorem, $H^*_{S^1}(S^2)=\ZZ[a,b]/(b^3+b^2a)$
with $\deg a=\deg b=2$. Take as before $E\to X$ a $S^1$ principal bundle, 
let $V=L^E\oplus\OOO$, and let $\pi:\PP(V)\to X$. Then
$\rho_E^*(a)=\pi^*c_1(L^E)$ and $\rho_E^*(b)=c_1(\OOO_{\PP(V)}(-1))$.

\subsection{Hitchin--Kobayashi correspondence and the moduli}
Suppose from now on that $X$ is a compact connected Riemann surface,
and take a holomorphic structure on the line bundle $L^E\to X$. 
Let $\phi:X\to\PP(V)$         
be a holomorphic map. The possible lifts of $\phi$ to a section
$\ophi\in H^0(V)$ are given by $H^0(\phi^*\OOO_{\PP(V)}(-1))$.
Now, since for any line bundle $L\to X$ we have
$\OOO_{\PP(V\otimes L)}(-1)=\OOO_{\PP(V)}(-1)\otimes \pi^*L$,
and since $\PP(V\otimes L)=\PP(V)$ canonically, by taking
$L=\phi^*\OOO_{\PP(V)}(1)$ we can lift $\phi$ to a nonzero section
$\ophi\in H^0(V\otimes L)=H^0(L^E\otimes L\oplus L)$, unique
up to $\CC^*$ (indeed, then we have $\phi^*\OOO_{\PP(V\otimes L)}(-1)=\OOO$).
Let $\ophi=(\ophi_0,\ophi_1)$ in the splitting
$H^0(L^E\otimes L\oplus L)=H^0(L^E\otimes L)\oplus H^0(L)$.
Then $\ophi_0$ and $\ophi_1$ have no common zeros (if they had
then the bundle $\phi^*\OOO_{\PP(V\otimes L)}(-1)$ would be of
positive degree; but by our choice of $L$ it has degree zero).

Of course, the converse is also true: since $\OOO$ is the unique 
line bundle of degree $0$ with a nonzero section, whenever
we have a line bundle $K\to X$ and sections $\ophi_0\in H^0(L^E\otimes K)$
and $\ophi_1\in H^0(K)$ with no common zero the induced section
$\phi\in\Gamma(\PP(V\otimes K))$ has
$\phi^*\OOO_{\PP(V\times K)}(-1)=\OOO$.
The complex gauge transformations which keep fixed the complex structure
of $E$ (that is, the holomorphic gauge transformations) are the constant
ones $\lambda\in\CC^*$, which send any section $\phi$ admitting
a lift $(\ophi_0,\ophi_1)$ to the section induced by 
$(\lambda\ophi_0,\ophi_1)$.

Let $S^jX$ denote the $j$-th symmetric product of $X$.
Let $\Delta^{p,q}\subset S^pX\times S^qX$ be the set of pairs
$s_p\in S^pX$ and $s_q\in S^qX$ with at least one common point.
Let $\AAA$ be the space of connections on $E$, 
$\SSS=\Gamma(\PP(L^E\oplus\OOO))$ and $\GGG^{\CC}=\Map(X,\CC^*)$,
acting on $\AAA\times\SSS$ as usual.
The preceeding discussion proves the following result.
\begin{lemma} Let $\deg(E)=d$. Then
$$\{(A,\Phi)\in\AAA\times\SSS\mid \ov{\partial}_A\Phi=0\}/\GGG^{\CC}
=\coprod_{p-q=d}S^pX\times S^qX\setminus \Delta^{p,q}.$$
\label{simpq}
\end{lemma}
Fix now $c\in\imag\RR$. We want to study which $\GGG^{\CC}$ orbits 
in $\AAA\times\SSS$ contain solutions to the equation
$\Lambda F_A+\mu(\Phi)=c$. Take a pair $(A,\Phi)\in\AAA\times\SSS$,
and suppose that $\Phi$ is contained neither in 
$\PP(0\oplus\OOO)\subset\PP(L^E\oplus\OOO)$ nor in
$\PP(L^E\oplus 0)\subset\PP(L^E\oplus\OOO)$. Then $(A,\Phi)$ is simple,
so we may apply the Hitchin--Kobayashi correspondence
(the usual metric in $\imag\RR$ is induced by the fundamental 
representation $S^1\to\U(1;\CC)$).
To check $c$-stability it is enough to consider the trivial
reduction $\sigma$ of the structure group of $E\times_{S^1}\CC^*$
to $\CC^*$ considered as a parabolic
subgroup of itself, and the antidominant characters
$\chi=\pm 1\in\imag\Lie(S^1)$. This choices give the constant sections
$g_{\sigma,\chi}=\pm 1\in\Map(X,\RR)=\imag \Lie \GGG$.
Now, by our assumptions on $\Phi$ and (\ref{mommap}), we must have
$$\int_{x\in X}\lambda(\Phi(x);-\imag)=
\int_{x\in X}\lambda(\Phi(x);\imag)=\Vol(X).$$
Consequently, $c$ stability amounts to the following two conditions:
$$\deg(E)+\Vol(X)>\Vol(X)\la c,\imag\ra
\mbox{ and }-\deg(E)+\Vol(X)>-\Vol(X)\la c,\imag\ra,$$
which are equivalent to this unique condition:
$$|\deg(E)-\Vol(X)\la c,\imag\ra|<\Vol(X).$$
The crucial point is that this condition does not depend on 
$A$ nor on $\Phi$ (we only made an assumption on $\Phi$ so that
$(A,\Phi)$ is simple). 

Dualising our description of the map $\rho_E^*$ above, we deduce
that fixing a homology class $B\in H_2(S^2_{S^1})$ and considering
the sections $\Phi\in\Gamma(\PP(V))$ such that ${\rho_E}_*\Phi_*[X]=B$
is the same as fixing two integers $(p,q)$ and considering
the sections $\Phi\in\Gamma(\PP(V))$ such that
$\deg E=p-q$ and $\deg\Phi^*\OOO_{\PP(V)}(-1)=-q$.
We will then denote $B=(p,q)$. If $0\neq q\neq p$ then, 
for any connection $A$ and any section $\Phi$ such that 
${\rho_E}_*\Phi_*[X]=(p,q)$, the pair $(A,\Phi)$ is simple.
Consequently, if we take any $c$ such that
$|\deg(E)-\Vol(X)\la c,\imag\ra|<\Vol(X)$ then we deduce
from lemma \ref{simpq} and the Hitchin--Kobayashi correspondence that
$$\MMM^{S^2,S^1}((p,q),c)=S^pX\times S^qX\setminus \Delta^{p,q}.$$

This is not only an identity of sets. It turns out that in this 
situation it is unnecessary to perturb the equations in order to obtain a 
smooth moduli, and the above equality is of smooth manifolds.
Observe for example that the virtual complex dimension for the moduli
is $p+q$. Indeed, $c_1^{S^1}(TS^2)=a-2b$ so
$$\la c_1^{S^1}(TS^2),B\ra+(n-1)(1-g)=(p-q)+2q=p+q$$
(since $n=1$), so it coincides with the actual dimension of 
$S^pX\times S^qX\setminus\Delta^{p,q}$.

\subsection{A nonzero invariant}
Suppose that $X=\CP^1$ and that $0\neq q\neq p$. We will compute in this               
case a nonzero $\oPhi$ invariant. (Observe that we can not hope to
obtain a nonzero $\Phi$ invariant, since all the equivariant cohomology
classes of $S^2$ have even degree, and the extended moduli space has
always odd dimension.)

Let $\MMM=\MMM^{S^2,S^1}((p,q),c)=S^pX\times S^qX\setminus \Delta^{p,q}$.
Observe that $S^nX=\CP^n$, the isomorphism being given by assigning
to $\{[\alpha_1:\beta_1],\dots,[\alpha_n:\beta_n]\}$
the class in $\CP^n$ of the coefficients of
$\prod_j (\alpha_jx-\beta_jy)$.
There is a universal bundle $\EE\to\MMM\times X$ and a universal section
$\ev\in\Gamma(\PP(\EE\otimes \LL\oplus \LL))$, where
$\LL\to\MMM\times X$ is another line bundle. We have
$$PD(\rho_{\EE}^*(a+b))=PD(c_1(\EE\otimes \LL))=\{(s_p,s_q,x)\in             
(S^pX\times S^qX\setminus \Delta^{p,q})\times X\mid
x\in s_p\}$$ and
$$PD(\rho_{\EE}^*b)=PD(c_1(\LL))=
\{(s_p,s_q,x)\in (S^pX\times S^qX\setminus 
\Delta^{p,q})\times X\mid x\in s_q\}.$$
On the other hand, if $\alpha\in H^*(\MMM\times X)$ is any class
then $PD(\alpha/[\pt])=PD(\alpha)\cap \MMM\times\{x\}$ for generic
$x\in X$. So $PD(c_1(\EE\otimes \LL)/[\pt])$ (resp. $PD(c_1(\LL)/[\pt])$)
is $H^p\times\CP^q\setminus\Delta^{p,q}$
(resp. $\CP^p\times H^q\setminus\Delta^{p,q}$)
where $H^p$ (resp. $H^q$) is a hyperplane of $S^pX=\CP^p$ 
(resp. $S^qX=\CP^q$).

We take as a compactification of our moduli $\bM=S^pX\times S^qX$.
The complementary $\bM\setminus\MMM$ has real codimension at least two.
The bundles $\EE,\LL$ extend naturally to bundles
$\oEE,\oLL\to\bM\times X$ and the section $\ev$ also extends
to a section $\oev\in\Gamma(\PP(\oEE\otimes\oLL\oplus\oLL))$.
By the theory of pseudocycles, the preceeding observations imply
that we can compute the invariants using this particular compactification
and without worrying whether it coincides with the compactification
with cusp THCs. But we then obtain that 
$$\oPhi^{\CP^1,S^2}_{(p,q),c}(\mid
\overbrace{a+b,\dots,a+b}^{p},\overbrace{b,\dots,b}^{q}\mid)=1.$$
(That is, we use the map $\omu_0$ in all the $p+q$ arguments.)

\appendix

\chapter{Some useful results}

\section{Vector bundles over fibre bundles}
\label{vbfb}            

Let $\pi^V:V\to F$ be a vector bundle. 
We denote $TV_v\subset TV$ the subbundle $\Ker d\pi^V$ of
vertical tangent vectors to $V$. Suppose that a compact
connected Lie group $K$ acts on $V$ linearly on the fibres.
Let $\nabla$ be a $K$-invariant connection on $V$. For any
vector field $\fX\in\Gamma(TF)$, let $\sigma_{\nabla}(\fX)\in\Gamma(TV)$
denote the lift of $\fX$ given by $\nabla$. 

\begin{definition}
The {\bf moment} of the action of $K$ on $V$ with respect to $\nabla$
is the map $\Omega=\Omega^V:\klie\to\End V$ defined as
$\Omega(s)=\sigma_{\nabla}(\fX_s^F)-\fX_s^V$
for any $s\in\klie$. 
\end{definition}

Note that, as defined, $\Omega(s):V\to TV_v$ is map
linear on the fibers of $V$.
We denote this by $\Omega(s)\in\Gamma(V;TV)^{\lin}$.
Using the canonical isomorphism
$TV_v\simeq V$ we regard $\Omega(s)\in\End V$.
The next lemma follows from an easy computation.

\begin{lemma}
For any $s\in\klie$ and any $\fX\in\Gamma(V)$ we have
$$L_{\fX^F_s}\fX=\nabla_{\fX^F_s}\fX+\Omega(s)\fX.$$
\label{Lieder}
\end{lemma}

\begin{example}
Let $(F,\omega_F)$ be a symplectic manifold, and assume that there is 
a line bundle $V\to F$ with a connection $\nabla$ whose curvature
is $\imag\omega_F$. Assume that $K$ acts on $F$ respecting $\omega_F$.
Now, if the action of $K$ on $F$ lifts to an action on $L$ we may
average and assume that $\nabla$ is $K$ invariant. Then, the resulting
moment $\Omega$ is a symplectic moment map for the action of $K$ on $F$.
Conversely, any symplectic moment map on $F$ gives rise to a lift of the
infinitesimal action of $\klie$ which leaves $\nabla$ invariant.
(See p. 244 in \cite{DoKr}.)
\end{example}

Let $\pi:E\to X$ be a $K$-principal bundle. Let $\VVV=E\times_K V$
and $\FFF=E\times_K F$. Then $\pvf:\VVV\to\FFF$ is a vector bundle.
Let $A$ be a connection on $E$. 

\subsection{A connection on $\VVV\to\FFF$}

We define a connection $\nv=\nv(A,\nabla)$ on $\VVV\to\FFF$
by giving the lift of any tangent field $\fX\in\Gamma(T\FFF)$
to a field $\sigma_{\nv}(\fX)\in\Gamma(T\VVV)$.

Consider the following diagram
$$\xymatrix{
K \ar[r] & E\times V \ar[r]^{q_{\VVV}} \ar[d]^{\id\times\pi^V} & 
\VVV \ar[d]^{\pvf} \\
K \ar[r] & E\times F \ar[r]^{q_{\FFF}} \ar[d]^{\pi_E} &
\FFF \ar[d]^{\pfx} \\
K \ar[r] & E \ar[r]^{\pi} & X.}$$

The three horizontal sequences of maps denote $K$-principal
bundles, and the two squares are cartesian. In particular, the
principal $K$ bundle $E\times F\to\FFF$ is isomorphic to the
pullback $\pfx^*E$. Consider on $E\times F\to\FFF$ the pullback
connection $A^{E\times F}=\pfx^*A$. This connection gives a lift
$\Sigma_A:T\FFF\to T(E\times F)$.
On the other hand, we have a map 
$$\id\times\sigma_{\nabla}:T(E\times F)\to T(E\times V).$$
We define
$$\sigma_{\nv}:=dq_{\VVV}(\id\times\sigma_{\nabla})
\Sigma_A:T\FFF\to T\VVV.$$

Alternatively, we can define the connection $\nv$ by giving 
the projection $$\rho_{\nv}:T\VVV\exh T\VVV_v,$$ where
$T\VVV_v=\Ker d\pvf$ is the field of vertical tangent vectors
of the fibration $\pvf:\VVV\to\FFF$. 
The connection $A$ induces a connection on the fibration $\pvx:\VVV\to X$,
which gives a projection $\rho_{A^{\VVV}}:T\VVV\exh\Ker d\pvx$.
On the other hand, the connection $\nabla$ gives a projection
$\rho_{\nabla}:TV\exh TV_v$. Since $\nabla$ is $K$-invariant,
so is the projection $\rho_{\nabla}$, and hence $\rho_{\nabla}$ 
extends globally to give another projection 
$\rho_{\nabla^{\FFF}}:TV_v\exh \Ker d\pvf$.
Then $\rho_{\nv}=\rho_{\nabla^{\FFF}}\circ\rho_{A^{\VVV}}$.

\subsection{The curvature of $\nv$}

Let $F_{\nabla}\in\Omega^2(F;\End V)$ be the curvature of $\nabla$.
Since $\nabla$ is $K$-equivariant, so is $F_{\nabla}$, and
hence $F_{\nabla}$ extends fibrewise to a map
$$F_{\nabla}:\Lambda^2 T\FFF_v\to\End \VVV.$$
Consider the projection $\rho_{A^{\FFF}}:T\FFF\to T\FFF_v$ given by
the connection $A$. 
We have $\rho_{A^{\FFF}}=1-\sigma_{A^{\FFF}} d\pfx$,
where $\sigma_{A^{\FFF}}:TX\to T\FFF$ is the lift given by $A$.
We then have $$\rho_{A^{\FFF}}^*F_{\nabla}\in\Omega^2(\FFF;\End \VVV).$$
Let $F_A\in\Omega^2(X;E\times_{\Ad}\klie)$ be the curvature of $A$.
Then 
$$\Omega(\pfx^* F_A)\in\Omega^2(\FFF;\End \VVV).$$

\begin{lemma}
The curvature $F_{\nv}$ of $\nv$ is 
$$F_{\nv}=\rho_{A^{\FFF}}^*F_{\nabla}-\Omega(\pfx^* F_A).$$
\label{curvatura}
\end{lemma}
\begin{pf}
Since the question is local on $X$, we may 
suppose that there is
an isomorphism $E\simeq X\times K$. Then there are induced isomorphisms
$\FFF\simeq X\times F$ and $\VVV\simeq X\times V$.
Using these isomorphisms we may write the map 
$dq_{\VVV}:T(X\times K\times V)\to T(X\times V)$ as
\begin{equation}
dq_{\VVV}(u,s,v)=(u,v+\fX_s^V)
\label{dosa}
\end{equation}
where $(u,s,v)\in TX\times\klie\times TV\simeq T(X\times K\times V)$.
We have the following diagram

$$\xymatrix{
TX\times\klie\times TV \ar[r]^{dq_{\VVV}} \ar[dd]_{\id\times\id\times d\pi^V}
& TX\times TV \ar[dd]_{d\pvf} \\ & \\ 
TX\times\klie\times TF \ar[r]^{dq_{\FFF}} 
                       \ar[dd]_{d\pi_E=\pi_{TX\times\klie}}
                       \ar@/_/@{.>}@<-1ex>[uu]
                         _{\id\times\id\times\sigma_{\nabla}} 
& TX\times TF \ar@/_1pc/@{.>}@<-1ex>[l]_{\Sigma_A}
              \ar[dd]_{d\pfx} 
              \ar@/_/@{.>}@<-1ex>[uu]_{
              \sigma_{\nv}}  \\ & \\
TX\times\klie \ar[r]^{d\pi}
                \ar@/_/@{.>}@<-1ex>[uu]_{\id\times 0}
& TX \ar@/_/@{.>}@<-1ex>[uu]_{\sigma_{A^{\FFF}}}
     \ar@/_1pc/@{.>}@<-1ex>[l]_{\sigma_A}.}$$

Take two fields $\fY_1,\fY_2\in\Gamma(TX\times TF)$. 
We identify the sections in $\Gamma(\VVV;T\VVV_v)^{\lin}$ linear
in the fibres of $\VVV$ with $\End\VVV$. We compute
\begin{align}
F_{\nv}(\fY_1,\fY_2) &= [\sigma_{\nv}\fY_1,\sigma_{\nv}\fY_2]
-\sigma_{\nv}[\fY_1,\fY_2] \notag \\
&=
[dq_{\VVV}(\id\times\sigma_{\nabla})\Sigma_A(\fY_1),
dq_{\VVV}(\id\times\sigma_{\nabla})\Sigma_A(\fY_2)]\notag\\
&-
dq_{\VVV}(\id\times\sigma_{\nabla})[\Sigma_A(\fY_1),
\Sigma_A(\fY_2)] \notag \\
&+
dq_{\VVV}(\id\times\sigma_{\nabla})[\Sigma_A(\fY_1),
\Sigma_A(\fY_2)] -
dq_{\VVV}(\id\times\sigma_{\nabla})\Sigma_A
[\fY_1,\fY_2]. \label{curva}
\end{align}
Let $\pi_3:TX\times\klie\times TF\to TX\times\klie\times TF$
be the map $\pi_3(u,s,v)=(0,0,v)$. The lift $\Sigma_A\fY_i$ is defined
by these properties
\begin{align}
& dq_{\FFF}(\Sigma_A\fY_i)=\fY_i \notag \\
& d\pi_E(\Sigma_A\fY_i)=\sigma_A(d\pfx(\fY_i))
\notag
\end{align}
and, on the other hand, $\sigma_{A^{\FFF}}=dq_{\FFF}(\id\times 0) \sigma_A$.
Now, since $\pi_3=1-(\id\times 0)d\pi_E$,
\begin{align}
dq_{\FFF}\pi_3\Sigma_A(\fY_i)
&= dq_{\FFF}\Sigma_A(\fY_i)-
dq_{\FFF} (\id\times 0) d\pi_E
\Sigma_A(\fY_i) \notag \\
&= \fY_i-dq_{\FFF} (\id\times 0)\sigma_A(d\pfx(\fY_i))\notag\\
&= \fY_i-\sigma_{A^{\FFF}} d\pfx(\fY_i)=\rho_{A^{\FFF}}(\fY_i).\notag
\end{align}
The map $dq_{\FFF}$ restricted to $\{0\}\times\{0\}\times TF$
sends any $(0,0,v)$ to $(0,v)$, so 
$\pi_3\Sigma_A(\fY_i)=(0,0,\rho_{A^{\FFF}}(\fY_i)).$
This implies that the first two summands in (\ref{curva})
are equal to
$$dq_{\VVV}(0,0,F_{\nabla}(\rho_{A^{\FFF}}(\fY_1),
\rho_{A^{\FFF}}(\fY_2)))=\rho_{A^{\FFF}}^* F_{\nabla}(\fY_1,\fY_2).$$
On the other hand, observe that 
\begin{align}
[\Sigma_A(\fY_1),\Sigma_A(\fY_2)]
&-\Sigma_A[\fY_1,\fY_2]=
\fX_{\pfx^* F_A (\fY_1,\fY_2)}^{E\times F} \notag\\
&= (0,\pfx^* F_A(\fY_1,\fY_2),-\fX_{\pfx^* F_A (\fY_1,\fY_2)}^F)
\in TX\times\klie\times TF. \notag
\end{align}
Using formula (\ref{dosa}) this implies that the last two
summands in (\ref{curva}) are equal to 
$$-\Omega(\pfx^* F_A).$$
\end{pf}

\subsection{The case $V=TF$}
\label{nablatf}
Suppose that $V=TF$ and that the action of $K$ on $V$ is the derivation of
the action of $K$ on $F$. Then $\nabla$ is a connection on $TF$. Using the 
connection $A$ to view $\pfx^*TX$ as a subbundle of $T\FFF$, we may write
\begin{equation}
T\FFF\simeq \VVV\oplus\pfx^*TX.
\end{equation}
Take on $TX$ a connection $\nx$. By means of the above splitting we
define a connection on $T\FFF$ as
$$\nf=\nf(A,\nabla,\nx)=
\nv\oplus \pfx^*\nx.$$

\begin{lemma}
The torsion of $\nf$ is
$$\Tor_{\nf}=\sigma_{A^{\FFF}}\circ d\pfx^*\Tor_{\nx}
-d\pfx^* F_{A^{\FFF}}+\rho_{A^{\FFF}}^*\Tor_{\nabla}.$$
\label{torsio}
\end{lemma}
\begin{remark}
The second summand should be interpreted as follows.
\begin{align}
d\pfx^* F_{A^{\FFF}}(\fY_1,\fY_2)
&=F_{A^{\FFF}}(d\pfx(\fY_1),d\pfx(\fY_2)) \notag\\
&=[\sigma_{A^\FFF}d\pfx(\fY_1),\sigma_{A^\FFF}d\pfx(\fY_2)]
-\sigma_{A^\FFF}[d\pfx(\fY_1),d\pfx(\fY_2)],\notag
\end{align}
which is a field of tangent vectors on $\FFF$. This field of 
tangent vectors may be written locally (using a trivialisation
$E\simeq X\times K$) as $\fX^F_{F_A}$.
\end{remark}
\begin{pf}
The statement is local on $X$, so we can suppose that there is
an isomorphism $E\simeq X\times K$. Then there are induced isomorphisms
$\FFF\simeq X\times F$ and $\VVV\simeq X\times V$.

Suppose first that $\fY_i=\sigma_{A^{\FFF}}\fZ_i$,
where $\fZ_1,\fZ_2\in\Gamma(TX)$. Clearly, $d\pfx(\fY_i)=\fZ_i$.
By definition
$$\nf_{\fY_1}\fY_2=\sigma_{A^{\FFF}}\nx_{\fZ_1}\fZ_2,$$
so that
\begin{align}
\Tor_{\nf}(\fY_1,\fY_2) &=
\nf_{\fY_1}\fY_2-\nf_{\fY_2}\fY_1-[\fY_1,\fY_2] \notag \\
&= \sigma_{A^{\FFF}}(\nx_{\fZ_1}\fZ_2-\nx_{\fZ_2}\fZ_1-[\fZ_1,\fZ_2])
-F_{A^{\FFF}}(\fZ_1,\fZ_2) \notag \\
&= \sigma_{A^{\FFF}}(\Tor_{\nx}(\fZ_1,\fZ_2))-F_{A^{\FFF}}(\fZ_1,\fZ_2).
\notag
\end{align}

Now suppose that $\fY_i=(0,\fZ_i)\in TX\times TF\simeq T\FFF$,
where $\fZ_1,\fZ_2\in\Gamma(TF)$.
$\nf_{\fY_1}\fY_2=\pi_F^*\nabla_{\fZ_1}\fZ_2$ and
$[\fY_1,\fY_2]=\pi_F^*[\fZ_1,\fZ_2]$. This implies that
$$\Tor_{\nf}(\fY_1,\fY_2)=\pi_F^*\Tor_{\nabla}(\fZ_1,\fZ_2).$$

Finally, suppose that $\fY_1=\sigma_{A^{\FFF}}\fZ_X$ and $\fY_2=(0,\fZ_F)$,
where $\fZ_F\in\Gamma(TF)$ and $\fZ_X\in\Gamma(TX)$.
By definition 
\begin{equation}
\Tor_{\nf}(\fY_1,\fY_2)=
\nf_{\fY_1}\fY_2-\nf_{\fY_2}\fY_1-[\fY_1,\fY_2]
\label{casmixte}
\end{equation}
The second summand in the left hand side is equal to zero:
$$\nf_{\fY_2}\fY_1=\sigma_{A^{\FFF}}\nx_{d\pfx\fY_2}\fZ_X=0.$$ 
On the other hand, if we view $\fY_2$ as a map
$\fY_2:\FFF\to T\FFF$, then
$$\nf_{\fY_1}\fY_2=\fY_2^*(d_{\fY_1}\fY_2-\sigma_{\nv}(\fY_1)).$$
We may write the lift $\sigma_A$ as 
$\sigma_A(\fY)=(\fY,s(\fY))$ for any field $\fY\in\Gamma(TX)$,
where $s\in\Omega^1(X;\klie)$. Then one easily sees that
$\Sigma_A(\fY_1)=(\fZ_X,s(\fZ_X),0)$ and that
$\sigma_{\nv}(\fY_1)=(\fZ_X,\fX^{TF}_{s(\fZ_X)})$.
In addition, $d_{\fY_1}\fY_2=(0,d_{\fX_{s(\fZ_X)}^F}\fZ_F)$,
so that $$\nf_{\fY_1}\fY_2=(0,\nabla_{\fX_{s(\fZ_X)}^F}\fZ_F)
+(0,\Omega(s(\fZ_X))\fZ_F).$$
Finally, the right hand side is equal, by lemma \ref{Lieder},
to $(0,[s(\fZ_X),\fZ_F])$, which is equal to $[\fY_1,\fY_2]$. Hence the 
whole expression in (\ref{casmixte}) vanishes.
\end{pf}

\section{Actions of compact groups on manifolds}

\subsection{Left and right actions}
We recall that an action of a group $K$ on $F$ is said to be on the left
(resp. on the right) if for any $g,g'\in K$ and any $x\in F$
we have $(gg')\cdot x=g\cdot(g'\cdot x)$ (resp.
$(gg')\cdot x=g'\cdot(g\cdot x)$), where $\cdot$ denotes the action.
The following lemma is rather elementary, although of crucial importance
in the computations.

\begin{lemma}
Suppose that a Lie group $K$ acts on $F$ on the left (resp. on the right).
Let $\fX_s$ be the field generated on $F$ by the infinitesimal action of 
an element $s\in\klie$. We then have $$[\fX_s,\fX_{s'}]=-\fX_{[s,s']}$$ 
(resp. $[\fX_s,\fX_{s'}]=\fX_{[s,s']}$) for any $s,s'\in\klie$.
\end{lemma}
\begin{pf}
See for example \cite{BeGeV} p. 208.
\end{pf}

\subsection{Stabiliser groups}
Let $F$ be a compact manifold and $K$ a compact Lie group. Assume
that $K$ acts smoothly on $F$ (say, on the left). Let $x\in F^K$ be a fixed
point of the action (that is, $k\cdot x=x$ for any $k\in K$).
Then there is a linear action of $K$ on the tangent space at $x$
$$\rho:K\to\GL(T_xF).$$
Take a $K$-invariant metric $g$ on $F$. The exponential map
$\exp^g_x$ going from a neighbourhood $U$ of $0\in T_xF$ to a neighbourhood
$V$ of $x\in F$ is a $K$-equivariant diffeomorphism. Hence a neighbourhood
of $x\in F^K\subset F$ is diffeomorphic to a neigbourhood of 
$0\in T_xF^K\subset T_xF$, the
fixed point set of the linear representation $\rho$. 
This proves the following lemma.

\begin{lemma} The fixed point set $F^K$ is a finite union
of disjoint smooth submanifolds of $F$. Furthermore, if $F_1,\dots,F_r$
are the connected components of $F^K$ and $x\in F_k$, then 
$$\dim F_k=\dim T_xF^K.$$
\end{lemma}

\begin{lemma} Suppose that $K$ is abelian. Then the set of subgroups
of $K$ arising as stabilisers of points in $F$ is a finite set.
\label{stabfinit}
\end{lemma}
\begin{pf}
Let $x\in F$, and let $S_x=\Stab_Kx$. Set $T=K/S_x$. 
For any $y\in F$, let $S_y'$ be the projection on $T$ of $\Stab_Ky$. 
It is easy to see that there exists a neighbourhood $T_\epsilon\subset T$ 
of the identity such that no nontrivial subgroup of $T$ is contained
in $T_\epsilon$. Hence, if there is a sequence $y_n\to x$
such that $S_y'\neq 1$, then $S_x'\neq 1$ as well, which is a contradiction.
This implies that there is a neighbourhood $U_x$ of $x$ such that
any $y\in U_y$ has stabiliser contained in $S_x$. 
On the other hand, the preceeding reasonement for the group $S_x$ implies
that a small neighbourhood of $x$ is $S_x$-equivariantly diffeomorphic to a 
neighbourhood of $0$ in the $S_x$ space $T_xF$.
But the set of stabilisers in $S_x$ of points in $T_xF$ is finite
(to see this, just split $T_xF$ in eigenspaces). 
Consequently, if $U_x$ is small enough
then the set of $K$ stabilisers of points in $U_x$ is finite.
Taking a finite subcovering of $\{U_x\mid x\in F\}$ (recall that we
are assuming $F$ to be compact) we get the desired result.
\end{pf}

For any connected component $F_j$ of the fixed point set $F^K$
we have an action of $K$ on the normal bundle $N_F(F_j)$ lifting
the identity on $F_j$. Fix any maximal torus $T\subset K$.
Since the group $K$ is compact, the weights
of the representation at $T_xF$ for $x\in F_j$ are the same for
all the points of $F_j$. Hence we may talk about the weights of the
representation of $K$ on $N_F(F_j)$ and they form a finite set.

\begin{definition}
The action of $K$ on $F$ is called {\bf almost-free} if the stabiliser
subgroup $\Stab_Kx$ of any $x\in F$ is connected.
\end{definition}

\begin{lemma}                                       
Suppose that $K=S^1$.  An action of $S^1$ on $F$ is almost-free if and only
if it is free at the complementary of the fixed point set $F^{S^1}$. In 
this case, the weights of the representation on the normal bundle of any 
connected component $F_j$ of $F^{S^1}$ are $1$ or $-1$.
\label{nomesmesmenysu}
\end{lemma}
\begin{pf}             
The only connected subgroups of $S^1$ are the identity $\{1\}\subset S^1$ 
and $S^1$, and the first claim follows. We now prove the second claim.
For any $x\in F_j$ the tangent space
splits $T_xF=\oplus V_\lambda$, where $S^1$ acts on $V_\lambda$
with weight $w(\lambda)$ and where $\dim_\RR V_\lambda$ is $1$ if
$w(\lambda)=0$ and it is $2$ otherwise. The trivial representations
span the tangent space $T_xF^{S^1}$.
Any $x\in V_\lambda\setminus\{0\}$ such that $w(\lambda)\neq 0$
has as stabiliser the cyclic group of $|w(\lambda)|$ elements.
So if not all $w(\lambda)$ belonged to $\{-1,0,1\}$ then there would
be some $v\in T_xF$ with stabiliser $\ZZ/m\ZZ$, where $m\geq 2$.
Finally, since $\exp_x^g$ is $S^1$ equivariant for a $S^1$-invariant
metric $g$, taking $v$ small enough we would get a point
$\exp_x^g(v)\in F$ whose stabiliser would be exactly $\ZZ/m\ZZ$.
Hence the action would not be almost-free.
\end{pf}

\begin{remark} An action of $K$ on $F$ which is free in $F\setminus F^K$
is called {\bf semi-free}.
\end{remark}

\begin{corollary}
Let $S^1$ act on $F$ nontrivially in an almost-free way, and let
$x$ belong to the fixed point set $F^{S^1}$. Then the induced
representation of $S^1$ on $T_xF$ has at least one weight equal
to $\pm 1$ (and those which are not equal to $\pm 1$ are equal to $0$).
\label{omesuomenysuozero}
\end{corollary}

\section{Principal bundles on a compact Riemann surface}

Let $K$ be a compact connected Lie group and let $X$ be a compact
connected Riemann surface. Let us write $P_K(X)$ for the
isomorphism classes of $K$-principal bundles over $X$.
            
Let $x_0\in X$ be a fixed point and let 
$X^*=X\setminus\{x_0\}$. Take any $E\in P_K(X)$. Since $X^*$
is homotopic to a bouquet of circles, the restriction of
$E$ to $X^*$ is trivial. On the other hand, $E$ restricts to a
trivial bundle on a small disk $D$ around $x_0$. 
Let us fix trivialisations of $E|_{X^*}$ and of $E|_D$.
The topological type of $E$ is described by the homotopy class of the 
transition function $\gamma(E):D\cap X^*\to K$ relating both 
trivialisations. Since $D\cap X^*\sim S^1$, we get an element
$\gamma(E)\in\pi_1(K)$.

The following lemma is easily proved.

\begin{lemma}
The map $\gamma:P_K(X)\to\pi_1(K)$ is a bijection.
\label{zerou}
\end{lemma}

We define a map $\eta:P_K(X)\to H_2(BK;\ZZ)$ as follows.
Let $E\in P_K(X)$ and let $\xi_0\in [X,BK]$ be its classifying map.
We set $\eta(E):=(\xi_0)_*([X])$.

\begin{lemma}
The map $\eta:P_K(X)\to H_2(BK;\ZZ)$ is a bijection.
\label{etabij}
\end{lemma}
\begin{pf}
Since $\pi_1(BK)=0$, the Leray-Serre spectral sequence applied
to $K\to EK\to BK$ tells us that
$$E^2_{p,q}=H_p(BK;H_q(K;\ZZ))\Rightarrow H_{p+q}(EK)=0.$$
From the fact that the sequence converges to zero we deduce that
the map $d^2:E^2_{2,0}\to E^2_{0,1}$ must be an isomorphism.
But we clearly have $E^2_{2,0}=H_2(BK;\ZZ)$ and $E^2_{0,1}=H_1(K;\ZZ)$,
so we obtain $d^2:H_2(BK;\ZZ)\stackrel{\sim}{\longrightarrow}H_1(K;\ZZ)$.
Using the canonical isomorphism $H_1(K;\ZZ)\simeq \pi_1(K)$
(since $K$ is a Lie group $\pi_1(K)$ is abelian) we may consider
$d^2\circ\eta$ as a map from $P_K(X)$ to $\pi_1(K)$.
Now, tracing the construction of the Leray-Serre spectral sequence
(see for example \cite{McCl}) one easily verifies that 
$d^2\circ \eta=\gamma$. Hence by lemma \ref{zerou} $\eta$ is a bijection.
\end{pf}

\section{(Co)homology classes in fibre bundles}
\label{pullingback}

In this subsection all the (co)homology groups will be taken, unless
otherwise stated, with coefficients in $\ZZ$. 

\begin{lemma}
Let $F$ will be a left $K$ space. Let $E\to X$ be a $K$ principal bundle on 
a topological space $X$, and let $\FFF=E\times_K F$.
There is a natural map $H_*(F)\to H_*(\FFF)$.
\label{triviallocsys}
\end{lemma}
\begin{pf}
We define the map by taking any point $x\in X$, trivialising
the fibre $E_x\simeq K$ and using the induced trivialisation
$\FFF_x\simeq F$ to map $H_*(F)\to H_*(\FFF_x)\to H_*(\FFF)$,
the lattest map given by the inclusion $\FFF_x\subset\FFF$.
This is well defined because the system of local coefficients
$\HHH_*(\FFF_x)$, where $x\in M$, is trivial. Indeed, since $K$ 
is connected, the action of $K$ on $H_*(F)$ is trivial. See p. 154 
in \cite{McCl} for a definition of a system of local coefficients.
\end{pf}

\begin{prop}
Let $E\to X$ be a $K$ principal bundle over a finite CW complex $X$. 
Let $c:X\to BK$
be the classifying map for $E$, and let $\phi:E\simeq c^*EK$
be any isomorphism (in other words, $\phi$ is a $K$ equivariant
map from $E$ to $EK$). Let $F$ be a finite CW complex with a left 
continous action of $K$, and let $\FFF=E\times_K F$ and $F_K=EK\times_K F$. 
The map $\phi\times\id:E\times F\to EK\times F$ descends to a map 
$\psi:\FFF\to F_K$. Let $g:E\to E$
be any gauge transformation, and let $g(\FFF):\FFF\to\FFF$ be
the induced gauge transformation. Then
$$\psi_*=\psi_*g(\FFF)_*:H_*(\FFF)\to H_*(F_K)$$
and similarly $\psi^*=g(\FFF)^*\psi^*:H^*(F_K)\to H^*(\FFF)$.
\label{indepframe}
\end{prop}
\begin{pf} For a very readable explanation of the techniques used
in this proof see chapter 7 in \cite{McCl}. We prove the statement
on homology. Let $\DGAlg$ be the category of differential graded 
algebras over $\ZZ$, and let $\Gamma$ be any object in $\DGAlg$.
Let $\DGMG$ (resp. $\DGGM$) be the category of differential right
(resp. left) modules over $\Gamma$ (these modules are in particular
modules over $\ZZ$). For any object $M$ of $\DGMG$ one can define 
the tensor product functor 
$M\otimes_{\Gamma}\cdot:\DGGM\to\DGGM$. This functor admits 
right derived functors $\Tor^{\Gamma}(M,\cdot)$. Suppose that
$\Lambda$ is another object of $\DGAlg$, $M'$ (resp. $N'$) is 
an object of $\DGML$ (resp. $\DGLM$) and that there are maps 
$f\in\Mor_{\DGAlg}(\Gamma,\Lambda)$ $g\in\Mor_{\ModZ}(M,M')$  
and $h\in\Mor_{\ModZ}(N,N')$ (here $\ModZ$ is the category of
$\ZZ$ modules) satisfying the obvious compatibility relations. 
Then there is a natural map 
$\Tor^f(g,h):\Tor^{\Gamma}(M,N)\to\Tor^{\Lambda}(M',N')$.
Furthermore, if $f,g,h$ induce isomorphisms in homology,
then $\Tor^f(g,h)$ is also an isomorphism (this is a theorem of
Moore; see corollary 7.6 in \cite{McCl}).

The actions $E\times K\to E$ and $K\times F\to F$ give $C_*(E)$ (resp. 
$C_*(F)$) a structure of differential right (resp. left)
module ofer $C_*(K)$ ($C_*(X)$ denotes the module of singular chains in 
$X$). A theorem of Moore (see theorem 7.16 in \cite{McCl}) asserts that 
there is a natural isomorphism
\begin{equation}
\theta_*:\Tor^{C_*(K)}(C_*(E),C_*(F))\stackrel{\simeq}{\longrightarrow}
H_*(E\times_K F).
\label{isoMo}
\end{equation}
Let now $g:E\to E$ be any gauge transformation. This induces a morphism
$g_*\in\Mor(C_*(E),C_*(E))$ as differential right module over $C_*(K)$
(thanks to equivariance under the action of $K$ on the right: the action
of the gauge group is on the left). 
Let $p:EK\to \pt$ be the map sending $EK$ to a point.
Since $EK$ is contractible, $p$ induces an isomorphism in homology.
Consider now the following diagram
$$\xymatrix{\Tor^{C_*(K)}(C_*(E),C_*(F)) \ar[d]_{\Tor^{\id}(g_*,\id)}
\ar[rr]^{\Tor^{\id}(\phi_*,\id)} 
&& \Tor^{C_*(K)}(C_*(EK),C_*(F)) \ar[d]^{\Tor^{\id}(p_*,\id)}\\
\Tor^{C_*(K)}(C_*(E),C_*(F)) \ar[rru]^{\Tor^{\id}(\phi_*,\id)}
&& \Tor^{C_*(K)}(C_*(\pt),C_*(F))}$$
Since $\psi_*$ corresponds to $\Tor^{\id}(\phi_*,\id)$ 
and $g(\FFF)_*$ to $\Tor^{\id}(g_*,\id)$ 
through the isomorphism (\ref{isoMo}), we have to prove that 
$$\Tor^{\id}(\phi_*,\id)=\Tor^{\id}(\phi_*,\id)\circ\Tor^{\id}(g_*,\id).$$
But since $p\phi=p\phi g$ we clearly have 
$$\Tor^{\id}(p_*,\id)\circ\Tor^{\id}(\phi_*,\id)=
\Tor^{\id}(p_*,\id)\circ\Tor^{\id}(\phi_*,\id)\circ\Tor^{\id}(g_*,\id)$$
and, since $\Tor^{\id}(p_*,\id)$ is an isomorphism, we are done.

The claim on cohomology is proved following exactly the same technique,
but using differential comodules over $C^*(K)$.
\end{pf}

\begin{corollary}
Given $E\to X$ and $F$ as above, there are canonical maps
$$\rho_*:H_*(\FFF)\to H_*(F_K)\qquad\mbox{ and }
\qquad \rho^*:H^*(F_K)\to H^*(\FFF).$$
\end{corollary}
\begin{pf}
We put $\rho_*=\psi_*$ and $\rho^*=\psi^*$, where $\psi$ is 
as in proposition \ref{indepframe}. The only extra choice
made in the definition of $\psi$ was a $K$ equivariant map 
$\phi:E\to EK$. Any two such maps $\phi_0,\phi_1$ 
are related by a gauge transformation: $\phi_0=\phi_1g$.
Then proposition \ref{indepframe} implies that the corresponding maps
$\psi_0,\psi_1:\FFF\to F_K$ induce the same maps in (co)homology:
${\psi_0}_*={\psi_1}_*$ and $\psi_0^*=\psi_1^*$.
\end{pf}

\begin{remark}
Observe that the map $\rho$ does not only depend on the bundle $\FFF$:
it depends also on the structure group $K$ (this is obvious!) and on the 
particular $K$ principal bundle $E\to X$ such that $\FFF=E\times_K F$.
When this is not clear from the context we will write $\rho(E)$ or
$\rho(K,E,F)$ instead of $\rho$.
\end{remark}

Finally, we quote the following trivial corollary of the definition
of the map $\rho^*$.

\begin{corollary}                           
Consider a map of $S^1$ principal bundles
$$\xymatrix{E \ar[r]^R \ar[d] & E' \ar[d] \\
X \ar[r] & X'.}$$
Then the map $\rho_E^*:H_{S^1}^*(F)\to H^*(E\times_{S^1}F)$
is equal to $R_F^*\rho_{E'}^*$, where
$\rho_{E'}^*:H_{S^1}^*(F)\to H^*(E'\times_{S^1}F)$
and where $R_F^*:H^*(E'\times_{S^1}F)\to H^*(E\times_{S^1}F)$
is given by the map $R_F:E\times_{S^1}F\to E'\times_{S^1}F$
induced by $R$.
\label{apobvi}
\end{corollary}

\end{document}